\documentclass[12pt]{article}

\usepackage[leqno]{amsmath}
\usepackage{amscd,latexsym,amsthm,amsfonts,amssymb,amsxtra,mathrsfs,bm}
\usepackage{fullpage}
\usepackage[mathscr]{eucal}
\usepackage{pictexwd,dcpic}
\usepackage[normalem]{ulem}
\usepackage{hyperref}
\usepackage{enumerate}
\usepackage[all]{xy}
\usepackage{color}

\renewcommand\theequation{\thesection.\arabic{equation}}

\newtheorem{thm}{Theorem}[section]
\newtheorem{cor}[thm]{Corollary}
\newtheorem{lem}[thm]{Lemma}

\newtheorem{prop}[thm]{Proposition}

\newtheorem {ques/conj}[thm]{Question/Conjecture}

\newtheorem{defn}[thm]{Definition}
\newtheorem{rmk}[thm]{Remark}

\DeclareMathOperator{\Tr}{Tr}
\DeclareMathOperator{\vol}{vol}
\DeclareMathOperator{\Nil}{Nil}

\DeclareMathOperator{\Norm}{Norm}
\DeclareMathOperator{\ad}{ad}

\DeclareMathOperator{\Supp}{Supp}

\newcommand{\BC}{{\mathbb {C}}}

\newcommand{\BR}{{\mathbb {R}}}

\newcommand{\CA}{{\mathcal {A}}}
\newcommand{\CB}{{\mathcal {B}}}
\newcommand{\CC}{{\mathcal {C}}}

\newcommand{\CF}{{\mathcal {F}}}

\newcommand{\CH}{{\mathcal {H}}}

\newcommand{\CK}{{\mathcal {K}}}
\newcommand{\CL}{{\mathcal {L}}}

\newcommand{\CN}{{\mathcal {N}}}
\newcommand{\CO}{{\mathcal {O}}}
\newcommand{\CP}{{\mathcal {P}}}

\newcommand{\CU}{{\mathcal {U}}}

\newcommand{\CX}{{\mathcal {X}}}
\newcommand{\CY}{{\mathcal {Y}}}

\newcommand{\Fa}{{\mathfrak {a}}}
\newcommand{\Fb}{{\mathfrak {b}}}

\newcommand{\Fg}{{\mathfrak {g}}}
\newcommand{\Fh}{{\mathfrak {h}}}

\newcommand{\Fn}{{\mathfrak {n}}}

\newcommand{\Fp}{{\mathfrak {p}}}

\newcommand{\Ft}{{\mathfrak {t}}}
\newcommand{\Fu}{{\mathfrak {u}}}

\newcommand{\RR}{{\mathrm {R}}}

\newcommand{\Ad}{{\mathrm{Ad}}}
\newcommand{\Aut}{{\mathrm{Aut}}}

\newcommand{\cont}{{\mathrm{cont}}}

\newcommand{\Cent}{{\mathrm{Cent}}}

\newcommand{\disc}{{\mathrm{disc}}}

\newcommand{\elli}{{\mathrm{ell}}}
\newcommand{\End}{{\mathrm{End}}}

\newcommand{\GL}{{\mathrm{GL}}}

\newcommand{\geom}{{\mathrm{geom}}}

\newcommand{\Hom}{{\mathrm{Hom}}}

\newcommand{\Id}{{\mathrm{Id}}}

\newcommand{\Ker}{{\mathrm{Ker}}}

\newcommand{\Mat}{{\mathrm{Mat}}}

\newcommand{\reg}{{\mathrm{reg}}}
\newcommand{\rs}{{\mathrm{rs}}}
\newcommand{\Res}{{\mathrm{Res}}}

\newcommand{\SL}{{\mathrm{SL}}}

\newcommand{\SO}{{\mathrm{SO}}}

\newcommand{\Sp}{{\mathrm{Sp}}}
\newcommand{\Stab}{{\mathrm{Stab}}}

\newcommand{\spec}{{\mathrm{spec}}}
\newcommand{\scusp}{{\mathrm{scusp}}}
\newcommand{\tr}{{\mathrm{tr}}}

\newcommand{\ind}{{\mathrm{ind}}}

\newcommand{\back}{\backslash}

\newcommand{\ssi}{\mathrm{ss}}
\newcommand{\fg}{\mathfrak{g}}
\newcommand{\cO}{\mathcal{O}}
\newcommand{\cA}{\mathcal{A}}
\newcommand{\cC}{\mathcal{C}}

\newcommand{\tm}{\widetilde{m}}
\newcommand{\tM}{\widetilde{M}}
\newcommand{\tP}{\widetilde{P}}
\newcommand{\tQ}{\widetilde{Q}}
\newcommand{\tB}{\widetilde{B}}
\newcommand{\tL}{\widetilde{L}}
\newcommand{\tH}{\widetilde{H}}
\newcommand{\tR}{\widetilde{R}}
\newcommand{\tS}{\widetilde{S}}
\newcommand{\tT}{\widetilde{T}}
\newcommand{\tG}{\widetilde{G}}
\newcommand{\tg}{\widetilde{g}}
\newcommand{\tpi}{\widetilde{\pi}}
\newcommand{\tw}{\widetilde{w}}
\newcommand{\tsigma}{\widetilde{\sigma}}
\newcommand{\Temp}{\mathrm{Temp}}
\newcommand{\Specsupp}{\mathrm{Supp}_{\mathrm{spec}}}

\def\oF{\overline{F}}

\def\C{\mathbb{C}}

\def\diag{{\rm diag}}

\def\th{{\widetilde{h}}}
\def\tn{{\widetilde{n}}}

\def\tX{{\widetilde{X}}}

\numberwithin{equation}{subsection}
\newcounter{keepeqno}
\newenvironment{num}
{\setcounter{keepeqno}{\value{equation}}%
	\begin{list}{(\theequation)}{\usecounter{equation}}%
		\setcounter{equation}{\value{keepeqno}}}
	{\end{list}}

\begin{document}
	
	\title{A local twisted trace formula for Whittaker induction of coregular symmetric pairs: the geometric side}

	\author{Rapha\"el Beuzart-Plessis,\; Chen Wan}

	\date{\today}

	\maketitle
	
	\begin{abstract}
		In this paper, we prove the geometric expansion of a local twisted trace formula for the Whittaker induction of any symmetric pairs	that are coregular. This generalizes the local (twisted) trace formula for reductive groups proved by Arthur \cite{A91} and Waldspurger \cite{WalFTLtordue}. We also prove a formula for the regular germs of quasi-characters associated to strongly cuspidal functions in terms of certain weighted orbital integrals. As a consequence of our trace formula and the formula for regular germs of quasi-characters, we prove a simple local trace formula of those models for strongly cuspidal test functions which implies a multiplicity formula for these models. We also present various applications of our trace formula and multiplicity formula, including a necessary condition for a discrete L-packet  to contain a representation with a unitary Shalika model (resp. a Galois model for classical groups) in terms of the associated Langlands parameter, and we also compute the summation of the corresponding multiplicities for certain discrete L-packets.

	\end{abstract}
	
	\tableofcontents

	\section{Introduction}
	Let $F$ be a local non-Archimedean field of characteristic 0, $G$ be a reductive group defined over $F$, $H\subset G$ be a unimodular subgroup and $\xi: H(F)\to \BC^\times$ be a smooth unitary character. Choosing Haar measures on $G(F)$ and $H(F)$ induces an invariant measure on $H(F)\backslash G(F)$ and we let $L^2(H(F)\backslash G(F),\xi)$ be the space of $\varphi: G(F)\to \CC^\times$ that transform by left multiplication by $H(F)$ according to the character $\xi$ (i.e. $\varphi(hg)=\xi(h)\varphi(g)$ for $(h,g)\in H(F)\times G(F)$) and whose norm is square-integrable on $H(F)\backslash G(F)$. Then, the natural action of $G(F)$ on $L^2(H(F)\backslash G(F),\xi)$ by right translation is a unitary representation and for $f\in C_{c}^{\infty}(G(F))$, we define by integration an operator $R(f)$ on $L^2(H(F)\backslash G(F),\xi)$. This operator is associated to the following kernel function (for simplicity we assume the center of $G$ is trivial in the introduction)
	$$\displaystyle K_f(x,y)=\int_{H(F)} f(x^{-1}hy) \xi(h)dh,\;\; x,y\in G(F).$$
	Formally, the trace of the operator $R(f)$ should be given by the integral of $K_f(x,x)$ over $x\in H(F)\backslash G(F)$. However, neither of these two expressions are well-defined in general. The goal of this paper is to define some canonical regularizations of the integral of $K_f$ over the diagonal for certain triples $(G,H,\xi)$ (essentially associated to symmetric varieties that we name ``coregular'') and to express the resulting distribution on $G(F)$ as a sum (or integral) of contributions naturally generalizing the weighted orbital integrals of Arthur \cite{A91}. This can be considered as the geometric side of a local trace formula for the corresponding unitary representations $L^2(H(F)\backslash G(F),\xi)$. We plan to develop in a subsequent paper a general spectral expansion for those trace formulas.
	
	In the so-called {\em group-case}, corresponding to $G=H\times H$ with $H$ embedded diagonally in the product, we recover the geometric side of Arthur local trace formula \cite{A91}. We actually also consider an enhancement of the previous setting where we fix an extra automorphism $\theta$ of the triple $(G,H,\xi)$ and we formally try to compute the trace of the composition $R(f)\circ \theta$. This can be more naturally formulated using the notion of twisted spaces due to Labesse. In the group-case again, we recover the geometric side of the local twisted trace formula due to Waldspurger \cite{WalFTLtordue}.
	
	Although not implied by our main results, the work of Waldspurger \cite{WalGGPI} on the local Gan-Gross-Prasad conjecture, whose main innovation was the development of a certain simple local trace formula, has been a main source of inspiration and motivation for this paper.

	We also present few applications of our general trace formula. Namely, specializing our geometric expression to a matrix coefficient of a supercuspidal or square-integrable representation, we obtain explicit integral formulas for multiplicities of certain models which generalize our previous results for Galois models \cite{BeuGalP} and the Shalika model \cite{BW}. This can then be further applied to establish necessary conditions, in terms of the associated Langlands parameters, for the distinction of discrete $L$-packets with respect to a unitary Shalika model or a Galois model for classical groups and we moreover compute the corresponding multiplicities of such packets under an extra assumption. In the case of Galois models, this confirms some consequences of a general conjecture made by D. Prasad \cite{Pras}.

	\subsection{Main results}
	\subsubsection{Whittaker induction of coregular symmetric varieties}
	
	Let $\iota$ be an involutive automorphism of $G$ and assume that $(G^\iota)^0\subset H\subset G^\iota$ where $G^\iota$ stands for the subgroup of fixed points and $(G^\iota)^0$ for its neutral component. In this situation, the quotient variety $H\backslash G$ is sometimes called a symmetric variety. In this paper, we impose an important condition on the variety $H\backslash G$ that we decided to name coregularity:
	
	\begin{defn}
		Let $X=H\backslash G$  be a homogeneous $G$-variety with $H$ reductive. We say that $X$ is \textbf{coregular} if there exists an non-empty open subset $U\subset X\times X$ such that for every $x\in U$, the stabilizer $G_x\subset G$ of $x$ for the diagonal action contains regular elements.
	\end{defn}
	
	In Section \ref{section coregular varieties} we give various alternative characterizations of coregular homogeneous $G$-varieties (including the case where $H$ is not reductive). Technically, the most important for us is the following property (where $G_{rs}\subset G$ denotes the open locus of regular semisimple elements and $D^G$, $D^H$ stand for the usual Weyl discriminants):
	
	\begin{center}
		A homogeneous $G$-variety $X=H\backslash G$ is coregular if and only if $H\cap G_{rs}$ is nonempty and the function $h\in H(F)\cap G_{rs}(F)\mapsto \frac{D^H(h)^2}{D^G(h)}$ is locally bounded on $H(F)$.
	\end{center}
	
	Examples of coregular symmetric varieties are the group case (that is $X=H^{diag}\backslash H\times H$), Galois symmetric varieties (i.e. homogeneous varieties of the form $X=H\backslash \Res_{E/F}H_E$ where $E/F$ is a quadratic extension and $\Res_{E/F}$ denotes Weil's restriction of scalars) or $\Sp_{2n}\backslash \GL_{2n}$. However, many other natural examples of homogeneous varieties such as $O_n\backslash \GL_{n}$, $\GL_n\times \GL_n\backslash \GL_{2n}$ or $\SO_n^{\diag}\backslash (\SO_n\times \SO_{n+1})$ (the so-called Gross-Prasad variety, that is not symmetric but at least spherical) are not coregular.
	
	In this paper, we will actually consider a slightly more general setting, essentially including all triples $(G,H,\xi)$ that are in a suitable sense ``Whittaker induction'' of a coregular symmetric pair $(M,H_0)$. More precisely, the most general triples $(G,H,\xi)$ that we can consider are constructed as follows. There exists an involution $\iota$ of $G$ as well as a parabolic subgroup $P\subset G$ that is $\iota$-split (recall that it means that $\overline{P}:=\iota(P)$ is opposite to $P$) and a semi-direct product decomposition $H=H_0\ltimes N$ where:
	\begin{itemize}
		\item $N$ is the unipotent radical of $P$ and $H_0$ is a subgroup of the Levi factor $M:=P\cap \overline{P}$;
		
		\item We have $H_0=(M^\iota)^0$ and the symmetric variety $H_0\backslash M$ is coregular;
		
		\item The restriction of the character $\xi$ to $N(F)$ is non-degenerate i.e. its orbit under the adjoint action of $M(F)$ is open in the $F$-vector space of all smooth characters $N(F)\to \BC^\times$;
		
		\item In the case where $P\neq G$, $H_0$ is precisely the neutral component of the stabilizer of $\xi$ in $M$.
	\end{itemize}
	Following \cite[Sect. 2.6]{SV}, we say that the pair $(H\backslash G,\xi)$ is the {\it Whittaker induction} of the symmetric (and coregular) variety $H_0\backslash M$. One example of such Whittaker induction is given by the triple $(\GL_{2n},\GL_n^{\diag}\ltimes \Mat_n, \psi\circ \Tr)$, where $\psi: F\to \BC^\times$ is a nontrivial character, which is related to so-called {\it Shalika models} of representations of $\GL_{2n}(F)$. In this particular case, our result on geometric expansions contains the main results of our previous work \cite{BW}. There is also a variant of this example for unitary groups that will be described in more details below, related to what we call unitary Shalika models.
	
	\subsubsection{Truncation on symmetric varieties}
	
	Fix a triple $(G,H,\xi)$ as in the previous paragraph. Our starting point will be to truncate in a meaningful way the (usually non-convergent) integral
	$$\displaystyle I(f)=\int_{H(F)\backslash G(F)} K_f(x,x)dx.$$
	For this we introduce a sequence of truncation functions $(\kappa_Y)_Y$ indexed by points $Y$ in a certain affine space \footnote{For the definition of our truncation functions, we do not need to assume $(G,H)$ is coregular. It works for all the symmetric varieties.}.
	
	More precisely, we fix from now on a special maximal compact subgroup $K\subset G(F)$ that is in good position with respect to the Levi subgroup $M$ as well as a minimal $\iota$-split parabolic subgroup $P_0\subset P$\footnote{Here, by a {\em minimal $\iota$-split} parabolic subgroup we mean a parabolic subgroup that is $\iota$-split and minimal for this property.}. We assume for simplicity that there is only one $H_0(F)$-conjugacy class of minimal $\iota$-split parabolic subgroups in $P$. (Otherwise we need to replace $P_0$ by a set of representatives of those conjugacy classes, there are always finitely many, or, even better, we should replace $P_0$ in the discussion that follows by the union of finitely many $P_0(F)$-orbits in $(P_0\cap H \backslash P_0)(F)$). Let $\CA_{P_0,\iota}$ be the subspace of the real vector space
	$$\displaystyle \CA_{P_0}:=\Hom(X^*(P_0),\BR)$$
	on which $\iota$ acts by $-\Id$. Let $\CA^{P,+}_{P_0}\subset \CA_{P_0}$ be the usual Weyl chamber associated to the parabolic subgroup $P_0\cap M$ of $M$ and $\CA_{P_0,\iota}^{P,+}$ be its projection to $\CA_{P_0,\iota}$. We also let $A_0$ be the maximal central split torus in $M_0=P_0\cap \iota(P_0)$ (a Levi factor of $P_0$) which is {\em $\iota$-split} (in the sense that $\iota(a)=a^{-1}$ for every $a\in A_0$) and denote by $H_{P_0,\iota}:P_0(F)\to \CA_{P_0,\iota}$ the composition of the usual Harish-Chandra morphism $P_0(F)\to \CA_{P_0}$ with the projection $\CA_{P_0}\to \CA_{P_0,\iota}$. Then, by the weak Cartan decomposition of \cite{DS} and \cite{BO}, we can find a compact subset $\omega_{P_0}\subset P_0(F)$ such that setting
	$$\displaystyle \mathcal{S}(\omega_{P_0})=\{x=pak\mid p\in \omega_{P_0}, k\in K, a\in A_0(F), H_{P_0,\iota}(a)\in \CA_{P_0,\iota}^{P,+} \}$$
	we have the decomposition
	$$\displaystyle G(F)=H(F)\mathcal{S}(\omega_{P_0}).$$
	Note the formal resemblance with the existence of Siegel domains in a global setting. Let ${}^- \CA_{P_0}\subset \CA_{P_0}$ be the cone defined by the negative simple weights with respect to $P_0$ and ${}^- \CA_{P_0,\iota}$ be its image in $\CA_{P_0,\iota}$. Then, for any $Y\in \CA_{P_0,\iota}^+$ that is ``sufficiently positive'', we denote by $\kappa_Y$ the characteristic function of the image in $H(F)\backslash G(F)$ of the set $\mathcal{S}(\omega_{P_0},Y)$ defined by
	$$\displaystyle \mathcal{S}(\omega_{P_0},Y):=\{x\in \mathcal{S}(\omega_{P_0})\mid H_{P_0,\iota}(x)\in Y+{}^-\CA_{P_0,\iota} \}$$
	where we have denoted also by $H_{P_0,\iota}$ the unique extension of $H_{P_0,\iota}$ to $G(F)$ that is $K$-invariant on the right.
	
	Although the family of truncation functions $(\kappa_Y)_Y$ a priori depends on the auxiliary choice of the compact subset $\omega_{P_0}$, it can be shown that it doesn't asymptotically in the following precise sense. For any pair of compact subsets $\omega_{P_0},\omega'_{P_0}\subset P_0(F)$ such that $G(F)=H(F)\mathcal{S}(\omega_{P_0})=H(F)\mathcal{S}(\omega_{P_0}')$, we can define as above two families of truncation functions $(\kappa_Y)_Y$ and $(\kappa'_Y)_Y$. Then, there exists $Y_+\in \CA_{P_0,\iota}$ such that $\kappa_Y=\kappa'_Y$ for every $Y\in Y_++\CA_{P_0,\iota}^+$. In particular, it makes sense to study the asymptotic behavior of the expression
	$$\displaystyle I_Y(f):=\int_{H(F)\backslash G(F)} K_f(x,x) \kappa_Y(x) dx$$
	when $Y\xrightarrow{P_0} \infty$ (where the latter notation means asymptotic along the filter generated by translates $Y_++\CA_{P_0,\iota}^+$ of the positive Weyl chamber). Moreover, the functions $\kappa_Y$ are so defined that the integrand in the above expression is compactly supported (see Lemma \ref{lem abs conv}).
	
	Finally, we can also suppress the dependence of our truncation process on the choice of $P_0$ (but not on that of $K$) as follows: for any other choice of a minimal $\iota$-split parabolic subgroup $P_0'\subset P$, there exists a natural affine isomorphism
	\begin{equation}\label{eq affine isom}
		\displaystyle \iota_{P_0,P_0',K}: \CA_{P_0,\iota}\simeq \CA_{P_0',\iota},
	\end{equation}
	such that as $Y\in \CA_{P_0,\iota}\xrightarrow{P_0} \infty$ we eventually have $\kappa_Y=\kappa_{Y'}$ where $Y'=\iota_{P_0,P_0',K}(Y)$. We emphasize that $\iota_{P_0,P_0',K}$ is not the most obvious isomorphism $\CA_{P_0,\iota}\simeq \CA_{P_0',\iota}$, namely the one induced by conjugation by an element $p\in P(F)$ such that $pP_0p^{-1}=P_0'$, which is not only affine but linear. Indeed, in general the map \eqref{eq affine isom} does not preserve the origins; a fact that is related to the existence of more than one $H(F)\cap K$-conjugacy class of minimal $\iota$-split parabolic subgroups $P_0\subset P$.
	
	Therefore, we can as well think of $Y$ as living in the inverse limit
	$$\displaystyle \CA_{X,K}=\varprojlim_{P_0} \CA_{P_0,\iota}$$
	where $P_0$ runs over the set of minimal $\iota$-split parabolic subgroups of $P$ and the transition maps are given by the affine isomorphisms \eqref{eq affine isom}. This is the point of view that we will adopt in the body of the paper.

	\subsubsection{The geometric expansion of a general local trace formula.}
	
	
	Let $\Gamma(H_0)$ (resp. $\Gamma_{\elli}(H_0)$) be the set of regular semisimple (resp. regular elliptic) conjugacy classes in $H_0(F)$. These two sets can be naturally equipped with measures, see Sections \ref{sec geometric expansion} and \ref{sec simple LTF} for details.
	
	For $t\in \Gamma(H_0)$, that we identify with a representative in $H(F)$, we denote by $H_{t}$, $G_t$, $M_{t}$, $N_t$ and $B_t=M_t N_t$ the neutral components of the centralizers of $t$ in $H$, $G$ $M$, $N$ and $P$ respectively. Then, for $t$ in general position $B_t$ is a Borel subgroup of $G_t$ and $\xi\mid_{N_t(F)}$ is a non-degenerate character on its unipotent radical (see Lemma \ref{lem generic character}, here we need to use the coregular assumption). For $f\in C_c^\infty(G(F))$ and $Y\in \CA_{X,K}$, we define the following expression
	$$\displaystyle J_Y(f)=\int_{\Gamma(H_0)}D^{H}(t)\xi(t) J_Y(t,f) dt$$
	where $J_Y(t,f)$ denotes some kind of ``weighted orbital integral''. More precisely, $J_Y(t,.)$ is a distribution of the form
	$$\displaystyle J_Y(t,f)=\int_{\CO_t} f(g) v_{\xi,\iota, Y}(g) dg$$
	where $\CO_t$ denotes the union of the (finitely many) {\em regular} $G(F)$-conjugacy classes with semisimple part $t$ and the function $g\mapsto v_{\xi,\iota, Y}(g)$ is a certain weight function. When $\xi=1$ (so that $P=G$ and, by the coregular assumption, $t$ is already regular in $G$ at least when it is in general position), this weight is very similar to the one appearing in the definition of Arthur's weighted orbital integrals as $v_{\xi,\iota, Y}(g^{-1}tg)$ is given by the volume of the convex hull of a certain family $(-H_{\overline{Q},\iota}(g)+Y_Q)_Q$ where $Q$ runs over the minimal $\iota$-split parabolic subgroups of $G$ containing $t$, $H_{\overline{Q},\iota}: G(F)\to \CA_{L,\iota}$ denotes the usual Harish-Chandra map for the parabolic subgroup $\overline{Q}=\iota(Q)=LN_{\overline{Q}}$ (chosen to be $K$-invariant on the right) composed with the projection $\CA_{\overline{Q}}\to \CA_{L,\iota}$ to the $\iota$-antifixed points in $\CA_L$ and $Y\mapsto Y_Q\in \CA_{Q,\iota}=\CA_{L,\iota}$ is the composition of the canonical isomorphism $\CA_{X,K}\simeq \CA_{P_0,\iota}$ with the natural projection $\CA_{P_0,\iota}\to \CA_{Q,\iota}$ for any minimal $\iota$-split parabolic subgroup $P_0\subset Q$ (it can be shown that the composition doesn't depend on $P_0$). In general, the precise definition looks like
	$$\displaystyle J_{Y}(t,f)=\int_{B_t(F)\back G(F)} \int_{N_t(F)}f(x^{-1}tu x)v_{B_t,\xi,\iota,Y}(x,u) dudx$$
	where we refer the reader to Section \ref{sec geometric expansion} for the definition of the weight $v_{B_t,\xi,\iota,Y}(x,u)$ when $\xi\neq 1$. Another important point is that, after Harish-Chandra, it is known that near singular point the typical order of growth of (weighted) orbital integrals is as the inverse of the square root of the Weyl discriminant. Therefore, our assumption on coregularity of the pair $(G,H)$ is what guarantees the absolute convergence of the expression defining $J_Y(f)$ above. Then, the aforementioned geometric expansion of the local trace formula for $(G,H,\xi)$ is contained in the following theorem. 
	
	\begin{thm}\label{thm2 introduction}
		Let $0<\epsilon<1$ and fix $f\in C_{c}^{\infty}(G(F))$. Then, for any $k>0$, we have
		$$\displaystyle \left\lvert I_{Y}(f)-J_{Y}(f)\right\rvert\ll N(Y)^{-k}$$
		for every $Y\in \cA_{X,K}$ with $d(Y)>\epsilon N(Y)$. Moreover, the function $Y\in \CA_{X,K}\mapsto J_Y(f)$ is a polynomial-exponential function in a suitable sense (see Section \ref{orthogonal sets}) and if the variety $X=H\backslash G$ is tempered (see Section \ref{S tempered varieties}), then the same statement holds for functions $f$ in the Harish-Chandra Schwartz space $\CC(G(F))$.
	\end{thm}
	
	In the above statement, $N(Y)$ stands for any norm on the affine space $\CA_{X,K}$ whereas, fixing a minimal $\iota$-split parabolic subgroup $P_0\subset P$ for convenience, the {\em depth} $d(Y)$ of $Y$ is defined by
	$$\displaystyle d(Y)=\min_{\alpha\in \Delta_0} \langle \alpha, Y-Y_0\rangle$$
	where $\Delta_0$ stands for the set of simple roots with respect to $P_0$. Therefore, in some loose sense, the above theorem describes the asymptotic behavior of $I_Y(f)$ as $Y$ goes to infinity in the direction of the positive Weyl chamber and ``sufficiently far from the walls''.
	
	As already mentioned, in the main body of the paper we actually prove a more general theorem of the above form for suitable {\em twisted} triples $(\tG,\tH,\xi)$. In particular, in the group case (i.e. when $H$ is diagonally embedded in $G=H\times H$) this recovers the geometric side of the twisted local trace formula \cite{WalFTLtordue}.

	\subsubsection{The case of strongly cuspidal functions and integral formulas for multiplicities}
	
	Most applications of our trace formula comes from a simple version obtained by specializing it to the case of {\em strongly cuspidal} test functions. More precisely, we recall following \cite{WalGGPI} that a function $f\in C_c^\infty(G(F))$ is said to be {\em strongly cuspidal} if for every proper parabolic subgroup $Q=LV\subset G$ we have
	$$\displaystyle \int_{V(F)} f(lu)du=0,\; \mbox{ for every } l\in L(F).$$
	It is then shown in {\it loc. cit.} that the regular semisimple weighted orbital integrals (in the sense of Arthur) of a strongly cuspidal function $f$ don't depend on any choice (except that of a Haar measure on $G(F)$) and that, correctly normalized by certain signs, they define a function
	$$\displaystyle \Theta_f: G_{\rs}(F)\to \BC$$
	which is $G(F)$-invariant by conjugation and a {\em quasi-character} in the following sense: for every semisimple element $x\in G(F)$, there exists an expansion
	$$\displaystyle \Theta_f(x\exp(X))=\sum_{\CO\in \Nil(\mathfrak{g}_x^*)} c_{f,\CO}(x) \widehat{j}(\CO,X),\;\; X\in \omega\cap \mathfrak{g}_{x,\rs}(F),$$
	where:
	\begin{itemize}
		\item $\omega\subset \mathfrak{g}_{x}(F)$ is a sufficiently small neighborhood of $0$ in the Lie algebra of $G_x(F)$;
		
		\item $\Nil(\mathfrak{g}_x^*)$ denotes the (finite) set of nilpotent coadjoint orbits in the dual $\mathfrak{g}_x(F)^*$ of $\mathfrak{g}_x(F)$;
		
		\item for $\CO\in \Nil(\mathfrak{g}_x^*)$, $\widehat{j}(\CO,.)$ stands for the unique smooth function on $\mathfrak{g}_{x,\rs}(F)$ that is locally integrable on $\mathfrak{g}_x(F)$ and represents the Fourier transform of the orbital integral over $\CO$ i.e. for every $\varphi\in C_c^\infty(\Fg_x(F))$ we have
		$$\displaystyle \int_{\mathfrak{g}_x(F)} \varphi(X) \widehat{j}(\CO,X)dX=\int_{\CO} \widehat{\varphi}(Y)dY$$
		where $dX$ is a Haar measure on $\Fg_x(F)$, $\widehat{\varphi}(Y)=\int_{\mathfrak{g}_x(F)} \varphi(X) \psi(\langle X,Y\rangle) dX$, $Y\in \Fg_x(F)^*$, denotes the Fourier transform of $\varphi$ (which depend in the auxiliary choice of a non-trivial additive character $\psi: F\to \BC^\times$) and $dY$ is the canonical Kirillov-Kostant measure on the coadjoint orbit $\CO$.
	\end{itemize}

	For $t\in H_{0,\rs}(F)$ in general position, the restriction $\xi_t:=\xi\mid_{N_{t}(F)}$ is a generic character of $N_{t}(F)$. We let $\CO_t\in \Nil(\Fg_{t}^{\ast})$ be the orbit associated to $\xi_{t}$\footnote{More precisely, $\CO_t$ is the unique nilpotent coadjoint orbit in $\Fg_t(F)^{\ast}$ containing an element $Y$ such that $\xi(\exp(X))=\psi(\langle Y,X\rangle)$ for every $X\in \Fn_t(F)$, the Lie algebra of $N_t(F)$.}. Then, for every strongly cuspidal test function $f\in C_c^\infty(G(F))$ we set
	$$\displaystyle I_{\geom}(f):=\int_{\Gamma_{\elli}(H_0)} D^{H}(t) c_{f,\CO_t}(t)\xi(t)dt.$$
	
	\begin{thm}
		Let $f\in C_{c}^{\infty}(G(F))$ be a strongly cuspidal function. Then,
		\begin{enumerate}
			\item We have
			$$\displaystyle \lim\limits_{Y\xrightarrow{P_0} \infty} I_Y(f)=I_{\geom}(f),$$
			in particular the limit exists.
			\item If moreover $f$ is a matrix coefficient of a supercuspidal representation $\pi$ of $G(F)$ and the dimension $m_{H,\xi}(\pi^\vee)$ of the space $\Hom_H(\pi^\vee,\xi)$ of $(H(F),\xi)$-equivariant linear forms on (the space of) the contragredient representation $\pi^\vee$ is finite, then the integral defining $I(f)$ is already convergent and we have
			$$\displaystyle I(f)=\frac{f(1)}{d(\pi)} m_{H,\xi}(\pi^{\vee})$$
			where $d(\pi)$ stands for the formal degree of $\pi$.
		\end{enumerate}
		Furthermore, if the pair $(G,H)$ is tempered then the same holds for strongly cuspidal test functions $f\in \CC(G(F))$ and matrix coefficients of square-integrable representations $\pi$ respectively.
	\end{thm}	
	
	As a corollary of the above theorem we can also obtain general integral formulas for the multiplicities $m_{H,\xi}(\pi)$. More precisely, for $\pi$ an irreducible representation of $G(F)$, it is known that the Harish-Chandra character $\Theta_\pi$ is also a quasi-character in the above sense. Therefore, we can define an expression $m_{\geom,H,\xi}(\pi)$ similar to $I_{\geom}(f)$ by formally replacing $\Theta_f$ by $\Theta_\pi$. Then, we have the following. (see Theorem \ref{thm multiplicity formula}).
	
	\begin{cor}
		\begin{enumerate}
			\item Assume that $\pi$ is supercuspidal and the multiplicity $m_{H,\xi^{-1}}(\pi)$ is finite. Then, we have
			\begin{equation}\label{eq multiplicities}
				\displaystyle m_{H,\xi^{-1}}(\pi)=m_{\geom,H,\xi}(\pi).
			\end{equation}
			
			\item If the pair $(G,H)$ is tempered, $\pi$ is square-integrable and the multiplicity $m_{H,\xi}(\pi)$ is finite, then the equality \eqref{eq multiplicities} also holds.
		\end{enumerate}
	\end{cor}
	
	In the case of Galois models or the Shalika model, the above corollary recovers one of the main result in \cite{BeuGalP} and \cite{BW} respectively. Actually for Galois models associated to classical groups, we can also deduce new results from the analog of the above corollary in certain twisted situations as explained in more details below.




	\subsubsection{An integral formula for regular germs of quasi-characters}
	
	One important technical result that we prove along the way to Theorem \ref{thm2 introduction}, and that may be of independent interest, is a certain explicit formula for some singular weighted orbital integrals of strongly cuspidal functions. More precisely, we are able to write down a set of measures on regular (but not necessarily semi-simple) conjugacy classes representing the distributions $f\mapsto c_{f,\CO}(x)$ for $x\in G(F)$ semisimple and $\CO\in \Nil(\Fg_x^{\ast})$ a regular nilpotent coadjoint orbit.
	
	More precisely, let us fix a Borel subgroup $B_x$ of $G_x$ with a Levi decomposition $B_x=T_xN_x$ as well as a generic character $\xi_x$ of $N_x(F)$, and we let $\CO_x\in \Nil(\Fg_{x}^{\ast})$ be the corresponding regular nilpotent coadjoint orbit (every regular nilpotent coadjoint orbit arises in this way). In Section 4 we will define a weighted function $v_{B_x,\xi_x}(g,u)$ for $g\in G(F)$ and $u\in N_{x}(F)$ {\em regular}. The next theorem expresses the regular germ of the quasi-character $\Theta_f$ in terms of certain weighted orbital integral (we refer the reader to Section 2-4 for various notation). It will be proved in Section 4.
	
	\begin{thm}
		For every strongly cuspidal function $f\in \mathcal{C}(G(F))$, we have
		\begin{equation*}
			\displaystyle c_{f,-\CO_x}(x)=\int_{B_x(F)\backslash G(F)} \int_{N_x(F)} f(g^{-1}x u g) v_{B_x,\xi_x}(u,g) du dg.
		\end{equation*}
	\end{thm}

	\subsubsection{The Galois model for classical groups}
	Let $E/F$ be a quadratic extension, $H$ be a reductive group defined over $F$, $\chi$ be a character of $H(F)$ and $G=Res_{E/F}H_E$. The model $(G,H,\chi)$ is the so-called Galois model. In \cite{Pras}, Prasad made a general conjectural regarding the multiplicity of Galois model. In this paper, we will study the case when $H$ is a classical group.
	
	Let $H$ be a quasi-split special orthogonal group or a symplectic group and $G=Res_{E/F}H_E$. If $H$ is the even special orthogonal group, let $H_0$ be a quasi-split special orthogonal group that is not a pure inner form of $H$ and such that $G=Res_{E/F}H_E=Res_{E/F}H_{0,E}$ (i.e. the determinants of the quadratic forms defining $H$ and $H_0$ belong to the same square class in $E^\times/(E^\times)^2$ but belong to different square classes in $F^\times/(F^\times)^2$). If $H=\Sp_{2n}$ or $\SO_{2n}$, let $\chi$ be the trivial character on $H$ (and $H_0$ if $H=\SO_{2n}$). If $H=\SO_{2n+1}$, let $\chi\in \{1,\eta_n\}$  where $\eta_n$ is the composition of the Spin norm character of $\SO_{2n+1}$ with the quadratic character $\eta_{E/F}$. 
	
	Our first result is a necessary condition for a discrete L-packet to be distinguished.
	
	\begin{thm}\label{thm Galois distinguished introduction}
		Let $H=\Sp_{2n},\SO_{2n}$ or $\SO_{2n+1}$, $G=Res_{E/F}H$, $\chi=1$ if $H=\Sp_{2n}$ or $\SO_{2n}$, and $\chi\in \{1,\eta_n\}$ if $H=\SO_{2n+1}$. Let $\Pi_\phi(G)$ be a discrete L-packet of $G(F)$ and $\Pi_\phi(G')$ be the endoscopic transfer of the L-packet to the general linear group $G'=\GL_a(E)$ (here $a=2n$ if $H=\SO_{2n}$ or $\SO_{2n+1}$ and $a=2n+1$ if $H=\Sp_{2n}$). Then the packet $\Pi_\phi(G)$ is distinguished (i.e. $m(\pi,\chi)\neq 0$ for some $\pi\in \Pi_\phi(G)$) only if $\Pi_\phi(G')$ is $(\GL_a(F),\chi')$-distinguished. Here $\chi'=1$ if $\chi=1$ and $\chi'=\eta_n':=\eta_{E/F}\circ \det$ if $\chi=\eta_n$.
	\end{thm}
	
	Our second result is to compute the summation of the multiplicities over certain discrete L-packets. Assume that $\Pi_{\phi}(G')$ is $(\GL_a(F),\chi')$-distinguished. By Theorem 4.2 of \cite{Mat}, $\Pi_\phi(G')$ is of the form 
	$$\Pi_{\phi}(G')=(\tau_1\times \cdots \times \tau_l)\times (\sigma_1\times \bar{\sigma}_{1})\times \cdots \times (\sigma_m\times \bar{\sigma}_{m})$$
	where 
	\begin{itemize}
		\item $\tau_i$ is a discrete series of $\GL_{a_i}(E)$ that is conjugate self-dual. Moreover, if $(H,\chi)=(\SO_{2n+1},\eta_n)$, $\tau_i$ is self-dual of symplectic type; otherwise, $\tau_i$ is self-dual of orthogonal type.
		\item $\sigma_j$ is a discrete series of $\GL_{b_i}(E)$ that is NOT conjugate self-dual. Moreover, if $(H,\chi)=(\SO_{2n+1},\eta_n)$, $\sigma_j$ is self-dual of symplectic type; otherwise, $\sigma_j$ is self-dual of orthogonal type.
		\item $\tau_i,\sigma_j$ are all distinct.
		\item $\sum_{i=1}^{l} a_i +2\sum_{j=1}^{m}b_j=a$.
	\end{itemize}
	We will consider the special case when $m=0$. The general case will be consider in our future paper. When $m=0$, $\Pi_\phi(G')$ appears discretely in the $L^2$ space of the Galois model $(\GL_a(E),\GL_a(F),\chi')$.
	
	\begin{thm}\label{thm Galois summation introduction}
		With the notation above, if $H$ is the symplectic group or the odd special orthogonal group, we have
		$$\sum_{\pi\in \Pi_\phi(G)} m(\pi,\chi)=2^{l-1}.$$
		If $H$ is the even special orthogonal group, we let $H_0$ be another even special orthogonal group as above. We use $m_0(\pi,\chi)$ to denote the multiplicity for the model $(G,H_0,\chi)$. Then we have
		$$\sum_{\pi\in \Pi_\phi(G)} m(\pi,\chi)+m_0(\pi,\chi)=2^{l-1}.$$
	\end{thm}
	
	\begin{rmk}
		By Theorem 1 of \cite{BeuGalP}, the above two theorems also hold if we replace $H$ (and $H_0$ if we are in the even orthogonal group case) by the non quasi-split classical group.
	\end{rmk}

	\subsubsection{The unitary Shalika model}
	Let $Z$ be a $E$-vector space of finite dimension $n\geqslant 1$. Let $Z^{*,c}$ be the conjugate-dual of $Z$ that is the space of $c$-linear forms on $Z$ (a similar notation will be applied later to other vector spaces). Set $V=Z\oplus Z^{*,c}$ and we equip with the nondegenerate Hermitian form
	$$\displaystyle h(v+v^*,w+w^*)=\langle v,w^*\rangle +\langle w,w^*\rangle^c,\;\;\; (v,v^*), (w,w^*)\in Z\oplus Z^{*,c}.$$
	Here $\langle .,.\rangle$ stands for the canonical pairing between $Z$ and $Z^{*,c}$. Let $G=U(V,h)$ be the unitary group associated to this Hermitian form. We define two maximal parabolic subgroups $Q$ and $\overline{Q}$ of $G$ as the stabilizers of the maximal isotropic subspaces $Z$ and $Z^{*,c}$ respectively. Then, $L=Q\cap \overline{Q}$ is a Levi component of $Q$ and restriction to $Z$ induces an isomorphism
	\begin{equation}
		\displaystyle L\simeq \Res_{E/F}GL(Z).
	\end{equation}
	Let $N$ be the unipotent radical of $Q$. Thus $Q=LN$ and restriction to $Z^{*,c}$ induces an isomorphism
	\begin{equation}
		\displaystyle N\simeq \left\{X\in \Hom(Z^{c,*},Z)\mid {}^T X^c=-X \right\}
	\end{equation}
	where ${}^T X^c$ denotes the transpose conjugate of $X$ (seen as a linear endomorphism $Z\to Z^{*,c}$ through the canonical identification $(Z^{*,c})^{*,c}=Z$). We will actually identify the right hand side above with the Lie algebra $\mathfrak{n}$ of $N$ in a way such that the above isomorphism becomes the exponential map.
	
	We henceforth choose two isomorphisms $W_+,W_-:Z\to Z^{*,c}$ satisfying ${}^T W_{\pm}^c=-W_{\pm}$ and such that the corresponding antihermitian forms on $Z$ are not equivalent (there are actually only two equivalence classes of antihermitian forms on $Z$). For $\epsilon\in \{ \pm\}$, we let $H_{0,\epsilon}\subset L\simeq \Res_{E/F}GL(Z)$ be the unitary group associated to $W_\epsilon$, that is the stabilizer of $W_\epsilon$ for the obvious action. Then, $H_{0,\epsilon}(F)$ coincides with the stabilizer in $L(F)$ of the character
	$$\displaystyle \xi_\epsilon:N(F)\to \mathbb{C}^\times,$$
	$$\displaystyle \exp(X)\mapsto \psi(\Tr(W_\epsilon X))\;\; (X\in \mathfrak{n}(F)).$$
	We will henceforth assume, as we may, that $W_{\pm}$ have been chosen so that $H_{0,+}$ is quasi-split.
	
	Set $H_\epsilon=H_{0,\epsilon}\ltimes N$. We extend $\xi_\epsilon$ to a character of $H_\epsilon(F)$ trivial on $H_{0,\epsilon}(F)$. We also fix a character $\chi$ of $E^1=\ker(N_{E/F})$ that we will consider as a character of $H_{0,\epsilon}(F)$ through composition with the determinant $\det: H_{0,\epsilon}(F)\to E^1$ . The model $(G,H_{\epsilon},\chi\otimes \xi_{\epsilon})$ is an analogue of the Shalika model for unitary groups and we will call it the unitary Shalika model. For a smooth irreducible representation $\pi$ of $G(F)$, we define the multiplicity
	$$\displaystyle m_\epsilon(\pi,\chi):=\dim(\Hom_{H_\epsilon(F)} (\pi,\chi\otimes \xi_\epsilon)).$$

	Our first result for the unitary Shalika model is that the multiplicity for the two models are equal to each other for all stable discrete series.
	
	\begin{thm}\label{thm multiplicity constant introduction}
		\begin{enumerate}
			\item Let $\pi$ be a finite length discrete series of $G(F)$ with central character $\chi^n$. If $\Theta_\pi$ is a stable distribution, then $m_+(\pi,\chi)=m_-(\pi,\chi)$.
			\item Let $\Pi_\phi(G)$ be a discrete L-packet of $G(F)$ with central character $\chi^n$. Then we have
			$$\sum_{\pi\in \Pi_\phi(G)} m_+(\pi,\chi)=\sum_{\pi\in \Pi_\phi(G)} m_-(\pi,\chi).$$
		\end{enumerate}
	\end{thm}
	
	Our second result for the unitary Shalika model is a necessary condition for a discrete L-packet to be distinguished. The character $\chi$ of $E^1$ induces a character $\chi'$ of $E^\times$ by $\chi'(x)=\chi(x/x^c)$. Let $\Pi_\phi(G)$ be a discrete L-packet of $G(F)$ and let $\Pi_\phi(G')$ be its base change to $G'(F)=\GL_{2n}(E)$. Then $\Pi_\phi(G')$ is an irreducible tempered representation. Let $H'(F)=\{\begin{pmatrix}h&0\\0&h \end{pmatrix}\begin{pmatrix}I_n&X\\0&I_n \end{pmatrix}|\; h\in \GL_n(E),X\in Mat_{n\times n}(E)\}$ be the Shalika subgroup and we define the character $\chi'\otimes \xi'$ on it to be
	$$\chi'\otimes \xi'(\begin{pmatrix}h&0\\0&h \end{pmatrix}\begin{pmatrix}I_n&X\\0&I_n \end{pmatrix})=\chi'(\det(h))\psi(\tr_{E/F}(\tr(X))).$$
	
	\begin{thm}\label{thm Shalika distinguished introduction}
		With the notation above, the packet $\Pi_\phi(G)$ is $(H_+,\chi\otimes \xi_+)$-distinguished (i.e. $m_+(\pi,\chi)\neq 0$ for some $\pi\in \Pi_\phi(G)$) only if $\Pi_\phi(G')$ is distinguished by the Shalika model $(H',\chi'\otimes \xi')$. 
	\end{thm}
	
	\begin{rmk}
		By Theorem \ref{thm multiplicity constant introduction}, we know that the packet $\Pi_\phi(G)$ is $(H_+,\chi\otimes \xi_+)$-distinguished if and only if it is $(H_-,\chi\otimes \xi_-)$-distinguished.
	\end{rmk}

	Our next result for the unitary Shalika model is to compute the summation of the multiplicities over certain discrete L-packets. Assume that $\Pi_\phi(G')$ is distinguished by the Shalika model $(H',\chi'\otimes \xi')$. By Corollary 1.1 of \cite{M}, $\Pi_\phi(G')$ is of the form ($\chi''$ is a character of $E^{\times}$ with $\chi'=(\chi'')^2$)
	$$\Pi_{\phi}(G')\otimes (\chi''\circ \det)^{-1}=(\tau_1\times \cdots \times \tau_l)\times (\sigma_1\times \sigma_{1}^{\vee})\times \cdots \times (\sigma_m\times \sigma_{m}^{\vee})$$
	where 
	\begin{itemize}
		\item $\tau_i$ is a discrete series of $\GL_{2a_i}(E)$ that is conjugate self-dual, self-dual and of symplectic type. In particular, $a_i$ is even.
		\item $\sigma_j$ is a discrete series of $\GL_{b_i}(E)$ that is conjugate self-dual, but NOT self-dual.
		\item $\tau_i,\sigma_j$ are all distinct.
		\item $\sum_{i=1}^{l} a_i +2\sum_{j=1}^{m}b_j=2n$.
	\end{itemize}
	We will consider the special case when $m=0$. The general case will be consider in our future paper. When $m=0$, $\Pi_\phi(G')$ appears discretely in the $L^2$-space of the Shalika model. 
	
	\begin{thm}\label{thm Shalika summation introduction}
		With the notation above, we have 
		$$\sum_{\pi\in \Pi_\phi(G)} m_+(\pi,\chi)=\sum_{\pi\in \Pi_\phi(G)} m_-(\pi,\chi)=2^{l-1}.$$
	\end{thm}

	The idea to prove our main theorems for the unitary Shalika model (resp. the Galois model for classical groups) is by comparing the simple trace formula of the unitary Shalika model (resp. the Galois model for classical groups) with the twisted simple trace formula for the Shalika model (resp. Galois model for general linear groups), we refer the reader to Section 8 (resp. Section 9) for details.

	In our next paper, we will prove the spectral side of the trace formula in the general case and we will use it to compute the multiplicity of all the discrete series for the Galois model for classical groups and the unitary Shalika model.

	\subsection{Organization of the paper}
	In Section 2, we introduce basic notations and conventions of this paper. This include some extended discussions of twisted weighted orbital integrals, germ expansions and the twisted local trace formula for strongly cuspidal functions. 
	
	In Section 3, we introduce the notion of coregular varieties and we will have an extended discussion of symmetric varieties.
	
	In Section 4, we will introduce certain ($\iota$-)weighted functions associated to singular semisimple elements and we will prove an integral formula of the regular germs of quasi-characters. We will also prove a descent formula for the $\iota$-weighted functions which will be used in later section.
	
	We prove a special case of the spectral expansion of the trace formula in Section 5 and in Section 6 we will prove the geometric expansion.
	
	In Section 7 we will discuss our first two applications of the trace formula, namely a simple trace formula for strongly cuspidal functions and a multiplicity formula.
	
	In Section 8 and 9 we will discuss another two applications of the trace formula. In Section 8 we will prove our theorems for the unitary Shalika models and in Section 9 we will prove our theorems for the Galois models of classical groups.
	
	In Appendix A we will prove some results regarding finitely generated convex sets and in Appendix B we will prove the Howe's conjecture for twisted weighted orbital integrals. The results in the two appendices will be used in Section 4 in our proof of the integral formula of the regular germs of quasi-characters.

	\subsection{Acknowledgement}
	We thank Rui Chen for pointing out that the base change of a character is always a square of another character. We thank Spencer Leslie and Yiannis Sakellaridis for a remark on the notion of coregular spherical varieties (Remark \ref{rmk coregular}). Finally, the first author is grateful to Yiannis Sakellaridis for many very inspiring discussions over the years on local trace formulas for spherical varieties that have partly inspired this work.
	
	The work of first author was funded by the European Union ERC Consolidator Grant, RELANTRA, project number 101044930. Views and opinions expressed are however those of the authors only and do not necessarily reflect those of the European Union or the European Research Council. Neither the European Union nor the granting authority can be held responsible for them. The work of the second author is partially supported by the NSF grant DMS-2000192 and DMS-2103720.

	\section{Groups and twisted spaces}\label{preliminaries}
	
	Throughout the paper, $F$ will be a non-Archimedean local field of characteristic zero with normalized absolute value $\lvert .\rvert$. Unless otherwise specified, all the groups and varieties that we will consider are implicitely supposed defined over $F$. We fix a non-trivial additive character $\psi:F\to \mathbb{C}^\times$ and, whenever convenient, we will also fix an algebraic closure $\overline{F}$ of $F$.
	
	For $V$ a real vector space, we write $V^*$ for its dual and we denote by $V_{\mathbb{C}}=V\otimes_{\mathbb{R}} \mathbb{C}$ its complexification. Moreover, $iV\subset V_{\mathbb{C}}$ will stand for the real subspace of purely imaginary vectors.
	
	For two complex valued functions $f$ and $g$ on a set $X$ with $g$ taking values in $\BR_{\geq 0}$, we write that
	$$
	f(x)\ll g(x),\;\; x\in X,
	$$
	and say that $f$ is essentially bounded by $g$, if there exists a constant $c>0$ such that for all $x\in X$, we have
	$$
	| f(x)| \leq cg(x).
	$$
	We say $f$ and $g$ are equivalent, which is denoted by
	$$f(x)\sim g(x)$$
	if $f$ is essentially bounded by $g$ and $g$ is essentially bounded by $f$.
	
	\subsection{Groups}
	
	In this section, we fix some notation relative to the datum of a linear algebraic group $G$ defined over $F$. First, we write $rk(G)$ for the (absolute) rank of $G$, that is $\dim(T)$ for any maximal torus $T\subset G$, $U_G$ for the unipotent radical of $G$ and we will denote Lie algebras by the corresponding gothic letter such as $\mathfrak{g}$ for $G$. The adjoint representations $G\to \GL(\fg)$ and $\fg\to \End(\fg)$ will be denoted by $g\mapsto \Ad_g$ and $X\mapsto \ad_X$ respectively. We also write $\mathfrak{g}^*$ for the dual of $\mathfrak{g}$. The exponential map, which is well-defined on some neighborhood of $0$ in $\mathfrak{g}(F)$ to $G(F)$, will be denoted by $\exp$. For every $g\in G$, we write $\Ad_g$ both for the adjoint action of $g$ on $G$ and $\mathfrak{g}$. We also denote by $X^*(G)$ the group of algebraic characters $G\to \mathbb{G}_m$ defined over $F$ and by $A_G$ the maximal central split torus in $G/U_G$. We set
	$$\displaystyle \cA_G=\Hom(X^*(G),\mathbb{R})=\Hom(X^*(A_G),\mathbb{R})$$
	and we let as usual $H_G:G(F)\to \cA_G$ be the homomorphism defined by $\langle H_G(g),\chi\rangle=\log \lvert \chi(g)\rvert$ for every $(g,\chi)\in G(F)\times X^*(G)$. We denote by $G_{\rs}$ and $\mathfrak{g}_{\rs}$ (resp. $G_{\reg}$ and $\fg_{\reg}$) the open subsets of regular semi-simple elements (resp. regular elements) in $G$ and $\mathfrak{g}$ respectively. The notation $\delta_G$ will stand for the modular character of $G(F)$ that is the character $\delta_G:G(F)\to \mathbb{R}_+^*$ defined by $\delta_G(g)=\lvert \det \Ad_g\mid_{\mathfrak{g}}\rvert$. For every semi-simple element $X\in \mathfrak{g}(F)$, we let
	$$\displaystyle D^G(X)=\lvert \det \ad_X\mid_{\mathfrak{g}/\mathfrak{g}_X}\rvert$$
	be its Weyl discriminant, where $\mathfrak{g}_X$ stands for the centralizer of $X$. The Weyl discriminant $D^G(g)$ for $g\in G$ semi-simple is defined in a similar way. We also define algebraic variants of the Weyl discriminant by
	$$\displaystyle D^G_{alg}(X)=\det \ad_X\mid_{\mathfrak{g}/\mathfrak{g}_X},\;\; D^G_{alg}(x)=\det 1-\Ad_x\mid_{\mathfrak{g}/\mathfrak{g}_x}$$
	for $X\in \fg_{\rs}$ and $x\in G_{\rs}$ regular semi-simple elements. Note that these extend to regular functions on $\fg$ and $G$ respectively.

	When $G$ is connected and $P\subset G$ is a parabolic subgroup, there is a natural splitting $\cA_P=\cA_G\oplus \cA_P^G$ and we can define as usual subsets $\Delta_P^G\subset (\cA_P^G)^*$, $\Delta_P^{G,\vee}\subset \cA_P^G$, that we call by abuse of terminology the sets of simple roots and coroots associated to the pair $P\subset G$, see e.g. \cite[\S 1.2]{LW}. When $G$ is moreover reductive and clear from the context, we will sometimes drop the superscript.

	\subsection{Twisted spaces}\label{Sect twisted spaces}
	
	In this paper, we will freely use the notion of twisted space due to Labesse as well as corresponding terminology. The main references are \cite{LW} and \cite{WalFTLtordue} but for the reader's convenience we recall most of the definitions here.
	
	A {\em twisted space} is a pair $(G,\tG)$ where $G$ is a group and $\tG$ is a set equipped with two commuting left and right actions
	\begin{equation}\label{eq1 groups}
		\displaystyle G\times \widetilde{G}\times G\to \widetilde{G},\;\; (g,\gamma,g')\mapsto g\gamma g'
	\end{equation}
	each making $\widetilde{G}$ into a principal $G$-homogeneous space and such that $\tG\neq \emptyset$. Similarly, a {\em twisted space over $F$} is a pair $(G,\tG)$ where $G$ is an algebraic group over $F$ and $\tG$ is an algebraic variety over $F$ equipped with two commuting regular actions as in \eqref{eq1 groups} both making $\tG$ into a principal $G$-torsor and such that $\tG(F)\neq \emptyset$. We moreover say that $(G,\tG)$ is {\em reductive} (resp. {\em connected}) if $G$ is so.
	
	Let $(G,\tG)$ be a twisted space over $F$ with $G$ linear. For every $\gamma\in \tG$, we denote by $\Ad_\gamma$ the unique automorphism of $G$ such that
	$$\displaystyle \gamma g=\Ad_\gamma(g)\gamma, \mbox{ for every } g\in G.$$
	We will denote by $\theta_{\tG}$ the outer automorphism of $G$ over $F$ (i.e. an element of $\Aut_F(G)/\Ad(G(F))$) associated to $\Ad_\gamma$ for any $\gamma\in \tG(F)$ (it is independent of $\gamma$). We will also write $\theta$ for $\theta_{\tG}$ when the twisted space $\tG$ is clear from the context. When $\theta_{\tG}$, or $\Ad_\gamma$, induces natural automorphisms on related objects these will invariantly be denoted by the same symbol. For example, $\theta_{\tG}$ induces an automorphism of $A_G$ and $\cA_G$. We write $A_{\tG}$ for the connected component of the subgroup of fixed points $A_G^\theta$. The following condition on $\theta_{\tG}$ will be assumed throughout:
	\begin{num}
		\item\label{eq1 twisted spaces} the outer automorphism $\theta_{\tG}$ is of finite order.
	\end{num}
	Note that if $G$ is reductive, this is equivalent to the restriction $\theta\mid_{Z(G)}$ to the center of $G$ being of finite order. We set
	$$\displaystyle \cA_{\tG}=\Hom(X^*(A_{\tG}),\mathbb{R}).$$
	Then, $\cA_{\tG}$ can naturally be identified with the subspace of $\theta_{\tG}$-invariants in $\cA_G$ and due to condition \eqref{eq1 twisted spaces} it admits a unique $\theta_{\tG}$-stable complement so that we have a canonical projection $\cA_G\to \cA_{\tG}$. We denote by
	$$\displaystyle H_{\tG}:G(F)\to \cA_{\tG}$$
	the composite of $H_G$ with that projection and by $\cA_{\tG,F}$ the lattice $H_{\tG}(G(F))$. We also set $i\cA_{\tG,F}^\vee=\Hom(\cA_{\tG,F},2i\pi \mathbb{Z})$, a subgroup of $i\cA_{\tG}^*$, and $i\cA_{\tG,F}^*=i\cA_{\tG}^*/i\cA_{\tG,F}^\vee$. We also denote by $\delta_{\tG}$ the ``modular character'' of $\tG(F)$ that is the function $\delta_{\tG}:\tG(F)\to \mathbb{R}^+$ defined by
	$$\displaystyle \delta_{\tG}(\gamma)=\lvert \det \Ad_\gamma\mid_{\mathfrak{g}}\rvert.$$
	
	For a subset $X\subset G$, the {\em normalizer} of $X$ in $\tG$ is the subset of $\gamma\in \tG$ such that $\Ad_{\gamma}(X)=X$. We denote by $\Norm_{\tG}(X)$ the normalizer of $X$ in $\tG$ and by $\Norm_{\tG(F)}(X)=\Norm_{\tG}(X)\cap \tG(F)$. Similarly, for a subset $X\subset \tG$ we write $\Norm_G(X)$ (resp. $Z_G(X)$) for the normalizer (resp. the centralizer) of $X$ in $G$ that is the subset of $x\in G$ such that $x^{-1}Xx=X$ (resp. $x^{-1}\gamma x=\gamma$ for every $\gamma\in X$) and we set $\Norm_{G(F)}(X)=\Norm_G(X)\cap G(F)$ (resp. $Z_{G(F)}(X)=Z_G(X)\cap G(F)$). When $\gamma\in \tG$, we simply write $Z_G(\gamma)$ for $Z_G(\{\gamma\})$ and we denote by $G_\gamma$ the neutral component of $Z_G(\gamma)$.
	
	We henceforth assume that $G$ is connected and reductive. A {\em twisted parabolic subspace} of $\tG$ is the normalizer $\tP=\Norm_{\tG}(P)$ of a parabolic subgroup $P\subset G$ satisfying $\tP(F)\neq \emptyset$ (or equivalently $\tP(\overline{F})\neq \emptyset$). Note that the parabolic subgroup $P$ is entirely determined by $\tP$ and that $(P,\tP)$ is a twisted space over $F$. If $\tP$ is a twisted parabolic subspace of $\tG$, a {\em Levi component} of $\tP$ is the normalizer $\tM=\Norm_{\tP}(M)$ in $\tP$ of a Levi component $M$ of $P$. Note that the condition $\tP(F)\neq \emptyset$ implies $\tM(F)\neq \emptyset$. A {\em twisted Levi subspace} of $\tG$ is a Levi component $\tM$ of some twisted parabolic subspace $\tP$ of $\tG$. Note that if $\tP\subset \tG$ is a parabolic subspace and $\tM\subset \tP$ is a Levi component of it, we have (canonically) $A_{\tM}=A_{\tP}$ and $\mathcal{A}_{\tM}=\mathcal{A}_{\tP}$.
	
	Let $\tM\subset \tG$ be a twisted Levi subspace. We denote by $\mathcal{P}(\tM)$ (resp. $\mathcal{F}(\tM)$) the set of twisted parabolic subspaces with Levi component $\tM$ (resp. containing $\tM$). For $\tQ\in \mathcal{F}(\tM)$, we have a natural decomposition
	$$\displaystyle \mathcal{A}_{\tM}=\mathcal{A}_{\tQ}\oplus \mathcal{A}_{\tM}^{\tQ}.$$
	We will also write $\mathcal{A}_{\tP}^{\tQ}$ for $\mathcal{A}_{\tM}^{\tQ}$ for every $\tP\in \mathcal{P}(\tM)$.
	
	For two parabolic subspaces $\tP\subset \tQ$, we denote by $\Delta^{\tQ,\vee}_{\tP}$ and $\Delta^{\tQ}_{\tP}$ the respective images of $\Delta^{Q,\vee}_{P}$ and $\Delta^Q_P$ by the natural projections $\CA_P^Q\to \CA_{\tP}^{\tQ}$ and $(\CA_P^Q)^*\to (\CA_{\tP}^{\tQ})^*$.

	If $M\subset G$ is a Levi subgroup (not necessarily corresponding to any twisted Levi subspace of $\tG$), we set
	$$\displaystyle W^G(M)=\Norm_{G(F)}(M)/M(F) \mbox{ and } W^{\tG}(M)=\Norm_{\tG(F)}(M)/M(F).$$
	Note that if $W^{\tG}(M)\neq \emptyset$ then $(W^G(M),W^{\tG}(M))$ is a twisted group.
	
	Two twisted parabolic subspaces $\tP$ and $\tQ$ of $\tG$ are called {\em opposite} if the corresponding parabolic subgroups $P$ and $Q$ of $G$ are so or, equivalently, if the intersection $\tP\cap \tQ$ is a common Levi component of $\tP$ and $\tQ$. If $\tM\subset \tG$ is a twisted Levi subspace and $\tP\in \mathcal{P}(\tM)$, there exists a unique $\tQ\in \mathcal{P}(\tM)$ which is opposite to $\tP$.
	
	There is also a notion of {\em twisted maximal torus}: it is a subvariety $\tT\subset \tG$ defined over $F$ for which there exists a Borel pair $(B,T)$, not necessarily defined over $F$, such that $\tT=\Norm_{\tG}(B)\cap \Norm_{\tG}(T)$ and $\tT(F)\neq \emptyset$. If $\tT\subset \tG$ is a twisted maximal torus, the torus $T\subset G$ is uniquely determined by $\tT$ and is defined over $F$. Moreover, the pair $(T,\tT)$ is a twisted space over $F$. We say that a twisted maximal torus $\tT\subset \tG$ is {\em elliptic} if $A_{\tT}=A_{\tG}$. For any twisted maximal torus $\tT\subset \tG$ we set
	$$\displaystyle W(G,\tT)=\Norm_{G(F)}(\tT)/T(F).$$
	
	An element $\gamma\in \tG$ is {\em semisimple} if $\Ad_\gamma$ normalizes a Borel pair $(B,T)$ (not necessarily defined over $F$). The subset of semisimple elements is denoted $\tG_{ss}$. Also, we say that $\gamma\in \tG$ (resp $\gamma\in \tG(F)$) is {\em regular semisimple} (resp. {\em regular elliptic}) if the neutral component $G_\gamma$ of its centralizer $Z_G(\gamma)$ is a torus (resp. a torus anisotropic modulo $A_{\tG}$). We denote by $\tG_{\rs}$ (resp. $\tG(F)_{\elli}$) the open subset of regular semi-simple elements in $\tG$ (resp. of regular elliptic elements in $\tG(F)$). We write $\Gamma_{\elli}(\tG)$ for the set of $G(F)$-conjugacy classes in $\tG(F)_{\elli}$ and for $\gamma\in \tG_{ss}(F)$ we define its {\em Weyl discriminant} by
	$$\displaystyle D^{\tG}(\gamma)=\lvert \det(1-\Ad_\gamma)\mid_{\mathfrak{g}/\mathfrak{g}_\gamma}\rvert$$
	where $\mathfrak{g}_\gamma$ stands for the Lie algebra of $G_\gamma$.
	
	We henceforth fix a minimal parabolic subgroup $P_{\min}$ of $G$ with a Levi decomposition $P_{\min}=M_{\min}U_{\min}$ and we let $\tP_{\min}=\Norm_{\tG}(P_{\min})$, $\tM_{\min}=\Norm_{\tP_{\min}}(M_{\min})$. Then, $\tP_{\min}$ is a minimal parabolic subspace of $\tG$ and $\tM_{\min}$ a Levi component of it. We denote by $\mathcal{L}(\tM_{\min})$ the set of twisted Levi subspaces of $\tG$ containing $\tM_{\min}$ and for every $\tM\in \mathcal{L}(\tM_{\min})$ we set
	$$\displaystyle \widetilde{W}^M=\Norm_{M(F)}(\tM_{\min}(F))/M_{\min}(F).$$
	We also fix a special maximal compact subgroup $K$ of $G(F)$ in good position relative to $M_{\min}$.
	
	\subsection{Log-norms and Harish-Chandra $\Xi$ function}\label{sect HCS space}
	
	In this paper we shall freely use the notion of log-norms on algebraic varieties over $F$ as defined in \cite[\S 1.2]{BeuGGP}, which are simple variants of the norms introduced by Kottwitz in \cite[\S 18]{KottHA}. For every algebraic variety $X$ over $F$, we will fix a log-norm $\sigma_X$ on it and, for $C>0$, we denote by $\textbf{1}_{\sigma_X\leq C}$ the characteristic function of the set
	$$\displaystyle \{x\in X(F)|\;\sigma_X(x)\leq C\}.$$
	In particular, we have log-norms $\sigma_G$ and $\sigma_{\tG}$ on $G$ and $\tG$ respectively that for simplicity we will both denote by $\sigma$. For any given base-point $\gamma_0\in \tG(F)$ we have $\sigma(g\gamma_0)\sim \sigma(g)$ for $g\in G(F)$. Moreover, it will be convenient to assume, as we may, that $\sigma$ is left and right $K$-invariant for some chosen special maximal compact subgroup $K\subset G(F)$.

	\begin{lem}\label{lem proper map}
		Let $W,Z$ be two algebraic varieties over $F$ and $f:W\to Z$ be a proper morphism. Then, we have
		\begin{equation}\label{eq lem log-norms2}
			\displaystyle \sigma_Z(f(x))\sim \sigma_W(x),\mbox{ for } x\in W(\overline{F}).
		\end{equation} 
	\end{lem}
	
	\begin{proof}
		By Chow's lemma \cite[\href{https://stacks.math.columbia.edu/tag/02O2}{Tag 02O2}]{stacks-project}, there exists for some $n\geqslant 0$ a closed subscheme $W'\hookrightarrow Z\times \mathbb{P}^n$ with a surjective regular morphism $\pi: W'\to W$ such that the following diagram commutes	
		$$\displaystyle \xymatrix{W \ar[rd]^{f} & W' \ar[l]^{\pi} \ar@{^{(}->}[r]  & Z\times \mathbb{P}^n \ar[ld]^{pr_1} \\ & Z & }$$
		where $pr_1: Z\times \mathbb{P}^n\to Z$ stands for the first projection. It is readily seen that
		$$\displaystyle \sigma_{\mathbb{P}^n}(y)\sim 1, \mbox{ for } y\in \mathbb{P}^n(\overline{F})$$
		and therefore
		$$\displaystyle \sigma_{Z\times \mathbb{P}^n}(z,y)\sim \sigma_{Z}(z)+\sigma_{\mathbb{P}^n}(y)\sim \sigma_{Z}(z),\mbox{ for } (z,y)\in Z(\overline{F})\times \mathbb{P}^n(\overline{F}).$$
		As $W'\hookrightarrow Z\times \mathbb{P}^n$ is a closed immersion, it follows that
		$$\displaystyle \sigma_{Z}(f(\pi(x'))) \ll \sigma_W(\pi(x')) \ll \sigma_{W'}(x')\sim \sigma_{Z \times \mathbb{P}^n}(x')\sim \sigma_{Z}(pr_1(x'))=\sigma_{Z}(f(\pi(x'))),\; x'\in W'(\overline{F}).$$
		Hence, $\sigma_{Z}(f(\pi(x'))) \sim \sigma_W(\pi(x'))$ for $x'\in W'(\overline{F})$. As $\pi$ is surjective, this implies \eqref{eq lem log-norms2}.
	\end{proof}

	We also denote by $\Xi^G$, or simply by $\Xi$, the basic spherical function of Harish-Chandra i.e. the normalized spherical matrix coefficient (with respect to some choice of special compact subgroup $K\subset G(F)$) of the unramified representation with trivial Satake parameter. Fixing a base-point $\gamma_0\in \tG(F)$, we also define a function $\Xi^{\tG}$ on $\tG(F)$ by
	$$\displaystyle \Xi^{\tG}(g\gamma_0)=\Xi^G(g),\mbox{ for } g\in G(F).$$
	Standard properties of $\Xi^G$ have obvious analogs for $\Xi^{\tG}$ e.g. we have (see \cite[proposition II.4.5]{WalPlanch}):
	\begin{num}
		\item\label{eq1 Xi function} Let $\tP=\tM U_P$ be a parabolic subspace of $\tG$. Then, for every $d>0$, there exists $d'>0$ such that
		$$\displaystyle \delta_{\tP}(\tm)^{1/2}\int_{U_P(F)} \Xi^{\tG}(\tm u) \sigma(\tm u)^{-d'} du\ll \Xi^{\tM}(\tm) \sigma(\tm)^{-d}, \mbox{ for } \tm\in \tM(F).$$
	\end{num}
	From the `doubling principle' \cite[lemme II.1.3]{WalPlanch} we also deduce:
	\begin{num}
		\item\label{eq2 Xi function} For every compact-open subset $\omega_{\tG}$ of $\tG(F)$ we have
		$$\displaystyle \int_{\omega_{\tG}} \Xi^{\tG}(x\gamma y)d\gamma\ll \Xi^G(x)\Xi^G(y), \mbox{ for } x,y\in G(F).$$
	\end{num}
	
	We let $\cC(\tG(F))$ be the {\em Harish-Chandra Schwartz space} of $\tG(F)$ i.e. the space of functions $f:\tG(F)\to \mathbb{C}$ that are left and right invariant by some compact-open subgroup of $G(F)$ and such that, for every $d>0$, we have
	\begin{equation*}
		\displaystyle \sup_{\gamma\in \tG(F)}\lvert f(\gamma)\rvert \Xi^{\tG}(\gamma)^{-1}\sigma(\gamma)^d <\infty.
	\end{equation*}
	For every compact-open subgroup $J\subset G(F)$, the subspace $\cC(J\backslash \tG(F)/J)$ of $J$-biinvariant functions is naturally a Fr\'echet space, the topology being associated to the seminorms defined by the above supremum for every $d>0$, and $\cC(\tG(F))=\bigcup_J \cC(J\backslash \tG(F)/J)$ is a strict LF space. Moreover, the subspace $C_c^\infty(\tG(F))$ of locally constant compactly supported functions is dense in $\cC( \tG(F))$. The Harish-Chandra Schwartz space $\cC(G(F))$ of $G(F)$ is defined similarly (it suffices to replace $\Xi^{\tG}$ and $\sigma_{\tG}$ by $\Xi^G$ and $\sigma_G$ in the definition). We denote by ${}^0 \cC(G(F))$ the subspace of {\it cusp forms} i.e. of functions $f\in \cC(G(F))$ such that for every proper parabolic subgroup $P=MU_P\subsetneq G$,
	$$\displaystyle \int_{U_P(F)} f(xu)du=0,\;\; \mbox{ for every } x\in G(F).$$
	
	\subsection{Measures}\label{Sect measures}
	
	Let $T$ be a torus (over $F$). We equip $T(F)$ with a Haar measure as follows: if $T$ is split we choose the unique Haar measure giving to the maximal compact subgroup $T(F)_c$ measure one, in general we endow $T(F)$ with the measure such that its quotient by the measure just defined on $A_T(F)$ gives $T(F)/A_T(F)=(T/A_T)(F)$ a total mass of one.
	
	Let $\tT\subset \tG$ be a twisted maximal torus. The neutral connected component $T^{\theta,0}$ of the subgroup $T^\theta$ of $\theta_{\tT}$-fixed points is a torus and therefore, $T^{\theta,0}(F)$ is already equipped with a measure as in the above discussion. Let $\tT(F)/(1-\theta)(T(F))$ be the quotient of $\tT(F)$ by the adjoint action of $T(F)$. We endow $\tT(F)/(1-\theta)(T(F))$ with the unique left and right $T(F)$-invariant measure such that, for every $\gamma\in \tT(F)$, the application
	\begin{equation*}
	\displaystyle T^{\theta,0}(F)\to \tT(F)/(1-\theta)(T(F))
	\end{equation*}
	$$\displaystyle t\mapsto \gamma t$$
	is locally measure preserving. Set $\tT_{\reg}(F):=\tT(F)\cap \tG_{\rs}(F)$. Then, $\tT_{\reg}(F)/(1-\theta)(T(F))$ is an open subset of $\tT(F)/(1-\theta)(T(F))$ (for the quotient topology). To simplify notation, we will write $\tT(F)_{/\theta}$ (resp. $\tT_{reg}(F)_{/\theta}$) for $\tT(F)/(1-\theta)(T(F))$ (resp. $\tT_{\reg}(F)/(1-\theta)(T(F))$). Similarly, we write $\tT_{/\theta}$ for the GIT quotient of $\tT$ by the adjoint action of $T$.
	
	We endow the real vector spaces $\mathcal{A}_G$ and $\mathcal{A}_{\tG}$ with the unique Haar measures for which the lattices $H_G(A_G(F))$ and $H_{\tG}(A_{\tG}(F))$ are of covolume one. Through the exponential map, $i\mathcal{A}_G^*$ and $i\mathcal{A}^*_{\tG}$ can be identified with the Pontryagin duals of $\mathcal{A}_G$ and $\mathcal{A}_{\tG}$. We equip them with the dual measures.
	
	Let $\mathcal{T}_{\elli}(\tG)$ (resp. $\mathcal{T}(\tG)$) be a set of representatives of the $G(F)$-conjugacy classes of elliptic twisted maximal tori (resp. twisted maximal tori) in $\tG$. We equip the set $\Gamma_{\elli}(\tG)$ (resp. $\Gamma(\tG)$) of regular elliptic conjugacy classes in $\tG(F)$ (resp. regular semisimple conjugacy classes in $\tG(F)$) with a measure which is characterized by:
	\begin{equation*}
		\displaystyle \int_{\Gamma_{\elli}(\tG)} \varphi(\gamma) d\gamma=\sum_{\tT \in \mathcal{T}_{\elli}(\tG)} \lvert W(G,\tT)\rvert^{-1} [T^\theta(F):T^{\theta,0}(F)]^{-1}\int_{\tT_{\reg}(F)_{/\theta}} \varphi(t)dt
	\end{equation*}
\begin{equation*}
	\displaystyle \left(\mbox{resp. } \int_{\Gamma(\tG)} \varphi(\gamma) d\gamma=\sum_{\tT \in \mathcal{T}(\tG)} \lvert W(G,\tT)\rvert^{-1} [T^\theta(F):T^{\theta,0}(F)]^{-1}\int_{\tT_{\reg}(F)_{/\theta}} \varphi(t)dt \right).
\end{equation*}
	for every ``reasonable'' function $\varphi$ on $\Gamma_{\elli}(\tG)$ (resp. on $\Gamma(\tG)$).
	
	These conventions also apply to the parabolic and Levi subgroups (resp. parabolic and Levi subspaces) of $G$ (resp. of $\tG$). In particular, the definition of the measures on $\Gamma_{\elli}(\tG)$ and $\Gamma(\tG)$ was chosen so that Weyl's integration formula for $\tG(F)$ (see \cite[\S 4.1]{WalFTLtordue}) takes the following forms:
	\begin{align}\label{Weyl integration formula}
		\displaystyle \int_{\tG(F)} f(\gamma) d\gamma & =\int_{\Gamma(\tG)} D^{\tG}(\gamma) \int_{G_\gamma(F)\backslash G(F)} f(g^{-1}\gamma g) dg d\gamma \\
	\nonumber	& =\sum_{\tM\in \mathcal{L}(\tM_{\min})} \lvert \widetilde{W}^M\rvert \lvert \widetilde{W}^G\rvert^{-1} \int_{\Gamma_{\elli}(\tM)} D^{\tG}(\gamma) \int_{G_\gamma(F)\backslash G(F)} f(g^{-1}\gamma g) dg d\gamma
	\end{align}
	for every $f\in L^1(\tG(F))$ where in the above formula, we have chosen a Haar measure on $G(F)$ from which we deduce a measure on $\tG(F)$ by translation by any element $\gamma\in \tG(F)$ and we put on the $F$-points of the torus $G_\gamma$ the canonical measure defined above.
	
	There is yet another description of the measure on $\Gamma(\tG)$ that will be useful to us. More precisely, for every $\gamma\in \tG_{\rs}(F)$ set $\tG_\gamma=\gamma G_\gamma$. Then, the pair $(G_\gamma,\tG_\gamma)$ is a twisted torus of a very special form, namely every element of $G_\gamma$ commutes with every element of $\tG_\gamma$, which is however rarely a maximal twisted torus of $\tG$ (unless the outer automorphism $\theta_{\tG}$ is trivial). Let $\mathcal{S}(\tG)$ be a set of representatives of twisted torus of the form $(G_\gamma,\tG_\gamma)$, $\gamma\in \tG_{\rs}(F)$, for the natural action by conjugation of $G(F)$. Then, there exist constants $(c_{\tS})_{\tS\in \mathcal{S}(\tG)}$ such that the measure on $\Gamma(\tG)$ takes the form:
	\begin{equation}\label{eq1 measures}
	\displaystyle \int_{\Gamma(\tG)} \varphi(\gamma) d\gamma=\sum_{\tS \in \mathcal{S}(\tG)} c_{\tS}\int_{\tS_{\reg}(F)} \varphi(s)ds
	\end{equation}
where, as before, we have set $\tS_{\reg}=\tS\cap \tG_{\rs}$ and the measure on $\tS_{\reg}(F)$ is the restriction of the translation of the natural Haar measure on $S(F)$ introduced above to $\tS(F)$. More precisely, we can take for $\mathcal{S}(\tG)$ the set of twisted tori $(T^{\theta,0},T^{\theta,0}t)$ where $\tT$ runs over the set of representatives $\mathcal{T}(\tG)$ fixed above and $t$ describes a set of representatives of the $T^{\theta,0}(F)$ orbits in $\tT(F)/\theta$. Then, the above integration formula follows readily from the definition of the measure on $\Gamma(\tG)$ with constants
$$\displaystyle c_{\tS}=\lvert W(G,\tT)\rvert^{-1} [T^\theta(F): T^{\theta,0}(F)]^{-1} \lvert T^{\theta,0}(F)\cap (1-\theta)(T(F))\rvert^{-1}.$$

	For $P\subset Q$ (resp. $\tP\subset \tQ$) two parabolic subgroups of $G$ (resp. two parabolic subspaces of $\tG$) we equip
	$$\displaystyle \mathcal{A}^Q_P=\mathcal{A}_P/\mathcal{A}_Q \mbox{ (resp. } \mathcal{A}_{\tP}^{\tQ}=\mathcal{A}_{\tP}/\mathcal{A}_{\tQ})$$
	with the quotient of the two Haar measures just defined.
	
	All other groups considered will be equipped with Haar measures whose normalization does not really matter. However, for some intermediate steps, it will be convenient to assume that for $P=MU_P$ a parabolic subgroup of $G$, the Haar measures are chosen so that we have the following integration formula:
	$$\displaystyle \int_{G(F)} f(g)dg=\int_{M(F)}\int_{U_P(F)} \int_{K} f(muk) dkdudm.$$
	
	Finally, for a Levi subspace $\tM$ of $\tG$, we endow $\tM(F)$ with the unique (biinvariant) measure such that for every $\gamma\in \tM(F)$ the bijection $m\in M(F)\mapsto \gamma m\in \tM(F)$ is measure-preserving.
	
	\subsection{Estimates}
	
	Let $\tT\subset \tG$ be a twisted maximal torus. In this section, we denote by $\theta=\theta_{\tT}$ the restriction of $\Ad_t$ to $T$ for any $t\in \tT$ (it does not depend on $t$). As in the previous section, we write $\tT_{\reg}(F)_{/\theta}=\tT_{\reg}(F)/(1-\theta)(T(F))$ for the quotient of $\tT_{\reg}(F)$ by the adjoint action of $T(F)$ and we denote by $\tT_{/\theta}=\tT/(1-\theta)T$ the categorical quotient of $\tT$ by the adjoint action of $T$.
	
	\begin{lem}\label{lem1 estimates}
		We have
		\begin{equation*}
			\displaystyle \inf_{t\in T(F)} \sigma_G(tg)\ll \sigma_{\tG}(g^{-1}\gamma g)+\lvert \log D^{\tG}(\gamma)\rvert
		\end{equation*}
		for $(g,\gamma)\in G(F)\times \tT_{\reg}(F)$.
	\end{lem}
	
	\begin{proof}
		Let $Y=\tT_{\reg}\times^T G$ be the quotient of $\tT_{\reg}\times G$ by the free action of $T$ given by $t\cdot (\gamma,g)=(t\gamma t^{-1},tg)$. Then, the regular map
		$$\displaystyle Y\to \tG_{\rs},\; [\gamma,g]\mapsto g^{-1}\gamma g$$
		is finite. Thus, by \cite[Proposition 18.1(1)]{KottHA}, we have
		\begin{equation}\label{eq1 estimates}
			\displaystyle \sigma_Y(\gamma,g)\sim \sigma_{\tG_{\rs}}(g^{-1}\gamma g)\sim \sigma_{\tG}(g^{-1}\gamma g)+\lvert \log D^{\tG}(\gamma)\rvert,\; \mbox{ for } [\gamma,g]\in Y(F).
		\end{equation}
		On the other hand, the regular map $Y\to T\backslash G$, $[\gamma,g]\mapsto Tg$, implies that
		\begin{equation}\label{eq2 estimates}
			\displaystyle \sigma_{T\backslash G}(g)\ll \sigma_Y(\gamma,g),\; \mbox{ for } [\gamma,g]\in Y(F).
		\end{equation}
		Finally, by \cite[Proposition 18.3]{KottHA}, we have
		\begin{equation}\label{eq3 estimates}
			\displaystyle \sigma_{T\backslash G}(g)\sim \inf_{t\in T(F)} \sigma(tg),\; \mbox{ for } g\in G(F).
		\end{equation}
		The lemma readily follows from the combination of \eqref{eq1 estimates}, \eqref{eq2 estimates} and \eqref{eq3 estimates}.
	\end{proof}
	
	For every positive function $f$ on $\tG(F)$ and $\gamma\in \tT_{\reg}(F)$ we set
	\begin{equation*}
		\displaystyle J_{\tG}(\gamma,f)=D^{\tG}(t)^{1/2} \int_{A_{\tT}(F)\backslash G(F)} f(g^{-1}\gamma g) dg
	\end{equation*}
	whether the integral is convergent or not. Note that this expression only depends on the image of $\gamma$ in $\tT_{\reg}(F)_{/\theta}$.
	
	\begin{prop}\label{prop conv WOI}
		For every $d>0$ there exists $d'>0$ such that the orbital integral $J_{\tG}(\gamma,\Xi^{\tG}\sigma_{\tG}^{-d'})$ is convergent for all $\gamma\in \tT_{\reg}(F)$ and we have
		\begin{equation*}
			\displaystyle \sup_{\gamma\in \tT_{\reg}(F)_{/\theta}} \sigma_{\tT_{/\theta}}(\gamma)^d J_{\tG}(\gamma,\Xi^{\tG}\sigma_{\tG}^{-d'})<\infty.
		\end{equation*}
	\end{prop}
	
	\begin{proof}
		Let $\tM\subset \tG$ be the centralizer of $A_{\tT}$. Then, $\tM$ is a twisted Levi subspace. Choose a parabolic subspace $\tP=\tM U_P\in \mathcal{P}(\tM)$. By the Iwasawa decomposition $G(F)=P(F)K$ and a standard Jacobian computation, up to a constant depending on measures, for every positive function $f$ on $\tG(F)$ we have
		\begin{equation*}
			\displaystyle J_{\tG}(\gamma,f)=J_{\tM}(\gamma,f_{\tP})
		\end{equation*}
		where $f_{\tP}$ is the function on $\tM(F)$ defined by
		\begin{equation*}
			\displaystyle f_{\tP}(\tm)=\delta_{\tP}(\tm)^{1/2} \int_K \int_{U_P(F)} f(k^{-1} \tm u k) du dk,\;\;\; \tm\in \tM(F).
		\end{equation*}
		Therefore, by \eqref{eq1 Xi function}, up to replacing $\tG$ by $\tM$ we may assume that $A_{\tT}=A_{\tG}$ i.e. that $\tT$ is elliptic in $\tG$. The statement of the proposition can also be readily reduced to the case where $A_{\tG}=1$ which we assume from now on. Then, as $\tT$ is elliptic the quotient $\tT(F)_{/\theta}$ is compact and we just need to show the existence of $d_0>0$ such that
		\begin{equation*}
			\displaystyle \sup_{\gamma\in \tT_{\reg}(F)_{/\theta}} J_{\tG}(\gamma,\Xi^{\tG}\sigma_{\tG}^{-d_0})<\infty.
		\end{equation*}
		Assume for one moment the following claim:
		\begin{num}
			\item\label{eq1 prop conv WOI} There exists $d_0>0$ such that for almost all $\gamma\in \tT_{\reg}(F)_{/\theta}$, the integral defining  $J_{\tG}(\gamma,\Xi^{\tG}\sigma_{\tG}^{-d_0})$ converges.
		\end{num}
		Then, we can conclude as in \cite[Corollary 2]{CloSchw} using Howe's conjecture for twisted groups \cite[Chap. 2, th\'eor\`eme 2.1]{WalFTLtordue}. Indeed, let $(\Omega_n)_{n\geqslant 1}$ be an increasing and exhaustive\footnote{Meaning that $\tG(F)=\bigcup_n \Omega_n$.} sequence of $K$-biinvariants compact subsets of $\tG(F)$ and set $f_n=\mathbf{1}_{\Omega_n}\Xi^{\tG} \sigma_{\tG}^{-d_0}$. Then, $(f_n)$ is an increasing sequence of functions in $C_c(K\backslash \tG(F)/K)$ converging pointwise to $\Xi^{\tG} \sigma_{\tG}^{-d_0}$ hence $J_{\tG}(\gamma,f_n)$ converges to $J_{\tG}(\gamma,\Xi^{\tG}\sigma_{\tG}^{-d_0})$ for all $\gamma\in \tT_{\reg}(F)_{/\theta}$ (whether the last integral is finite or infinite). However, by \cite[Chap. 1, 4.2 (1)]{WalFTLtordue}, the functions $\gamma\in \tT_{\reg}(F)_{/\theta}\mapsto J_{\tG}(\gamma,f_n)$ are locally constant and bounded whereas by \cite[Chap. 2, th\'eor\`eme 2.1]{WalFTLtordue} (``Howe's conjecture'' for the twisted group $\tG$) the vector space they span is finite dimensional. It follows that the function $\gamma\in \tT_{\reg}(F)_{/\theta}\mapsto J_{\tG}(\gamma,\Xi^{\tG}\sigma_{\tG}^{-d_0})$ has the same properties (i.e. it is locally constant and bounded) and this proves the proposition.
		
		It remains to show \eqref{eq1 prop conv WOI}. Set $T^\theta$ for the subgroup of $\theta$-fixed points in $T$ (recall that $\theta=\theta_{\tT}$). Let $\gamma\in \tT_{\reg}(F)$ and let $(T^\theta)'$ be the inverse image of $\tG_{\rs}$ by the morphism $t\in T^\theta\mapsto \gamma t$. Then, the map $t\in (T^\theta)'(F) \mapsto \gamma t((1-\theta)(T(F)))\in \tT_{\reg}(F)_{/\theta}$ is a local homeomorphism. Therefore, by Fubini, we just need to check that for every compact-open subset $\omega\subset (T^\theta)'(F)$, the integral
		\begin{equation}\label{eq3 prop conv WOI}
			\displaystyle \int_{\omega} \int_{G(F)} \Xi^{\tG}(g^{-1}\gamma t g) \sigma_{\tG}(g^{-1}\gamma tg)^{-d_0} dg dt
		\end{equation}
		converges. First, we show that
		\begin{equation}\label{eq2 prop conv WOI}
			\displaystyle \sigma_{G}(g)\ll \sigma_{\tG}(g^{-1}\gamma t g), \mbox{ for } (g,t)\in G(F)\times \omega.
		\end{equation}
		The morphism 
		\begin{equation}
			(T^\theta)'\times T^\theta\backslash G\to \tG_{\rs}, \;\;\; (t,g)\mapsto g^{-1}\gamma t g
		\end{equation}
		is finite \'etale. Therefore, we have
		\begin{equation*}
			\displaystyle \sigma_{T^\theta\backslash G}(g)+\sigma_{(T^\theta)'}(t)\ll \sigma_{\tG_{\rs}}(g^{-1}\gamma t g)\sim \sigma_{\tG}(g^{-1}\gamma tg)+\lvert \log D^{\tG}(\gamma t)\rvert
		\end{equation*}
		for $(g,t)\in G(F)\times (T^\theta)'(F)$. On the other hand, since $\omega$ is compact, we have $\sigma_{(T^\theta)'}(t)\sim 1$ and $\lvert \log D^{\tG}(\gamma t)\rvert\sim 1$ for $t\in \omega$. Combining this with the previous inequality, gives
		\begin{equation*}
			\displaystyle \sigma_{T^\theta\backslash G}(g)\ll \sigma_{\tG}(g^{-1}\gamma tg), \mbox{ for } (g,t)\in G(F)\times \omega.
		\end{equation*}
		Moreover, by \cite[Proposition 18.3]{KottHA} we have
		\begin{equation*}
			\displaystyle \sigma_{T^\theta\backslash G}(g)\sim \inf_{t\in T^\theta(F)} \sigma_G(tg)\sim \sigma_G(g), \mbox{ for } g\in G(F),
		\end{equation*}
		(Recall that the twisted torus $\tT$ is elliptic and $A_{\tG}=1$ hence $T^\theta(F)$ is compact.) and this implies \eqref{eq2 prop conv WOI}.
		
		We now consider the integral \eqref{eq3 prop conv WOI}. By \eqref{eq2 prop conv WOI} it is essentially bounded by
		$$\displaystyle \int_{\omega} \int_{G(F)} \Xi^{\tG}(g^{-1}\gamma t g) \sigma_{G}(g)^{-d_0} dgdkdt$$
		which can be rewritten as
		$$\displaystyle \int_{\omega} \int_K \int_{G(F)} \Xi^{\tG}(g^{-1}k^{-1}\gamma t kg) \sigma_{\tG}(g)^{-d_0} dgdkdt.$$
		Since the map $(T^\theta)'(F)\times T^\theta(F)\backslash G(F) \ni (t,g)\mapsto g^{-1}\gamma t g\in \tG(F)$ is a local $F$-analytic isomorphism, the last expression above is bounded up to a (multiplicative) constant by
		$$\displaystyle \int_{\omega_{\tG}} \int_{G(F)} \Xi^{\tG}(g^{-1}\tg g) \sigma_{\tG}(g)^{-d_0} dgd\tg$$
		for some compact-open subset $\omega_{\tG}$ of $\tG(F)$. Furthermore, by \eqref{eq2 Xi function}, we have $\int_{\omega_{\tG}} \Xi^{\tG}(g^{-1}\tg g) d\tg\ll \Xi^G(g)^2$ for $g\in G(F)$ and the integral \eqref{eq3 prop conv WOI} is therefore bounded up to a constant by
		$$\displaystyle \int_{G(F)} \Xi^G(g)^2 \sigma_G(g)^{-d_0} dg$$
		which is well-known to converge for $d_0$ sufficiently large, see \cite[lemme II.1.5]{WalPlanch}.
	\end{proof}
	
	\subsection{Quasi-characters}\label{Section qc}
	
	Following \cite[\S 1.6]{WalGGPIV}, by a {\em quasi-character} on $\tG(F)$ we mean a function $\Theta: \tG_{\rs}(F)\to \C$ such that for every semisimple element $x\in \tG_{\ssi}(F)$, there is a local expansion 
	\begin{equation}\label{eq1qc}
		\displaystyle \Theta(x\exp(X))=\sum_{\mathcal{O}\in \Nil(\mathfrak{g}_x^*)} c_{\Theta,\mathcal{O}}(x) \widehat{j}_\psi(\mathcal{O},X)
	\end{equation}
	valid for $X\in \mathfrak{g}_{x,\rs}(F)$ sufficiently close to $0$ and where
	\begin{itemize}
		\item $\Nil(\mathfrak{g}_x^*)$ stands for the set of nilpotent $G_x(F)$-coadjoint orbits in $\mathfrak{g}^*_x(F)$;
		
		\item $c_{\Theta,\mathcal{O}}(x)\in \mathbb{C}$ for every $\mathcal{O}\in \Nil(\mathfrak{g}_x^*)$;
		
		\item For $\mathcal{O}\in \Nil(\mathfrak{g}_x^*)$, $\widehat{j}_\psi(\mathcal{O},.)$ is the unique locally integrable function on $\mathfrak{g}_x(F)$ which is locally constant on $\mathfrak{g}_{x,\rs}(F)$ and such that
		$$\displaystyle \int_{\mathfrak{g}_x(F)} \varphi(X) \widehat{j}_\psi(\mathcal{O},X) dX=\int_{\mathcal{O}} \widehat{\varphi}(Z)dZ$$
		for every $\varphi\in C_c^\infty(\mathfrak{g}_x(F))$, where $\widehat{\varphi} \in C_c^\infty(\mathfrak{g}_x^*(F))$ denotes the Fourier transform
		$$\displaystyle Y\in \mathfrak{g}_x^*(F)\mapsto \widehat{\varphi}(Y)=\int_{\mathfrak{g}_x(F)} \varphi(X) \psi(\langle X,Y\rangle) dX$$
		and $dZ$ is the Kirillov-Kostant $G_x(F)$-invariant measure on $\mathcal{O}$ deduced from the canonical symplectic form on $\mathcal{O}$ and the self-dual measure on $F$ associated to $\psi$ (see \cite{HCH}).
	\end{itemize}
	
	For $x\in \tG_{\ssi}(F)$, we denote by $\Nil_{\reg}(\mathfrak{g}_x^*)\subset \Nil(\mathfrak{g}_x^*)$ the subset of regular nilpotent coadjoint orbits.
	
	\begin{lem}\label{lem bounded germ}
		Let $\Theta$ be a quasi-character on $\tG(F)$. The function
		$$\displaystyle x\in \tG_{ss}(F)\mapsto D^{\tG}(x)^{1/2} \max_{\mathcal{O}\in \Nil_{\reg}(\mathfrak{g}_x^*)} \lvert c_{\Theta,\mathcal{O}}(x)\rvert$$
		is locally bounded.
	\end{lem}
	
	\begin{proof}
		Let $x\in \tG_{ss}(F)$ be a semisimple element. By \cite[Lemma 3.2]{HCH}, we have
		$$\displaystyle \widehat{j}_\psi(\cO,tX)=\lvert t\rvert^{-\dim(\cO)/2}\widehat{j}_\psi(\cO,X)$$
		for every $\cO\in \Nil(\fg_x^*)$, $X\in \fg_{x,\rs}(F)$ and $t\in F^{\times,2}$. Moreover, for $X\in \fg_{x,\rs}(F)$ sufficiently close to $0$, we have
		$$\displaystyle D^{\tG}(x\exp(X))=D^{\tG}(x) D^{G_x}(X) \mbox{ and } D^{G_x}(tX)=\lvert t\rvert^{\delta(G_x)}D^{G_x}(X)$$
		for every $t\in F^\times$ where we have set $\delta(G_x)=\dim(G_x)-rk(G_x)$. As for every $\cO\in \Nil(\fg_x^*)$ we have $\dim(\cO)\leqslant \delta(G_x)$ with equality if and only if $\cO$ is regular, we deduce from the expansion \eqref{eq1qc} that for every $X\in \fg_{x,\rs}(F)$ we have
		\begin{equation*}
			\displaystyle \lim\limits_{t\in F^{\times,2}, t\to 0} D^{\tG}(x\exp(tX))^{1/2}\Theta(x\exp(tX))=D^{\tG}(x)^{1/2}D^{G_x}(X)^{1/2}\sum_{\cO\in \Nil_{\reg}(\fg_x^*)} c_{\Theta,\cO}(x) \widehat{j}_\psi(\cO,X).
		\end{equation*}
		Since the functions $\widehat{j}_\psi(\cO,.)$, $\cO\in \Nil(\fg_x^*)$, are linearly independent this implies
		\begin{equation*}
			\displaystyle D^{\tG}(x)^{1/2} \max_{\mathcal{O}\in \Nil_{\reg}(\mathfrak{g}_x^*)} \lvert c_{\Theta,\mathcal{O}}(x)\rvert\leqslant c(G_x) \limsup\limits_{y\in \tG_{\rs}(F), y\to x} D^{\tG}(y)^{1/2} \lvert \Theta(y)\rvert
		\end{equation*}
		where $c(G_x)>0$ is a constant that depends only on the isomorphism class of $G_x$. By \cite[Corollary 6.3]{HCH}, the function $(D^{\tG})^{1/2} \Theta$ is locally bounded and the lemma follows as there are only finitely many isomorphism classes of centralizers $G_x$ for $x\in \tG_{ss}(F)$.
	\end{proof}
	
	\subsection{Representations}\label{Sect rep of G}
	
	In this paper, by a {\em representation} of $G(F)$ we mean a pair $(\pi,V_\pi)$ where $V_\pi$ is a complex vector space and $\pi: G(F)\to GL(V_\pi)$ is a smooth representation of $G(F)$ on $V_\pi$. Most of the time we will omit the space $V_\pi$ and just write $\pi$ for a representation of $G(F)$. For $\lambda\in \mathcal{A}_{G,\C}^*$, we denote by $\pi\mapsto \pi_\lambda$, where $\pi_\lambda(g)=e^{\langle \lambda, H_G(g)\rangle}\pi(g)$, the twisting operation by $\lambda$ on representations of $G(F)$.
	
	Let $\pi$ be a representation of $G(F)$. We denote by $\pi^\vee$ the smooth contregredient of $\pi$ realized in the usual way on the space $V_\pi^\vee$ of smooth functionals on $V_\pi$. We denote by $\langle .,.\rangle$ the canonical pairing on $V_\pi\times V_\pi^\vee$. We say that $\pi$ is {\em tempered} if it is of finite length and for every $(v,v^\vee)\in V_\pi\times V_\pi^\vee$ there exists a constant $C>0$ such that
	\begin{equation*}
		\displaystyle \lvert \langle \pi(g)v,v^\vee\rangle \rvert\leqslant C \Xi^G(g),\; \mbox{ for every } g\in G(F).
	\end{equation*}
	
	We write $\Pi_2(G)$ (resp. $\Temp(G)$) for the set of isomorphism classes of unitary square-integrable (resp. tempered)  irreducible representations of $G(F)$. If $P=MU$ is a parabolic subgroup of $G$ and $\sigma$ is a representation of $M(F)$, we let $I_P^G(\sigma)$ be the smooth normalized parabolic induction of $\sigma$ to $G(F)$. When $\sigma\in \Temp(M)$, we write $I_M^G(\sigma)$ for the isomorphism class of $I_P^G(\sigma)$ where $P\in \mathcal{P}(M)$ (it does not depend on this choice). Define $\Temp_\ind(G)$ as the set of isomorphism classes of representations of $G(F)$ of the form $I_M^G(\sigma)$ where $M$ is a Levi subgroup of $G$ and $\sigma\in \Pi_2(M)$. According to Harish-Chandra, every $\pi\in \Temp(G)$ can be embedded in $I_M^G(\sigma)$ for such a pair $(M,\sigma)$ which is moreover unique up to conjugacy by $G(F)$. Thus, we get a map
	\begin{equation*}
		\displaystyle pr_{G}:\Temp(G)\to \Temp_{\ind}(G).
	\end{equation*}
	
	We equip $\Temp_{\ind}(G)$ with a topology that can be described as follows. Let $M$ be a Levi subgroup of $G$ and $\sigma\in \Pi_2(M)$. Then, the set
	\begin{equation*}
		\displaystyle \mathcal{O}=\{I_M^G(\sigma_\lambda)\mid \lambda\in i\cA_M^* \}
	\end{equation*}
	is a connected component of $\Temp_{\ind}(G)$ and the topology on $\mathcal{O}$ is the quotient topology inherited from $i\cA_M^*$ via the surjection
	\begin{equation}\label{eq1rep}
		\displaystyle \lambda\in i\cA_M^*\mapsto I_M^G(\sigma_\lambda)\in \mathcal{O}.
	\end{equation}
	We say that a function $z:\Temp_{\ind}(G)\to \mathbb{C}$ is {\em smooth} if for every pair $(M,\sigma)$ as before, the composition of $z$ with the map \eqref{eq1rep} gives a $C^\infty$ function on $i\cA_M^*$ in the usual sense. We denote by $C^\infty(\Temp_{\ind}(G))$ the vector space of smooth functions on $\Temp_{\ind}(G)$. It is an algebra for pointwise multiplication. Moreover, by the description of the image by Fourier transform of the Harish-Chandra Schwartz space $\cC(G(F))$ \cite{WalPlanch}, there exists an action
	\begin{equation}\label{eq2rep}
		\displaystyle C^\infty(\Temp_{\ind}(G))\times \cC(G(F))\to \cC(G(F)),\;\; (z,f)\mapsto z\star f
	\end{equation}
	of $C^\infty(\Temp_{\ind}(G))$ on $\cC(G(F))$ which is characterized by
	\begin{equation}\label{eq3rep}
		\displaystyle \pi(z\star f)=z(\pi) \pi(f)
	\end{equation}
	for every $(\pi,z,f)\in \Temp_{\ind}(G)\times C^\infty(\Temp_{\ind}(G))\times \cC(G(F))$. See also \cite{SZ} for a different approach where $C^\infty(\Temp_{\ind}(G))$ is shown to coincide with the so-called {\em tempered Bernstein center} of $G(F)$.
	
	The outer automorphism $\theta$ of $G(F)$ induces a bijection $\theta: \Temp_{\ind}(G)\to \Temp_{\ind}(G)$. We denote by $\Temp_{\ind}(G)^\theta$ the subset of fixed points.
	
	\subsection{Twisted representations}\label{Section twrep}
	
	A (smooth) {\em representation} of the twisted space $\widetilde{G}(F)$ is a pair $(\pi,\widetilde{\pi})$ where $\pi$ is a representation of $G(F)$ and $\widetilde{\pi}$ is a map $\widetilde{G}(F)\to GL(V_\pi)$ satisfying
	$$\widetilde{\pi}(g\gamma g')=\pi(g) \tpi(\gamma) \pi(g'),\mbox{ for every } (g,\gamma,g')\in G(F)\times \tG(F)\times G(F).$$
	Most of the time, we will simply refer to a representation of $\tG(F)$ by the map $\tpi$, the underlying representation $(\pi,V_\pi)$ of $G(F)$ being understood. Note that if $\tpi$ is a representation of $\tG(F)$ then so is $c\tpi$ for every $c\in \C^\times$. Moreover, a representation $\pi$ of $G(F)$ extends to a representation $(\pi,\tpi)$ of $\tG(F)$ (although not uniquely) if and only if its isomorphism class is fixed by the outer automorphism $\theta$. 
	
	Let $\tpi$ be a representation of $\tG(F)$. We say that $\tpi$ is {\em $G$-irreducible} if $\pi$ is irreducible in the usual sense i.e. if there is no nontrivial $G(F)$-invariant subspace of $V_\pi$. We also say that $\tpi$ is {\em admissible} (resp. {\em tempered}) if $\pi$ is so. We denote by $\tpi^\vee$ the {\em smooth contragredient} of $\tpi$ that is the representation of $\tG(F)$ on the space $V_\pi^\vee$ of smooth functionals on $V_\pi$ characterized by
	\begin{equation*}
		\displaystyle \langle \tpi(\gamma)v,\tpi^\vee(\gamma)v^\vee\rangle=\langle v,v^\vee\rangle, \;\mbox{ for } (\gamma,v,v^\vee)\in \tG(F)\times V_\pi\times V_\pi^\vee.
	\end{equation*}
	
	Assume that $\pi$ is of finite length. For every $f\in C_c^\infty(\tG(F))$, we define as usual an operator $\tpi(f)$ on $V_\pi$ characterized by
	\begin{equation}\label{eq1 rep}
		\displaystyle \langle \tpi(f)v,v^\vee\rangle=\int_{\tG(F)} f(\gamma) \langle\tpi(\gamma)v,v^\vee\rangle d\gamma,\;\mbox{ for } (v,v^\vee)\in V_\pi\times V_\pi^\vee.
	\end{equation}
	These operators are of finite rank and, according to \cite[Theorem 3]{Clo}, there exists a quasi-character $\Theta_{\tpi}$ on $\tG(F)$ (in the sense of Section \ref{Section qc}), called the {\em Harish-Chandra character of $\tpi$}, such that
	\begin{equation}\label{eq2 rep}
		\displaystyle \Tr \tpi(f)=\int_{\tG(F)} f(g) \Theta_{\tpi}(g) dg, \mbox{ for every } f\in C_c^\infty(\tG(F)).
	\end{equation}
	For ease of notation, we will denote by
	$$\displaystyle c_{\tpi,\cO}(x)=c_{\Theta_{\tpi},\cO}(x),\mbox{ for every } x\in \tG_{\ssi}(F) \mbox{ and } \cO\in \Nil(\mathfrak{g}_x),$$
	the various coefficients of the germs expansions of $\Theta_{\tpi}$. If moreover $\tpi$ is tempered, the definition \eqref{eq1 rep} of the operator $\tpi(f)$ still makes sense and the formula \eqref{eq2 rep} is still valid for $f\in \cC(\tG(F))$ (the integral being absolutely convergent).

	Let $\tP=\tM U$ be a parabolic subspace of $\tG$ and $\tsigma$ be a representation of $\tM(F)$. We denote by $I_{\tP}^{\tG}(\tsigma)$ the {\em normalized parabolic induction of $\tsigma$} i.e. the representation of $\tG(F)$ on the space of smooth functions $e:G(F)\to V_\sigma$ satisfying
	$$\displaystyle e(mug)=\delta_P(m)^{1/2} \sigma(m)e(g)$$
	for every $(m,u,g)\in M(F)\times U(F)\times G(F)$ the action of $\tG(F)$ being given by
	$$\displaystyle (I_{\tP}^{\tG}(\gamma,\tsigma)e)(g)=\delta_{\tP}(\gamma')^{1/2} \tsigma(\gamma') e(g')$$
	for $(\gamma,g)\in \tG(F)\times G(F)$ where $g\gamma=\gamma' g'$ is any decomposition with $(\gamma',g')\in \tM(F)\times G(F)$ (the right hand side is readily seen to be independent of this decomposition). Note that the underlying representation $I_P^G(\sigma)$ of $G(F)$ is the usual normalized parabolic induction of $\sigma$.
	
	Let $M\subset G$ be a Levi subgroup and $\sigma\in \Pi_2(M)$. We set
	$$\displaystyle \Norm_{\tG(F)}(\sigma)=\{\tn\in \Norm_{\tG(F)}(M)\mid \sigma\circ \Ad_{\tn}\simeq \sigma \},$$ 
	$$\Norm_{G(F)}(\sigma)=\{n\in \Norm_{G(F)}(M)\mid \sigma\circ \Ad_{n}\simeq \sigma \}$$
	and
	$$\displaystyle W^{\tG}(\sigma)=\Norm_{\tG(F)}(\sigma)/M(F),\;\; W^G(\sigma)=\Norm_{G(F)}(\sigma)/M(F).$$
	Assume that $W^{\tG}(\sigma)\neq \emptyset$. Then, $W^{\tG}(\sigma)$ is a torsor under $W^G(\sigma)$ both for left and right multiplication i.e. the pair $(W^G(\sigma), W^{\tG}(\sigma))$ is a twisted space. Let $P\in \mathcal{P}(M)$. As in \cite[Chap 1, \S 2.8]{WalFTLtordue}, and making auxilliary choices (including a regularization of standard intertwining operators), we can define for each $\tw\in W^{\tG}(\sigma)$ an extension of $I_P^G(\sigma)$ to a representation $\widetilde{I_P^G(\sigma)}(\tw,.)$ of $\tG(F)$. This extension depends on the auxilliary choices only up to multiplication by a nonzero scalar. Let $W_0^G(\sigma)$ be the distinguished subgroup of elements $w\in W_0^G(\sigma)$ such that for each $\tw\in W^{\tG}(\sigma)$ the representations $\widetilde{I_P^G(\sigma)}(\tw,.)$ and $\widetilde{I_P^G(\sigma)}(\tw w,.)$ are the same up to a scalar. The {\em twisted $R$-group of $(M,\sigma)$} is the quotient $R^{\tG}(\sigma)=W^{\tG}(\sigma)/W^G_0(\sigma)$. We also denote by $R^G(\sigma)=W^G(\sigma)/W_0^G(\sigma)$ the corresponding $R$-group so that $(R^G(\sigma),R^{\tG}(\sigma))$ is again a twisted space. To every $r\in R^{\tG}(\sigma)$, we associate the twisted representation $\widetilde{I_P^G(\sigma)}(r,.)=\widetilde{I_P^G(\sigma)}(\tw_r,.)$ where $\tw_r\in W^{\tG}(\sigma)$ is some chosen lift of $r$. This representation still depends, up to a scalar, on various choices but, henceforth, we will always assume that all such choices have been made and we will denote its isomorphism class by $\widetilde{I_M^G(\sigma)}(r,.)$ (which, as the notation suggests, does not depend on $P$, at least up to a scalar). Note, however, that the isomorphism class of the twisted representation $\widetilde{I_M^G(\sigma)}(r,.)\otimes \widetilde{I_M^G(\sigma)}(r,.)^\vee$ of $\tG(F)\times \tG(F)$ is completely canonical and independent of all the choices involved.
	
	
	Let $E(\tG)$ be the set of $G(F)$-conjugacy classes of triples $(M,\sigma,r)$ where $(M,\sigma)$ is as above and $r\in R^{\tG}(\sigma)$ is such that the character of $\widetilde{I_M^G(\sigma)}(r,.)$ (which, again, is only well-defined up to a scalar) is nonzero. For $\tau\in E(\tG)$ represented by a triple $(M,\sigma,r)$, we will write $\tpi_\tau$ for the twisted representation $\widetilde{I_M^G(\sigma)}(r,.)$. Actually, for $\tau,\tau'\in E(\tG)$ the representations $\tpi_\tau$ and $\tpi_{\tau'}$ are isomorphic if and only if $\tau=\tau'$ (this follows e.g. from \cite[Chap. 1, proposition 2.9]{WalFTLtordue}) and we will also sometimes identify $E(\tG)$ with the set of isomorphism classes $\{\tpi_\tau\mid \tau\in E(\tG) \}$. Note that for every $\tpi\in E(\tG)$ the isomorphism class of the underlying representation $\pi$ belongs to $\Temp_{\ind}(G)^\theta$. 
	
	
	Each $\tw\in W^{\tG}(\sigma)$ induces an automorphism of $\cA_M$ (induced from $\Ad_{\tn}$ for any lifting $\tn\in \tG(F)$ of $\tw$). Let $W^{\tG}_{\reg}(\sigma)$ be the subset of $\tw\in W^{\tG}(\sigma)$ such that $\cA_M^{\tw}=\cA_{\tG}$. Following \cite[\S 2.11]{WalFTLtordue}, we define $E_{\disc}(\tG)$ (resp. $E_{\elli}(\tG)$) to be the subset of triples $\tau=[M,\sigma,r]\in E(\tG)$ such that $W_0^G(\sigma)r\cap W_{\reg}^{\tG}(\sigma)\neq \emptyset$ (resp. $W_0^G(\sigma)=\{1 \}$ and $r\in W_{\reg}^{\tG}(\sigma)$). We also introduce the further subset $E_2(\tG)$ of triples $\tau=[M,\sigma,r]\in E(\tG)$ such that $W^{\tG}(\sigma)=W_{\reg}^{\tG}(\sigma)$. By \cite[lemme 2.11]{WalFTLtordue}, we have $E_2(\tG)\subset E_{\elli}(\tG)\subset E_{\disc}(\tG)$.
	
	\begin{rmk}
		The set $E(\tG)$, $E_{\disc}(\tG)$ and $E_{\elli}(\tG)$ do not exactly coincide with the ones defined in \cite[Chap. 1, \S 2.9]{WalFTLtordue} but correspond rather to the sets denoted by $E(\tG)/conj$, $E_{\disc}(\tG)/conj$ and $E_{\elli}(\tG)/conj$ in {\em loc. cit.}
	\end{rmk}
	
	
	There is a natural action of $i\cA_{\tG}^*$ on $E(\tG)$ given by $\lambda\cdot [M,\sigma,r]=[M,\sigma_\lambda,r]$ \footnote{Identifying $E(\tG)$ with a set of isomorphism classes of tempered representations of $\tG(F)$ as before, this action is also sending $\tpi$ to its ``twist'' by $\lambda$ but this twist is only well-defined up to a scalar (it requires the choice of an extension to $\tG(F)$ of the unramified character associated to $\lambda$ e.g. through the choice of a base-point).}. This action factors through $i\cA_{\tG,F}^*$ and preserves the subsets $E_{\disc}(\tG)$, $E_{\elli}(\tG)$ and $E_2(\tG)$. Let $J\subset G(F)$ be a compact-open subgroup. Then, we have:
	\begin{num}
		\item\label{eq3 rep} the subset $E_{\disc}(\tG)^J$ of triples $\tau\in E_{\disc}(\tG)$ such that the representation $\pi_{\tau}$ admits nonzero $J$-invariant vectors is finite modulo the action of $i\cA_{\tG}^*$;
	\end{num}
	(see \cite[Chap. 2, Proposition 2.2]{WalFTLtordue} for the case of $E_{\elli}(\tG)$ the proof being entirely similar for $E_{\disc}(\tG)$).
	
	We equip $E_{\disc}(\tG)$ with the unique measure such that for every $\tau\in E_{\disc}(\tG)$, the twisting map $\lambda\in i\cA_{\tG}^*\mapsto \lambda\cdot \tau$ is locally measure preserving. Thus, denoting by $E_{\disc}(\tG)/i\cA_{\tG}^*$ the set of orbits in $E_{\disc}(\tG)$ under the action of $i\cA_{\tG}^*$, for every sufficiently nice function $\varphi:E_{\disc}(\tG)\to \C$\footnote{In practice, we will only consider functions $\varphi$ that are supported in a finite number of $i\cA_{\tG}^*$-orbits and such that for every $\tau\in E_{\disc}(\tG)$, $\lambda\in i\cA_{\tG}^*\mapsto \varphi(\lambda\cdot \tau)$ is continuous (even $C^\infty$).} we have
	\begin{equation*}
		\displaystyle \int_{E_{\disc}(\tG)} \varphi(\tau) d\tau=\sum_{\tau\in E_{\disc}(\tG)/i\cA_{\tG}^*} \lvert \Stab(i\cA_{\tG,F}^*,\tau)\rvert^{-1} \int_{i\cA_{\tG,F}^*} \varphi(\lambda\cdot \tau) d\lambda
	\end{equation*}
	where we have denoted by $\Stab(i\cA_{\tG,F}^*,\tau)$ the stabilizer of $\tau$ in $i\cA_{\tG,F}^*$.
	
	For $\tau=[M,\sigma,r]\in E_{\disc}(\tG)$, we set (following \cite[Sect. 2.11]{WalFTLtordue})
	\begin{equation}\label{eq4 rep}
		\displaystyle \iota(\tau)= \lvert R^{G}(\sigma)_r\rvert^{-1} \lvert W_0^G(\sigma)\rvert^{-1} \sum_{\tw \in W^G_0(\sigma)r\cap W^{\tG}_{\reg}(\sigma)} \epsilon_\sigma(\tw)\lvert \det (1-\tw)\mid_{\cA_M^{\tG}} \rvert^{-1}
	\end{equation}
	where $\cA_M^{\tG}=\cA_M/\cA_{\tG}$, $R^{G}(\sigma)_r$ denotes the centralizer of $r$ in $R^G(\sigma)$ and the $\epsilon_\sigma(\tw)$ are certain signs defined in {\it loc. cit.} In the particular case where $\tau\in E_{\elli}(\tG)$ this simplifies to
	$$\displaystyle \iota(\tau)=D(\tau):=\lvert R^{G}(\sigma)_r\rvert^{-1} \lvert \det (1-r)\mid_{\cA_M^{\tG}} \rvert^{-1}.$$

	
	\begin{lem}\label{lem rep}
		Let $\tpi\in E(\tG)$. Then, $\tpi\in E_2(\tG)$ if and only if $\{\pi_\lambda\mid \lambda\in i\cA_{\tG}^*\}$ is a connected component of $\Temp_{\ind}(G)^\theta$.
	\end{lem}
	
	\begin{proof}
		Let $M\subset G$ be a Levi subgroup and $\sigma\in \Pi_2(M)$. Then, $\pi=I_M^G(\sigma)\in \Temp_{\ind}(G)^\theta$ if and only if $W^{\tG}(\sigma)\neq \emptyset$. Assume this is the case and set $\pi_\lambda=I_M^G(\sigma_\lambda)$ for every $\lambda\in i\cA_M^*$. Then, it suffices to show that
		$$\displaystyle \pi\otimes i\cA_{\tG}^*:=\{\pi_\lambda\mid \lambda\in i\cA_{\tG}^*\}$$
		is a connected component of $\Temp_{\ind}(G)^\theta$ if and only if $W^{\tG}(\sigma)=W^{\tG}_{\reg}(\sigma)$. This, in turn, is an easy consequence of the following claim:
		\begin{num}
			\item There exists a neighborhood $U\subset i\cA_M^*$ of $0$ such that for every $\lambda\in U$, $\pi_\lambda\in \Temp_{\ind}(G)^\theta$ if and only if there exists $\tw\in W^{\tG}(\sigma)$ such that $\tw \lambda=\lambda$.
		\end{num}
		To prove the claim, we first observe that, for $\lambda\in i\cA_M^*$, $\pi_\lambda\in \Temp_{\ind}(G)^\theta$ if and only if there exists $\tw\in W^{\tG}(M)$ such that $\sigma_\lambda\circ \Ad_{\tw}\simeq \sigma_\lambda$. Moreover, we can find a sufficiently small $W^G(\sigma)$-invariant neighborhood $U\subset i\cA_M^*$ of $0$ such that:
		\begin{itemize}
			\item For every $\lambda,\mu\in U$, we have $\sigma_\lambda\simeq \sigma_\mu$ if and only if $\lambda=\mu$;
			\item for every $\tw\in W^{\tG}(M)\setminus W^{\tG}(\sigma)$ and $\lambda\in U$, we have $\sigma_\lambda\circ \Ad_{\tw}\notin \sigma\otimes U$.
		\end{itemize}
		It follows that, for $\lambda\in U$, we have $\pi_\lambda\in \Temp_{\ind}(G)^\theta$ if and only if there exists $\tw\in W^{\tG}(\sigma)$ such that $\sigma_\lambda\circ \Ad_{\tw} \simeq \sigma_\lambda$ or equivalently, since $\sigma_\lambda\circ \Ad_{\tw} \simeq\sigma_{\tw^{-1}\lambda}$, $\tw \lambda=\lambda$.
		
	\end{proof}

	We can extend \eqref{eq2rep} to an action of $C^\infty(\Temp_{\ind}(G))$ on $\cC(\tG(F))$ as follows. Choose $\gamma\in \tG(F)$ and set, for every $f\in \cC(\tG(F))$, $f_\gamma(g)=f(g\gamma)$ ($g\in G(F)$). This function belongs to $\cC(G(F))$ and for $(z,f)\in C^\infty(\Temp_{\ind}(G))\times \cC(\tG(F))$, we define $z\star f\in \cC(\tG(F))$ by
	\begin{equation*}
		\displaystyle (z\star f)(g\gamma):=(z\star f_\gamma)(g),\; \mbox{ for } g\in G(F).
	\end{equation*}
	As the endomorphism $z\star$ commutes with right translations, this definition is easily seen to be independent on the choice of $\gamma$. Moreover, we have
	\begin{equation}
		\displaystyle (zz')\star f=z\star(z'\star f)
	\end{equation}
	and
	\begin{equation}
		\displaystyle \tpi(z\star f)=z(\pi) \tpi(f)
	\end{equation}
	for every $(z,z')\in C^\infty(\Temp_{\ind}(G))\times C^\infty(\Temp_{\ind}(G))$, $f\in \cC(\tG(F))$ and $\tpi\in E(\tG)$.

	\subsection{Orthogonal sets}\label{orthogonal sets}

	Let $(G,\tG)$ be a twisted space. We briefly recall the notion of  $(\tG,\tM)$-families from \cite{LW}.
	
	Let $\tM$ be a Levi subspace of $\tG$. Two parabolic subspaces $\tP, \widetilde{Q}\in \mathcal{P}(\tM)$ are said to be {\em adjacent} if the intersection $\Delta_{\tP}^\vee\cap -\Delta_{\widetilde{Q}}^\vee$ is a singleton $\{\alpha_{\tP,\tQ}^\vee\}$. If this is the case, the hyperplane $\{X\in i\mathcal{A}_{\tM}^*\mid \langle \alpha_{\tP,\tQ}^\vee,X\rangle=0 \}$ is called {\em the wall separating} $\tP$ and $\widetilde{Q}$.
	
	By definition {\em $(\tG,\tM)$-orthogonal set} is a family $\mathcal{Y}=(Y_{\tP})_{\tP\in \mathcal{P}(\tM)}$ of points in $\mathcal{A}_{\tM}$ such that for every adjacent parabolic subspaces $\tP, \widetilde{Q}\in \mathcal{P}(\tM)$, we have
	$$\displaystyle Y_{\tP}-Y_{\tQ}\in \mathbb{R}\alpha_{\tP,\tQ}^\vee$$
	where $\Delta_{\tP}^\vee\cap -\Delta_{\widetilde{Q}}^\vee=\{ \alpha_{\tP,\tQ}^\vee\}$. We further say that $\mathcal{Y}$ is {\em positive} if
	$$\displaystyle Y_{\tP}-Y_{\tQ}\in \mathbb{R}_{> 0}\alpha_{\tP,\tQ}^\vee$$
	for every pair of adjacent parabolic subspaces $\tP, \widetilde{Q}\in \mathcal{P}(\tM)$.
	
	For any $(\tG,\tM)$-orthogonal set $\mathcal{X}=(X_{\tP})_{\tP\in \mathcal{P}(\tM)}$,  we set
	\begin{equation*}
		\displaystyle d(\mathcal{X})=\min_{\tP\in \mathcal{P}(\tM)} \min_{\alpha\in \Delta_{\tP}} \alpha(X_{\tP}), \;\; N(\mathcal{X})=\max_{\tP\in \mathcal{P}(\tM)} \max_{\alpha\in \Delta_{\tP}} \lvert \alpha(X_{\tP})\rvert
	\end{equation*}
	that we shall call the {\em depth} and the {\em norm} of $\mathcal{X}$ respectively.
	
	Let $\mathcal{Y}=(Y_{\tP})_{\tP\in \mathcal{P}(\tM)}$ be a $(\tG,\tM)$-orthogonal set. For $\tQ=\tL U_Q\in \mathcal{F}(\tM)$, we denote by $Y_{\tQ}$ the projection to $\cA_{\tL}$ of $Y_{\tP}$ for any $\tP\in \mathcal{P}(\tM)$ such that $\tP\subset \tQ$ (this projection does not depend on the choice of $\tP$). To $\mathcal{Y}$ we associate functions $\Gamma^{\tQ}_{\tL}(.,\mathcal{Y})$ on $\cA_{\tL}^{\tQ}$ and complex numbers $v^{\tQ}_{\tL}(\mathcal{Y})\in \mathbb{C}$ for every $\tL\in \mathcal{L}(\tM)$ and $\tQ\in \mathcal{F}(\tL)$ as follows:
	$$\displaystyle \Gamma^{\tQ}_{\tL}(H,\mathcal{Y})=\sum_{\tP\in \mathcal{F}(\tL), \tP\subset \tQ} (-1)^{a^{\tQ}_{\tP}} \widehat{\tau}_{\tP}^{\tQ}(H-Y_{\tP}),\;\;\; H\in \cA_{\tL}^{\tQ},$$
	and
	$$\displaystyle v^{\tQ}_{\tL}(\mathcal{Y})=\int_{\cA_{\tL}^{\tQ}} \Gamma^{\tQ}_{\tL}(H,\mathcal{Y}) dH.$$
	Here $\widehat{\tau}_{\tP}^{\tQ}$ denotes the characteristic function of the cone in $\cA$ characterized by (where $\varpi_\alpha$ is the fundamental weight associated to $\alpha$)
	\begin{equation*}
		\displaystyle \widehat{\tau}_{\tP}^{\tQ}(H)=1 \Leftrightarrow \varpi_\alpha(H)>0,\; \forall \alpha\in \Delta_{\tP}^{\tQ}.
	\end{equation*}
	When $\tQ=\tG$, we will sometimes drop the superscript $\tQ$.	If $\mathcal{Y}$ is positive, $v^{\tQ}_{\tL}(\mathcal{Y})$ is simply the volume of the convex hull of $(Y_{\tP})_{\tP\in \mathcal{P}(\tL), \tP\subset \tQ}$. Once again, we will sometimes drop the superscript when $\tQ=\tG$. We will also use $\tau_{\tP}^{\tQ}$ to denote the characteristic function of the cone in $\cA$ characterized by
	\begin{equation*}
		\displaystyle \tau_{\tP}^{\tQ}(H)=1 \Leftrightarrow \alpha(H)>0,\; \forall \alpha\in \Delta_{\tP}^{\tQ}.
	\end{equation*}
	
	
	Let $K$  be a special compact subgroup of $G(F)$. Using the Iwasawa decomposition $G(F)=P(F)K$, for every parabolic subspace $\tP\subset \tG$, we can extend the homomorphism $H_{\tP}$ to a map $G(F)\to \cA_{\tP}$. Then, for every Levi subspace $\tM\subset \tG$ and $g\in G(F)$, the family $\mathcal{H}_{\tM}(g)=(-H_{\tP}(g))_{\tP\in \mathcal{P}(\tM)}$ is a positive $(\tG,\tM)$-orthogonal set and we define
	$$\displaystyle v_{\tM}^{\tQ}(g)=v_{\tM}^{\tQ}(\mathcal{H}_{\tM}(g)),\;\; \mbox{ for } \tQ \in \mathcal{F}(\tM).$$

	Let $\Lambda\subset \CA_{\tM,\mathbb{Q}}:=X_*(A_{\tM})\otimes_{\mathbb{Z}}\mathbb{Q}$ be a $\mathbb{Z}$-lattice. We say that a $(\tG,\tM)$-orthogonal set $\mathcal{Y}=(Y_{\tP})_{\tP\in \mathcal{P}(\tM)}$ is {\it $\Lambda$-rational} if for every $\tP\in \mathcal{P}(\tM)$, we have $Y_{\tP}\in \Lambda$ and we say that it is {\it rational} if it is $\Lambda$-rational for some lattice $\Lambda$. We denote by $\CC_\Lambda(\tG,\tM)$ (resp. $\CC_{\mathbb{Q}}(\tG,\tM)$) the set of all $\Lambda$-rational (resp. rational) $(\tG,\tM)$-orthogonal sets. Then, a function $\CY\in \CC_{\mathbb{Q}}(\tG,\tM)\mapsto f(\CY)\in \mathbb{C}$ is said to be a {\it unitary polynomial-exponential} if for every lattice $\Lambda\subset \CA_{\tM,\mathbb{Q}}$ we can find a family of polynomial functions $Q_{\mu,\Lambda,\tP}\in \mathbb{C}[\CA_{\tM}]$ for $\tP\in \mathcal{P}(\tM)$ and $\mu\in \widehat{\Lambda}:=\Hom(\Lambda,\mathbb{S}^1)$ that are almost all equal to $0$ and such that
	$$\displaystyle f(\CY)=\sum_{\tP\in \mathcal{P}(\tM)}\sum_{\mu\in \widehat{\Lambda}} Q_{\mu,\Lambda,\tP}(Y_{\tP})\mu(Y_{\tP})$$
	for every $\mathcal{Y}=(Y_{\tP})_{\tP\in \mathcal{P}(\tM)}\in \CC_\Lambda(\tG,\tM)$. Moreover, we say that a unitary polynomial-exponential function $f$ is {\it of degree at most $r$} if the polynomials $Q_{\mu,\Lambda,\tP}$ are of degree at most $r$ for every lattice $\Lambda\subset \CA_{\tM,\mathbb{Q}}$, $\tP\in \mathcal{P}(\tM)$ and $\mu\in \widehat{\Lambda}$.

	\subsection{Weighted orbital integrals}\label{Section WOI}
	
	Let $\tM$ be a Levi subspace of $\tG$, $\gamma\in \tM(F)\cap \tG_{rs}(F)$ and $\tQ\in \mathcal{F}(\tM)$. For $f\in \cC(\tG(F))$, we define the {\em twisted weighted orbital integral}
	\begin{equation*}
		\displaystyle \Phi_{\tM}^{\tQ}(\gamma,f)=\int_{G_\gamma(F)\backslash G(F)} f(g^{-1}\gamma g) v_{\tM}^{\tQ}(g)dg
	\end{equation*}
	as well as its normalized version
	\begin{equation*}
		\displaystyle J_{\tM}^{\tQ}(\gamma,f)=D^{\tG}(\gamma)^{1/2} \Phi_{\tM}^{\tQ}(\gamma,f).
	\end{equation*}

	The above integral is absolutely convergent. More precisely, for $\tT\subset \tM$ a maximal twisted torus, we claim:
	\begin{num}
		\item\label{eq1 WOI} There exist $p>0$ and, for every $d>0$, a continuous semi-norm $\nu_d$ on $\cC(\tG(F))$ such that
		$$\displaystyle \left\lvert J_{\tM}^{\tQ}(\gamma,f)\right\rvert\leqslant \nu_d(f) (1+\lvert \log D^{\tG}(\gamma)\rvert)^p \sigma_{\tT_{/\theta}}(\gamma)^{-d}$$
		for every $\gamma\in \tT_{\reg}(F)$ and $f\in \cC(\tG(F))$.
	\end{num}
	Indeed, there exists $p>0$ such that $v_{\tM}^{\tQ}(g)\ll \sigma_G(g)^{p}$ for $g\in G(F)$. As $v_{\tM}^{\tQ}$ is left invariant by $T(F)$, by Lemma \ref{lem1 estimates} this implies $v_{\tM}^{\tQ}(g)\ll (1+\lvert \log D^{\tG}(\gamma)\rvert)^p \sigma_{\tG}(g^{-1}\gamma g)^{p}$ for $(g,\gamma)\in G(F)\times \tT_{\reg}(F)$. The claim is now a straightforward consequence of Proposition \ref{prop conv WOI}.
	
	\vspace{2mm}

	Now consider the case where $\tG=\tH\times \tH$ where $\tH$ is a connected twisted reductive space over $F$ (with underlying reductive group $H$). Let $\tM_H$ be a Levi subspace of $\tH$. Then, $\tM=\tM_H\times \tM_H$ is a Levi subspace of $\tG$. Let $\gamma\in \tM_H(F)\cap \tH_{\rs}(F)$ and $f_1,f_2\in \cC(\tH(F))$. We set (following \cite[Chap. 1, \S 4.8]{WalFTLtordue})
	
	\begin{equation*}
		\displaystyle J^{\tH}_{\tM_H}(\gamma,f_1,f_2)=\int_{H_\gamma(F)\backslash H(F)\times H_\gamma(F)\backslash H(F)} f_1(x^{-1}\gamma x) f_2(y^{-1}\gamma y) v_{\tM_H}^{\tH}(x,y) dxdy
	\end{equation*}
	where the ``weight'' $v_{\tM_H}^{\tH}(x,y)$ is the volume associated to the positive $(\tH,\tM_H)$-orthogonal set
	$$\displaystyle \tP_H\in \mathcal{P}(\tM_H)\mapsto H_{\widetilde{\overline{P}}_H}(y)-H_{\tP_H}(x).$$
	(Here $\widetilde{\overline{P}}_H$ denotes the unique parabolic subspace opposite to $\tP_H$ such that $\widetilde{\overline{P}}_H\cap \tP_H=\tM$.)
	
	By the descent formulas of \cite[Chap. 1, lemme 5.4]{WalFTLtordue}, we have
	\begin{equation}\label{eq1WOI}
		\displaystyle J_{\tM_H}^{\tH}(\gamma,f_1,f_2)=\sum_{\tL_1,\tL_2\in \mathcal{L}(\tM_H)} d_{\tM_H}^{\tH}(\tL_1,\tL_2) J_{\tM_H}^{\widetilde{\overline{Q}}_1}(\gamma,f_1) J_{\tM_H}^{\widetilde{Q}_2}(\gamma,f_2)
	\end{equation}
	where $\tQ_1$, $\tQ_2$ are certain parabolic subspaces in $\mathcal{P}(\tL_1)$, $\mathcal{P}(\tL_2)$ respectively and $d_{\tM_H}^{\tH}(\tL_1,\tL_2)$ is a certain real numbers which is zero unless $\cA_{\tM_H}^{\tH}=\cA_{\tL_1}^{\tH}\oplus \cA_{\tL_2}^{\tH}$ and moreover $d_{\tM_H}^{\tH}(\tH,\tM_H)=1$.

	\subsection{Twisted weighted characters}\label{Section WC}
	
	Let $\tM$ be a Levi subspace of $\tG$, $\widetilde{R}\in \mathcal{F}(\tM)$ and $\tpi$ be a tempered representation of $\tM(F)$. First assume that $\tpi$ is in ``general position'' (more precisely, this means that $\pi$ is in some open-dense subset of the family $\{\pi_{\lambda}\mid \lambda\in i\cA_{\tM}^* \}$). Then, we define as in \cite[Chap.1, \S 2.7]{WalFTLtordue} a {\em weighted character}
	\begin{equation*}
		\displaystyle f\in \cC(\tG(F))\mapsto J_{\tM}^{\widetilde{R}}(\tpi,f):= \Tr(\mathcal{M}_{\tM}^{\widetilde{R}}(\pi)I_{\tP}^{\tG}(\tpi,f))
	\end{equation*}
	where $\tP$ is any chosen parabolic subspace in $\mathcal{P}(\tM)$ (the distribution $J_{\tM}^{\widetilde{R}}(\tpi,.)$ does not depend on this choice) and $\mathcal{M}_{\tM}^{\widetilde{R}}(\pi)$ is the operator on (the space of) $I_{P}^{G}(\pi)$ associated to the $(\tG,\tM)$-family of operators $(\mathcal{M}(\pi;\Lambda,\tQ))_{\tQ\in \mathcal{P}(\tM)}$ defined as in {\em loc. cit.}. Similarly, for $f_1,f_2\in \cC(G(F))$ we set
	\[\begin{aligned}
		\displaystyle J^{\tG}_{\tM}(\tpi,f_1,f_2)=\Tr\left(\mathcal{M}_{\tM}^{\tG}(\pi^\vee\otimes \pi) I_{\tP\times \tP}^{\tG\times \tG}(\tpi^\vee\otimes \tpi, f_1\otimes f_2)\right),
	\end{aligned}\]
	where this time the operator $\mathcal{M}_{\tM}^{\tG}(\pi^\vee\otimes \pi)$ is associated to the $(\tG,\tM)$-family
	$$\displaystyle \tQ\in \mathcal{P}(\tM) \mapsto \mathcal{M}(\pi^\vee\otimes \pi; \Lambda, \tQ)=\mathcal{M}(\pi^\vee; \Lambda, \widetilde{\overline{Q}})\otimes \mathcal{M}(\pi; \Lambda, \tQ)$$
	of operators on $I_{P\times P}^{G\times G}(\pi^\vee \otimes \pi)$. The genericity assumption on $\tpi$ is necessary for the above $(\tG,\tM)$-families to be well-defined. However, the definitions of $J_{\tM}^{\widetilde{R}}(\tpi,f)$ and $J^{\tG}_{\tM}(\tpi,f_1,f_2)$ extend to every tempered representation $\tpi$ thanks to the following property (see \cite[Chap. 1, proposition 2.7]{WalFTLtordue}):
	\begin{num}
		\item\label{eq1 WC} The operator valued functions $\lambda\mapsto \mathcal{M}_{\tM}^{\widetilde{R}}(\pi_\lambda)$ and $\lambda\mapsto \mathcal{M}_{\tM}^{\tG}(\pi^\vee_{\lambda}\otimes \pi_{\lambda})$ , a priori only well-defined on an dense open subset of $i\cA_{\tM}^*$, extend to smooth functions on all of $i\cA_{\tM}^*$.
	\end{num} 
	
	Finally, by the descent formula of \cite[Chap 1, lemme 5.4]{WalFTLtordue}, for every $f_1,f_2\in \cC(\tG(F))$ we have
	\begin{equation}\label{eq2 WC}
		\displaystyle J_{\tM}^{\tG}(\tpi,f_1,f_2)=\sum_{\tL_1,\tL_2\in \mathcal{L}(\tM)} d_{\tM}^{\tG}(\tL_1,\tL_2) J_{\tM}^{\widetilde{\overline{Q}}_1}(\tpi,f_1) J_{\tM}^{\widetilde{Q}_2}(\tpi,f_2)
	\end{equation}
	where $\tQ_1$, $\tQ_2$ and $d_{\tM}^{\tG}(\tL_1,\tL_2)$ are as in \eqref{eq1WOI}.

	\subsection{Twisted strongly cuspidal functions}\label{Sect strongly cuspi}
	
	We say that a function $f\in \cC(\widetilde{G}(F))$ is {\em strongly cuspidal}, if for every  parabolic subspace $\widetilde{P}=\widetilde{M}U_P$ of $\tG$ and $x\in G(F)$, the function defined by
	\begin{equation}\label{eq1 scf}
		\displaystyle {}^xf_{(\tP)}(\tm):=\delta_{\tP}(\tm)^{1/2}\int_{U_P(F)} f(x^{-1}\widetilde{m}ux) du,\; \mbox{ for } \widetilde{m}\in \widetilde{M}(F),
	\end{equation}
	is identically zero. By a change of variable, this last condition is equivalent to
	\begin{equation}\label{eq1bis scf}
		\displaystyle \int_{U_P(F)} f(x^{-1}u^{-1}\widetilde{m}ux) du=0,\; \mbox{ for every } \widetilde{m}\in \widetilde{M}(F)\cap \tG_{\rs}(F) \mbox{ and } x\in G(F).
	\end{equation}
	We denote by $\CC_{\scusp}(\tG(F))\subseteq \CC(\tG(F))$ the subspace of strongly cuspidal functions.
	
	Let $f\in \cC_{\scusp}(\widetilde{G}(F))$. Let $\tM$ be a Levi subspace and $\gamma\in \tM(F)\cap \tG_{rs}(F)$. For $\tQ=\tL U_Q\in \mathcal{F}(\tM)$, the weight $v_{\tM}^{\tQ}$ is left invariant by $U_Q(F)$. Hence, by \eqref{eq1bis scf}, we have
	\begin{equation}\label{eq2 scf}
		\displaystyle J_{\tM}^{\tQ}(\gamma,f)=0 \mbox{ unless } \tQ=\tG.
	\end{equation}
	Then, by the same argument as for \cite[lemme 5.2]{WalGGPI}, it follows that the weighted orbital integral $\Phi_{\tM}(\gamma,f)$ does not depend on the choice of $K$.
	
	We define a function $\Theta_f$ on $\tG_{rs}(F)$ by
	\begin{equation*}
		\displaystyle \Theta_f(\gamma)=(-1)^{a_{G_\gamma}-a_{\tG}}\Phi_{\tM(\gamma)}^{\tG}(\gamma,f),\;\;\; \gamma\in \tG_{\rs}(F),
	\end{equation*}
	where $\tM(\gamma)$ stands for the centralizer of $A_{G_\gamma}$ in $\tG$ (it is the minimal Levi subspace containing $\gamma$), $a_{G_\gamma}=\dim(A_{G_\gamma})$ and $a_{\tG}=\dim(A_{\tG})$. It is proved in \cite[proposition 1.7]{WalGGPIV} that if $f$ is compactly supported then $\Theta_f$ is a {\em quasi-character} on $\tG(F)$ in the sense of Section \ref{Section qc}. We extend this result to every strongly cuspidal function $f\in \cC_{\scusp}(\tG(F))$ in Section \ref{Section TLTF strongly cuspidal} (see Corollary \ref{cor TLTF srongly cuspidal}). For ease of notation, for every $x\in \tG_{\ssi}(F)$, we set
	$$\displaystyle c_{f,\cO}(x)=c_{\Theta_f,\cO}(x),\;\;\; \cO\in \Nil(\mathfrak{g}_x),$$
	for the coefficients of the germ expansion of $\Theta_f$ near $x$.
	
	Let again $f\in \cC_{\scusp}(\tG(F))$. Let $\tM\subset \tG$ be a twisted Levi subspace, $\tpi$ be a tempered representation of $\tM(F)$ and $\tQ\in \mathcal{F}(\tM)$. By \cite[lemme 1.13]{WalGGPIV}, we have
	\begin{equation}\label{eq3 scf}
		\displaystyle J_{\tM}^{\tQ}(\tpi,f)=0 \mbox{ unless } \tQ=\tG.
	\end{equation}
	Still by \cite[lemme 1.13]{WalGGPIV}, we also have
	\begin{equation}\label{eq4 scf}
		\displaystyle J_{\tM}^{\tG}(\tpi,f)=0 \mbox{ if } \tpi \mbox{ is properly parabolically induced (e.g. if } \tpi\in E(\tM)\setminus E_{\elli}(\tM)).
	\end{equation}
	On the other hand, for $\tpi\in E_{\elli}(\tM)$, we set
	\begin{equation*}
		\displaystyle \widehat{\Theta}_f(\tpi)=(-1)^{a_{\tM}-a_{\tG}}J_{\tM}^{\tG}(\tpi,f), \mbox{ for } \tpi\in E_{\elli}(\tM).
	\end{equation*}
	
	Recall that in Section \ref{Section twrep}, we have defined an action of $C^\infty(\Temp_{\ind}(G))$ on $\cC(\tG(F))$. We denote by $C^\infty(\Temp_{\ind}(G))^\theta$ the subspace of $\theta$-invariant functions in $C^\infty(\Temp_{\ind}(G))$.
	
	\begin{lem}\label{lem strongly cuspidal function}
		Let $f\in \cC(\tG(F))$ and $z\in C^\infty(\Temp_{\ind}(\tG))^\theta$. Then, if $f$ is strongly cuspidal so is $z\star f$.
	\end{lem}
	
	\begin{proof}
		Let $\tP=\tM U_P$ be a proper parabolic subspace of $\tG$. By \eqref{eq1 Xi function}, for every $x\in G(F)$ the function ${}^xf_{(\tP)}$ defined by the integral \eqref{eq1 scf} belongs to $\cC(\tM(F))$. For $\tsigma$ a tempered representation of $\tM(F)$, we set
		$$\displaystyle I_{\tP}^{\tG}(\tsigma,f)(x,y)= \int_{\tM(F)} \delta_{\tP}(\tm)^{1/2} \int_{U_P(F)}f(x^{-1}\tm u y) du\, \tsigma(\tm) d\tm,\mbox{ for } (x,y)\in G(F)\times G(F).$$
		This operator-valued function is the kernel of the operator $I_{\tP}^{\tG}(\tsigma,f)$ in the sense that
		$$\displaystyle (I_{\tP}^{\tG}(\tsigma,f)e)(x)=\int_{P(F)\backslash G(F)} I_{\tP}^{\tG}(\tsigma,f)(x,y)e(y)dy$$
		for every $e\in I_P^G(\sigma)$ and $x\in G(F)$. Note that
		\begin{equation}\label{eq1 sp char scf}
			\displaystyle I_{\tP}^{\tG}(\tsigma,f)(x,x)=\tsigma({}^xf_{(\tP)}),\mbox{ for } x\in G(F).
		\end{equation}
		Let $d\geqslant 1$ be the order of the outer automorphism $\theta=\theta_{\tM}$ of $M(F)$ and set, for $\sigma\in \Temp_{\ind}(M)$,
		$$\displaystyle \sigma\langle \theta\rangle:=\sigma\oplus \sigma^{\theta}\oplus\ldots\oplus \sigma^{\theta^{d-1}}.$$
		It is clear that $\sigma\langle \theta\rangle$ extends to a twisted representation of $\tM(F)$ and we will denote by $\widetilde{\sigma\langle \theta\rangle}$ one such extension. It is well-known, and this follows e.g. from the Harish-Chandra-Plancherel formula \cite{WalPlanch}, that a function $f'\in \cC(M(F))$ is zero if and only if $\sigma(f')=0$ for every $\sigma\in \Temp_{\ind}(M)$. This implies a similar equivalence for $\tM$: a function $f'\in \cC(\tM(F))$ is zero if and only if $\widetilde{\sigma\langle \theta\rangle}(f')=0$ for every $\sigma\in \Temp_{\ind}(M)$. Therefore, from \eqref{eq1 sp char scf} and the definition of a strongly cuspidal function, we see that $f$ is strongly cuspidal if and only if $I_{\tP}^{\tG}(\widetilde{\sigma\langle \theta\rangle},f)(x,x)=0$ for every proper parabolic subspace $\tP=\tM U_P$, every $\sigma\in \Temp_{\ind}(M)$ and every $x\in G(F)$. The lemma is a direct consequence of this characterization since for $z\in C^\infty(\Temp_{\ind}(\tG))^\theta$ and every $\tP$ and $\sigma$ as before, since $z(I_M^G(\sigma^{\theta^i}))=z(I_M^G(\sigma))$ for all $i$ (by $\theta$-invariance of $z$), we have
		$$\displaystyle I_{\tP}^{\tG}(\widetilde{\sigma\langle \theta\rangle},z\star f)=z(I_M^G(\sigma))I_{\tP}^{\tG}(\widetilde{\sigma\langle \theta\rangle},f).$$
	\end{proof}

	\subsection{Twisted local trace formula for strongly cuspidal functions}\label{Section TLTF strongly cuspidal}

	The twisted local trace formula of \cite[Chap. 1, th\'eor\`eme 5.1]{WalFTLtordue} is an equality of distributions
	\begin{equation}\label{eq LTTF}
		\displaystyle J^{\tG}_{\spec}(f_1,f_2)=J^{\tG}_{\geom}(f_1,f_2)
	\end{equation}
	where $f_1,f_2\in C_c^\infty(\tG(F))$ and
	\begin{equation}\label{eq2 LTTF}
		\displaystyle J^{\tG}_{\spec}(f_1,f_2)= \sum_{\tM\in \mathcal{L}(\tM_{\min})} \lvert \widetilde{W}^M\rvert \lvert \widetilde{W}^G\rvert^{-1} (-1)^{a_{\tM}-a_{\tG}}  \int_{E_{\disc}(\tM)} \iota(\tau) J^{\tG}_{\tM}(\tpi_{\tau},f_1,f_2) d\tau,
	\end{equation}
	\begin{equation}\label{eq3 LTTF}
		\displaystyle J^{\tG}_{geom}(f_1,f_2)=\sum_{\tM\in \mathcal{L}(\tM_{\min})} \lvert \widetilde{W}^M\rvert \lvert \widetilde{W}^G\rvert^{-1} (-1)^{a_{\tM}-a_{\tG}} \int_{\Gamma_{\elli}(\tM)} J_{\tM}^{\tG}(\gamma,f_1,f_2) d\gamma.
	\end{equation}
	We refer the reader to Section \ref{Section twrep} for the definition of $\iota(\tau)$ as well as of the measure on $E_{\disc}(\tM)$ and to Sections \ref{Section WOI} and \ref{Section WC} for the definitions of $J_{\tM}^{\tG}(\gamma,f_1,f_2)$ and $J^{\tG}_{\tM}(\tpi_{\tau},f_1,f_2)$ respectively.
	
	\begin{rmk}
		Despite the notation, the distributions $J^{\tG}_{\spec}$ and $J^{\tG}_{\geom}$ depend on the choice of the pair $(M_{\min},K)$ (at least up to conjugacy). They also depend, incidentally, on the choice of the Haar measure on $G(F)$.
	\end{rmk}
	
	\begin{prop}\label{prop LTF HCS}
		The expressions \eqref{eq2 LTTF} and \eqref{eq3 LTTF} are both absolutely convergent for $(f_1,f_2)\in \cC(\tG(F))^2$ and they define continuous bilinear forms on $\mathcal{C}(\tG(F))\times \mathcal{C}(\tG(F))$. In particular, the identity \eqref{eq LTTF} extends by continuity to all $f_1,f_2\in \cC(\tG(F))$.
	\end{prop}
	
	\begin{proof}
		The same argument as in the non-twisted case \cite[p.189]{Art94} applies here noticing that \eqref{eq1 WOI} gives the required twisted analog of the estimates (5.7) of {\em loc. cit.}.
	\end{proof}
	
	Let $f_1,f_2\in \mathcal{C}(\tG(F))$ and assume that $f_1$ is strongly cuspidal. By the descent formulas \eqref{eq1WOI}, \eqref{eq2 WC} as well as the vanishing \eqref{eq2 scf}, \eqref{eq3 scf}, \eqref{eq4 scf} we then have
	\begin{equation*}
		\displaystyle J_{\tM}^{\tG}(\gamma,f_1,f_2)=(-1)^{a_{\tM}-a_{\tG}}D^{\tG}(\gamma)^{1/2}\Theta_{f_1}(\gamma) J_{\tG}(\gamma,f_2)
	\end{equation*}
	and
	\begin{equation*}
		\displaystyle J_{\tM}^{\tG}(\tpi_\tau,f_1,f_2)=\left\{\begin{array}{ll}
			(-1)^{a_{\tM}-a_{\tG}}\widehat{\Theta}_{f_1}(\tpi^\vee_\tau) J_{\tG}(\tpi_\tau,f_2) \mbox{ if } \tau\in E_{\elli}(\tM); \\
			\\
			0 \mbox{ otherwise,}
		\end{array} \right.
	\end{equation*}
	for every $\gamma\in \Gamma_{\elli}(\tM)$ and $\tau\in E_{\disc}(\tM)$. Thus, in this case the distributions $J^{\tG}_{\spec}$ and $J^{\tG}_{\geom}$ can be rewritten as
	\[\begin{aligned}
		\displaystyle J^{\tG}_{\spec}(f_1,f_2) & = \sum_{\tM\in \mathcal{L}(\tM_{\min})} \lvert \widetilde{W}^M\rvert \lvert \widetilde{W}^G\rvert^{-1}  \int_{E_{\elli}(\tM)} D(\tau) \widehat{\Theta}_{f_1}(\tpi^\vee_\tau) J_{\tG}(\tpi_\tau,f_2) d\tau 
	\end{aligned}\]
	and
	\[\begin{aligned}
		\displaystyle J^{\tG}_{geom}(f_1,f_2) & =\sum_{\tM\in \mathcal{L}(\tM_{\min})} \lvert \widetilde{W}^M\rvert \lvert \widetilde{W}^G\rvert^{-1}  \int_{\Gamma_{\elli}(\tM)} D^{\tG}(\gamma)^{1/2}\Theta_{f_1}(\gamma) J_{\tG}(\gamma,f_2) d\gamma \\
		& =\int_{\tG(F)} \Theta_{f_1}(\gamma) f_2(\gamma) d\gamma
	\end{aligned}\]
	where the last equality follows from the Weyl integration formula \eqref{Weyl integration formula}. Moreover, by definition of the Harish-Chandra characters $\Theta_{\tpi_\tau}$, the spectral side $J^{\tG}_{\spec}(f_1,f_2)$ can be further rewritten as
	\[\begin{aligned}
		\displaystyle J^{\tG}_{\spec}(f_1,f_2) & = \sum_{\tM\in \mathcal{L}(\tM_{\min})} \lvert \widetilde{W}^M\rvert \lvert \widetilde{W}^G\rvert^{-1}  \int_{E_{\elli}(\tM)} \int_{\tG(F)} D(\tau) \widehat{\Theta}_{f_1}(\tpi^\vee_\tau) \Theta_{\tpi_\tau}(\gamma) f_2(\gamma)d\gamma d\tau.
	\end{aligned}\]
	The above expression being absolutely convergent (note that, by \eqref{eq3 rep}, the support of $\tau\mapsto \widehat{\Theta}_{f_1}(\tpi^\vee_\tau)$ in $E_{\elli}(\tM)$ is contained in a finite union of orbits under the action of $i\cA_{\tM}^*$), from \eqref{eq LTTF} we get the identity
	\[\begin{aligned}
		\displaystyle \int_{\tG(F)} \Theta_{f_1}(\gamma) f_2(\gamma) d\gamma=\int_{\tG(F)}\sum_{\tM\in \mathcal{L}(\tM_{\min})} \lvert \widetilde{W}^M\rvert \lvert \widetilde{W}^G\rvert^{-1}  \int_{E_{\elli}(\tM)}  D(\tau) \widehat{\Theta}_{f_1}(\tpi^\vee_\tau) \Theta_{\tpi_\tau}(\gamma)  d\tau f_2(\gamma)d\gamma
	\end{aligned}\]
	for every $f_1,f_2\in \cC(\tG(F))$ with $f_1$ strongly cuspidal. Fixing $f_1$ and varying $f_2$, we deduce:
	
	\begin{prop}\label{prop LTTF strongly cuspidal}
		Let $f\in \mathcal{C}(\tG(F))$ be a strongly cuspidal function. Then, for every $\gamma\in \tG_{\rs}(F)$, we have
		\begin{equation}
			\displaystyle \Theta_{f}(\gamma)=\sum_{\tM\in \mathcal{L}(\tM_{\min})} \lvert \widetilde{W}^M\rvert \lvert \widetilde{W}^G\rvert^{-1}  \int_{E_{\elli}(\tM)}  D(\tau) \widehat{\Theta}_{f}(\tpi^\vee_\tau) \Theta_{\tpi_\tau}(\gamma)  d\tau
		\end{equation}
		where the right hand side is absolutely convergent.
	\end{prop}
	
	As a corollary, we can now show the following extension of \cite[proposition 1.7]{WalGGPIV}.
	
	\begin{cor}\label{cor TLTF srongly cuspidal}
		For $f\in \cC(\tG(F))$ strongly cuspidal, the function $\Theta_f$ is a quasi-character.
	\end{cor}
	
	\begin{proof}
		This is a direct consequence of Proposition \ref{prop LTTF strongly cuspidal} combined with the following facts (valid for every $\tM\in \mathcal{L}(\tM_{\min})$):
		\begin{itemize}
			\item[$\bullet$] For every $\tau\in E_{\elli}(\tM)$, $\Theta_{\tpi_\tau}$ is a quasi-character \cite[Theorem 3]{Clo};
			
			\item[$\bullet$] The function $\tau\in E_{\elli}(\tM)\mapsto \widehat{\Theta}_{f}(\tpi^\vee_\tau)$ is supported on a finite number of orbits under the action of $i\cA_{\tM}^*$ (see \eqref{eq3 rep});
			
			\item[$\bullet$] For every $i\cA_{\tM}^*$-orbit $\Omega\subset E_{\elli}(\tM)$ and compact subset $\mathcal{K}\subset \tG(F)$, the vector space spanned by the restrictions
			$$\displaystyle \{\Theta_{\tpi_\tau}\mid_{\mathcal{K}_{\rs}}\mid \tau\in \Omega \},$$
			where $\mathcal{K}_{\rs}:=\mathcal{K}\cap \tG_{\rs}(F)$, is of finite dimension (this follows e.g. from the induction formula \cite[lemme 1.12]{WalGGPIV}).
			
		\end{itemize}
	\end{proof}
	
	
	
	\subsection{Spectral localization of strongly cuspidal functions}
	
	Let $f\in \cC(\tG(F))$ and choose a base-point $\gamma_0\in \tG(F)$. Put $f_{\gamma_0}(g)=f(g\gamma_0)$ for every $g\in G(F)$. Note that $f_{\gamma_0}\in \cC(G(F))$. We define the {\em spectral support} of $f$, henceforth denoted by $\Specsupp(f)$ to be the support of the operator-valued function
	\begin{equation*}
		\displaystyle \Temp_{\ind}(G)\ni \pi\mapsto \pi(f_{\gamma_0})\in \End(V_\pi).
	\end{equation*}
	Note that $\Specsupp(f)$ does not depend on the choice of $\gamma_0$: changing the base-point replaces $f_{\gamma_0}$ by one of its right translates which acts non trivially on the same tempered representations as $f_{\gamma_0}$.
	
	\begin{prop}\label{prop pesudo coefficient}
		Let $\tau\in E_2(\tG)$ (see Section \ref{Section twrep} for the definition of $E_2(\tG)$) and $\omega$ be a compact neighborhood of $\pi_\tau$ in $\Temp_{\ind}(G)$ (see Section \ref{Sect rep of G} for the topology on $\Temp_{\ind}(G)$). Then, there exists a strongly cuspidal function $f\in \cC(\tG(F))$ such that
		\begin{equation}\label{eq-2 pseudocoeff}
			\displaystyle \Specsupp(f)\subset \omega
		\end{equation}
		and for every $\tau'\in E(\tG)$ we have
		\begin{equation}\label{eq-1 pseudocoeff}
			\displaystyle \Tr \tpi_{\tau'}(f)=\left\{ \begin{array}{ll}
				0 \mbox{ if } \tau'\neq \lambda\cdot \tau \mbox{ for every } \lambda\in i\cA_{\tG}^*, \\ \\
				1 \mbox{ if } \tau'=\tau.
			\end{array}
			\right.
		\end{equation}
		Moreover, if $f\in \cC(\tG(F))$ is such a strongly cuspidal function, we have
		\begin{equation}\label{eq0 pseudocoeff}
			\displaystyle \Theta_f(\gamma)=\lvert \Stab(i\cA_{\tG,F}^*,\tau)\rvert^{-1} D(\tau)\int_{i\cA_{\tG,F}^*} \Tr \tpi_{\lambda\cdot \tau}(f) \Theta_{\tpi_{\lambda\cdot \tau}^\vee}(\gamma) d\lambda
		\end{equation}
		for every $\gamma\in \tG_{\rs}(F)$ where $\Stab(i\cA_{\tG,F}^*,\tau)$ stands for the stabilizer of $\tau$ in $i\cA_{\tG,F}^*$ (for the action by twisting).
	\end{prop}
	
	\begin{proof}
		For simplicity of notation, let us set $\tpi=\tpi_\tau$, $\pi=\pi_\tau$ as well as $\tpi_\lambda=\tpi_{\lambda\cdot \tau}$ and $\pi_\lambda=\pi_{\lambda\cdot \tau}=(\pi_\tau)_\lambda$ for every $\lambda\in i\cA_{\tG}^*$. By Lemma \ref{lem rep}, up to shrinking $\omega$ we may assume that it is $\theta$-stable and that
		\begin{equation}\label{eq1 pseudocoeff}
			\displaystyle \omega\cap \Temp_{\ind}(G)^\theta\subset\{\pi_\lambda\mid \lambda\in i\cA_{\tG}^* \}.
		\end{equation}
		Let $S$ be the finite set of $\tau'\in E(\tG)$ such that $\pi_{\tau'}=\pi_{\tau}$. By definition of $E_2(\tG)$, we have $S\subset E_2(\tG)$. Let $S_0\subset S$ be a subset such that for every $\tau'\in S$ there exists an unique $\tau'_0\in S_0$ as well as $\lambda\in i\cA_{\tG}^*$ (not necessarily unique) such that $\tau'= \lambda\cdot \tau'_0$. We may and will assume that $\tau\in S_0$. Moreover, by \eqref{eq1 pseudocoeff}, we have:
		\begin{num}
			\item\label{eq1bis pseudocoeff} for $\tau'\in E(\tG)$ if $\pi_{\tau'}\in \omega$ then there exist $\lambda\in i\cA_{\tG}^*$ and $\tau'_0\in S_0$ such that $\tau'=\lambda\cdot \tau'_0$.
		\end{num}

		Let $G(F)^1$ be the kernel of the homomorphism $H_G:G(F)\to \mathcal{A}_G$. By the orthogonality relations \cite[th\'eor\`eme 7.3]{WalFTLtordue} between elliptic twisted characters, the restrictions of the twisted characters $\Theta_{\tpi_{\tau'}}$, for $\tau'\in S$, to the elliptic locus $\tG(F)_{\elli}$ are linearly independent. More precisely, fixing $\gamma\in \tG(F)_{\elli}$ and since elements of $S_0$ all have different orbits under $i\cA_{\tG}^*$, the restrictions of the twisted characters $\Theta_{\tpi_{\tau'}}$, for $\tau'\in S_0$, to $\tG(F)_{\elli}\cap G(F)^1\gamma$ are linearly independent. Thus, we can find a function $f_0\in C_c^\infty(\tG(F))$ supported in $\tG(F)_{\elli}\cap G(F)^1\gamma$ such that
		\begin{equation}\label{eq2 pseudocoeff}
			\displaystyle \Tr \tpi_{\tau'}(f_0)=\left\{\begin{array}{ll}
				0 \mbox{ if } \tau'\neq \tau, \\ \\
				1 \mbox{ if } \tau'=\tau
			\end{array} \right.
		\end{equation}
		for every $\tau'\in S_0$. Note that, since $f_0$ is supported in $\tG(F)_{\elli}$, it is a strongly cuspidal function. Moreover, since $f_0$ is supported in a unique coset modulo $G(F)^1$, for every $\tau'\in E(\tG)$ and $\lambda\in i\cA_{\tG}^*$, $\Tr \tpi_{\lambda\cdot \tau'}(f_0)$ is equal (up to a non-zero multiplicative constant which depends on how we normalized $\tpi_{\lambda\cdot \tau'}$) to $\Tr \tpi_{\tau'}(f_0)$. In particular, by \eqref{eq2 pseudocoeff}, we also have
		\begin{equation}\label{eq2bis pseudocoeff}
			\displaystyle \Tr \tpi_{\lambda\cdot \tau'}(f_0)=0
		\end{equation}
		for every $\tau'\in S_0\setminus \{\tau \}$ and $\lambda\in i\cA_{\tG}^*$.
		
		Let now $z\in C^\infty(\Temp_{\ind}(G))^\theta$ be a $\theta$-invariant $C^\infty$ function on $\Temp_{\ind}(G)$ which is supported in $\omega$ and such that $z(\pi)=1$ (such a function certainly exists). Using the action of $C^\infty(\Temp_{\ind}(G))^\theta$ on $\cC(\tG(F))$ defined in Section \ref{Section twrep}, we set $f=z\star f_0\in \cC(\tG(F))$. By Lemma \ref{lem strongly cuspidal function}, $f$ is strongly cuspidal. On the other hand, by the spectral characterization of the action of $C^\infty(\Temp_{\ind}(G))^\theta$ on $\cC(\tG(F))$, $f$ clearly satisfies condition \eqref{eq-2 pseudocoeff}. Similarly, \eqref{eq-1 pseudocoeff} follows from the combination of \eqref{eq1bis pseudocoeff}, \eqref{eq2 pseudocoeff} and \eqref{eq2bis pseudocoeff}. Finally, the equality \eqref{eq0 pseudocoeff} is an immediate consequence of Proposition \ref{prop LTTF strongly cuspidal}, remembering that the restriction of the measure on $E_{\elli}(\tG)$ to the orbit $\{\lambda\cdot \tau\mid \lambda\in i\cA_{\tG}^* \}$ is equal to $\lvert \Stab(i\cA_{\tG,F}^*,\tau)\rvert^{-1}$ times the pushforward of the measure on $i\cA_{\tG,F}^*$ by the map
		$$\displaystyle \lambda\in i\cA_{\tG,F}^*\mapsto \lambda\cdot \tau.$$
	\end{proof}
	
	\section{Spherical spaces}

	\subsection{Coregular varieties}\label{section coregular varieties}
	
	Let $G$ be a connected reductive group over $F$ and $H\subset G$ be a closed subgroup. We let $X=H\backslash G$ be the corresponding homogenous variety. We let $TX$, $T^*X$ be the tangent and cotangent bundles of $X$ respectively. Both are naturally equipped with a right action of $G$.
	
	Let $\CB$ be the flag variety of $G$. Recall that the variety $X$ is called {\it spherical} if $H$ has an open orbit in $\CB$ or, equivalently, if $G$ has an open orbit in $X\times \CB$ for the diagonal action.
	
	In the proposition below, by the {\it generic stabilizer} of a $G$-variety $Y$ we mean a conjugacy class of closed subgroups $S\subset G$ such that for some dense open subset $U\subset Y$, the stabilizer of every $y\in U$ is conjugated to $S$. Generic stabilizers do not always exist but they do in the cases considered in the proposition below by the references cited in the proof, namely \cite{Knop93} and \cite{Knop94}.
	
	\begin{prop}
		Assume that $X=H\backslash G$ is quasi-affine and that $H$ is connected. Then, the following assertions are equivalent:
		\begin{enumerate}[(i)]
			\item The generic stabilizer of $T^*X$ contains regular elements;
		\end{enumerate}
		\begin{enumerate}[(i')]
			\item The generic stabilizer of $T^*X$ contains regular semisimple elements;
		\end{enumerate}
		\begin{enumerate}[(i)]
			\setcounter{enumi}{1}
			\item The generic stabilizer of $X\times \CB$ contains regular elements;
		\end{enumerate}
		\begin{enumerate}[(i')]
			\setcounter{enumi}{1}
			\item The generic stabilizer of $X\times \CB$ contains regular semisimple elements;
		\end{enumerate}
		\begin{enumerate}[(i)]
			\setcounter{enumi}{2}
			\item We have $H\cap G_{\rs}\neq \emptyset$ and the function $h\in H\cap G_{\rs}\mapsto \frac{D^{H}_{alg}(h)^2}{D^{G}_{alg}(h)}$ extends to a regular function on $H$.
		\end{enumerate}
		\begin{enumerate}[(i')]
			\setcounter{enumi}{2}
			\item We have $H\cap G_{\rs}\neq \emptyset$ and the function $h\in H(F)\cap G_{\rs}(F)\mapsto \frac{D^H(h)^2}{D^G(h)}$ is locally bounded on $H(F)$ (i.e. it is bounded on the intersection of $G_{\rs}(F)$ with any compact subset of $H(F)$).
		\end{enumerate}
		Moreover, the above assertions imply that $X$ is spherical. If furthermore $H$ is reductive, then the above conditions are also equivalent to:
		\begin{enumerate}[(i)]
			\setcounter{enumi}{3}
			\item The generic stabilizer of $TX$ in $G$ contains regular elements;
			\item The generic stabilizer of $X\times X$ in $G$ for the diagonal action contains regular elements.
		\end{enumerate}
	\end{prop}
	
	\begin{rmk}
		The above proposition does not hold without the assumption that $H$ is connected as the example of $X=O(2)\backslash \GL_2$ shows. Indeed, for $X=O(2)\backslash \GL_2$, conditions (i), (i'), (ii), (ii') are satisfied but not (iii) and (iii'). We also believe that (i), (ii), (ii'), (iii) and (iii') are still equivalent when $X$ is not necessarily quasi-affine (but still assuming that $H$ is connected) but that (i') is strictly stronger (e.g. take $X=\CB$).
	\end{rmk}

	We will say that the variety $X$ is {\it coregular}, or that the pair $(G,H)$ is {\it coregular}, if the equivalent conditions (i)-(iii) (or (i)-(v) if $H$ is reductive) of the above proposition are satisfied.
	
	\begin{proof}
		
		Pick a Borel subgroup $B\subset G_{\overline{F}}$ with unipotent radical $N$ and let $P(X)$, $U(X)$ be the respective stabilizers of the generic $B$ and $N$ orbits in $X$. In other words, there exists an open dense subset $\CU\subset X_{\overline{F}}$ such that $xp\in xB$ (resp. $xu\in xN$) for every $(x,p)\in \CU\times P(X)$ (resp. $(x,u)\in \CU\times U(X)$) and the subgroups $P(X)$, $U(X)$ are maximal for these properties.
		
		Let $L(X)\subset P(X)$ be a Levi factor and set $S(X)=L(X)\cap U(X)$. By \cite[Korollar 2.9]{Knop93}, we know that $U(X)$ is a normal subgroup of $P(X)$ and the quotient $A_X:=P(X)/U(X)=L(X)/S(X)$ is a torus. Moreover, by the local structure theorem of \cite[Theorem 2.3, Proposition 2.4]{Knop94}, there exists a locally closed subvariety $\Sigma\subset X_{\overline{F}}$ which is $L(X)$-stable, on which the $L(X)$-action factors through the quotient $L(X)\to A_X$ and on which the resulting $A_X$-action is free, such that the $P(X)$-action induces an open embedding:
		\begin{equation}\label{coreg eq1}
			\displaystyle \Sigma\times^{L(X)} P(X)\hookrightarrow X.
		\end{equation}
		Since $P(X)=L(X)B$, it follows that the generic stabilizer of $X\times \CB$ exists and is given by the conjugacy class of $S(X)\cap B$. On the other hand, by the construction of \cite[\S 3]{Knop94} there exists a $L(X)$-equivariant embedding
		$$\displaystyle \Sigma\times \Fa_X^*\hookrightarrow T^*X$$
		whose image intersects every generic $G$-orbit in $T^*X$ \cite[Theorem 3.2, Lemma 3.1]{Knop94} and whose composition with the moment map $T^*X\to \Fg^*$ is the second projection $\Sigma\times \Fa_X^*\to \Fa_X^*$ (followed by the natural inclusions $\Fa_X^*\subset \mathfrak{l}(X)^*\subset \mathfrak{g}^*$). As the centralizer of a generic element in $\Fa_X^*$ is $L(X)$ \cite[Lemma 2.1]{Knop94}, this shows that the generic stabilizer of $T^*X$ exists and is the same as that of $\Sigma$ in $L(X)$, i.e. $S(X)$.
		
		As every conjugacy class, over $\overline{F}$, in $L(X)$ meets $L(X)\cap B$ and $S(X)$ is normal in $L(X)$, every element of $S(X)$ is $G$-conjugated to an element of $S(X)\cap B$ and this shows $(i)\Leftrightarrow (ii)$, $(i')\Leftrightarrow (ii')$. Moreover, if $S(X)$ contains a regular element of $G$ then it contains a regular semisimple one. Indeed, if $g\in S(X)$ is $G$-regular and $g=su$ is its Jordan decomposition, then $u$ is a regular unipotent element of the connected centralizer $Z_G(s)^0$. However, $u$ belongs to the Levi subgroup $L(X)\cap Z_G(s)^0$ of $Z_G(s)^0$ and therefore $Z_G(s)^0\subset L(X)$. If this is so, an element of the form $st$ for $t\in Z_G(s)^0_{der}$ in general position will be regular semisimple in $G$ and this proves the claim as $L(X)_{der}\subset S(X)$ implies $Z_G(s)^0_{der}\subset S(X)$.
		
		Thus, if (i) and (ii) are satisfied, $S(X)$ contains a regular semisimple element and so does $S(X)\cap B$. This proves the equivalence between (i), (i'), (ii) and (ii'). 
		
		Assume now that (i') is satisfied i.e. that there exists $h\in S(X)$ which is regular semisimple in $G$. Then, $T_G=Z_G(h)^0\subset L(X)$ acts transitively on all the connected components of the subvariety of fixed points $X^h$. As $\Sigma$ is a connected subvariety of $X^h$ this shows that $\Sigma$ is actually homogeneous under $L(X)$ and it follows, by the open embedding \eqref{coreg eq1}, that $X$ is spherical. Up to conjugacy, we may assume that the canonical base point $x_0=H1$ of $X$ belongs to $\Sigma$ i.e. that $HB$ is open in $G$. Then by \eqref{coreg eq1}, choosing a splitting of the surjection $\Fa_{L(X)}\twoheadrightarrow \Fa_X$, we have a direct sum decomposition
		\begin{equation*}
			\displaystyle \Fg=\Fh\oplus \Fa_X\oplus \Fn(X)
		\end{equation*}
		where $\Fn(X)$ denotes the nilradical of $\Fp(X)$. Note that this decomposition is stable under the adjoint action of $T_H=Z_H(h)^0$ as the latter is a maximal torus of $H$ contained in $L(X)$. Therefore, for $t\in T_H\cap G_{rs}$ we have
		\begin{equation*}
			\displaystyle D^G_{alg}(t)=D^H_{alg}(t)\det(1-\Ad_t\mid_{\Fn(X)})
		\end{equation*}
		from which it follows that
		\begin{equation*}
			\displaystyle \frac{D^H_{alg}(t)^2}{D^G_{alg}(t)}=\frac{D^G_{alg}(t)}{\det(1-\Ad_t\mid_{\Fn(X)})^2}=\delta_{P(X)}^{alg}(t)^{-1}D^{L(X)}_{alg}(t)
		\end{equation*}
		which implies (iii).
		
		It is clear that (iii) implies (iii').
		
		Assume now that $(G,H)$ satisfies (iii') and let us show that (i') is also satisfied. First, we make a reduction to the case where the generic stabilizer $S(X)$ of $T^*X$ is a torus. Indeed, since $(T^* X)(F)$ is Zariski dense in $T^* X$, up to conjugating we may assume that $S(X)$ is the stabilizer of a point $p=(x,\xi)\in (T^* X)(F)$, hence is defined over $F$, and even that $x=x_0$ is the canonical base-point of $X=H\backslash G$ (so that, in particular, $S(X)\subset H$). Let $T_S\subset S(X)$ be a maximal torus and let $H'=Z_H(T_S)$, $G'=Z_G(T_S)$ be the centralizers of $T_S$ in $H$ and $G$ respectively. We need to show that $G'$ is a torus (i.e. that $T_S$ contains regular semi-simple elements of $G$). Let $X'=H'\backslash G'$ be the connected component of the subvariety of fixed points $X^{T_S}$ containing $x$. Then, $X'$ is a homogeneous $G'$-variety and $T^*(X')$ is a connected component of the subvariety of fixed points $(T^* X)^{T_S}$ in the cotangent bundle. We claim that:
		\begin{num}
			\item\label{coreg eq2} The pair $(G',H')$ also satisfies condition (iii') i.e. $H'\cap G'_{rs}\neq \emptyset$ and the function $h\in H'(F)\cap G'_{\rs}(F)\mapsto \frac{D^{H'}(h)^2}{D^{G'}(h)}$ is locally bounded on $H'(F)$.
		\end{num}
		Indeed, $H'$ contains a maximal torus of $H$ hence regular semisimple elements of $G$ by assumption but such elements are a fortiori also regular semisimple in $G'$. Moreover, for $h\in H'(F)$ and $t\in T_S(F)$ we have
		$$\displaystyle D^H(ht)=D^{H'}(h) \lvert \det(1-\Ad_{ht}\mid_{\Fh/\Fh'})\rvert,\; D^G(ht)=D^{G'}(h) \lvert \det(1-\Ad_{ht}\mid_{\Fg/\Fg'})\rvert$$
		and for each $h_0\in H'(F)$ we can find $t\in T_S(F)$ as well as an open neighborhood $U\subset H'(F)$ of $h_0$ such that $h\mapsto \lvert \det(1-\Ad_{ht}\mid_{\Fh/\Fh'})\rvert$ is bounded from below and $h\mapsto \lvert\det(1-\Ad_{ht}\mid_{\Fg/\Fg'})\rvert$ is bounded from above on $U$. By the assumption that $(G,H)$ satisfies condition (iii'), this shows that the function $h\mapsto \frac{D^{H'}(h)^2}{D^{G'}(h)}$ is bounded on $U$ hence the function is locally bounded everywhere.
		
		Let $B\subset G$ be a Borel subgroup containing $T_S$ (not necessarily defined over $F$) and set $B'=Z_B(T_S)$, a Borel subgroup of $G'$ containing a maximal torus $T$ of $G$, and let $T_H\subset H'$ be a maximal torus. Taking $T_S$-invariants of the embedding \eqref{coreg eq1}, we see that $X'$ contains an open subset $B'$-equivariantly isomorphic to $\Sigma\times^{T} B'$. Furthermore, in a neighborhood of $0\in \Ft_H(F)$, the functions $X\mapsto D^{H'}(e^X)$ and $X\mapsto D^{G'}(e^X)$ are products of $\dim(H')-\dim(T_H)$ and $\dim(G')-\dim(T)$ absolute values of linear forms respectively. Thus, \eqref{coreg eq2} implies that
		\begin{equation*}
			\displaystyle \dim(G')-\dim(T)\leqslant 2(\dim(H')-\dim(T_H))
		\end{equation*}
		or equivalently
		\begin{equation}\label{coreg eq3}
			\displaystyle \dim(X')\leqslant \frac{\dim(G')+\dim(T)}{2}-\dim(T_H)=\dim(B')-\dim(T_H).
		\end{equation}
		However, as $\Sigma\times^{T} B'$ is open in $X'$, we also have
		\begin{equation}\label{coreg eq4}
			\displaystyle \dim(X')=\dim(\Sigma)+\dim(B')-\dim(T).
		\end{equation}
		Combining \eqref{coreg eq3} with \eqref{coreg eq4}, we obtain that
		\begin{equation*}
			\displaystyle \dim(\Sigma)\leqslant \dim(T)-\dim(T_H).
		\end{equation*}
		However, as $T/T_S$ acts freely on $\Sigma$ and $T_H$ contains $T_S$ this last inequality is only possible if $T_S=T_H$. But then, by the assumption that $H\cap G_{rs}\neq \emptyset$ and since $H$ is connected, this implies that $T_S$ and hence also $S(X)$ contains a regular semisimple element i.e. (i') is verified. This proves that (iii')$\Rightarrow$ (i') and therefore that the conditions (i), (i'), (ii), (ii'), (iii) and (iii') are all equivalent.
		
		It remains to show that these are also equivalent to (iv) and (v) when $H$ is reductive. The equivalence (iv)$\Leftrightarrow$ (v) follows from Luna's slice theorem \cite{Luna} applied to the diagonal $G$-orbit in $X\times X$ and noting that the normal bundle to the diagonal in $X\times X$ is isomorphic to $TX$. On the other hand, we have
		$$\displaystyle TX=\Fg/\Fh\times^H G,\;\; T^*X=\Fh^\perp\times^H G$$
		where $\Fh^\perp$ stands for the orthogonal of $\Fh$ in $\Fg^*$. As both $H$ and $G$ are reductive, the adjoint representation of $H$ on $\Fh$ is isomorphic to the coadjoint action of $H$ on $\Fh^\perp$ and this shows that (v)$\Leftrightarrow$(i).
	\end{proof}
	
	It is clear from the above discussion that if $(G,H)$ is coregular then $H_{rs}\subset G_{rs}$. However, the opposite direction is not true in general. For example, when $(G,H)=(\GL_3,\SL_2)$, we have $H_{rs}\subset G_{rs}$ but the pair is not coregular. The next lemma shows that in the case of symmetric pairs, the coregular condition is equivalent to $H_{rs}\subset G_{rs}$.
	
	\begin{lem}\label{regular symmetric pair}
		A symmetric pair $(G,H)$ is coregular if and only if $H_{rs}\subset G_{rs}$.
	\end{lem}
	
	\begin{proof}
		The ``only if" direction is obvious, we will only prove the other direction. So assume that $H_{rs}\subset G_{rs}$. Let $T_H\subset H$ be a maximal torus and let $T=Z_G(T_H)$ be its centralizer in $G$ (a maximal torus of $G$ by the assumption that $H_{\rs}\subset G_{\rs}$). Pick a cocharacter $\lambda\in X_*(T_H)$ that is $G$-regular and let $B\subset G$ be the Borel subgroup consisting of elements $b_in G$ such that $\lambda(t)b\lambda(t)^{-1}$ has a limit when $t\to 0$. Let $\iota$ be the involution of $G$ such that $H=(G^\iota)^\circ$. Obviously, $B$ is $\iota$-stable and contains $T$. Moreover, $B_H:=B\cap H$ is a Borel subgroup of $H$ containing the maximal torus $T_H$.
		
		Let $\Sigma^+$ to denote the set of positive roots of $T$ with respect to $B=TN$ and for each $\alpha\in \Sigma$, let  $X_{\alpha}\in \Fn(F)$ denote a nonzero element in the root subspace corresponding to $\alpha$. Then we can decompose $\Sigma$ into a union of three subsets $\Sigma=\Sigma_1\cup \Sigma_2\cup \Sigma_3$ where
		$$\Sigma_1=\{\alpha\in \Sigma|\;\iota(\alpha)\neq \alpha\},\;\Sigma_2=\{\alpha\in \Sigma|\;\iota(\alpha)=\alpha,\;\iota(X_\alpha)=X_\alpha\},$$
		$$\Sigma_3=\{\alpha\in \Sigma|\;\iota(\alpha)=\alpha,\;\iota(X_\alpha)=-X_\alpha\}.$$
		
		Note that $\iota$ induces an involution without fixed points of $\Sigma_1$ and that, denoting by $\Sigma_1/\iota$ the set of $\iota$-orbits in $\Sigma_1$, the set of positive roots $\Sigma_H^+$ of $T_H$ with respect to $B_H$ is in bijection with the set $\Sigma_1/\iota\cup \Sigma_2$ by the map
		$$\displaystyle \alpha\in \Sigma_1/\iota\cup \Sigma_2\mapsto \alpha\in \Sigma_H^+.$$
		In particular, we have
		$$\displaystyle D^G_{alg}(t)=\prod_{\alpha\in \Sigma^+} (1-\alpha(t))(1-\alpha(t)^{-1}),\;\; D^H_{alg}(t)=\prod_{\alpha\in \Sigma_1/\iota\cup \Sigma_2} (1-\alpha(t))(1-\alpha(t)^{-1}) \mbox{ for } t\in T_H.$$
		From these identities, we see that it suffices to show that $\Sigma_3=\emptyset$.

We will prove this by contradiction. Assume $\Sigma_3$ is non-empty and let $\alpha\in \Sigma_3$. The complements of $T_H\cap H_{\rs}$ in $T_H$ can be described as a union of divisors as follows
$$\displaystyle T_H\setminus (T_H\cap G_{\rs})=\bigcup_{\beta\in \Sigma_1/\iota\cup \Sigma_2^+} D_\beta$$
where $D_\beta:=\{t\in T_H\mid \beta(t)=1 \}$. Note that each of these divisors $D_\beta$ is actually a finite disjoint union of translates by the subtorus $T_\beta=\Ker(\beta\mid_{T_H})^0$ and that the subtori $(T_\beta)_{\beta}$ are two by two distinct (as $\Sigma_H^+$ is reduced). Since $T_H\setminus (T_H\cap G_{\rs})$ also contains the divisor $D_\alpha=\{t\in T_H\mid \alpha(t)=1 \}$ and $T_H\cap G_{\rs}=T_H\cap H_{\rs}$ by assumption, there must exist a $\beta\in \Sigma_1\cup \Sigma_2$ such that $T_\beta=\Ker(\alpha\mid_{T_H})^0$ and so $D_\alpha\subset D_\beta$ i.e. $\beta\mid_{T_H}=k \alpha\mid_{T_H}$ for some integer $k\in \mathbb{Z}$. Note that this last equality is equivalent to
$$\displaystyle \beta+\iota(\beta)=k(\alpha+\iota(\alpha))=2k\alpha.$$

We now distinguish two cases. First, if $\beta\in \Sigma_2$ then $\iota(\beta)=\beta$ and the above identity becomes $\beta=k\alpha$ hence $\alpha=\beta\in \Sigma_2$ (as the system of positive roots $\Sigma^+$ is reduced) which is a contradiction. On the other hand, if $\beta\in \Sigma_1$, denoting by $\alpha^\vee$ the coroot associated to $\alpha$, we have
$$\displaystyle 2\langle \alpha^\vee,\beta\rangle=\langle \alpha^\vee,\beta+\iota(\beta)\rangle=2k\langle \alpha^\vee,\alpha\rangle=4k.$$
Hence, $\langle \alpha^\vee,\beta\rangle= 2k$ and therefore $k=1$. Thus, we have $2\alpha=\beta+\iota(\beta)$ and $\langle \alpha^\vee,\beta\rangle=2$. Since $\beta\neq \iota(\beta)$ this also implies that the length of $\beta$ is twice the length of $\alpha$ and thus $\langle \beta^\vee,\alpha\rangle=\langle \alpha^\vee,\beta\rangle/2=1$. Finally, denoting by $s_\beta$ the simple reflection corresponding to $\beta$, we deduce that $s_\beta \alpha=\alpha-\beta=(\iota(\beta)-\beta)/2$ is a root of $T$ in $G$. But that is a contradiction since this root would be $\iota$-antiinvariant and there is no such root.
	\end{proof}

	\begin{defn}
		Assume now given twisted spaces $(G,\tG)$ and $(H,\tH)$ with an embedding $\tH\hookrightarrow \tG$ that is compatible with the inclusion $H\subset G$. Following the above discussion, we will say that the pair $(\tG,\tH)$ is {\it coregular} if the following condition is satisfied:
		\begin{num}
			\item We have $\tH\cap \tG_{rs}\neq \emptyset$ and the function $h\in \tH(F)\cap \tG_{rs}(F)\mapsto \frac{D^{\tH}(h)^2}{D^{\tG}(h)}$ is locally bounded on $\tH(F)$.
		\end{num}
	\end{defn}

	\begin{rmk}\label{rmk coregular}
		\begin{itemize}
			\item There is yet another characterization of the coregular spherical varieties which reads as follows: if $X=H\backslash G$ is a quasi-affine homogeneous spherical variety then $X$ is coregular if and only if $\mathrm{rk}(X)=\mathrm{rk}(G)-\mathrm{rk}(H)$ where $\mathrm{rk}(G)$, $\mathrm{rk}(H)$, $\mathrm{rk}(X)$ denote the (absolute) ranks of $G$, $H$ and $X$ respectively. (We recall that the rank of a $G$-variety is, by definition, the rank of the torus $T_X=B/B_xU$ where $B=TU\subset G$ is a Borel subgroup and $x\in X$ a point in general position e.g., in the case of a spherical variety, a point in the open $B$-orbit.)
			
			Indeed, if $\mathrm{rk}(X)=\mathrm{rk}(G)-\mathrm{rk}(H)$ and $B\subset G$ is a Borel subgroup with $HB$ open, $B\cap H$ contains a maximal torus $T_H$ of $H$ and we have an isomorphism of $T_H$-representations $\fg/\Fb\simeq \Fh/\Fb\cap \Fh$. However, as a $T_H$-representation, $\fg/\Fb$ is isomorphic to the dual of $\Fu$ (the nilradical of $\Fb$) and $\Fh/\Ft_H$ does not contain the trivial representation of $T_H$. Therefore, all restrictions of the roots of $G$ to $T_H$ are non-trivial which implies that $T_H$ contains regular semi-simple element. As $T_H$ is contained in the generic stabilizer of $X\times \mathcal{B}$ this proves point (ii') of the above proposition.
			
			Conversely, assume that $X$ is coregular. Let again $B$ be a Borel subgroup with $HB$ open, then $H\cap B$ contains a $G$-regular semi-simple element $h$ (again by characterization (ii')), hence the maximal torus $T=G_h$ of $G$ is included in $B$ and therefore also the maximal torus $T_H=H_h$ of $H$ and this shows that the universal Cartan $T_X$ of $X$ is a finite quotient of $T/T_H$ hence $\mathrm{rk}(X)=\mathrm{rk}(G)-\mathrm{rk}(H)$.

			\item It is not true that $(\tG,\tH)$ is coregular if and only if $(G,H)$ is so. For example, let $H$ be connected reductive and take $\tG=(H\times H)\iota$ where $\iota(h_1,h_2)=(h_2,h_1)$ and $\tH=H$ with the embedding $\tH\hookrightarrow \tG$ given by $h\mapsto (h,h)\iota$. Then, the pair $(G,H)$ is always coregular whereas the pair $(\tG,\tH)$ is coregular if and only if for every $h\in H(F)$, $\det(1+\Ad_h)\neq 0$.
				\end{itemize}
	\end{rmk}
	
	\subsection{Tempered varieties}\label{S tempered varieties}
	
	We continue to consider the setting at the end of the last section: $(G,\tG)$ is a connected reductive twisted space over $F$ and $(H,\tH)$ is a closed connected twisted subspace of it. We also assume that $H$ is unimodular (this implies that $X=H\backslash G$ is quasi-affine).
	
	Following \cite[\S 2.7]{BeuGalP}, we say that the pair $(\tG,\tH)$ is {\it tempered} if it satisfies the following condition:
	\begin{num}
		\item There exists $d>0$ such that the integral
		$$\displaystyle \int_{\tH(F)} \Xi^{\tG}(h) \sigma_{\tG}(h)^{-d} dh$$
		is convergent.
	\end{num}
	Note that the pair $(\tG,\tH)$ is tempered if and only if $(G,H)$ is so. Moreover, by {\it loc. cit.}, a pair $(G,H)$ is tempered in the above sense if and only if $L^2(H(F)\backslash G(F))$ is tempered as a unitary representation of $G(F)$. (But we will not need this fact.)
	
	\begin{lem}\label{lem coregular and tempered}
		Assume that the pair $(\tG,\tH)$ is coregular and tempered. Then, the function 
		$$\displaystyle h\in \tH(F)\cap \tG_{rs}(F)\mapsto \frac{D^{\tH} (h)}{D^{\tG} (h)^{1/2}}$$
		is globally bounded.
	\end{lem}
	
	\begin{proof}
		Let $\tT\subset \tH$ be a maximal twisted torus. It is enough to show that the function 
		$$\displaystyle t\in \tT(F)\cap \tG_{rs}(F)\mapsto \frac{D^{\tH}(t)}{D^{\tG}(t)^{1/2}}$$
		is globally bounded. 
		
		Set $M=Z_G(A_{\tT})$ and $\widetilde{M}=Z_{\tG}(A_{\tT})=M \tT$. Then, $\tM$ is the minimal Levi subspace of $\tG$ containing $\tT$. For each $\tP\in \CF(\tM)$, we let $\tT_P^+$ be the subset of those $t\in \tT(F)$ such that all the eigenvalues of the restriction of $\Ad_t$ to $\mathfrak{p}(\overline{F})$ are of absolute value $\geq 1$. Then, we have a partition
		$$\displaystyle \tT(F)=\bigsqcup_{\tP\in \CF(\tM)} \tT_{P}^{+}$$
		and it is enough to show that, for any fixed $\tP\in \CF(\tM_{T})$, the function 
		$$\displaystyle t\in \tT_{P}^{+}\mapsto \frac{D^{\tH}(t)}{D^{\tG}(t)^{1/2}}$$
		is bounded.

		
		Let $\tP\in \CF(\tM)$ and assume that $\tT_P^+\neq \emptyset$.
		Let $\tL$ be the unique Levi factor of $\tP$ containing $\tM$ and set
		$$\displaystyle P_H=P\cap H,\; L_H=L\cap H,\; A_P^+=\{a\in A_L(F)\mid \lvert \alpha(a)\rvert\geqslant 1\; \forall \alpha\in \Delta_P \} \mbox{ and } A_{P,H}^+=A_P^+\cap H(F).$$
		Then, $\tT_P^+$ is right invariant by the monoid $A_{P,H}^+$ and the quotient $\tT_P^+/A_{P,H}^+$ is compact. Moreover, for $t\in \tT_P^+$, the Lie algebras $\mathfrak{p}(\oF)$ and $\mathfrak{p}_H(\oF)$ (resp. $\mathfrak{l}(\oF)$ and $\mathfrak{l}_H(\oF)$) are the maximal subspaces of $\mathfrak{g}(\oF)$ and $\mathfrak{h}(\oF)$ where all the eigenvalues of $\Ad_t$ are of absolute value $\geq 1$ (resp. $=1$). Setting $\tP_H=\tP\cap \tH$ and $\tL_H=\tL \cap \tH$, it follows that
		$$\displaystyle D^{\tG}(t)=D^{\tL}(t)\delta_{\tP}(t),\;D^{\tH}(t)=D^{\tL_H}(t)\delta_{\tP_H}(t)$$
		for every $t\in \tT_P^+$. Therefore,
		\begin{equation}\label{lem coregular and tempered eq1}
			\displaystyle \frac{D^{\tH} (t)}{D^{\tG} (t)^{1/2}}=\frac{D^{\tL_H}(t)}{D^{\tL}(t)^{1/2}}\cdot \frac{\delta_{\tP_H}(t)}{\delta_{\tP}(t)^{1/2}},\mbox{ for } t\in \tT_{P}^{+}.
		\end{equation}
		Since $(\tG,\tH)$ is coregular, this in particular entails that  the function
		$$\displaystyle t\in \tT_P^+ \mapsto D^{\tL_H}(t)D^{\tL}(t)^{-1/2}$$
		is locally bounded. Because $\tT_{P}^{+}/A_{P,H}^+$ is compact and $D^{\tL_H}$, $D^{\tL}$ are both $A_L(F)$-invariant, we deduce that the function $t\in \tT_P^+\mapsto D^{\tL_H}(t)D^{\tL}(t)^{-1/2}$ is globally bounded. Therefore, by \eqref{lem coregular and tempered eq1} it only remains to check that the function
		$$\displaystyle t\in \tT_P^+\mapsto \frac{\delta_{\tP_H}(t)}{\delta_{\tP}(t)^{1/2}}$$
		is bounded. Again because $\tT_P^+/A_{P,H}^+$ is compact, it is equivalent to work with the function
		$$\displaystyle a\in A_{P,H}^+\mapsto \frac{\delta_{P_H}(a)}{\delta_P(a)^{1/2}}.$$

		
		Let $J_H\subset H(F)$ be a compact-open subgroup. Then, we have
		\begin{equation}
			\displaystyle \vol(J_H aJ_H)\sim \delta_{P_H}(a),\mbox{ for } a\in A_{P,H}^+.
		\end{equation}
		Moreover, we can assume that $J_H\cap A_{L_H}(F)=A_{L_H}^c$ is the maximal compact subgroup of $A_{L_H}(F)$ and that the cosets $J_HaJ_H$, $a\in A_{P,H}^+/A_{L_H}^c$ are disjoint. As (see \cite[lemme II.1.1]{WalPlanch})
		$$\displaystyle \delta_{P}(a)^{-1/2}\ll \Xi^{G}(a), \mbox{ for } a\in A_P^+,$$
		and $(G,H)$ is tempered, we can find $d>0$ such that
		\[\begin{aligned}
			\displaystyle \sum_{a\in A_{P,H}^+/A_{L_H}^c} \frac{\delta_{P_H}(a)}{\delta_{P}(a)^{1/2}} \sigma(a)^{-d} & \ll \sum_{a\in A_{P,H}^+/A_{L_H}^c}\Xi^{G}(a)\sigma(a)^{-d}\vol(J_HaJ_H) \\
			& \leqslant \int_{H(F)} \Xi^G(h)\sigma(h)^{-d} dh<\infty.
		\end{aligned}\]
		Since $A_{P,H}^{+}/A_{L_H}^c$ is a finitely generated monoid and $\delta_{P_H}$, $\delta_{P}$ are characters on it, the above estimate implies that $a\in A_{P,H}^+\mapsto \frac{\delta_{P_H}(a)}{\delta_{P}(a)^{1/2}}$ is bounded and the lemma is proved.
		
	\end{proof}

	\subsection{Symmetric pairs}
	
	In this paper, by a {\it symmetric pair} (over $F$) we will mean a pair $(G,\iota)$ where $G$ is a connected linear group over $F$ and $\iota$ an involutive automorphism of $G$ defined over $F$. Let $(G,\iota)$ be a symmetric pair. We denote by $G^\iota$ the closed subgroup of $\iota$-fixed points, by $X^*(G)_\iota$ the subgroup of characters $\chi\in X^*(G)$ that are trivial on $G^\iota$ and by $A_{G,\iota}$ the neutral component of the subgroup $\{a\in A_G\mid \iota(a)=a^{-1} \}$. We also set
	$$\displaystyle \cA_{G,\iota}^*=X^*(A_{G,\iota})\otimes \mathbb{R}=X^*(G)_\iota\otimes \mathbb{R},\;\;\; \cA_{G,\iota}=X_*(A_{G,\iota})\otimes \mathbb{R}=\Hom(X^*(G)_\iota,\mathbb{R}).$$
	Then, $\cA_{G,\iota}^*$ (resp. $\cA_{G,\iota}$) can be identified with the subspace of $\iota$-antiinvariant vectors in $\cA_G^*$ (resp. in $\cA_G$). We also denote by $H_{G,\iota}: G(F)\to \cA_{G,\iota}$ the composition of $H_G$ and of the natural projection $\cA_G\to \cA_{G,\iota}$.
	
	Assume from now on that $G$ is reductive and connected. Recall that a parabolic subgroup $P\subset G$ is said to be {\em $\iota$-split} if $\iota(P)$ is opposite to $P$ and that a Levi subgroup $M\subset G$ is said to be {\em $\iota$-split} if there exists a $\iota$-split parabolic subgroup $P$ such that $M=P\cap \iota(P)$.
	
	We will denote by $\CP_\iota$ and $\CL_\iota$ the sets of all $\iota$-split parabolic subgroups and $\iota$-split Levi subgroups of $G$ respectively. We will also write $[\CP_\iota]=\CP_\iota/\sim$ where, for $P,P'\in \CP_\iota$, $P\sim P'$ if $P$ and $P'$ are $G(F)$-conjugate (or, equivalently, $G(\overline{F})$-conjugate). For $P\in \CP_\iota$, we set $\overline{P}=\iota(P)$ for the unique $\iota$-split parabolic subgroup opposite to $P$, $M_P=P\cap \iota(P)\in \CL_\iota$ for its unique $\iota$-split Levi factor and we denote by $[P]$ its image in $[\CP_\iota]$.
	
	By a {\it minimal $\iota$-split parabolic subgroup} of $G$, we mean a parabolic subgroup that is $\iota$-split, defined over $F$ and minimal for these properties. We denote by $\CP^{\min}_\iota\subseteq \CP_\iota$ the subset of minimal $\iota$-split parabolic subgroups. By \cite[Proposition 4.9]{HW}, all minimal $\iota$-split parabolic subgroups are conjugated under $G(F)$ i.e. the image $[\CP^{\min}_\iota]$ of $\CP^{\min}_\iota$ in $[\CP_\iota]$ is a singleton. For $M\in\CL_\iota$, we set $\mathcal{P}_\iota(M)=\CP(M)\cap \CP_\iota$, $\mathcal{F}_\iota(M)=\CF(M)\cap \CP_\iota$ and $\mathcal{L}_\iota(M)=\CL(M)\cap \CL_\iota$.

	For every $P\in \CP_\iota$, we set $A_{P,\iota}=A_{M_P,\iota}$, $\CA_{P,\iota}=\CA_{M_P,\iota}$ and we denote by $H_{P,\iota}:P(F)\to \CA_{P,\iota}$ the composition of the projection $P(F)\to M_P(F)$ with $H_{M_P,\iota}$. Then, for $P,Q\in \CP_\iota$ with $P\subset Q$ we have the decomposition $\cA_{P,\iota}=\cA_{P,\iota}^Q\oplus \cA_{Q,\iota}$ where $\cA_{P,\iota}^{Q}=\cA_{P,\iota}\cap \cA_P^{Q}$ and we denote by $\Delta_{P,\iota}^Q$, $\Delta_{P,\iota}^{Q,\vee}$ the respective projections of $\Delta_{P}^Q$ and $\Delta_P^{Q,\vee}$ to $\cA_{P,\iota}^{Q,*}$ and $\cA_{P,\iota}^Q$. When $Q=G$, we will sometimes drop the superscript and when $M=M_P$ we will sometimes write $\CA_{M,\iota}^Q$ for $\CA_{P,\iota}^Q$. For $P\in \CP_\iota$, we also set
	$$\displaystyle A_{P,\iota}^+=\{a\in A_P(F)\mid \langle \alpha, H_{P,\iota}(a)\rangle\geq 0\; \forall \alpha\in \Delta_{P,\iota} \}.$$
	
	Let $P_0\in \CP^{\min}_\iota$ and set $M_0=M_{P_0}=P_0\cap \iota(P_0)$, $A_{0,\iota}=A_{M_0,\iota}$ and $\cA^G_{0,\iota}=\cA^G_{M_0,\iota}$. It is known that the set $\Sigma_{0,\iota}\subseteq \cA^{G,*}_{0,\iota}$ of nonzero weights for the adjoint action of $A_{0,\iota}$ on $\mathfrak{g}$ is a root system and that the subset of weights $\Sigma_{0,\iota}^+\subset \Sigma_{0,\iota}$ appearing in $\mathfrak{p}_0$ forms a system of positive roots with associated set of simple roots $\Delta_{0,\iota}=\Delta_{P_0,\iota}$ see \cite[\S 5]{HW}. We will denote by $W_{0,\iota}$ the Weyl group of this root system. We also set
	$$\displaystyle \cA_{P_0,\iota}^{+}=\{X\in \cA_{P_0,\iota}\mid \langle \alpha, X\rangle\geq 0\; \forall \alpha\in \Delta_{0,\iota} \},\; {}^-\cA_{P_0,\iota}=\{X\in \cA_{P_0,\iota}\mid \langle \varpi_\alpha, X\rangle\leq 0\; \forall \alpha\in \Delta_{0,\iota} \}.$$

	There is a bijection between $\mathcal{F}_\iota(M_0)$ and the collection of parabolic subsets of $\Sigma_{0,\iota}$ obtained by sending $Q\in \mathcal{F}_\iota(M_0)$ to the set $\Sigma_{0,Q,\iota}$ of nonzero weights of $A_{0,\iota}$ in $\mathfrak{q}$. Furthermore, for every $\iota$-split parabolic subgroup $P\supset P_0$, the subset $\Delta_{0,\iota}^P:=\Delta_{P_0,\iota}^P$ coincides with the set of simple roots $\Sigma_{0,P,\iota}\cap -\Sigma_{0,P,\iota}\cap \Delta_{0,\iota}$ and elements of $\Delta_{0,\iota}^{P,\vee}$ are {\em positively proportional} to the coroots associated to $\Delta_{0,\iota}^P$.
	
	Let $P,Q\in \CP_\iota$ with $[P]=[Q]$ and choose $\gamma\in G(F)$ such that $\gamma P\gamma^{-1}=Q$ and $\gamma M_{P}\gamma^{-1}=M_{Q}$. Then, $\Ad_\gamma$ induces isomorphisms
	$$\displaystyle I_{P,Q}: \CA_{P,\iota}\simeq \CA_{Q,\iota},\;\;\; I_{P,Q}: A_{P,\iota}\simeq A_{Q,\iota}$$
	sending respectively $\Delta_{P,\iota}$ to $\Delta_{Q,\iota}$, $A_{P,\iota}^+$ to $A_{Q,\iota}^+$ and which is independent of the choice of $\gamma$. Moreover, it is readily seen that the element $\iota(\gamma)$ still conjugates the pair $(P,M_P)$ to $(Q,M_{Q})$ from which it follows that $\gamma\iota(\gamma)^{-1}\in M_Q(F)$.

	Let us further fix a special maximal compact subgroup $K\subset G(F)$ that we use to extend $H_{P,\iota}$ to a right $K$-invariant map $G(F)\to \cA_{P,\iota}$ for every $P\in \CP_\iota$ by mean of the Iwasawa decomposition $G(F)=P(F)K$. Then, to every $P,Q\in \CP_\iota$ with $[P]=[Q]$ we associate a point $Y^K_{Q,P}\in \cA_{Q,\iota}$ as follows. Pick $\gamma\in G(F)$ such that $\gamma P\gamma^{-1}=Q$ and $\gamma M_P\gamma^{-1}=M_Q$. Then, recalling that $\gamma\iota(\gamma)^{-1}\in M_Q(F)$, we set
	$$\displaystyle Y^K_{Q,P}:=H_{\overline{Q},\iota}(\gamma)-\frac{1}{2}H_{M_Q,\iota}(\gamma \iota(\gamma)^{-1})$$
	where $\overline{Q}=\iota(Q)$ is the $\iota$-split parabolic subgroup opposite to $Q$.

	\begin{lem}\label{lem YQP}
		The element $Y^K_{Q,P}\in \cA_{Q,\iota}$ so constructed doesn't depend on the choice of $\gamma$ (i.e. it only depends on $P$, $Q$ and $K$).
	\end{lem}
	
	\begin{proof}
		Because the normalizer of the pair $(Q,M_Q)$ in $G$ is equal to $M_Q$, for any other element $\gamma'\in G(F)$ satisfying $\gamma' P(\gamma')^{-1}=Q$, $\gamma' M_P(\gamma')^{-1}=M_Q$, there exists $m\in M_Q(F)$ such that $\gamma'=m \gamma$. Then, it follows that $$\displaystyle H_{\overline{Q},\iota}(\gamma')=H_{M_Q,\iota}(m)+H_{\overline{Q},\iota}(\gamma) \mbox{ and } H_{M_Q,\iota}(\gamma' \iota(\gamma')^{-1})=2H_{M_Q,\iota}(m)+H_{M_Q,\iota}(\gamma\iota(\gamma)^{-1}),$$
		hence
		$$\displaystyle H_{\overline{Q},\iota}(\gamma')-\frac{1}{2}H_{M_Q,\iota}(\gamma'\iota(\gamma')^{-1})=H_{\overline{Q},\iota}(\gamma)-\frac{1}{2}H_{M_Q,\iota}(\gamma \iota(\gamma)^{-1}).$$
	\end{proof}

	\subsection{Symmetric varieties}
	
	Let $(G,\iota)$ be a symmetric pair with $G$ reductive and connected as in the previous section. We set $H=G^\iota$ and let $X=H\backslash G$ be the corresponding {\em symmetric variety}.
	
	For every $M\in \CL_\iota$, we set $H_M=H\cap M$ and $X_M=H_M\backslash M$. Note that $X_M$ is the symmetric variety associated to the symmetric pair $(M,\iota\mid_M)$ and that it is naturally a closed subvariety of $X$. Any character $\chi\in X^*(M)_\iota$ is by definition trivial on $H$ and therefore descends to a regular map $X_M\to \mathbb{G}_m$ that we shall denote by the same letter. We define a {\em Harish-Chandra map}
	$$\displaystyle H_{M,\iota}: X_M(F)\to \CA_{M,\iota}=\Hom(X^*(M)_\iota,\mathbb{R})$$
	by $\langle \chi, H_{M,\iota}(x)\rangle=\log \lvert \chi(x)\rvert$ for every $x\in X_M(F)$ and $\chi\in X^*(M)_{\iota}$. We also have a left action of $A_M$ on $X_M$ given by $a\cdot x=xa$ for $(a,x)\in A_M\times X_M$ which commutes with the (right) $M$-action and satisfies $H_{M,\iota}(a\cdot x)=H_{M,\iota}(a)+H_{M,\iota}(x)$ for every $(a,x)\in A_M(F)\times X_M(F)$. 
	
	For $P_0\in \CP^{\min}_\iota$ with $M_0=M_{P_0}$, we set
	$$\displaystyle X_{P_0}^+=\{x\in X_{M_0}(F)\mid \langle \alpha,H_{M_0,\iota}(x)\rangle \geqslant 0\; \forall \alpha\in \Delta_{P_0,\iota} \}.$$
	If moreover $P\in \CP_\iota$ is such that $P_0\subset P$ and $C\geq 0$ we define
	$$\displaystyle X_{P_0}^+(\geq C,P)=\{x\in X_{P_0}^+\mid \langle \alpha,H_{M_0,\iota}(x)\rangle \geqslant C\; \forall \alpha\in \Delta_{P_0,\iota}\setminus \Delta^P_{P_0,\iota} \}.$$
	
	Set
	$$\displaystyle \cA_X:=\varprojlim_{P_0\in \CP_\iota^{\min}} \cA_{P_0,\iota}$$
	where the transition maps are given by the isomorphisms $I_{P_0,P_0'}:\CA_{P_0,\iota}\simeq \CA_{P_0',\iota}$ for $P_0,P_0'\in \CP_\iota^{\min}$. Then, $\cA_X$ is a real vector space equipped with canonical isomorphisms $\cA_X\simeq \cA_{P_0,\iota}$ for every $P_0\in \CP_\iota^{\min}$. The images by this isomorphism of the cones $\CA_{P_0,\iota}^+$ and ${}^- \CA_{P_0,\iota}$ don't depend on the choice of $P_0$, we will denote them by $\CA_X^+$, ${}^- \CA_X$ respectively. Moreover, for every $P\in \CP_\iota$ we can choose $P_0\in \CP_\iota^{\min}$ such that $P_0\subset P$ and we get an embedding $\CA_{P,\iota}\hookrightarrow \CA_X$ given by the composition of the natural inclusion $\CA_{P,\iota}\subset \CA_{P_0,\iota}$ with the isomorphism $\CA_{P_0,\iota}\simeq \cA_X$. This embedding actually does not depend on the choice of $P_0$ as can readily be checked. 
	
	Let $P,Q\in \CP_\iota$ with $[P]=[Q]$ and choose $\gamma\in G(F)$ such that $\gamma P\gamma^{-1}=Q$, $\gamma M_P\gamma^{-1}=M_Q$. Then, we have $\gamma\in HM_P$: indeed, as already argued $\iota(\gamma)^{-1}\gamma\in M_P$ but this element is also in the neutral component of the subvariety of $\iota$-antiinvariant elements in $M_P$ hence there exists $m\in M_P$ such that $\iota(\gamma)^{-1}\gamma=\iota(m)^{-1}m$ or equivalently $\gamma\in Hm$. It now readily follows that $X_{M_Q}=X_{M_P}\gamma^{-1}$ and
	\begin{equation}\label{eq HM0iota}
		\displaystyle H_{M_Q,\iota}(x\gamma^{-1})=\Ad_\gamma(H_{M_P,\iota}(x))-\frac{1}{2}H_{M_Q,\iota}(\gamma\iota(\gamma)^{-1})=I_{P,Q}(H_{M_P,\iota}(x))+\frac{1}{2}H_{M_Q,\iota}(\gamma\iota(\gamma)^{-1})
	\end{equation}
	for every $x\in X_{M_P}(F)$. (Recall that $\gamma\iota(\gamma)^{-1}\in M_Q(F)$.)

	\subsection{Neighborhoods of infinity}
	
	Recall the following {\em weak Cartan decomposition} from \cite{DS} and \cite{BO}: for every $P_0\in \CP_\iota^{\min}$ we can find a compact subset $\CK\subset G(F)$ such that
	\begin{equation}\label{Weak Cartan}
		\displaystyle X(F)=X_{P_0}^+ \CK.
	\end{equation}
	
	Let $P\in \CP_\iota$ and set $M=M_P$. Choose $P_0\in \CP_\iota^{\min}$ with $P_0\subset P$ and a compact subset $\CK$ satisfying the equality \eqref{Weak Cartan}. Then, following \cite{DelPlanch}, we define a {\it neighborhood of $\infty_P$} in $X(F)$ to be a subset of the latter containing
	\begin{equation*}
		\displaystyle X_{P_0}^+(\geq C,P)\CK
	\end{equation*}
	for some large enough constant $C>0$. This notion actually only depends on the class $[P]$ in $[\CP_\iota]$, and in particular not on the auxilliary choices of $P_0$ and $\CK$. Indeed, using \eqref{eq HM0iota} this readily reduces to showing the following: if $\CK'\supset \CK$ is a bigger compact subset then for every $C\geq 0$ we can find $C'\geq 0$ with $X_{P_0}^+(\geq C,P)\CK\supset X_{P_0}^+(\geq C',P)\CK'$. This, in turn, is a consequence of the following lemma.
	
	\begin{lem}\label{lemma weak Cartan}
		Let $P_0\in \CP_\iota^{\min}$, $\CK\subset G(F)$ be a compact subset and set $M_0=M_{P_0}$. Then, there exists $d>0$ such that for every $x,y\in X_{P_0}^+$, $x\CK\cap y\CK\neq \emptyset$ implies $\lVert H_{M_0,\iota}(x)-H_{M_0,\iota}(y)\rVert\leq d$.
	\end{lem}
	
	\begin{proof}
		Set $A_{0,\iota}=A_{M_0,\iota}$ and recall from \cite[Proposition 4.7 (iii)]{HW} that $X_{M_0}(F)/A_{0,\iota}(F)$ is compact. It follows that we can find a compact subset $\Omega_0\subset X_{M_0}(F)$ such that
		\begin{equation}\label{eq1 lemma weak Cartan}
			X_{P_0}^+\subset \Omega_0 A_{P_0,\iota}^+.
		\end{equation}
		Let $X^*(M_0)_\iota^+$ be the subset of {\it dominant} weights $\chi\in X^*(M_0)_\iota$ i.e. such that $\langle \alpha^\vee, \chi\rangle\geq 0$ for every $\alpha^\vee\in \Delta_{P_0}^\vee$. Then, for every $\chi\in X^*(M_0)_\iota^+$ there exists a nonzero regular function $f_{2\chi}\in F[X]$ such that $f_{2\chi}(xp_0)=f(x)\chi(p_0)^2$ for every $(x,p_0)\in X\times P_0$\footnote{Indeed, since $\overline{P}_0=\iota(P_0)=M_0\overline{N}_0$ is opposite to $P_0$, there exists a nonzero regular function $\varphi_\chi\in F[G]$ such that $\varphi_\chi(\overline{u}p)=\chi(p)$ for $(\overline{u},p)\in \overline{N}_0\times P_0$ and it suffices to take $f_{2\chi}(x)=\varphi_\chi(\iota(x)^{-1}x)$.}. Moreover, up to scaling $f_{2\chi}$ we may assume that $f_{2\chi}(x)=\chi(x)^2$ for every $x\in X_{M_0}$. Let $V_\chi\subset F[X]$ be the $G$-submodule generated by $f_{2\chi}$ for the action by right translation $R$. Then the weights of $A_{M_0}$ in $V$ are of the form $2\chi-\sum_{\alpha\in \Delta_{P_0}}n_\alpha \alpha$ where $n_\alpha\in \mathbb{N}$. From this and \eqref{eq1 lemma weak Cartan} it follows that for every compact subset $L\subset V$ we can find $c_\chi^L>0$ such that
		\begin{equation}\label{eq2 lemma weak Cartan}
			\displaystyle \lvert f'(x)\rvert\leq c^L_\chi \lvert\chi(x)\rvert^2 \mbox{ for every } (f',x)\in L\times X_{P_0}^+.
		\end{equation}
		We will apply this to $L=R(\CK')f_{2\chi}$ where $\CK'=\CK\CK^{-1}$, setting $c_\chi=c_\chi^L$ for simplicity. Indeed, for $x,y\in X_{P_0}^+$ such that $x\CK\cap y\CK\neq \emptyset$ we can find $f'\in L$ such that $f_{2\chi}(x)=f'(y)$. Thus, applying \eqref{eq2 lemma weak Cartan} we get
		$$\displaystyle \lvert \chi(x)\rvert^2=\lvert f_{2\chi}(x)\rvert=\lvert f'(y)\rvert\leq c_\chi \lvert \chi(y)\rvert^2$$
		i.e. $\langle \chi,H_{M_0,\iota}(x)\rangle\leq \langle \chi,H_{M_0,\iota}(y)\rangle+\frac{1}{2}\log(c_\chi)$. By symmetry, we also have the inequality with $x$, $y$ permuted. As this holds for every $\chi\in X^*(M_0)_\iota^+$ and $X^*(M_0)_\iota^+$ generates $\CA_{M_0,\iota}^*$ this gives the desired result. 
	\end{proof}
	
	We shall denote by $\CN(\infty_P)$ the collection of all neighborhoods of $\infty_P$ in $X(F)$. The set $\CN(\infty_P)$ is stable by finite intersections and translations by elements of $G(F)$. By a {\em basis} of $\CN(\infty_P)$ we mean a subset $\CN'\subset \CN(\infty_P)$ such that every element $\Omega\in \CN(\infty_P)$ contains at least one $\Omega'\in \CN'$.

	We define similarly the notion of {\em neighborhood of $\infty^M_P$} in $X_M(F)$ as follows. By the weak Cartan decomposition \eqref{Weak Cartan} applied to the symmetric variety $X_M$, we can find a compact subset $\CK_M\subseteq M(F)$ such that $X_M(F)=X_{P_0\cap M}^+ \CK_M$. Then, by definition, a neighborhood of $\infty^M_P$ in $X_M(F)$ is a subset of the latter containing $X_{P_0}^+(\geq C,P)\CK_M$ for a suitable $C>0$. Once again, using the above lemma, we can show that this notion is independent on the choices of $P_0$ and $\CK_M$. We will denote by $\CN(\infty^M_{P})$ the collection of all neighborhoods of $\infty^M_P$ in $X_M(F)$. Note that $\CN(\infty^M_{P})$ admits a basis consisting of (left) $A_{P,\iota}^+$-invariant subsets (e.g. the family of subsets $X_{P_0}^+(\geq C,P)\CK_M$ would do).
	
	Let us now fixed a special maximal compact subgroup $K\subseteq G(F)$. Let $\overline{P}=MU_{\overline{P}}\in \CP^\iota$ be the parabolic subgroup opposite to $P$ with respect to $M$. Then, every $\gamma\in G(F)$ admits an Iwasawa decomposition $\gamma=m_{\overline{P}}(\gamma)u_{\overline{P}}(\gamma)k_{\overline{P}}(\gamma)$ with $m_{\overline{P}}(\gamma)\in M(F)$, $u_{\overline{P}}(\gamma)\in U_{\overline{P}}(F)$ and $k_{\overline{P}}(\gamma)\in K$.
	
	\begin{lem}\label{lem1 neighborhoods}
		Let $\CK\subset G(F)$ be a compact subset. Then, we can find $\Omega_P^M\in \CN(\infty^M_{P})$ such that
		$$\displaystyle x\gamma K=xm_{\overline{P}}(\gamma)K$$
		for every $(x,\gamma)\in \Omega_P^M\times \CK$.
	\end{lem}
	
	\begin{proof}
		For every neighborhood $\Omega\in \CN(\infty^M_{P})$, we can find another one $\Omega'\in \CN(\infty^M_{P})$ such that $\Omega'm_{\bar{P}}(\gamma)\subset \Omega$ for every $\gamma\in \CK$. It follows that we may assume that $\CK\subset U_{\bar{P}}(F)$ and $m_{\overline{P}}(\gamma)=1$ for every $\gamma\in \CK$.
		
		The lemma is then a variant of the wavefront Lemma \cite[Corollary 5.3.2]{SV}. Indeed, let us fix $P_0=M_0U_0\in \CP_\iota^{\min}$ with $P_0\subset P$ as well as representatives $x_1,\ldots,x_n$ for the $M_0(F)$-orbits in $X_{M_0}(F)$. Set, for $C\geq 0$,
		$$\displaystyle M_0^+(\geq C,P)=\{m_0\in M_0(F)\mid \langle \alpha, H_{M_0}(m_0)\rangle \geq 0 \; \forall \alpha\in \Delta_{P_0},\;  \langle \alpha, H_{M_0}(m_0)\rangle \geq C \; \forall \alpha\in \Delta_{P_0}\setminus \Delta_{P_0}^P\}.$$
		Then, there exists a compact $\CK_M\subset M(F)$ such that the subsets
		$$\displaystyle \bigsqcup_i x_i M_0^+(\geq C,P) \CK_M,\;\; C>0,$$
		form a basis of neighborhoods of $\infty^M_P$ in $X_M(F)$. Fix $1\leq i\leq n$ and set
		$$\displaystyle \CK_{\overline{U}}:=\bigcap_{k\in \CK_M}k\CK k^{-1}\cap U_{\overline{P}}(F), \;\; \CK':=\bigcup_{k\in \CK_M}k^{-1}Kk,$$
		two compact subsets of $U_{\bar{P}}(F)$ and $G(F)$ respectively. It suffices to show that
		$$\displaystyle x_i M_0^+(\geq C,P) \CK_{\overline{U}}\subset x_i M_0^+(\geq C,P) \CK'$$
		for $C$ large enough. Let $J_{P_0}\subset P_0(F)$ be a compact-open subgroup small enough so that $m_0^{-1}J_{P_0}m_0\subset \CK'$ for every $m_0\in M_0^+$. For every compact-open subgroup $J_{\overline{U}}\subset U_{\bar{P}}(F)$, provided $C$ is large enough we have $m_0\CK_{\overline{U}}m_0^{-1}\subset J_{\overline{U}}$ for every $m_0\in M_0^+(\geq C,P)$. Therefore, it only remains to check that $J_{\overline{U}}$ can be chosen such that $x_iJ_{\overline{U}}\subset x_i J_{P_0}$ but this follows from the fact that $x_iP_0$, which is the image of $HP_0$ by the natural projection $G\to X$, is open in $X$ (so that $x_i J_{P_0}$ contains a neighborhood of $x_i$ in $X(F)$) since $P_0$ is $\iota$-split.
	\end{proof}
	
	A consequence of the previous lemma is that for every $\Omega^M_P\in \CN(\infty_P^M)$ we have $\Omega_P:=\Omega^M_PK\in \CN(\infty_P)$ and moreover that, if $\Omega_P^M$ is sufficiently small, the natural surjection $\Omega^M_P\to \Omega_P/K$ descends to a map $\Omega_P^M/K_M\to \Omega_P/K$ where $K_M$ denotes the image of $K\cap P(F)$ by the natural surjection $P(F)\to M(F)$. We recall the following result from \cite[Theorem 2]{DelPlanch}:
	
	\begin{num}
		\item\label{eq1 neighborhoods} If $\Omega_P^M$ is sufficiently small, the map $\Omega_P^M/K_M\to \Omega_P/K$ is a bijection.
	\end{num}
	
	In particular, if $\Omega_P^M$ is sufficiently small and $A_{P,\iota}^+$-invariant, there exists a map $H_{\overline{P},\iota}: \Omega_P/K\to \CA_{M,\iota}$ and a left action of $A_{P,\iota}^+$ on $\Omega_P/K$ characterized by
	$$\displaystyle H_{\overline{P},\iota}(xK)=H_{M,\iota}(x) \mbox{ and } a\cdot(xK)=(a\cdot x)K$$
	for every $x\in \Omega_P^M$ and $a\in A_{P,\iota}^+$. For simplicity, we will henceforth assume that such a choice of $\Omega_P^M$ has been made for every parabolic $P\in \CP_\iota$ so that if $\Omega_P\in \CN(\infty_P)$ is sufficently small, $H_{\overline{P},\iota}(x)$ and $a\cdot x$ are well-defined for every $x\in \Omega_P/K$ and $a\in A_{P,\iota}^+$. These satisfy
	\begin{equation}\label{eq2 neighborhoods}
		\displaystyle H_{\overline{P},\iota}(a\cdot x)=H_{M,\iota}(a)+H_{\overline{P},\iota}(x)
	\end{equation}
	for every $(a,x)\in A_{P,\iota}^+\times \Omega_P/K$. 
	
	\begin{lem}\label{lem2 neighborhoods}
		Let $P,Q\in \CP_\iota$ with $[P]=[Q]$. Then, for $\Omega\in \CN(\infty_P)=\CN(\infty_Q)$ sufficiently small we have
		$$\displaystyle H_{\overline{Q},\iota}(x)=I_{P,Q}(H_{\overline{P},\iota}(x))+Y^K_{P,Q} \mbox{ and } I_{P,Q}(a)\cdot x=a\cdot x$$
		for every $(a,x)\in A_{P,\iota}^+\times \Omega/K$.
	\end{lem}
	
	\begin{proof}
		Let us choose $\gamma\in G(F)$ such that $\gamma P\gamma^{-1}=Q$, $\gamma M_P\gamma^{-1}=M_Q$. Then, the map $\Omega^{M_P}_P\mapsto \Omega_{P}^{M_P} \gamma^{-1}$ induces a bijection $\CN(\infty_P^{M_P})\simeq \CN(\infty_Q^{M_Q})$ and thus we may assume that there exists a small enough $\Omega_P^{M_P}\in \CN(\infty_P^{M_P})$ such that $\Omega\subset \Omega_P^{M_P}K$ and $\Omega\subset \Omega_Q^{M_Q}K$ where we have set $\Omega_Q^{M_Q}:=\Omega_P^{M_P}\gamma^{-1}$. Thus, an element $x\in \Omega/K$ can both be written as $x=x_PK$ and $x=x_QK$ for $x_P\in \Omega_P^{M_P}$, $x_Q\in \Omega_Q^{M_Q}$. Then, since $x_P\gamma^{-1}\in \Omega_Q^{M_Q}$, provided $\Omega^{M_P}_P$ has been chosen sufficiently small, by Lemma \ref{lem1 neighborhoods} we have
		$$\displaystyle x=(x_P\gamma^{-1})\gamma K=(x_P\gamma^{-1})m_{\overline{Q}}(\gamma)K.$$
		Together with \eqref{eq1 neighborhoods} we get $x_QK_{M_Q}=(x_P\gamma^{-1})m_{\overline{Q}}(\gamma)K_{M_Q}$. By \eqref{eq HM0iota} and Lemma \ref{lem YQP} this implies that
		$$\displaystyle H_{\overline{Q},\iota}(x)=H_{M_Q,\iota}(x_P\gamma^{-1}m_{\overline{Q}}(\gamma))=I_{P,Q}(H_{M_P,\iota}(x_P))+\frac{1}{2} H_{M_Q,\iota}(\gamma \iota(\gamma)^{-1})+H_{\overline{Q},\iota}(\gamma)=I_{P,Q}(H_{P,\iota}(x))+Y_{P,Q}^K$$
		and this shows the first equality of the lemma. For the second one, we notice that if $\Omega^{M_P}_P$ has been chosen sufficiently small and $A_{P,\iota}^+$-invariant, we have
		\[\begin{aligned}
			\displaystyle a\cdot x & =x_PaK=x_Pa\gamma^{-1} \gamma K=x_Pa\gamma^{-1} m_{\overline{Q}}(\gamma) K \\
			& =x_P\gamma^{-1}m_{\overline{Q}}(\gamma)\gamma a\gamma^{-1}K=I_{P,Q}(a)\cdot (x_P\gamma^{-1}m_{\overline{Q}}(\gamma) K)= I_{P,Q}(a)\cdot x
		\end{aligned}\]
		for every $a\in A_{P,\iota}^+$.
	\end{proof}

	\subsection{The map $H_X$}\label{sec HX}
	
	In this subsection, we continue to fix a special maximal compact subgroup $K\subset G(F)$. Recall that for every $P,Q\in \CP_\iota$ with $[P]=[Q]$, we have introduced an element $Y^K_{P,Q}\in \cA_{Q,\iota}$. For two minimal $\iota$-split parabolic subgroups $P_0,P_0'$, we introduce the following affine isomorphism
	$$\displaystyle I^K_{P_0,P_0'}:\cA_{P_0,\iota}\simeq \cA_{P_0',\iota},\;\; I^K_{P_0,P_0'}(H)=I_{P_0,P_0'}(H)+Y^K_{P_0,P_0'}.$$
	These isomorphisms compose well (i.e. for any third $P_0''\in \CP_\iota^{\min}$ we have $I^K_{P_0,P_0''}=I^K_{P_0',P_0''}I^K_{P_0,P_0'}$) and we can introduce the real affine space
	$$\displaystyle \cA_{X,K}:=\varprojlim_{P_0} \cA_{P_0,\iota}$$
	where the transition maps are this time given by the (affine) isomorphisms $I^K_{P_0,P_0'}$. Note that the space of translations of $\cA_{X,K}$ is $\cA_X$ and for every $P_0\in \CP_\iota^{\min}$ there is an affine isomorphism $\cA_{X,K}\simeq \cA_{P_0,\iota}$ compatible with the identification $\cA_X\simeq \cA_{P_0,\iota}$.
	
	Let $P\in \CP_\iota$ be a $\iota$-split parabolic subgroup (not necessarily minimal) and choose $P_0\in \CP_\iota^{\min}$ with $P_0\subset P$. Then, the composition of the isomorphism $\cA_{X,K}\simeq \cA_{P_0,\iota}$ with the projection $\cA_{P_0,\iota}\twoheadrightarrow \cA_{P,\iota}$ is independent on the choice of $P_0$ and will be denoted
	$$\displaystyle proj_P: \cA_{X,K}\to \cA_{P,\iota}$$
	or simply $Y\mapsto Y_P$.
	
	In the following, we fix a norm $\lVert .\rVert$ on the real vector space $\CA_X$ that we transfer to $\CA_{X,K}$ through the choice of (an arbitrary) base-point.
	
	\begin{prop}\label{prop map HX}
		There exists a $K$-invariant map $H_X:X(F)\to \CA_{X,K}$ satisfying the following conditions: for every $P\in \CP_\iota$, there exists a small enough neighborhood $\Omega_P$ of $\infty_P$ in $X(F)$ such that:
		\begin{enumerate}
			\item For every $x\in \Omega_P$, we have $proj_P(H_X(x))=H_{\overline{P},\iota}(x)$;
			
			\item For every $(a,x)\in A_{P,\iota}^+\times \Omega_P/K$, $H_X(a\cdot x)=H_{M_P,\iota}(a)+H_X(x)$.
			
			\item $1+\lVert H_X(x)\rVert\sim \sigma_X(x)$ for $x\in X(F)$.
			
			\item For every $P_0\in \CP_\iota^{\min}$, we can find $Y_{P_0,\iota}^-\in \CA_{P_0,\iota}$ such that $H_X(x)_{P_0}\in Y_{P_0,\iota}^- + \CA_{P_0,\iota}^+$ for every $x\in X(F)$.
		\end{enumerate}
	\end{prop}
	
	\begin{proof}
		Let $P,Q\in \CP_\iota$ be such that $[P]=[Q]$. Then, it readily follows from Lemma \ref{lem2 neighborhoods} that a $K$-invariant map $H_X: X(F)\to \CA_{X,K}$ satisfies conditions 1 and 2 for $P$ if and only if it satisfies the same conditions for $Q$. Similarly, for $P_0,P_0'\in \CP_\iota^{\min}$, condition 4 holds for $P_0$ if and only if it holds for $P_0'$. Therefore, fixing $P_0\in \CP_\iota^{\min}$, it suffices to show the existence of a $K$-invariant map $H_X: X(F)\to \CA_{X,K}$ satisfying conditions 1.-4. for every parabolic subgroup $P\in \CP_\iota$ with $P_0\subset P$. We will call such parabolics standard and we will denote by $\CP_\iota^{std}$ the subset of them. We will also use the identification $\cA_{X,K}=\cA_{P_0,\iota}$. 
		
		For each $P\in \CP_\iota^{std}$ we can find a neighborhood $\Omega_P\in \CN(\infty_P)$ such that:
		\begin{itemize}
			\item For each $P\in \CP_\iota^{std}$, $\Omega_P$ is $A_{P,\iota}^+\times K$-stable and is small enough that the map $H_{\overline{P},\iota}:\Omega_P\to \cA_{P,\iota}$ as well as the action of $A_{P,\iota}^+$ on $\Omega_P/K$ are well-defined;
			
			\item For each $P,Q\in \CP_\iota^{std}$, $\Omega_P\cap \Omega_Q\subseteq \Omega_{P\cap Q}$.
		\end{itemize}
		Then, for each $P\in \CP_\iota^{std}$, we set
		$$\displaystyle \omega_P:=\Omega_P\setminus \bigcup_{Q\subsetneq P} \Omega_Q.$$
		From the second bullet point above, it follows that we have a partition in $K$-invariant subsets
		\begin{equation}\label{eq1 prop HX}
			\displaystyle X(F)=\bigsqcup_{P\in \CP_\iota \\ P_0\subset P} \omega_P.
		\end{equation}
		We define a map $H_X:X(F)\to \cA_{X,K}$ by $H_X(x)=H_{\overline{P},\iota}(x)$ for $x\in \omega_P$. Clearly $H_X$ is $K$-invariant.
		
		Let $P\in \CP_\iota$ be standard and let us check that $H_X$ satisfies conditions 1. and 2. Let $x\in \Omega_P$. By definition of the partition \eqref{eq1 prop HX} there exists a standard $Q\in \CP^{std}_\iota$ with $Q\subset P$ such that $x\in \omega_Q$. Since, by definition of $H_{\overline{P},\iota}$ and $H_{\overline{Q},\iota}$, we have $proj_P H_{\overline{Q},\iota}(x)=H_{\overline{P},\iota}(x)$, condition 1. is immediate. Also, since, by our choice of neighborhoods $(\Omega_P)_P$, $\omega_Q/K$ is invariant by $A_{Q,\iota}^+$, hence also by $A_{P,\iota}^+$, from \eqref{eq2 neighborhoods} we deduce that for every $a\in A_{P,\iota}^+$ we have
		$$\displaystyle H_X(a\cdot x)=H_{\overline{Q},\iota}(a\cdot x)=H_{M_Q,\iota}(a)+H_{\overline{Q},\iota}(x)=H_{M_P,\iota}(a)+H_X(x)$$
		and this proves condition 2.
		
		Let us now check condition 3. First, since $\lVert H_{\overline{P},\iota}(x)\rVert\ll \sigma_X(x)$ for every $P\in \CP_\iota^{std}$ and $x\in \Omega_P$, it follows from the above definition of $H_X$ that we have
		$$\displaystyle \lVert H_X(x)\rVert\ll \sigma_X(x),\;\; \mbox{ for } x\in X(F).$$
		Thus, we just need to prove the converse inequality. By the weak Cartan decomposition \eqref{Weak Cartan}, it suffices to check it for $x\in X_{P_0}^+$. Let $C>0$ that will be assumed large enough in what follows. Let $x\in X_{P_0}^+$ and let $P\in \CP_\iota^{std}$ be such that
		$$\displaystyle \Delta_{P_0,\iota}\setminus \Delta_{P_0,\iota}^P=\{\alpha\in \Delta_{P_0,\iota}\mid \langle \alpha, H_{M_0,\iota}(x)\rangle\geq C \}.$$
		Then, provided $C$ is large enough, we have $x\in \Omega_P$. Hence, by property 1. we have $proj_P H_X(x)=H_{\overline{P},\iota}(x)$. On the other hand, it is easy to see that $\sigma_X(x)\ll 1+\lVert H_{\overline{P},\iota}(x)\rVert$. Hence, $\sigma_X(x)\ll 1+\lVert H_{X}(x)\rVert$ and we are done.
		
		It only remains to prove that $H_X$ satisfies condition 4. Let us fix a weak Cartan decomposition like \eqref{Weak Cartan}. Then, by definition of neighborhoods of $\infty_P$, there exists $C>0$ such that $X_{P_0}^+(\geq C,P)\CK \subset \Omega_P$ for every $P\in \CP_\iota^{std}$. Then, for $P\subset Q$ we have $X_{P_0}^+(\geq C,P)\subseteq X_{P_0}^+(\geq C,Q)$ and the subsets
		$$\displaystyle X_{P_0}^+(\geq C,P)\setminus \bigcup_{Q\subsetneq P} X_{P_0}^+(\geq C,Q)$$
		are relatively compact modulo $A_{P,\iota}^+$. It follows that we can find compact subsets $\omega_P'\subset \Omega_P$ such that
		$$\displaystyle X(F)=\bigcup_{P\in \CP_\iota^{std}} A_{P,\iota}^+\omega_P'.$$
		By property 2., we have $H_X(A_{P,\iota}^+\omega_P')=H_{P,\iota}(A_{P,\iota}^+)+H_X(\omega_P')$ for each $P\in \CP_\iota^{std}$. Since, by property 3., $H_X(\omega_P')\subset \CA_{X,K}$ is relatively compact and $H_{P,\iota}(A_{P,\iota}^+)\subset \CA_{P,\iota}^+\subset \CA_{P_0,\iota}^+$, property 4. follows.
	\end{proof}
	
	\begin{prop}\label{prop2 HX}
		Let $H_X:X(F)\to \CA_{X,K}$ be the map as in the previous proposition. Let $P\in \CP_\iota$ and set $M=M_P$. Then, there exists $c>0$ such that for every $(a,x)\in A_{P,\iota}^+\times X_{M}(F)$, we can find $Q\in \CF_\iota(P)$ such that
		\begin{equation}\label{eq1 prop2 HX}
			\displaystyle proj_Q H_X(ax)=H_{M_Q,\iota}(ax),
		\end{equation}
		
		\begin{equation}\label{eq2 prop2 HX}
			\displaystyle \lVert H_X(ax)-proj_Q H_X(ax)\rVert\leqslant c \sigma_X(x),
		\end{equation}
		and
		\begin{equation}\label{eq3 prop2 HX}
			\displaystyle \lVert H_{M_Q,\iota}(ax)-H_{M_P,\iota}(ax)  \rVert\leqslant c \sigma_X(x).
		\end{equation}
	\end{prop}
	
	\begin{proof}
		We prove this by induction on $\dim(A_P)$. So we assume that the statement holds for $P$ replaced by any parabolic $R\in \CP_\iota$ with $P\subsetneq R$. Let $\Omega^{M}_P\in \CN(\infty^{M}_P)$ be a small enough $A_{P,\iota}^+$-stable neighborhood of $\infty^{M}_P$ in $X_M(F)$ such that the first and second points of the previous proposition are satisfied for $\Omega_P=\Omega^M_P K$ and $H_{\overline{P},\iota}(x)=H_{M,\iota}(x)$ for every $x\in \Omega_P^M$. Then, there exists a constant $c_1>0$ such that for every $(a,x)\in A_{P,\iota}^+\times X_{M}(F)$ satisfying
		$$\displaystyle \langle \alpha, H_{M,\iota}(a)\rangle \geq c_1\sigma_X(x), \; \mbox{ for every } \alpha\in \Delta_{P,\iota},$$
		we have $ax\in \Omega^M_P$. (This follows e.g. from using a weak Cartan decomposition for $X_{M}(F)$.) Moreover, there also exists a constant $c'>c_1$ such that for any such $a$ and $x$ we can find $a'\in A_{P,\iota}^+$ with $\sigma(a')\leq c'\sigma_X(x)$, $a'x\in \Omega^M_P$ and $a\in a'A_{P,\iota}^+$. From this and Proposition \ref{prop map HX}, we get
		$$\displaystyle proj_P H_X(ax)=H_{\overline{P},\iota}(ax)=H_{M,\iota}(ax)=H_{M,\iota}(a)+H_{M,\iota}(x)$$
		and
		$$\displaystyle H_X(ax)=H_X(a(a')^{-1}a'x)=H_{M,\iota}(a(a')^{-1})+H_X(a'x)=H_{M,\iota}(a)-H_{M,\iota}(a')+H_X(a'x).$$
		Hence,
		$$\displaystyle \lVert H_X(ax)-proj_P H_X(ax)\rVert= \lVert H_X(a'x)-H_{M,\iota}(a')-H_{M,\iota}(x)\rVert.$$
		Since $\lVert H_{M,\iota}(a')\rVert\ll \sigma_X(x)$, $\lVert H_X(a'x)\rVert\ll \sigma_X(x)$ and $\lVert H_{M,\iota}(x)\rVert\ll \sigma_X(x)$, we get that in this case both \eqref{eq1 prop2 HX} and \eqref{eq2 prop2 HX} holds for $Q=P$ and a suitable constant $c>0$ (note that \eqref{eq3 prop2 HX} is trivial when $P=Q$). 
		
		
		It remains to treat the case where there exists $\alpha\in \Delta_{P,\iota}$ such that $\langle \alpha, H_{M,\iota}(a)\rangle < c_1\sigma_X(x)$. Let $R\in \CP_\iota$ be the unique $\iota$-split parabolic subgroup such that $R\supset P$ and $\Delta_{P,\iota}^{R}=\{\alpha \}$. Then, for every $(a,x)\in A_{P,\iota}^+\times X_{M}(F)$ with $\langle \alpha, H_{M,\iota}(a)\rangle < c_1\sigma_X(x)$ we can find $a_R\in A_{R,\iota}^+$ such that $\sigma(aa_R^{-1})\ll \sigma_X(x)$. Then, writing $ax=a_R(a_R^{-1}ax)$ where $a_R^{-1}ax\in X_M(F)\subset X_{M_R}(F)$, by the induction hypothesis there exists a constant $c_R>0$ as well as $Q\in \CP_\iota$, $Q\supset R$ such that $proj_Q H_X(ax)=H_{M_Q,\iota}(ax)$ and
		$$\displaystyle \lVert H_X(ax)-proj_Q H_X(ax)\rVert\leq c_R \sigma_X(a_R^{-1}ax),\;\displaystyle \lVert H_{M_Q,\iota}(ax)-H_{M_R,\iota}(ax)  \rVert\leq c_R \sigma_X(a_R^{-1}ax).$$
		Since $\sigma_X(a_R^{-1}ax)\ll \sigma_X(x)$ and 
		$\lVert H_{M_R,\iota}(ax)-H_{M_P,\iota}(ax)| \rVert \ll \sigma_X(x),$
		this again gives \eqref{eq1 prop2 HX}, \eqref{eq2 prop2 HX} and \eqref{eq3 prop2 HX} for a suitable constant $c>0$ and the proposition is proved.

	\end{proof}
	
	\subsection{Twisted symmetric pairs}\label{sect twisted symmetric pairs}
	
	We define a {\it twisted symmetric pair} (over $F$) to be a triple $(G,\tG,\iota)$ where $(G,\tG)$ is a linear twisted space, $(G,\iota)$ is a symmetric pair both defined over $F$ and we have extended $\iota$ to an involutive automorphism of $\iota: \tG\to \tG$ (still defined over $F$) with $\iota(g_1\gamma g_2)=\iota(g_1)\iota(\gamma) \iota(g_2)$ for every $(\gamma,g_1,g_2)\in \tG\times G\times G$. We will usually refer to twisted symmetric pairs by $(\tG,\iota)$, the underlying group $G$ being implicit.
	
	Let $(\tG,\iota)$ be a twisted symmetric pair. We denote by $A_{\tG,\iota}$ the neutral component of the subgroup $\{a\in A_{\tG}\mid \iota(a)=a^{-1} \}$ and we set
	$$\displaystyle \cA_{\tG,\iota}^*=X^*(A_{\tG,\iota})\otimes \mathbb{R},\;\;\; \cA_{\tG,\iota}=X_*(A_{\tG,\iota})\otimes \mathbb{R}.$$
	Then, $\cA_{\tG,\iota}^*$ (resp. $\cA_{\tG,\iota}$) can be identified with the subspace of $\iota$-antiinvariant vectors in $\cA_{\tG}^*$ (resp. in $\cA_{\tG}$). We also denote by $H_{\tG,\iota}: G(F)\to \cA_{\tG,\iota}$ the composition of $H_{\tG}$ with the natural projection $\cA_{\tG}\to \cA_{\tG,\iota}$.
	
	We assume from now on that $G$ is connected and reductive. Let $H=G^\iota$, $\tH=\tG^\iota$ be the subvarieties of $\iota$-fixed points in $G$ and $\tG$ respectively. Then, $(H,\tH)$ is a reductive twisted space over $F$ and we will always assume that $\tH(F)\neq \emptyset$.
	
	A parabolic subspace $\widetilde{P}\subset \widetilde{G}$ is called {\em $\iota$-split} if the underlying parabolic subgroup $P\subset G$ is $\iota$-split or equivalently if $\iota(\tP)$ is a parabolic subspace opposite to $\tP$. Similarly, a Levi subspace $\tM\subset \tG$ is said to be {\em $\iota$-split} if there exists a $\iota$-split parabolic subspace $\tP$ such that $\tM=\tP \cap \iota(\tP)$. For $\tM$ a $\iota$-split Levi subspace, we denote by $\mathcal{P}_\iota(\tM)$ (resp. $\mathcal{F}_\iota(\tM)$, resp. $\mathcal{L}_\iota(\tM)$) the set of $\iota$-split parabolic subspaces having $\tM$ as a Levi component (resp. containing $\tM$, resp. the set of $\iota$-split Levi subspaces containing $\tM$).

	We equip $\mathcal{A}_{\tM,\iota}$ with the unique Haar measure for which the lattice $H_{\tM,\iota}(A_{\tM}(F))$ is of covolume one. We will also write $\mathcal{A}_{\tP,\iota}$ for $\mathcal{A}_{\tM,\iota}$ for every $\tP\in \mathcal{P}_\iota(\tM)$ and we denote by $H_{\tP,\iota}: P(F)\to \cA_{\tP,\iota}$ the composition of the projection $P(F)\to M(F)$ with $H_{\tM,\iota}$.
	
	We will denote by $\widetilde{\CP}_\iota$ and $\widetilde{\CL}_\iota$ the sets of all $\iota$-split parabolic subspaces and $\iota$-split Levi subspaces of $\tG$ respectively. 
	
	Let $\tM$ be a $\iota$-split Levi subspace of $\tG$, $\tQ\in \mathcal{F}_\iota(\tM)$ and set $\cA^{\tQ}_{\tM,\iota}=\cA_{\tM,\iota}/\cA_{\tQ,\iota}$. We equip this space with the quotient of the Haar measures on $\cA_{\tM,\iota}$ and $\cA_{\tQ,\iota}$. For $\tP\in \mathcal{P}_\iota(\tM)$ with $\tP\subset \tQ$, we denote by
	$$\displaystyle \Delta^{\tQ,\vee}_{\tP,\iota}\subset \cA^{\tQ}_{\tM,\iota} \mbox{ and } \Delta^{\tQ}_{\tP,\iota}\subset \cA^{\tQ,*}_{\tM,\iota}$$
	the images of $\Delta^{Q,\vee}_{P,\iota}$ and $\Delta_{P,\iota}^Q$ by the natural projections
	$$\displaystyle \cA^Q_{M,\iota}\twoheadrightarrow \cA^{\tQ}_{\tM,\iota} \mbox{ and } \cA^{Q,*}_{M,\iota}\twoheadrightarrow \cA^{\tQ,*}_{\tM,\iota}$$
	respectively. These form basis of $\cA^{\tQ}_{\tM,\iota}$ and $\cA^{\tQ,*}_{\tM,\iota}$ respectively and we write $\widehat{\Delta}^{\tQ}_{\tP,\iota}\subseteq \cA^{\tQ,*}_{\tM,\iota}$ for the basis dual to $\Delta^{\tQ,\vee}_{\tP,\iota}$. We denote by $\tau_{\tP,\iota}^{\tQ}$, $\widehat{\tau}_{\tP,\iota}^{\tQ}$ the characteristic functions of the cone in $\cA$ characterized by 
	$$\displaystyle \tau_{\tP,\iota}^{\tQ}(H)=1\iff \langle \alpha,H\rangle>0,\; \forall \alpha\in  \Delta^{\tQ}_{\tP,\iota},\; \widehat{\tau}_{\tP,\iota}^{\tQ}(H)=1\iff \langle \varpi,H\rangle>0,\; \forall \varpi\in  \widehat{\Delta}^{\tQ}_{\tP,\iota}$$ respectively. When $\tQ=\tG$ we will sometimes drop the superscript $\tQ$.
	
	Recall that $\CP_\iota^{\min}$ stands for the set of all minimal $\iota$-split parabolic subgroups of $G$. We have
	\begin{num}
		\item\label{eqminimaliotasplit} For every $P_0\in \CP_\iota^{\min}$, $\tP_0:=Norm_{\tG}(P_0)$ is a $\iota$-split parabolic subspace of $\tG$ i.e. $\tP_0(F)\neq \emptyset$.
	\end{num}
	Indeed, we just need to check the existence of an element $\gamma\in \tG(F)$ such that the parabolic subgroups $P_0$ and $\Ad_\gamma(P_0)$ are in the same conjugacy class. However, for $\gamma\in \tH(F)$ the parabolic subgroup $\Ad_\gamma(P_0)$ is also $\iota$-split minimal and by \cite[Proposition 4.9]{HW} all minimal $\iota$-split parabolic subgroups are in the same conjugacy class.
	
	Let $\widetilde{\CP}_\iota^{\min}\subset \widetilde{\CP}_\iota$ be the subset of minimal elements of $\widetilde{\CP}_\iota$ (for the inclusion relation). Then, by \eqref{eqminimaliotasplit}, the map $P_0\mapsto \tP_0$ gives a bijection $\CP_\iota^{\min}\simeq \widetilde{\CP}_\iota^{\min}$.
	
	Set $\tM_0=\tP_0\cap \sigma(\tP_0)$. It is easy to see that the automorphism $\theta$ of $\cA_{M_0,\iota}^{G,*}$ preserves the root system $\Sigma_{0,\iota}$ as well as its subset of simple roots $\Delta_{0,\iota}$.

	Let $X=H\backslash G$ be the homogeneous symmetric variety associated to $(G,\iota)$. Then, there exists a unique regular map $X\times \tG\to X$, $(x,\gamma)\mapsto x\gamma$ such that $(Hg)\gamma=H\Ad_\gamma^{-1}(g)$ for every $(g,\gamma)\in G\times \tH$. Note that we have
	$$\displaystyle ((xg_1)\gamma)g_2=x(g_1\gamma g_2), \mbox{ for every } (x,g_1,g_2,\gamma)\in X\times G\times G\times \tG.$$
	We will usually write $\tX$ to mean $X$ equipped with this ``twisted action'' of $\tG$. This twisted action naturally induces an automorphism $\theta$ of the real vector space $\CA_X$ and we set
	$$\displaystyle \CA_{\tX}:=\CA_X^\theta.$$
	Then, for any $P_0\in \CP_\iota^{\min}$, the canonical isomorphism $\CA_X\simeq \CA_{P_0,\iota}$ induces an isomorphism $\CA_{\tX}\simeq \CA_{\tP_0,\iota}$. We will also write $\CA_{\tX}^+$, ${}^- \CA_{\tX}$ for the respective images of $\CA_X^+$, ${}^- \CA_X$ by the natural projection $\cA_X\to \CA_{\tX}$ and $\phi_{\tX}$ for the characteristic function of the cone ${}^- \CA_{\tX}$.
	
	
	Let $K\subset G(F)$ be a special maximal compact subgroup and let $\CA_X(1-\theta)$ be the kernel of the natural projection  $\cA_X\to \CA_{\tX}$. We set
	$$\displaystyle \CA_{\tX,K}:= \CA_{X,K}/\CA_X(1-\theta).$$
	It is an affine space with direction $\CA_{\tX}$ and for every $P_0\in \CP_\iota^{\min}$ the (affine) isomorphism $\CA_{X,K}\simeq \CA_{P_0,\iota}$ induces an isomorphism $\CA_{\tX,K}\simeq \CA_{\tP_0,\iota}$. Moreover, for every $\tP\in \widetilde{\CP}_\iota$ there is a natural affine projection
	\begin{equation}\label{eqprojAX}
		\displaystyle \CA_{\tX,K}\to \CA_{\tP,\iota}
	\end{equation}
	that can be described as the composition of $\CA_{\tX,K}\simeq \CA_{\tP_0,\iota}$ with the projection $\CA_{\tP_0,\iota}\to \CA_{\tP,\iota}$ for any $\tP_0\in \widetilde{\CP}_\iota$ with $\tP_0\subset \tP$. For every $Y\in \CA_{\tX,K}$, we will denote by $Y_{\tP}$ its image by the projection \eqref{eqprojAX}.
	
	If $H_X: X(F)\to \CA_{X,K}$ is a map as in Proposition \ref{prop map HX}, we will usually write $H_{\tX}:\tX(F)\to \CA_{\tX,K}$ for the composition of $H_X$ with the natural projection $\CA_{X,K}\to \CA_{\tX,K}$.
	
	\subsection{Orthogonal sets}\label{orthogonal sets twist}

	Let $(\tG,\iota)$ be a twisted symmetric pair. In \cite[\S 2.8.2]{BeuGalP}, we have introduced notions of $(G,M,\iota)$-families extending in an obvious way Arthur's definition of $(G,M)$-families in the context of symmetric pairs. This actually exactly corresponds to Arthur's theory applied to the root system $\Sigma_{0,\iota}$. There is a similar combinatorics for twisted groups as developed in \cite{LW} which can be applied to any automorphism of a root system preserving a positive system. In particular, starting from the pair $(\Sigma_{0,\iota},\theta)$ there is a corresponding notion of $(\tG,\tM,\iota)$-orthogonal sets that we now briefly describe.
	
	Let $\tM$ be a $\iota$-split Levi subspace of $\tG$. Two parabolic subspaces $\tP, \widetilde{Q}\in \mathcal{P}_\iota(\tM)$ are said to be {\em $\iota$-adjacent} if the intersection $\Delta_{\tP,\iota}^\vee\cap -\Delta_{\widetilde{Q},\iota}^\vee$ is a singleton $\{\alpha_{\tP,\tQ}^\vee\}$. If this is the case, the hyperplane $\{X\in i\mathcal{A}_{\tM,\iota}^*\mid \langle \alpha_{\tP,\tQ}^\vee,X\rangle=0 \}$ is called {\em the wall separating} $\tP$ and $\widetilde{Q}$.
	
	By definition {\em $(\tG,\tM,\iota)$-orthogonal set} is a family $\mathcal{X}=(X_{\tP,\iota})_{\tP\in \mathcal{P}_\iota(\tM)}$ of points in $\mathcal{A}_{\tM,\iota}$ such that for every $\iota$-adjacent parabolic subspaces $\tP, \widetilde{Q}\in \mathcal{P}_\iota(\tM)$, we have
	$$\displaystyle X_{\tP,\iota}-X_{\tQ,\iota}\in \mathbb{R}\alpha_{\tP,\tQ}^\vee$$
	where $\Delta_{\tP,\iota}^\vee\cap -\Delta_{\widetilde{Q},\iota}^\vee=\{ \alpha_{\tP,\tQ}^\vee\}$. We further say that $\mathcal{X}$ is {\em positive} if
	$$\displaystyle X_{\tP,\iota}-X_{\tQ,\iota}\in \mathbb{R}_{\geqslant 0}\alpha_{\tP,\tQ}^\vee$$
	for every pair of $\iota$-adjacent parabolic subspaces $\tP, \widetilde{Q}\in \mathcal{P}_\iota(\tM)$.
	
	As in Subsection \ref{orthogonal sets}, we define the {\it depth} and the {\it norm} of a $(\tG,\tM,\iota)$-orthogonal set $\mathcal{X}=(X_{\tP})_{\tP\in \mathcal{P}(\tM)}$ by
	\begin{equation*}
		\displaystyle d(\mathcal{X})=\min_{\tP\in \mathcal{P}(\tM)} \min_{\alpha\in \Delta_{\tP}} \alpha(X_{\tP}) \mbox{ and } N(\mathcal{X})=\max_{\tP\in \mathcal{P}(\tM)} \max_{\alpha\in \Delta_{\tP}} \lvert \alpha(X_{\tP})\rvert
	\end{equation*}
respectively. Note that $\mathcal{X}$ is positive if and only if $d(\mathcal{X})\geq 0$.

	Let $\mathcal{X}=(X_{\tP,\iota})_{\tP\in \mathcal{P}_\iota(\tM)}$ be a $(\tG,\tM,\iota)$-orthogonal set. For $\tQ=\tL U_Q\in \mathcal{F}_\iota(\tM)$, we denote by $X_{\tQ,\iota}$ the projection to $\cA_{\tL,\iota}$ of $X_{\tP,\iota}$ for any $\tP\in \mathcal{P}_\iota(\tM)$ such that $\tP\subset \tQ$ (this projection does not depend on the choice of $\tP$). To $\mathcal{X}$ we associate functions $\Gamma^{\tQ}_{\tL,\iota}(.,\mathcal{X})$ on $\cA_{\tL,\iota}^{\tQ}$ and complex numbers $v^{\tQ}_{\tL,\iota}(\mathcal{X})\in \mathbb{C}$ for every $\tL\in \mathcal{L}_\iota(\tM)$ and $\tQ\in \mathcal{F}_\iota(\tL)$ as follows:
	$$\displaystyle \Gamma^{\tQ}_{\tL,\iota}(H,\mathcal{X})=\sum_{\tP\in \mathcal{F}_\iota(\tL), \tP\subset \tQ} (-1)^{a^{\tQ}_{\tP,\iota}} \widehat{\tau}_{\tP,\iota}^{\tQ}(H-X_{\tP,\iota}),\;\;\; H\in \cA_{\tL,\iota}^{\tQ}$$
	and
	$$\displaystyle v^{\tQ}_{\tL,\iota}(\mathcal{X})=\int_{\cA_{\tL,\iota}^{\tQ}} \Gamma^{\tQ}_{\tL,\iota}(H,\mathcal{X}) dH.$$
	If $\mathcal{X}$ is positive, $v^{\tQ}_{\tL,\iota}(\mathcal{X})$ is simply the volume of the convex hull of the family $(X_{\tP,\iota})_{\tP\in \mathcal{P}_\iota(\tL), \tP\subset \tQ}$. Once again, we will sometimes drop the superscript when $\tQ=\tG$. 
	
	
	Let $K$ be a special compact subgroup of $G(F)$. Using the Iwasawa decomposition $G(F)=P(F)K$, for every $\iota$-split parabolic subspace $\tP\subset \tG$, we can extend the homomorphism $H_{\tP,\iota}$ to a map $G(F)\to \cA_{\tP,\iota}$. Then, for every $\iota$-split Levi subspace $\tM\subset \tG$ and $g\in G(F)$, the family $\mathcal{H}_{\tM,\iota}(g)=(-H_{\tP,\iota}(g))_{\tP\in \mathcal{P}_\iota(\tM)}$ is a positive $(\tG,\tM,\iota)$-orthogonal set and we define
	$$\displaystyle v_{\tM,\iota}^{\tQ}(g)=v_{\tM,\iota}^{\tQ}(\mathcal{H}_{\tM,\iota}(g)),\;\; \mbox{ for } \tQ \in \mathcal{F}_\iota(\tM).$$

	Let $Y\in \CA_{\tX,K}$ and $\tM\in \widetilde{\CL}_\iota$. Then, the family $\mathcal{Y}_{\tM}:=(Y_{\tP})_{\tP\in \CP_\iota(\tM)}$ is $(\tG,\tM,\iota)$-orthogonal set. Indeed, this follows from the fact that if $\tP,\tQ\in \mathcal{P}_{\iota}(\tM)$ are $\iota$-adjacent then $Y^K_{Q,P}\in \mathbb{R} \alpha^\vee_{\tP,\tQ}$ as can be directly checked on the definition.
	
	Whenever convenient, we will also fix a minimal $\iota$-split parabolic subspace $\tP_0\subset \tG$ to define the {\it depth} and {\it norm} of an element $Y\in \CA_{\tX,K}$ by
	$$\displaystyle d(Y)=\min_{\alpha\in \Delta_{\tP_0,\iota}} \alpha(Y_{\tP_0,\iota}) \mbox{ and } N(Y)=\max_{\alpha\in \Delta_{\tP_0,\iota}} \lvert \alpha(Y_{\tP_0,\iota})\rvert$$
	respectively. We note that for every $\iota$-split Levi subspace $\tM\subset \tG$, there exist constants $c_1,c_2>0$ such that
	$$\displaystyle d(Y)-c_1\leq d\left(\mathcal{Y}_{\tM}\right) \mbox{ and } N\left( \mathcal{Y}_{\tM}\right)\leq N(Y)+c_2$$
	for every $Y\in \CA_{\tX,K}$.

	\subsection{$\iota$-weighted orbital integrals}\label{Section iota WOI}
	
	Let $\tM$ be a $\iota$-split Levi subspace of $\tG$, $\gamma\in \tM(F)\cap \tG_{rs}(F)$ and $\tQ\in \mathcal{F}_\iota(\tM)$. For $f\in \cC(\tG(F))$, we define the {\em $\iota$-twisted weighted orbital integral}
	\begin{equation*}
		\displaystyle \Phi_{\tM,\iota}^{\tQ}(\gamma,f)=\int_{G_\gamma(F)\backslash G(F)} f(g^{-1}\gamma g) v_{\tM,\iota}^{\tQ}(g)dg
	\end{equation*}
	as well as its normalized version
	\begin{equation*}
		\displaystyle J_{\tM,\iota}^{\tQ}(\gamma,f)=D^{\tG}(\gamma)^{1/2} \Phi_{\tM,\iota}^{\tQ}(\gamma,f).
	\end{equation*}
	By the same argument as in \eqref{eq1 WOI}, the above integral is absolutely convergent, and for $\tT\subset \tM$ a maximal twisted torus, we have:
	\begin{num}
		\item\label{eq1 WOI iota} There exist $p>0$ and, for every $d>0$, a continuous semi-norm $\nu_d$ on $\cC(\tG(F))$ such that
		$$\displaystyle \left\lvert J_{\tM,\iota}^{\tQ}(\gamma,f)\right\rvert\leqslant \nu_d(f) (1+\lvert \log D^{\tG}(\gamma)\rvert)^p \sigma_{\tT_{/\theta}}(\gamma)^{-d}$$
		for every $\gamma\in \tT_{\reg}(F)$ and $f\in \cC(\tG(F))$.
	\end{num}

	\section{Harmonic analysis for certain singular conjugacy classes}
	
	In this section, we fix a connected reductive twisted space $\tG$ as well as a semisimple element $x\in \tG_{ss}(F)$ and a regular nilpotent coadjoint orbit $\CO\in \Nil_{\reg}(\fg_x^*)$. Our goal is to establish an alternative description of the distribution
	$$\displaystyle f\in \CC_{\scusp}(\tG(F))\mapsto c_{f,\CO}(x)$$
	that was defined in Subsection \ref{Sect strongly cuspi}.
	
	 Note that since $G_x$ admits a regular coadjoint orbit, it is quasi-split. We fix once and for all a Borel subgroup $B_x$ of $G_x$ with a Levi decomposition $B_x=T_xN_x$. Let $Y\in \CO$ whose stabilizer in $N_x$ is trivial. It defines a generic character $\xi$ of $N_x(F)$ by the formula
	 \begin{equation}\label{eq1introchap4}
	 \displaystyle \xi(\exp(X))=\psi(\langle X,Y\rangle),\;\; X\in \mathfrak{n}_x(F).
	 \end{equation}
	 We then say that the generic character $\xi$ is {\it associated} to the nilpotent coadjoint orbit $\cO$. The formula we will give for the distribution $f\mapsto c_{f,\cO}(x)$ will be in term of this character $\xi$. Note that there are more than one generic characters of $N_x(F)$ associated to $\cO$ but they form a unique $T_x(F)$-orbit and the resulting description of $c_{f,\cO}(x)$ will be easily seen to be independent of the choice of $\xi$. Conversely, given a generic character $\xi: N_x(F)\to \mathbb{C}^\times$, we can always find a regular nilpotent element $Y\in \fg_x^*(F)$ such that $\xi$ is given by formula \eqref{eq1introchap4} and moreover such $Y$ is unique up to $N_x(F)$-conjugacy. We will then say that the (nilpotent) coadjoint orbit of $Y$ is {\it associated} to the generic character $\xi$ and we will dentoe it by $\cO_\xi$.

	\subsection{The function $\Gamma_{B_x}(.,\mathcal{X})$}
	Let $A_x\subset T_x$ be the maximal split subtorus, $\CA_x=X_*(A_x)\otimes \RR$, $M(x)=Z_G(A_x)$ and $\tM(x)=M(x)x$, a Levi subspace of $\tG$. We let $\mathcal{P}_{B_x}(\tM(x))$ (resp. $\mathcal{F}_{B_x}(\tM(x))$) be the set of parabolic subspaces $\tP\in \mathcal{P}(\tM(x))$ (resp. $\tQ\in \mathcal{F}(\tM(x))$) such that $P\cap G_x=B_x$ (resp. $Q\supset B_x$). Let also
	$$\displaystyle W_x=Norm_{G_x(F)}(T_x)/T_x$$
	be the Weyl group of $T_x$ in $G_x$.
	
	\begin{lem}
		\begin{enumerate}[(i)]
			\item We have $A_{\tM(x)}=A_x$.
			\item There is a natural embedding of $W_x$ into the Weyl group $W(G,\tM(x))=Norm_{G(F)}(\tM(x))/M(x)(F)$ and we have a partition
			\begin{equation}\label{partition Wx}
				\displaystyle \mathcal{P}(\tM(x))=\bigsqcup_{w\in W_x} w\mathcal{P}_{B_x}(\tM(x)).
			\end{equation}
		\end{enumerate}
	\end{lem}
	
	\begin{proof}
		\begin{enumerate}[(i)]
			\item Since $A_x$ centralizes with $M(x)$ and $x$, it centralizes $\tM(x)$ and therefore $A_x\subset A_{\tM(x)}$. On the other hand, every element of $A_{\tM(x)}$ centralizes $x$ (so that $A_{\tM(x)}\subset G_x$) and $T_x$. It follows that $A_{\tM(x)}\subset Z_{G_x}(T_x)=T_x$ and finally $A_{\tM(x)}\subset A_x$.
			
			\item Since every element of $Norm_{G_x}(T_x)$ centralizes $x$ and normalizes $A_x$, $Norm_{G_x}(T_x)$ is contained in the normalizer of $Z_{G}(A_x)x=\tM(x)$ i.e. $Norm_{G_x}(T_x)\subset Norm_{G}(\tM(x))$. Moreover, $Norm_{G_x}(T_x)\cap M(x)$ is equal to $T_x$ because the centralizer of $A_x$ in $G_x$ is $T_x$. This explains the ``natural'' embedding $W_x\hookrightarrow W(G,\tM(x))$.
			
			We have
			$$\displaystyle w\mathcal{P}_{B_x}(\tM(x))=\{\tP\in \mathcal{P}(\tM(x))\mid P\cap G_x=wB_x \}$$
			so that \eqref{partition Wx} is just the partition corresponding to the fibers of the map $\mathcal{P}(\tM(x))\to \mathcal{P}^{G_x}(T_x)$, $\tP\mapsto P\cap G_x$.
		\end{enumerate}
	\end{proof}
	
	By the above lemma, we have a containment of set of roots 
	$$\displaystyle \Sigma(A_x,G_x)\subset \Sigma(A_x,G)=\Sigma(A_{\tM(x)},G).$$
	Thus, for every $\alpha\in \Sigma(A_x,G_x)$ there are a priori two associated coroots $\alpha^\vee_1,\alpha^\vee_2\in \cA_x$. Namely, we can either see $\alpha$ as a root of $A_x$ in $G_x$ and consider the corresponding coroot $\alpha_1^\vee\in X_*(A_x)$ or we can view $\alpha$ as a root of $A_{\tM(x)}$ in $G$ and consider the corresponding coroot $\alpha_2^\vee\in \cA_x$. It turns out that $\alpha_1^\vee$ and $\alpha_2^\vee$ are always positively proportional. This can be seen as follows. Take a maximal split torus $A_x\subset A_{\min}\subset G$ and fix an inner product on $\cA_{\min}:=\cA_{A_{\min}}$ which is invariant under the action of the Weyl group $W_{\min}=\Norm_{G(F)}(A_{\min}(F))/\Cent_{G(F)}(A_{\min})$. By restriction to $\cA_x\subset \cA_{\min}$, this gives an inner product on $\cA_x$, hence an identification $\cA_x\simeq \cA_x^*$ such that for every $\alpha\in \Sigma(A_{\tM(x)},G)$, $\alpha^\vee_2$ is positively proportional to $\alpha$. This inner product is still $W(G,\tM(x))$-invariant, hence $W_x$-invariant by the second point of the above lemma. It follows that, for every $\alpha\in \Sigma(A_x,G_x)$, the identification $\cA_x\simeq \cA_x^*$ also sends $\alpha^\vee_1$ to a positive multiple of $\alpha$. Since in what follows, the coroots will only matter up to a positive scalar, we will not really have to distinguish between $\alpha_1^\vee$ and $\alpha_2^\vee$. However, to fix ideas, when there is an ambiguity we will always use $\alpha_2^\vee$ instead of $\alpha_1^\vee$.
	
	Let $\mathcal{X}=(X_{\tP})_{\tP\in \mathcal{P}(\tM(x))}$ be a $(\tG,\tM(x))$-orthogonal set in $\CA_{x}=\CA_{\tM(x)}$. For $\tQ\in \mathcal{F}_{B_x}(\tM(x))$ and $H\in \CA_x$, we set (the function $\widehat{\tau}_{\tP}^{\tQ}$ is defined in Section \ref{orthogonal sets})
	\begin{equation*}
		\displaystyle \Gamma_{B_x}^{\tQ}(H,\mathcal{X})=\sum_{\tP\in \mathcal{F}_{B_x}(\tM(x)), \tP\subset \tQ} (-1)^{a_{\tP}^{\tQ}}\widehat{\tau}_{\tP}^{\tQ}(H-X_{\tP}).
	\end{equation*}
	
	For the next lemma, we recall, for two parabolic subspaces $\tQ,\tR\in \CF(\tM(x))$ with $\tQ\subset \tR$, the function
	$$\displaystyle \Gamma_{\tQ}^{\tR}: \CA_{\tQ} \times \CA_{\tQ}\to \mathbb{R}$$
	that can be defined by
	$$\displaystyle \Gamma_{\tQ}^{\tR}(H,X)=\sum_{\tQ\subset \tP\subset \tR} (-1)^{a_{\tP}^{\tR}} \tau_{\tQ}^{\tP}(H) \widehat{\tau}_{\tP}^{\tR}(H-X).$$
	Then, provided $X\in \CA_{\tQ}^+$, we have \cite[lemme 1.8.3]{LW}
	\begin{equation}\label{eq Gammafunction}
	\displaystyle \Gamma_{\tQ}^{\tR}(H,X)=\tau_{\tQ}^{\tR}(H) \phi_{\tQ}^{\tR}(H-X)
	\end{equation}
where $\phi_{\tQ}^{\tR}$ is the characteristic function of those $Y\in \CA_{\tQ}$ such that $\langle \varpi, Y\rangle \leq 0$ for every $\varpi\in \widehat{\Delta}_{\tQ}^{\tR}$ (in other words it is the characteristic function with support the closure of the support of the function $X\mapsto \widehat{\tau}_{\tQ}^{\tR}(-X)$).

	\begin{lem}\label{splitting formula}
		For two $(\tG,\tM(x))$-orthogonal sets $\mathcal{X}$ and $\mathcal{Y}$, we have 
		$$\displaystyle \Gamma_{B_x}^{\tR}(H,\mathcal{X}+\mathcal{Y})=\sum_{\tQ\in \CF_{B_x}(\tM(x)),\tQ\subset \tR} \Gamma_{B_x}^{\tQ}(H,\mathcal{X})\Gamma_{\tQ}^{\tR}(H-X_{\tQ},Y_{\tQ}).$$
	\end{lem}
	
	\begin{proof}
		The proof follows from the same argument as Lemma 1.8.6 of \cite{LW}.
	\end{proof}

	\begin{prop}\label{prop partition gammaBx functions}
		For every $H\in \CA_x$ and $\tR\in \mathcal{F}_{B_x}(\tM(x))$ , we have (the function $\tau_{\tQ}^{\tR}$ is defined in Section \ref{orthogonal sets})
		\begin{equation}\label{ident1 GammaBx}
			\displaystyle \sum_{\tR\supset \tQ\in \mathcal{F}_{B_x}(\tM(x))} \Gamma_{B_x}^{\tQ}(H,\mathcal{X}) \tau^{\tR}_{\tQ}(H-X_{\tQ})=1.
		\end{equation}
		Moreover, if $\mathcal{X}$ is {\em positive} the function $H\mapsto \Gamma_{B_x}(H,\mathcal{X})$ is the characteristic function of the set of $H\in \CA_x$ such that
		$$\displaystyle \varpi_\alpha(H-X_{\tP})\leqslant 0$$
		for every $\tP\in \mathcal{P}_{B_x}(\tM(x))$ and $\alpha\in \Delta_{\tP}$.
	\end{prop}
	
	\begin{proof}
		Let $\tR\in \mathcal{F}_{B_x}(\tM(x))$. By definition of $\Gamma_{B_x}^{\tQ}(.,\mathcal{X})$, for $H\in \CA_x$, we have
		\[\begin{aligned}
			\displaystyle & \sum_{\tR\supset \tQ\in \mathcal{F}_{B_x}(\tM(x))} \Gamma_{B_x}^{\tQ}(H,\mathcal{X}) \tau^{\tR}_{\tQ}(H-X_{\tQ}) \\
			& = \sum_{\tR \supset \tQ\in \mathcal{F}_{B_x}(\tM(x))} \sum_{\mathcal{F}_{B_x}(\tM(x))\ni \tP\subset \tQ}(-1)^{a_{\tP}^{\tQ}} \widehat{\tau}_{\tP}^{\tQ}(H-X_{\tP})\tau^{\tR}_{\tQ}(H-X_{\tQ}) \\
			& =\sum_{\tP\in \mathcal{F}_{B_x}(\tM(x))} \sum_{\substack{\mathcal{F}(\tM(x))\ni \tQ \\ \tP\subset \tQ\subset \tR}}(-1)^{a_{\tP}^{\tQ}} \widehat{\tau}_{\tP}^{\tQ}(H-X_{\tP})\tau^{\tR}_{\tQ}(H-X_{\tQ}).
		\end{aligned}\]
		Moreover, by \cite[proposition 1.7.1, lemme 2.9.2]{LW} the inner sum
		$$\displaystyle \sum_{\substack{\mathcal{F}(\tM(x))\ni \tQ \\ \tP\subset \tQ\subset \tR}}(-1)^{a_{\tP}^{\tQ}} \widehat{\tau}_{\tP}^{\tQ}(H-X_{\tP})\tau^{\tR}_{\tQ}(H-X_{\tQ})$$
		equals $1$ if $\tP=\tR$ and $0$ otherwise. The identity \eqref{ident1 GammaBx} follows.
		
		Assume now that $\mathcal{X}$ is positive. Fix $\tP\in \CP_{B_x}(\tM(x))$. For $\tP'\in \CP_{B_x}(\tM(x))$, we denote by $\phi_{\tP',\tP}$ the characteristic function of the set of $H\in \CA_x$ such that for every $\alpha\in \Delta_{\tP'}$ we have
		\begin{equation*}
			\varpi_\alpha(H)\leqslant 0 \mbox{ if } \alpha\in \Sigma_{\tP}^+
		\end{equation*}
		and
		\begin{equation*}
			\varpi_\alpha(H)>0 \mbox{ if } \alpha\in \Sigma^-_{\tP}.
		\end{equation*}
		Then, we have
		\begin{align}\label{eq1 proof GammaBxX}
			\displaystyle \Gamma_{B_x}(.,\CX) & =\sum_{\tQ\in \CF_{B_x}(\tM(x))} (-1)^{a_{\tQ}^{\tG}}\widehat{\tau}_{\tQ}(.-X_{\tQ}) \\
			\nonumber & =\sum_{\tP'\in \CP_{B_x}(\tM(x))} \sum_{\substack{\tP'\subset \tQ \\ \tP\cap \tQ=\tP\cap \tP'}} (-1)^{a_{\tQ}^{\tG}} \widehat{\tau}_{\tQ}(.-X_{\tQ}) \\
			\nonumber & =\sum_{\tP'\in \CP_{B_x}(\tM(x))}(-1)^{\lvert \Delta_{\tP'}\cap \Sigma_{\tP}^-\rvert} \phi_{\tP',\tP}(.-X_{\tP'})
		\end{align}
		where the last identity follows from \cite[lemme 1.7.4, lemme 2.9.2]{LW}. From the above we deduce that for every $H\in \CA_x$ satisfying $\Gamma_{B_x}(H,\CX)\neq 0$ there exists $\tP'\in \CP_{B_x}(\tM(x))$ such that $\phi_{\tP',\tP}(H-X_{\tP'})=1$ which, since $\CX$ is positive, further implies
		$$\displaystyle \varpi_\alpha(H)\leqslant \varpi_\alpha(X_{\tP'})\leqslant \varpi_\alpha(X_{\tP})$$
		for every $\alpha\in \Delta_{\tP}$. Thus, as $\tP\in \CP_{B_x}(\tM(x))$ was arbitrary, we see that $\Supp (\Gamma_{B_x}(.,\CX))$ is included in the subset of those $H\in \CA_x$ such that
		\begin{equation}\label{eq2 proof GammaBxX}
			\displaystyle \varpi_\alpha(H-X_{\tP})\leqslant 0,\; \forall \tP\in \CP_{B_x}(\tM(x)),\; \forall \alpha\in \Delta_{\tP}.
		\end{equation}
		Conversely, assume that $H\in \CA_x$ satisfies the inequalities \eqref{eq2 proof GammaBxX}. Then, for a chosen $\tP\in \CP_{B_x}(\tM(x))$, we have $\phi_{\tP,\tP}(H-X_{\tP})=1$ whereas for $\tP\neq \tP'\in \CP_{B_x}(\tM(x))$, as $\Delta_{\tP'}\cap \Sigma_{\tP}^-\neq \emptyset$, we have $\phi_{\tP',\tP}(H-X_{\tP'})=0$. From identity \eqref{eq1 proof GammaBxX} this readily implies that $\Gamma_{B_x}(H,\CX)=1$. This gives the last part of the proposition.

	\end{proof}

\begin{cor}\label{cor gammaBx functions}
Let $\CX,\CY$ be two positive $(\tG,\tM(x))$-orthogonal sets. Then, for every $\tQ,\tR\in \mathcal{F}_{B_x}(\tM(x))$ and $Y\in \CA_x$ we have
$$\displaystyle \Gamma_{B_x}^{\tQ}(H,\mathcal{X}) \tau^{\tR}_{\tQ}(H-X_{\tQ})\Gamma_{B_x}^{\tR}(H,\mathcal{X}+\mathcal{Y})=\Gamma_{B_x}^{\tQ}(H,\mathcal{X}) \tau^{\tR}_{\tQ}(H-X_{\tQ})\phi_{\tQ}^{\tR}(H-X_{\tQ}-Y_{\tQ}).$$
In other words, on the support of the function $\Gamma_{B_x}^{\tQ}(.,\mathcal{X}) \tau^{\tR}_{\tQ}(.-X_{\tQ})$, $\Gamma_{B_x}^{\tR}(.,\mathcal{X}+\mathcal{Y})$ is equal to $\phi_{\tQ}^{\tR}(.-X_{\tQ}-Y_{\tQ})$.
\end{cor}

\begin{proof}
This follows from the combination of Lemma \ref{splitting formula} with the identities \eqref{ident1 GammaBx} and \eqref{eq Gammafunction}, noting that each of the functions $\Gamma_{B_x}^{\tQ}(.,\mathcal{X}) \tau^{\tR}_{\tQ}(.-X_{\tQ})$, $\tQ,\tR\in \mathcal{F}_{B_x}(\tM(x))$, are characteristic functions by Proposition \ref{prop partition gammaBx functions}.
\end{proof}

	To simplify the notation, we will use $\Delta_{x}$ (resp. $\Delta_{x}^\vee$) to denote the set of roots $\Delta_{B_x}\subset X^*(A_x)$ (resp. of coroots $\Delta_{B_x}^\vee\subset X_*(A_x)$).
	
	\begin{prop}\label{prop aleternative description GammaBxX}
		Assume that $\CX$ is positive. Then, $\Gamma_{B_x}(.,\CX)$ is the characteristic function of
		$$\displaystyle Conv\{X_{\tP}\mid \tP\in \CP_{B_x}(\tM(x)) \}+{}^- \CA_{B_x}+\CA_{\tG}$$
		where $Conv\{X_{\tP}\mid \tP\in \CP_{B_x}(\tM(x)) \}$ denotes the convex hull of the finite set $\{X_{\tP}\mid \tP\in \CP_{B_x}(\tM(x)) \}$ whereas ${}^- \CA_{B_x}$ stands for the closed cone generated by $-\Delta_{x}^\vee$.
	\end{prop}
	
	\begin{proof}
		Set
		$$\displaystyle \CC_{B_x}(\CX):=Conv\{X_{\tP}\mid \tP\in \CP_{B_x}(\tM(x)) \}+{}^- \CA_{B_x}+\CA_{\tG}.$$
		This is obviously a closed convex subset of $\CA_x$. Moreover, its set of extreme points is contained in $\{X_{\tP}\mid \tP\in \CP_{B_x}(\tM(x)) \}$. Recall that for every closed convex subset $\CC\subset \CA_x$, denoting by $Ext(\CC)$ its set of extreme points and, for $X\in Ext(\CC)$, by $\CC_X$ the cone centered at $X$ generated by $\CC$, that is $\CC_X=\{X+t(Y-X)\mid Y\in \CC, t\geq 0 \}$, we have
		$$\displaystyle \CC=\bigcap_{X\in Ext(\CC)} \CC_X.$$
		According to the previous proposition $\Gamma_{B_x}(.,\CX)$ is the characteristic function of
		$$\displaystyle \{H\in \CA_x\mid \varpi_\alpha(H-X_{\tP})\leqslant 0\; \forall \tP\in \CP_{B_x}(\tM(x)),\; \forall \alpha\in \Delta_{\tP} \}.$$
		Therefore, it suffices to show that for every $\tP\in \CP_{B_x}(\tM(x))$ we have
		$$\displaystyle \CC_{B_x}(\CX)_{X_{\tP}}= \{H\in \CA_x\mid \varpi_\alpha(H-X_{\tP})\leqslant 0\; \forall \alpha\in \Delta_{\tP} \}(=X_{\tP}+{}^- \CA_{\tP}+\CA_{\tG}).$$
		As $\CX$ is positive, and ${}^- \CA_{B_x}\subset {}^- \CA_{\tP}$, the inclusion $\CC_{B_x}(\CX)_{X_{\tP}}\subseteq X_{\tP}+{}^- \CA_{\tP}+\CA_{\tG}$ is clear. On the other hand, for every $\alpha\in \Delta_{\tP}$ we either have:
		\begin{itemize}
			\item $\alpha^\vee$ is positively proportional to an element of $\Delta_{x}^\vee$;
			
			\item there exists $\tP'\in \CP_{B_x}(\tM(x))$ such that $\Sigma_{\tP}^+\cap \Sigma_{\tP'}^-=\{\alpha \}$ in which case $X_{\tP'}-X_{\tP}\in \BR_{>0}\alpha^\vee$.
		\end{itemize}
		This implies that, in both cases, $\CC_{B_x}(\CX)_{X_{\tP}}$ is invariant by translation by $\BR_{\leq 0}\alpha^\vee$. As this holds for all $\alpha\in \Delta_{\tP}$, this gives the reverse inclusion $X_{\tP}+{}^- \CA_{\tP}+\CA_{\tG}\subseteq \CC_{B_x}(\CX)_{X_{\tP}}$ and therefore $\CC_{B_x}(\CX)_{X_{\tP}}= X_{\tP}+{}^- \CA_{\tP}+\CA_{\tG}$.
	\end{proof}


	
	
	\subsection{The weight $v_{B_x,\xi}(u,g)$}
	
	Let $N_{x,\reg}=N_x\cap G_{x,\reg}$ be the open subset of regular elements in $N_x$ and $T_{x,c}\subset T_x(F)$ be the maximal compact subgroup. We equip $T_{x,c}$ with the Haar measure of total mass $1$ and we also fix a log-norm $\sigma_{x,\reg}: G_{x,reg}(F)\to \mathbb{R}_{\geqslant 1}$ on $G_{x,reg}(F)$ (see Section \ref{sect HCS space}). Set $r=\dim(\CA_x)-a_{\tG}$.

	
	\begin{lem}\label{lem def vBxpsi}
		For any $u\in N_{x,\reg}(F)$ and any positive $(\tG,\tM(x))$-orthogonal set $\mathcal{X}$, the iterated integral
		\begin{equation}\label{eq1 lem def vBxpsi}
			\displaystyle \int_{T_x(F)/A_{\tG}(F)} \int_{T_{x,c}} \xi(a^{-1}t^{-1}uta)\; dt\;  \Gamma_{B_x}(H_{T_x}(a),\mathcal{X})\; da
		\end{equation}
		is absolutely convergent in that order and will be denoted by
		\begin{equation*}
			\displaystyle \tilde{v}_{B_x,\xi}(u,\mathcal{X}):=\int_{T_x(F)/A_{\tG}(F)}^* \xi(a^{-1}ua) \Gamma_{B_x}(H(a),\mathcal{X})\; da.
		\end{equation*}
		Moreover, there exists a constant $C>0$ such that for every $u\in N_{x,reg}(F)$ and every positive $(\tG,\tM(x))$-orthogonal set $\mathcal{X}$, we have (where $N(\CX)$ denotes the norm of $\CX$ defined in Section \ref{orthogonal sets})
		\begin{equation*}
			\displaystyle \left\lvert \tilde{v}_{B_x,\xi}(u,\mathcal{X})\right\rvert\leqslant C(\sigma_{x,\reg}(u)+N(\mathcal{X}))^r.
		\end{equation*}
	\end{lem}

	\begin{proof}
		The inner integral over $T_{x,c}$ in \eqref{eq1 lem def vBxpsi} is obviously convergent. Let $N_{x,der}$ be the derived subgroup of $N_x$ and let $N_x/N_{x,der}=\bigoplus_{\alpha\in \Delta_x} (N_x/N_{x,der})_\alpha$ be the isotypic decomposition with respect to the adjoint action of $A_x$. We fix a norm $\lVert .\rVert$ on the $F$-vector space $N_x(F)/N_{x,der}(F)$ and for every $u\in N_x$ and $\alpha\in \Delta_x$, let us we denote by $u_\alpha$ the projection of $u$ to $(N_x/N_{x,der})_\alpha$.  Then, since $\xi$ is a generic character, there exists $C_1>0$ such that for all $a\in T_x(F)$ and $u\in N_{x,reg}(F)$ we have
		$$\displaystyle \int_{T_{x,c}} \xi(a^{-1}t^{-1}uta) dt\neq 0\Rightarrow \langle \alpha, H_{T_x}(a)\rangle \geqslant \log \lVert u_\alpha\rVert-C_1\;\text{for all}\; \alpha\in \Delta_{x}.$$
		On the other hand, there exists $C_2>0$ such that $\log \lVert u_\alpha\rVert-C_1\geqslant -C_2 \sigma_{x,\reg}(u)$ for all $(u,\alpha)\in N_{x,reg}(F)\times \Delta_x$. Combining this  with Proposition \ref{prop aleternative description GammaBxX}, we see that, for $u\in N_{x,reg}(F)$ and $\mathcal{X}$ a positive $(\tG,\tM(x))$-orthogonal set, the function 
		\begin{equation}\label{eq2 lem def vBxpsi}
			\displaystyle a\in T_x(F)/A_{\tG}(F)T_{x,c}\mapsto  \Gamma_{B_x}(H_{T_x}(a),\mathcal{X}) \int_{T_{x,c}} \xi(a^{-1}t^{-1}uta)\; dt 
		\end{equation}
		is supported in the compact subset (where we identify $T_x(F)/T_{x,c}$ with a subset of $\CA_x$ via the map $H_{T_x}$)
		\begin{equation}\label{eq3 lem def vBxpsi}
			\displaystyle \left(Conv\{X_{\tP}\mid \tP\in \mathcal{P}_{B_x}(\tM(x)) \}+{}^- \CA_{B_x}+\CA_{\tG}\right) \cap\left\{H\in \CA_x\mid \langle \alpha,H\rangle\geqslant -C_2\sigma_{x,\reg}(u),\; \forall \alpha\in \Delta_x \right\}
		\end{equation}
		of $\CA_x/\CA_{\tG}$. Since the function \eqref{eq2 lem def vBxpsi} is also obviously bounded by $1$, the lemma follows up to noticing the existence of $C_3>0$ such that the subset \eqref{eq3 lem def vBxpsi} is contained in $\mathbb{B}(C_3(\sigma_{x,\reg}(u)+N(\mathcal{X})))+\CA_{\tG}$ for any $u\in N_{x,reg}(F)$ and for any positive $(\tG,\tM(x))$-orthogonal set $\mathcal{X}$. Here for $R>0$, we use $\mathbb{B}(R)$ to denote the ball of radius $R$ centered at $0$ in $\CA_x$ for a given norm.
	\end{proof}

	
	
	\begin{lem}\label{lem exponential polynomial}
		There exists $C>0$ and, for every $u\in N_{x,reg}(F)$, a unique unitary polynomial-exponential function $v_{B_x,\xi}(u,.)$ on $\CC_{\mathbb{Q}}(\tG,\tM(x))$ such that for every rational $(\tG,\tM(x))$-orthogonal set $\mathcal{X}\in \CC_{\mathbb{Q}}(\tG,\tM(x))$ with $d(\mathcal{X})\geqslant C\sigma(u)$ (we refer the reader to Section \ref{orthogonal sets} for various notation), we have
		$$\displaystyle v_{B_x,\xi}(u,\mathcal{X})=\widetilde{v}_{B_x,\xi}(u,\mathcal{X}).$$
		Moreover, as $u$ varies, the set of those unitary polynomial-exponential functions $\{v_{B_x,\xi}(u,.)|\; u\in N_{x,reg}(F)\}$ spans a finite dimensional vector space and there exists $C'>0$ such that for every $u\in N_{x,reg}(F)$ and $\mathcal{X}\in \CC_{\mathbb{Q}}(\tG,\tM(x))$ we have
		\begin{equation*}
			\displaystyle \left\lvert v_{B_x,\xi}(u,\mathcal{X})\right\rvert\leqslant C'(\sigma_{x,reg}(u)+N(\mathcal{X}))^r.
		\end{equation*}
	\end{lem}
	
	\begin{proof}
		Before proving the lemma, we need some preparation. For every $u\in N_{x,reg}(F)$ and $\tQ\in \CF_{B_x}(\tM(x))$ with Levi decomposition $Q=L_QU_Q$ (where $M(x)\subset L_Q$), there is a unique decomposition $u=u^Qu_Q$ where $u^Q\in L_{Q}(F)$, $u_Q\in U_{Q}(F)$ and we set
		$$\displaystyle \xi^{c,u,\tQ}_x(t):=\int_{T_{x,c}}\xi(t^{-1}t_c^{-1}u^Qt_ct) dt_c,\; \mbox{ for } t\in T_x(F).$$
		Then, these functions satisfy:
		\begin{itemize}
			\item For every $u\in N_{x,reg}(F)$ and $\tQ\in \CF_{B_x}(\tM(x))$, $\xi^{c,u,\tQ}$ is invariant by translation by $A_{\tQ}(F)$;
			
			\item There exists $C_1>0$ such that for every $u\in N_{x,reg}(F)$, $\tQ\in \CF_{B_x}(\tM(x))$, $t\in T_x(F)$ and $(\tG,\tM(x))$-orthogonal set $\mathcal{X}$ satisfying $d(\mathcal{X})\geqslant C_1\sigma(u)$ the condition
			$$\displaystyle \Gamma_{B_x}^{\tQ}(H(t),\CX)\tau_{\tQ}(H(t)-X_{\tQ})\neq 0$$
			implies $\xi^{c,u,\tG}(t)=\xi^{c,u,\tQ}(t)$.
		\end{itemize}
		The first bullet point is obvious. Let's prove the second bullet point. Pick $C_1>0$ and let $\tQ$, $u$, $t$, $\mathcal{X}$ be as above satisfying
		\begin{equation}\label{eq-3 proof bulet point}
			\displaystyle d(\mathcal{X})\geqslant C_1\sigma(u),
		\end{equation}
		\begin{equation}\label{eq-2 proof bulet point}
			\displaystyle \Gamma_{B_x}^{\tQ}(H(t),\CX)\tau_{\tQ}(H(t)-X_{\tQ})\neq 0.
		\end{equation}
		Then, we will show that provided $C_1$ is large enough, we have
		\begin{equation}\label{eq-1 proof bulet point}
			\displaystyle	\xi^{c,u,\tG}(t)=\xi^{c,u,\tQ}(t).
		\end{equation}
		By Proposition \ref{prop aleternative description GammaBxX} (applied to $\tL_Q$ instead of $\tG$), the condition \eqref{eq-2 proof bulet point} is equivalent to
		\begin{equation*}
			\displaystyle H^{\tQ}(t)\in Conv\{X_{\tP}^{\tQ}\mid \tB_x\subset \tP\subset \tQ \}+{}^- \CA_{B_x}^{\tQ} \mbox{ and } H_{\tQ}(t)\in X_{\tQ}+\CA_{\tQ}^+
		\end{equation*}
		where $H^{\tQ}(t)$, $H_{\tQ}(t)$ denote the respective projections of $H(t)$ onto $\CA_{\tM(x)}^{\tQ}$, $\cA_{\tQ}$ and we have set $\tB_x=xB_x$. Hence, it implies that
		\begin{equation}\label{eq1 proof bulet point}
			\displaystyle H(t)=H^{\tQ}(t)+H_{\tQ}(t)\in Conv(X_{\tP}\mid \tB_x\subset \tP\subset \tQ)+\CA_{\tQ}^++{}^- \CA_{B_x}^{\tQ}.
		\end{equation}
		On the other hand, we have
		$$\displaystyle \xi^{c,u,\tG}(t)=\int_{T_{x,c}} \xi(t^{-1}t_c^{-1}u^Qt_ct) \xi(t^{-1}t_c^{-1}u_Qt_ct)dt_c.$$
		Thus, it suffices to show that, when $C_1$ is sufficiently large, we have
		$$\displaystyle \xi(t^{-1}t_c^{-1}u_Qt_ct)=1,\;\; \forall t_c\in T_{x,c}.$$
		There exists $C_2>0$ such that this last condition is implied by the inequalities
		\begin{equation*}
			\displaystyle \langle \alpha,H(t)\rangle\geqslant C_2\sigma(u),\; \mbox{ for every } \alpha\in \Delta_x\setminus \Delta_x^{\tQ}.
		\end{equation*}
		where $\Delta_x^{\tQ}=\Delta_{B_x}^{Q\cap B_x}$.	However, as every $\alpha\in \Delta_x\setminus \Delta^{\tQ}_x$ takes non-negative values on ${}^- \CA_{B_x}^{\tQ}$ and on $\CA_{\tQ}^+$, \eqref{eq1 proof bulet point} implies
		\begin{equation}\label{eq4 proof bulet point}
			\displaystyle \langle \alpha, H(t)\rangle\geqslant \min_{\tB_x\subset \tP\subset \tQ} \langle \alpha, X_{\tP}\rangle \geqslant d(\CX),\; \mbox{ for every } \alpha\in \Delta_{x}\setminus \Delta_{x}^{\tQ}.
		\end{equation}
		Therefore, taking $C_1\geqslant C_2$ gives the required identity.

		Let us now prove the lemma. Fix a lattice $\Lambda\subset \CA_{x,\mathbb{Q}}$. Then, we can find a constant $C_2>0$ and for every $u\in N_x(F)$, a $\Lambda$-rational orthogonal set $\mathcal{X}_u=(X_{u,\tP})_{\tP\in \CP(\tM(x))}\in \CC_{\Lambda}(\tG,\tM(x))$ such that $d(\CX_u)\geq C_1\sigma(u)$ and $N(\CX_u)\leq C_2\sigma(u)$. Obviously, it suffices to show that for every $u\in N_{x,reg}(F)$, the function 
		$$\displaystyle \mathcal{Y}\in \CC_{\mathbb{Q}}(\tG,\tM(x)) \mapsto \tilde{v}_{B_x,\xi}(u,\mathcal{X}_u+\mathcal{Y})$$
		coincides, for $\mathcal{Y}$ positive, with a unitary polynomial-exponential function. Applying the splitting formula of Lemma \ref{splitting formula} as well as the two bullet points above, we obtain 
		\[\begin{aligned}
			\displaystyle & \tilde{v}_{B_x,\xi}(u,\mathcal{X}_u+\mathcal{Y}) =\sum_{\tQ\in \CF_{B_x}(\tM(x))} \int_{T_x(F)/A_{\tG}(F)} \xi^{c,u,\tG}(t) \Gamma^{\tQ}_{B_x}(H(t),\mathcal{X}_u) \Gamma_{\tQ}(H(t)-X_{u,\tQ},Y_{\tQ}) dt \\
			& =\sum_{\tQ\in \CF_{B_x}(\tM(x))} \int_{T_x(F)/A_{\tG}(F)} \xi^{c,u,\tQ}(t) \Gamma^{\tQ}_{B_x}(H(t),\mathcal{X}_u) \Gamma_{\tQ}(H(t)-X_{u,\tQ},Y_{\tQ}) dt \\
			& =\sum_{\tQ\in \CF_{B_x}(\tM(x))} \int_{T_x(F)/A_{\tQ}(F)} \xi^{c,u,\tQ}(t) \Gamma^{\tQ}_{B_x}(H(t),\mathcal{X}_u) \int_{A_{\tQ}(F)/A_{\tG}(F)} \Gamma_{\tQ}(H_{\tQ}(at)-X_{u,\tQ},Y_{\tQ}) da dt
		\end{aligned}\]
		for every positive $(\tG,\tM(x))$-orthogonal set $\CY$ and where the second equality is based on the fact that since $Y_{\tQ}\in \CA_{\tQ}^+$, $\Gamma_{\tQ}(H-X_{u,\tQ},Y_{\tQ})\neq 0$ implies $\tau_{\tQ}(H-X_{u,\tQ})\neq 0$. For $\tQ\in \CF_{B_x}(\tM(x))$, the function
		$$\displaystyle t\in T_x(F)/A_{\tQ}(F)\mapsto \xi^{c,u,\tQ}(t) \Gamma^{\tQ}_{B_x}(H(t),\mathcal{X}_u)$$
		is compactly supported so that in the above expression the integral over $T_x(F)/A_{\tQ}(F)$ can be written as a finite sum. On the other hand for every fixed $t\in T_x(F)$, the function
		$$\displaystyle Y_{\tQ}\in \CA_{\tQ,\mathbb{Q}} \mapsto \int_{A_{\tQ}(F)/A_{\tG}(F)}  \Gamma_{\tQ}(H_{\tQ}(at)-X_{u,\tQ},Y_{\tQ}) da$$
		is a unitary polynomial-exponential function and the set of these functions, as $t\in T_x(F)$ and $X_{u,\tQ}\in \Lambda$ vary, spans a finite dimensional vector space. This shows the lemma except for the last estimate.
		
		By the above computation, we have
		$$\displaystyle v_{B_x,\xi}(u,\CX_u+\CY)=\sum_{\tQ\in \CF_{B_x}(\tM(x))} \int_{T_x(F)/A_{\tQ}(F)} \xi^{c,u,\tQ}(t) \Gamma^{\tQ}_{B_x}(H(t),\mathcal{X}_u) \int_{A_{\tQ}(F)/A_{\tG}(F)} \Gamma_{\tQ}(H_{\tQ}(at)-X_{u,\tQ},Y_{\tQ}) da dt$$
		for every $u\in N_{x,\reg}(F)$ and $\CY\in \CC_{\mathbb{Q}}(\tG,\tM(x))$. However, the integral
		$$\displaystyle \int_{A_{\tQ}(F)/A_{\tG}(F)} \lvert \Gamma_{\tQ}(H_{\tQ}(at)-X_{u,\tQ},Y_{\tQ})\rvert da$$
		is essentially bounded by $(1+N(\CX_u)+N(\CY))^{a^{\tG}_{\tQ}}$ whereas by a similar reasoning as in the proof of Lemma \ref{lem def vBxpsi}, the integral
		$$\displaystyle \int_{T_x(F)/A_{\tQ}(F)} \left\lvert\xi^{c,u,\tQ}(t) \Gamma^{\tQ}_{B_x}(H(t),\mathcal{X}_u)\right\rvert dt$$
		is essentially bounded by $(N(\CX_u)+\sigma_{x,\reg}(u))^{r-a^{\tG}_{\tQ}}$. Since $N(\CX_u)\ll \sigma(u)$ this shows that $\lvert v_{B_x,\xi}(u,\CX_u+\CY)\rvert\ll (\sigma_{x,\reg}(u)+N(\CX_u+\CY))^r$ and the lemma is proved.
	\end{proof}

	For $g\in G(F)$, applying the above definition to the $(\tG,\tM(x))$-orthogonal set $\mathcal{Y}(g)=(H_{\tP}(g))_{\tP\in \mathcal{P}(\tM(x))}$, we define the weight
	\begin{equation*}
		\displaystyle v_{B_x,\xi}(u,g)=v_{B_x,\xi}(u,\mathcal{Y}(g)).
	\end{equation*}
	It satisfies the relation
	\begin{equation}\label{eq equivariance vBxpsi}
		\displaystyle v_{B_x,\xi}(u,bg)=v_{B_x,\xi}(b^{-1}ub,g) \mbox{ for every } (u,b,g)\in N_{x,\reg}(F)\times B_x(F)\times G(F).
	\end{equation}

	\subsection{A formula of regular germs for quasi-characters}

	\begin{thm}\label{thm: formula for the germ}
		For every strongly cuspidal function $f\in \mathcal{C}(\tG(F))$, we have \footnote{$-\CO_\xi$ is the same as $\CO_{\xi^{-1}}$}
		\begin{equation*}
			\displaystyle c_{f,-\CO_\xi}(x)=\int_{B_x(F)\backslash G(F)} \int_{N_x(F)} f(g^{-1}x u g) v_{B_x,\xi}(u,g) du dg.
		\end{equation*}
	\end{thm}
	
	Let us remark that thanks to \eqref{eq equivariance vBxpsi}, the expression in the right-hand side of the above theorem makes sense formally. We will check its absolute convergent in the next subsection.
	
	\subsection{Some estimates}
	
	In this subsection we prove some estimates that in particular imply the convergence of the right-hand side of Theorem \ref{thm: formula for the germ}. 
	
	Let $S$ be the connected center of $G_x$ (a torus) and set $\tS=Sx$. Let $\tS'$ the open subset of those $s\in \tS$ such that $G_{s}=G_x$.
	
	Recall that $r=\dim(A_x)-\dim(A_{\tG})$ and $\sigma_{x,\reg}$ denotes a log-norm on $N_{x,\reg}=N_x\cap G_{x,\reg}$. We fix log-norms $\sigma_{\reg}$ and $\sigma_{\tS'}$ on $\tG_{\reg}(F)$ and $\tS'(F)$ respectively.
	
	\begin{lem}\label{lem norm descent property}
		We have inequalities
		\begin{equation}
			\displaystyle \inf_{b\in B_x(F)} (\sigma_{x,\reg}(bub^{-1})+\sigma(bg))\ll \sigma_{\reg}(g^{-1}sug)+\sigma_{\tS'}(s),
		\end{equation}
		and
		\begin{equation}
			\displaystyle \inf_{b\in B_x(F)} (\sigma(bub^{-1})+\sigma(bg))\ll \sigma(g^{-1}sug)+\sigma_{\tS'}(s)
		\end{equation}
		for $(s,u,g)\in \tS'(F)\times N_{x,\reg}(F)\times G(F)$.
	\end{lem}
	
	\begin{proof}
		Let $\CN_x\subset G_x$ be the unipotent cone and $\CN_x\times^{G_x} G$ be the quotient of $\CN_{x}\times G$ by the free action of $G_x$ given by $g_x\cdot (u,g)=(g_xug_x^{-1},g_xg)$. Then, the regular map
		$$\displaystyle \CN_{x}\times^{G_x} G\times \tS'\to \tG\times \tS',\;\; (u,g,s)\mapsto (g^{-1}sug,s)$$
		is a closed embedding with image the subset of those $(\gamma,s)\in \tG\times \tS'$ such that the semisimple part of $\gamma$ is in the same geometric conjugacy class as $s$. Let $\CN_{x,\reg}=\CN_x\cap G_{x,\reg}$ be the open subset of regular unipotent elements. Then, the previous map restricts to a closed embedding $\CN_{x,\reg}\times^{G_x} G\times \tS'\to \tG_{\reg}\times \tS'$. Furthermore, the natural map $\CN_{x}\times^{B_x} G\to \CN_{x}\times^{G_x} G$ is proper and $N_x$ (resp. $N_{x,\reg}$) is a closed subset of $\CN_x$ (resp. $\CN_{x,\reg}$). It follows that the two regular maps
		$$\displaystyle N_{x}\times^{B_x} G\times \tS' \to \tG\times \tS' \mbox{ and } N_{x,\reg}\times^{B_x} G\times \tS'\to \tG_{\reg}\times \tS'$$
		are proper (the second one being actually a closed embedding). By Lemma \ref{lem proper map}, we have
		\begin{equation*}
			\displaystyle \sigma_{N_{x}\times^{B_x} G}(u,g)\ll \sigma(g^{-1}sug)+\sigma_{\tS'}(s),\;\; \sigma_{N_{x,\reg}\times^{B_x} G}(u,g)\ll \sigma_{\reg}(g^{-1}sug)+\sigma_{\tS'}(s)
		\end{equation*}
		for $(u,s,g)\in N_{x,\reg}(F)\times \tS'(F)\times G(F)$. It remains to check that
		$$\displaystyle \sigma_{N_{x}\times^{B_x} G}(u,g)\sim \inf_{b\in B_x(F)} (\sigma(bub^{-1})+\sigma(bg)) \mbox{ and } \sigma_{N_{x,\reg}\times^{B_x} G}(u,g)\sim \inf_{b\in B_x(F)} (\sigma_{x,\reg}(bub^{-1})+\sigma(bg))$$
		for $(u,g)\in N_{x,\reg}(F)\times G(F)$ i.e. that the two natural projections $N_x\times G\to N_{x}\times^{B_x} G$ and $N_{x,\reg}\times G\to N_{x,\reg}\times^{B_x} G$ have the norm descent property. Since both are pullbacks of the projection $G\to B_x\backslash G$, it suffices to check that the latter has the norm descent property.
		
		Choose $P=M(x)U\in \CP_{B_x}(M(x))$ a parabolic subgroup with Levi $M(x)$ such that $P\cap G_x=B_x$. As $P\backslash G$ is proper, we just need to check that $P\to B_x\backslash P$ has the norm descent property. Let $\pi: B_x\backslash P\to T_x\backslash M(x)$ be the natural map. According to Kottwitz, $M(x)\to T_x\backslash M(x)$ already has the norm descent property. Thus, for every $p\in P(F)$, we can find $m\in M(x)(F)$ such that $p\in T_x(F)U(F)m$ and
		$$\displaystyle \sigma(m)\sim \sigma_{T_x\backslash M(x)} (\pi(B_xp))\ll \sigma_{B_x\backslash P} (B_xp).$$
		
		Choose $C>0$ large enough such that for all $p\in P(F)$, there exists $p'=u'm'\in B_x(F)p$ with $u'\in U(F)$ and $m'\in M(x)(F)$ such that 
		$$\sigma(m')\leq C\sigma_{B_x\backslash P} (B_xp).$$
		Fix another constant $C_1>0$ large enough (with respect to $C$). If $\sigma_{N_x\backslash U}(u')\leq C_1\sigma_{B_x\backslash P} (B_xp)$, then since $U\to N_x\backslash U$ has the norm descent property (this is because this quotient map admits a regular section), there exists $n\in N_x(F)$ such that
		$$\sigma(nu'm')\leq 2C_1\sigma_{B_x\backslash P} (B_xp).$$
		This implies that
		$$\inf_{b\in B_x(F)}\sigma(bp)\leq 2C_1\sigma_{B_x\backslash P} (B_xp).$$
		
		If $\sigma_{N_x\backslash U}(u')> C_1\sigma_{B_x\backslash P} (B_xp)$, since $N_x\backslash U\rightarrow B_x\backslash P$ is a closed embedding and since $\sigma(m')\leq C\sigma_{B_x\backslash P} (B_xp)$, we have 
		$$\sigma_{B_x\backslash P}(p)=\sigma_{B_x\backslash P}(u'm')\geq \frac{1}{2}\sigma_{B_x\backslash P}(u')=\frac{1}{2}\sigma_{N_x\backslash U}(u')$$ 
		and 
		$$\inf_{b\in B_x(F)}\sigma(bp)\leq \inf_{n\in N_x(F)}\sigma(nu'm')\leq  2\sigma_{N_x\backslash U}(u').$$
		This implies that
		$$\inf_{b\in B_x(F)}\sigma(bp)\leq 4\sigma_{B_x\backslash P} (B_xp).$$
		As a result we have proved that the map $P\to B_x\backslash P$ has the norm descent property and this finishes the proof of the lemma.
	\end{proof}
	
	\begin{cor}
		We have
		$$\displaystyle \lvert v_{B_x,\xi}(u,g)\rvert\ll (\sigma_{\reg}(g^{-1}sug)+\sigma_{\tS'}(s))^r$$
		for every $u\in N_{x,reg}(F)$, $g\in G(F)$ and $s\in \tS'(F)$.
	\end{cor}
	
	\begin{proof}
		According to Lemma \ref{lem exponential polynomial}, we have $\lvert v_{B_x,\xi}(u,g)\rvert\ll(\sigma_{x,\reg}(u)+\sigma(g))^r$. Combining this with the equation \eqref{eq equivariance vBxpsi} and the previous lemma, we obtain
		$$\displaystyle \lvert v_{B_x,\xi}(u,g)\rvert\ll\inf_{b\in B_x(F)} (\sigma_{x,\reg}(bub^{-1})+\sigma(bg))^r\ll  (\sigma_{\reg}(g^{-1}sug)+\sigma_{\tS'}(s))^r$$
		for $u\in N_{x,reg}(F)$, $g\in G(F)$ and $s\in \tS'(F)$.
	\end{proof}

	\begin{prop}\label{estimate 1}
		Let $r_0>0$. Then, for every $f\in \mathcal{C}(\tG(F))$ and every $d>0$, we have
		\begin{equation}\label{eq0 estimate}
			\displaystyle D^{\tG}(s)^{1/2}\int_{B_x(F)\backslash G(F)}\int_{N_x(F)} \lvert f(g^{-1}sug)\rvert  \sigma_{\reg}(g^{-1}sug)^{r_0} dudg\ll_{f,d} \sigma(s)^{-d} \sigma_{\tS'}(s)^{r_0}
		\end{equation}
		for $s\in \tS'(F)$. In particular, the integral in Theorem \ref{thm: formula for the germ} is absolutely convergent.
	\end{prop}
	
	\begin{proof}
		Let $K\subset G(F)$ be a compact-open subgroup. First, we show that
		\begin{equation}\label{eq1 estimate}
			\displaystyle \vol(gKg^{-1}\cap N_x(F))^{-1}\int_{gKg^{-1}\cap N_x(F)} \sigma_{\reg}(g^{-1}sukg)^{r_0}dk\ll (\sigma(g^{-1}sug)+\sigma_{\tS'}(s))^{r_0}
		\end{equation}
		for every $(g,s,u)\in G(F)\times \tS'(F)\times N_x(F)$. Since the left-hand side of the inequality is invariant by the transformation $(g,s,u)\mapsto (bg,s,bub^{-1})$, by Lemma \ref{lem norm descent property} it suffices to establish that 
		\begin{equation*}
			\displaystyle \vol(gKg^{-1}\cap N_x(F))^{-1}\int_{gKg^{-1}\cap N_x(F)} \sigma_{\reg}(g^{-1}sukg)^{r_0}dk\ll (\sigma(g)+\sigma(u)+\sigma_{\tS'}(s))^{r_0}
		\end{equation*}
		for $(g,s,u)\in G(F)\times \tS'(F)\times N_x(F)$. Note that $\sigma_{\reg}(g^{-1}sug)\ll \sigma(g)+\sigma_{\reg}(su)\ll \sigma(g)+\sigma_{\tS'}(s)+\sigma_{x,\reg}(u)$. Therefore, we are reduced to show
		\begin{equation}\label{eq2 estimate}
			\displaystyle \vol(gKg^{-1}\cap N_x(F))^{-1}\int_{gKg^{-1}\cap N_x(F)} \sigma_{x,\reg}(uk)^{r_0}dk\ll (\sigma(g)+\sigma(u))^{r_0} \mbox{ for } u\in N_x(F).
		\end{equation}
		
		Let $\Delta_x$ be the set of simple roots of $A_x$ in $B_x$ and for $\alpha\in \Delta_x$, let $\mathfrak{n}_{x,\alpha}\subset \mathfrak{n}_x$ be the corresponding root subspace. Then, we have a natural projection $\mathfrak{n}_x\to \mathfrak{n}_{x,\alpha}$ and for $u\in N_x$, we denote by $\log(u)_\alpha$ the image of $\log(u)$ in $\mathfrak{n}_{x,\alpha}$ where $\log:N_x\to \mathfrak{n}_x$ denotes the logarithmic map (a regular morphism). Fix an ultrametric norm $\lVert .\rVert$ on $\mathfrak{n}_x(F)$ and set $v(.)=-\log \lVert .\rVert$. Then, we have
		\begin{equation*}
			\displaystyle \sigma_{x,reg}(u)^{r_0}\sim \sigma(u)^{r_0}+\sum_{\alpha\in \Delta_x} \max(1,v(\log(u)_\alpha))^{r_0} \mbox{ for } u\in N_{x,\reg}(F).
		\end{equation*}
		Thus, to show \eqref{eq2 estimate} it suffices to bound the integral
		\begin{equation}\label{eq3 estimate}
			\displaystyle \vol(gKg^{-1}\cap N_x(F))^{-1}\int_{gKg^{-1}\cap N_x(F)} \max(1,v(\log(u)_\alpha+\log(k)_\alpha))^{r_0}dk
		\end{equation}
		by a constant times $(\sigma(u)+\sigma(g))^{r_0}$. For this, we remark that there exists $C>0$ such that the image of $gKg^{-1}\cap N_x(F)$ in $\mathfrak{n}_{x,\alpha}(F)$ contains the ball
		$$\displaystyle \mathbb{B}(C\sigma(g)):=\{X\in \mathfrak{n}_{x,\alpha}(F)\mid v(X)\geq C\sigma(g) \}$$
		for every $g\in G(F)$. Therefore, since $\sigma(uk)\ll \sigma(u)+\sigma(g)$ for every $g\in G(F)$, $u\in N_x(F)$ and $k\in gKg^{-1}$, the desired estimate for \eqref{eq3 estimate} follows from the elementary inequality
		$$\displaystyle \vol(\mathbb{B}(R))^{-1}\int_{\mathbb{B}(R)} \max(1,v(X+Y))^{r_0} dY\ll R^{r_0} \mbox{ for every } R\geq 1, X\in \mathfrak{n}_{x,\alpha}(F),$$
		and this ends the proof of \eqref{eq1 estimate}.
		
		From \eqref{eq1 estimate} applied to some compact-open subgroup $K\subset G(F)$ leaving $f$ invariant in the right, we get that the left-hand side of \eqref{eq0 estimate} is essentially bounded by
		\begin{equation*}
			\displaystyle \sigma_{\tS'}(s)^{r_0}D^{\tG}(s)^{1/2}\int_{B_x(F)\backslash G(F)}\int_{N_x(F)} \lvert f(g^{-1}sug)\rvert  \sigma(g^{-1}sug)^{r_0} dudg.
		\end{equation*}
		Note that the function $\sigma^{r_0}\lvert f\rvert$ belongs to the Harish-Chandra Schwartz space $\CC(G(F))$. Therefore, up to replacing $f$ by this function, it suffices to show that for every $d>0$ we have
		\begin{equation}\label{eq4 estimate}
			\displaystyle D^{\tG}(s)^{1/2}\int_{B_x(F)\backslash G(F)}\int_{N_x(F)} \lvert f(g^{-1}sug) \rvert dudg\ll_d \sigma(s)^{-d}
		\end{equation}
		for $s\in \tS'(F)$.
		
		Pick a parabolic subspace $\tP=\tM(x)U\in \CP_{B_x}(\tM(x))$ with $P\cap G_x=B_x$ as well as a compact subgroup $K\subset G(F)$ such that $G(F)=P(F)K$. Then, by the usual change of variable the last integral above is equal to
		$$\displaystyle D^{\tM(x)}(s)^{1/2} \delta_{\tP}(s)^{1/2} \int_{T_x(F)\backslash M(x)(F)}  \int_{K} \int_{U(F)} \lvert f(k^{-1}m^{-1}smuk) \rvert du dk dm.$$
		Thus, since the function $m\in \tM(x)(F)\mapsto \delta_{\tP}(m)^{1/2} \int_K \int_{U(F)} \lvert f(k^{-1} mu k) \rvert du$ is Harish-Chandra Schwartz \cite[Proposition II.4.5]{WalPlanch}, the estimate \eqref{eq4 estimate} is now a consequence of \cite[Lemma 2.9.2]{BeuGalP}.
	\end{proof}
	
	\subsection{Definition of a sequence of test functions}
	
	As a preparation for the proof of Theorem \ref{thm: formula for the germ} we introduce a sequence of functions $\phi_n\in C_c^\infty(G_x(F))$ as follows, the construction being inspired from \cite{MW}. 
	
	First, let $\Xi$ be the unique element of $\overline{\Fb}_x^\perp(F)\subset \mathfrak{g}^*(F)$ (where $\overline{B}_x=T_x\overline{N}_x$ denotes the Borel opposite to $B_x$) such that for every $X\in \Fn_x(F)$ we have $\xi(\exp(X))=\psi(\langle \Xi, X\rangle)$. Then, $\Xi\in \CO$ (by definition of the generic character $\xi$) and, denoting by $\Fg_{x,\Xi}$ the centralizer of $\Xi$ in $\Fg_x$, we have $\Fn_x\cap \Fg_{x,\Xi}=0$. Moreover, the image of $\Fn_x(F)$ in the quotient $\Fg_x(F)/\Fg_{x,\Xi}(F)$ is maximal isotropic with respect to the bicharacter
	\begin{equation}\label{bicharacter}
		\displaystyle (X,Y)\in \Fg_x(F)/\Fg_{x,\Xi}(F)\mapsto \psi(\langle \Xi,[X,Y]\rangle).
	\end{equation}
	Let $L\subset \fg_x(F)$ be a lattice such that:
	\begin{itemize}
		\item The image $L^\xi$ of $L$ in $\fg_x(F)/\fg_{x,\Xi}(F)$ is self-dual with respect to the bicharacter \eqref{bicharacter} i.e. $L^\xi=\{X\in  \fg_x(F)/\fg_{x,\Xi}(F)\mid \psi(\langle \Xi,[X,Y]\rangle)=1 \;\forall Y\in L^\xi\}$;
		
		\item The preimage of $L^\xi$ in $\Fn_x(F)$ is $\Fn_x(F)\cap L$.
	\end{itemize}
	We then choose an integer $n_0>0$ large enough such that:
	\begin{num}
		\item\label{properties n0}
		\begin{itemize}
			\item The exponential map $\exp: \Fg_x(F)\to G_x(F)$ is well-defined on $\varpi^{n_0}L$;
			
			\item For every $n\geq n_0$, $m\geq n_0$, $X\in \varpi^n L$ and $Y\in \varpi^m L$, we have $e^Xe^Y\in \exp(X+Y+\varpi^{m+n-n_0}L)$;
			
			\item $[L,L]\subset \varpi^{-n_0}L$;
			
			\item The restriction of $\xi$ to $\varpi^{n_0}L$ is trivial;
			
			\item For every $n\geq n_0$, $m\in \mathbb{Z}$, $Y\in \varpi^n L$, $X\in \varpi^m L$, we have
			$$\displaystyle \Ad_{e^Y}(X)-X-[Y,X]\in \varpi^{2n+m-n_0}L.$$
		\end{itemize}
	\end{num} 
	That the last point above is satisfied for $n_0$ large enough is a consequence of the series expansion $\Ad_{e^Y}X=\sum_{k\geq 0} \frac{\ad_Y^k(X)}{k!}$ (valid for $Y$ small enough).
	
	For every integer $n\geqslant n_0$, we set
	\begin{itemize}
		\item $a_n=(2\rho^\vee_x)(\varpi)^n$ where $\rho^\vee_x$ denotes half the sum of the positive coroots of $A_x$ with respect to $B_x$,;
		
		\item $L_n=(Ad_{a_n})^{-1} \varpi^n L$, $K_n=\exp(L_n)$ and $K'_n=\exp(\varpi^n L)$;
		
		\item $\xi_n: K_n\to \C^\times$, $\xi'_n: K'_n\to \C^\times$ the locally constant functions defined by $\xi_n(\exp(X))=\psi(\langle \Xi, X\rangle )$ and $\xi_n'(\exp(Y))=\psi(\langle \Xi,\varpi^{-2n}Y\rangle)$ for every $X\in L_n$ and $Y\in \varpi^n L$ respectively.
	\end{itemize}
	Note that, by the second condition on $n_0$, $K_n$ and $K'_n$ are compact-open subgroups of $G_x(F)$. Moreover, we have $K'_n=a_nK_na_n^{-1}$ and $\xi'_n(a_n ka_n^{-1})=\xi_n(k)$ for every $k\in K_n$ which follows from the fact that $\Ad_{a_n}^*\Xi=\varpi^{-2n}\Xi$. From the last condition on $n_0$, we also deduce that the function $\xi_n$ (resp. $\xi'_n$) is $K_n$-invariant (resp. $K'_n$-invariant) by conjugation\footnote{Actually, for $n$ large enough and if the residue characteristic is different from $2$, it can be shown that $\xi_n$ and $\xi'_n$ are characters of $K_n$ and $K'_n$ respectively. But we will not need this fact in the sequel.}
	
	We fix Haar measures on $\Fg_x(F)$ and $\Fn_x(F)$ compatibly with the measures on $G_x(F)$ and $N_x(F)$ i.e. such that exponential maps are locally measure preserving. Using the additive character $\psi$, we can identify $\Fn_x^\perp(F)$ with the Pontryagin dual of $\Fg_x(F)/\Fn_x(F)$, and we endow $\Fn_x^\perp(F)$ with the dual of the quotient measure on the latter. This is the only invariant measure on $\Fn_x^\perp(F)$ such that for every lattice $\Lambda\subset \Fg_x(F)$ we have $\vol(\Lambda^\perp\cap \Fn_x^\perp(F))=\vol(\Lambda\cap \Fn_x(F))\vol(\Lambda)^{-1}$.
	
	For $n\geq n_0$, we let $\varphi_n\in C_c^\infty(\mathfrak{g}_x(F))$ be the function that is equal to $X\mapsto \vol(L_n)^{-1}\psi(\langle \Xi,X\rangle )$ on $L_n$ and is equal to zero outside. We define similarly $\phi_n\in C_c^\infty(G_x(F))$ as the function that coincides with $\vol(K_n)^{-1}\xi_n$ on $K_n$ and is equal to zero outside. Note that, $\phi_n(\exp(X))=\varphi_n(X)$ for every $X\in L_n$ and that the Fourier transform $\widehat{\varphi}_n$ is the characteristic function of the coset $-\xi+L_n^\perp$. Here the Fourier transform is normalized as in Subsection \ref{Section qc} and for any $\cO_F$-lattice $\mathcal{L}\subset \mathfrak{g}(F)$, we denote by $\mathcal{L}^\perp$ the dual lattice in $\mathfrak{g}^*(F)$ with respect to $\psi$, that is
	$$\displaystyle \mathcal{L}^\perp=\{Y\in \mathfrak{g}^*(F)\mid \psi(\langle Y,X\rangle)=1,\; \forall X\in \mathcal{L}\}.$$

	\begin{lem}\label{lem lattice L_n}
		For $n$ large enough, the following hold.
		\begin{enumerate}[(i)]
			\item $\Xi+(L_n^\perp\cap \Fn_{x}^{\perp}(F))$ is invariant under the conjugation of $K_n\cap N_x(F)$.
			\item For $u\in N_x(F)$, if $u^{-1} (\Xi+(L_n^\perp\cap \Fn_{x}^{\perp}(F))) u \cap  \Xi+(L_n^\perp\cap \Fn_{x}^{\perp}(F)) \neq \emptyset$, then $u\in K_n\cap N_x(F)$.
			
			\item For every $f\in C_c^\infty(\fg_x(F))$, we have
			\begin{equation*}
				\displaystyle \int_{\fg_x(F)} \varphi_n(X) f(X)dX=\vol(K_n)^{-1}\vol(K_n\cap N_x(F))^{-1}\int_{K_n}\int_{\Xi+\mathfrak{n}_x^\perp(F)} \mathbf{1}_{\Xi+L_n^\perp}(Y) \widehat{f}(kYk^{-1})dYdk.
			\end{equation*}
			
			\item For every $\cO'\in \Nil(\widehat{\Fg}_x)$, the coadjoint orbital integral of $\widehat{\varphi}_n$ on $\cO'$ (normalized using the Kirillov-Kostant measure as in Subsection \ref{Section qc}) is given by
			$$\displaystyle \int_{\mathcal{O}'}\widehat{\varphi}_n(Y)dY=\left\{ \begin{array}{ll}
				1 & \mbox{ if } \mathcal{O}'=-\mathcal{O}, \\
				0 & \mbox{ otherwise.}
			\end{array}\right.$$
		\end{enumerate}
	\end{lem}
	
	\begin{proof}
		\begin{enumerate}[(i)]
			\item We have
			$$\displaystyle \Xi+(L_n^\perp\cap \Fn_{x}^{\perp}(F))=\left(\Xi+\Fn_x^\perp(F)\right)\cap \left(\Xi+L_n^\perp\right).$$
			Furthermore, since $\Xi$ restricts to a $N_x(F)$-invariant character on $\Fn_x(F)$, $\Xi+\Fn_x^\perp(F)$ is $N_x(F)$-invariant whereas, since the function $\xi_n$ is $K_n$-invariant, $\Xi+L_n^\perp$ is $K_n$-invariant. The claim follows.

			\item Let $u\in N_x(F)$ and set $X=\log(u)\in \Fn_x(F)$. After conjugating everything by $a_n$, the statement is equivalent to
			$$\displaystyle u^{-1} (\Xi+(\varpi^nL^\perp\cap \Fn_{x}^{\perp})) u \cap  \Xi+(\varpi^nL^\perp\cap \Fn_{x}^{\perp})\neq \emptyset \Rightarrow X\in \varpi^n L.$$
			By the theory of Kostant section, the map 
			$$\displaystyle N_x\times \Xi+\Fn_{x}^{\perp}\to  \left(\Xi+\Fn_{x}^{\perp}\right)\times \left( \Xi+ \Fn_{x}^{\perp}\right),\; (n,X)\mapsto (X,n^{-1}Xn)$$
			is a closed embedding. Hence, for any $m>0$ we can choose $n$ large enough such that if
			$$\displaystyle u^{-1} (\Xi+(\varpi^nL^\perp\cap \Fn_{x}^{\perp})) u \cap  \Xi+(\varpi^nL^\perp\cap \Fn_{x}^{\perp})\neq \emptyset,$$
			then $X\in \varpi^m L$. Now let $k$ be the largest integer such that $X\in \varpi^k L$. We know that $k\geq m$ and we need to show $k\geq n$. There exists some absolute constant $C\in \mathbb{N}$ such that
			$$u^{-1}\Xi u\in \Xi-[X,\Xi]+\varpi^{2k-C} L^{\perp},\; u^{-1} (\varpi^nL^\perp) u \subset \varpi^{n} L^{\perp}+ \varpi^{n+k-C} L^{\perp}.$$
			If $n>k$, once we choose $m$ to be large enough (with respect to $C$), the above relations imply that
			$$u^{-1} (\Xi+(\varpi^nL^\perp\cap \Fn_{x}^{\perp}(F))) u \subset  \Xi+[X,\Xi]+(\varpi^{k+1}L^\perp\cap \Fn_{x}^{\perp}(F)).$$
			Since the image of $L$ in $\fg_x/\fg_{x,\Xi}$ is self-dual, we know that $[X,\Xi]\in \varpi^{k}L^\perp\cap \Fn_{x}^{\perp}(F)$ and $[X,\Xi]\notin \varpi^{k+1}L^\perp\cap \Fn_{x}^{\perp}(F)$. This implies that 
			$$\displaystyle u^{-1} (\Xi+(\varpi^nL^\perp\cap \Fn_{x}^{\perp}(F))) u \cap  \Xi+(\varpi^{k+1}L^\perp\cap \Fn_{x}^{\perp}(F))=\emptyset$$
			which is a contradiction. Hence we must have $k\geq n$ and this proves the lemma.
			
			\item Let $D$ be the distribution on $\Fg_x(F)$ defined by
			$$\displaystyle D(f)=\int_{K_n}\int_{\Xi+\mathfrak{n}_x^\perp(F)} \mathbf{1}_{\Xi+L_n^\perp}(Y) \widehat{f}(kYk^{-1})dYdk,\;\; f\in C_c^\infty(\Fg_x(F)).$$
			Then, it has the following properties:
			\begin{enumerate}[(a)]
				\item It is $K_n$-invariant: $D({}^kf)=D(f)$ for every $(k,f)\in K_n\times C_c^\infty(\Fg_x(F))$;
				
				\item It is supported in $\Ad(K_n)(L_n+\Fn_x(F))$: this follows from the Fourier inversion formula
				$$\displaystyle \int_{\Xi+\mathfrak{n}_x^\perp(F)} \mathbf{1}_{\Xi+L_n^\perp}(Y) \widehat{f}(kYk^{-1})dY=\int_{L_n+\mathfrak{n}_x(F)} f(k^{-1}Xk) \psi(\langle \Xi,X\rangle) dX;$$
				
				\item It is $(L_n,\xi_n)$-equivariant: $D(L(X)f)=\psi(\langle \Xi,X\rangle)D(f)$ for every $(X,f)\in L_n\times C_c^\infty(\Fg_x(F))$ (this is a consequence of the same Fourier inversion formula and the fact that the restriction of $\psi(\langle \Xi,.\rangle)$ to $L_n$ is $K_n$-invariant).
			\end{enumerate}
			We claim:
			\begin{num}
				\item\label{prop0 lem lattice L_n} For $n$ large enough, every distribution on $\Fg_x(F)$ satisfying the properties $(a)$, $(b)$ and $(c)$ above is proportional to the distribution
				$$\displaystyle f\in C_c^\infty(\Fg_x(F))\mapsto \int_{\Fg_x(F)} \varphi_n(X) f(X)dX.$$
			\end{num}
			
			Indeed, every distribution $D$ verifying $(a)$, $(b)$, $(c)$ is represented by a function $F\in C^\infty(\fg_x(F))$ which is $K_n$-invariant, satisfies $F(Y+X)=\psi(\langle \Xi,X\rangle )F(Y)$ for $(X,Y)\in L_n\times \Fg_x(F)$ and is supported in $\Ad(K_n)(L_n+\Fn_x(F))$. It then suffices to show that such a function is necessarily supported in $L_n$ which would be a consequence of the following property: for every $X\in \Fn_x(F)\setminus (L_n\cap \Fn_x(F))$, we can find $k\in K_n$ such that $\Ad_k X-X\in L_n$ and $\psi(\langle \Xi,\Ad(k)X-X\rangle)\neq 1$. Conjugating everything by $a_n$, this property can be restated as:
			
			\begin{num}
				\item\label{prop1 lem lattice L_n} Provided $n$ is large enough, for every $X\in \Fn_x(F)\setminus (L\cap \Fn_x(F))$, we can find $k\in K'_n$ such that $\Ad(k)X-X\in L$ and $\psi(\varpi^{-n}\langle \Xi,\Ad_kX-X\rangle)\neq 1$.
			\end{num}
			
			Indeed, let $X\in \Fn_x(F)\setminus (L\cap \Fn_x(F))$ and set $-k=val_L(X)<0$. Set $m=\max(n,k+2n_0)$. Then, for every $Y\in \varpi^{m}L$, by the last and third points of \eqref{properties n0} respectively, we have
			$$\displaystyle \Ad_{e^Y}(X)-X \in [Y,X]+\varpi^{2m-k-n_0}L\subset [Y,X]+\varpi^{n+n_0}L$$
			and
			$$\displaystyle [Y,X]\in \varpi^{m-k-n_0}L\subset \varpi^{n_0} L.$$
			Thus, $\Ad_{e^Y}(X)-X\in L$ and, by the fourth point of \eqref{properties n0}, 
			$$\displaystyle \psi(\varpi^{-n}\langle \Xi,\Ad_{e^Y}(X)-X\rangle )=\psi(\varpi^{-n}\langle \Xi,[Y,X]\rangle )$$
			for every $Y\in \varpi^m L$. Since $\exp(\varpi^m L)\subset K'_n$ it therefore suffices to find $Y\in \varpi^m L$ such that $\psi(\varpi^{-n}\langle \Xi,[Y,X]\rangle)\neq 1$ (provided $n$ is large enough). However, if $\psi(\varpi^{-n}\langle \Xi,[Y,X]\rangle)=1$ for all $Y\in \varpi^m L$ then, since the lattice $L^\xi\subset \Fg_x(F)/\Fg_{x,\Xi}(F)$ is self-dual, the image of $X$ in $\Fg_x(F)/\Fg_{x,\Xi}(F)$ belongs to $\varpi^{n-m}L^\xi$. As the preimage of $L^\xi$ in $\Fn_x(F)$ is $\Fn_x(F)\cap L$, this would imply $Y\in \varpi^{n-m}L$ hence $n-m\leq -k$ or equivalently $n\leq m-k=\max(n-k,2n_0)$. This last inequality is obviously false for $n>2n_0$ so that the claim \eqref{prop1 lem lattice L_n} is satisfied for such a $n$.
			
			This shows \eqref{prop0 lem lattice L_n}. As a consequence, we can find a constant $c$ such that
			\begin{equation*}
				\displaystyle \int_{\fg_x(F)} \varphi_n(X) f(X)dX=c\int_{K_n}\int_{\Xi+\mathfrak{n}_x^\perp(F)} \mathbf{1}_{\Xi+L_n^\perp}(Y) \widehat{f}(kYk^{-1})dYdk
			\end{equation*}
			for every $f\in C_c^\infty(\Fg_x(F))$. Plugging in $f=\overline{\varphi_n}$, we have $\widehat{f}=\mathbf{1}_{\Xi+L_n^\perp}$ and we obtain
			$$\displaystyle \vol(L_n)^{-1}=c\vol(K_n)\vol(L_n^\perp\cap \Fn_x(F)^\perp).$$
			By our choice of measures, we have
			$$\displaystyle \vol(L_n)\vol(L_n^\perp\cap \Fn_x(F)^\perp)=\vol(L_n\cap \Fn_x(F))=\vol(K_n\cap N_x(F))$$
			and therefore $c=\vol(K_n)^{-1}\vol(K_n\cap N_x(F))^{-1}$ as claimed.
			
			\item This follows from the computation in the middle of p.437 of \cite{MW}.
		\end{enumerate}
	\end{proof}

	\subsection{Application of the local trace formula}
	
	We now start the proof of Theorem \ref{thm: formula for the germ} which will be finished in the next subsection. The proof for general reductive twisted spaces is basically the same as in the untwisted case, i.e. when $\tG=G$, the only difference is to replace the local trace formula in \cite{A91} (resp. Howe's conjecture for weighted orbital integrals in \cite[Lemma 8.2]{Art94} \cite{BarM}) by the local twisted trace formula in \cite{WalFTLtordue} (resp. Howe's conjecture for twisted weighted orbital integrals in Appendix B). Hence to simplify notation, we will only write the proof when $\tG=G$. For further simplification of notation, we will also assume that the split center $A_G$ is trivial.

	We need to recall some material on the local trace formula from \cite{A91}. Fix a minimal Levi subgroup $M_{min}$ of $G$ as well as $P_{min}\in \mathcal{P}(M_{min})$ and a special maximal compact subgroup $K\subset G(F)$ in good position relative to $M_{min}$. We set $W=Norm_{G(F)}(M_{min})/M_{min}(F)$. Let $T\in \CA_{\min,\mathbb{Q}}$. For $P\in \mathcal{P}(M_{min})$, we set $T_P=w_PT$ where $w_P\in W$ is the unique element such that $w_PP_{min}w_P^{-1}=P$ and, for $g_1,g_2\in G(F)$, we define a $(G,M_{min})$-orthogonal set by
	$$\displaystyle \mathcal{Y}_P(g_1,g_2,T)=T_P+H_P(g_1)-H_{\overline{P}}(g_2),\;\; P\in \mathcal{P}(M_{min}),$$
	where $\overline{P}$ denotes the parabolic subgroup opposite to $P$ (with respect to $M_{\min}$). For $M\in \mathcal{L}(M_{\min})$, we set
	$$\displaystyle v_M(g_1,g_2,T)=\int_{A_M(F)} \Gamma_M(H(a), \mathcal{Y}(g_1,g_2,T)) da.$$
	
	We also denote by $M_{min,\leqslant T}$ the set those $m\in M_{min}(F)$ such that
	$$\displaystyle 0\leqslant \langle \alpha, H_{M_{min}}(m)\rangle \mbox{ and } \langle \varpi_\alpha, H_{M_{min}}(m)\rangle \leqslant \langle \varpi_\alpha, T\rangle, \mbox{ for every } \alpha\in \Delta_{min},$$
	and we let $u(.,T)$ be the characteristic function of the subset $KM_{min,\leqslant T}^+ K$ of $G(F)$.
	
	For $f\in \mathcal{C}(G(F))$ and $\varphi\in C_c^\infty(G(F))$, we define
	$$\displaystyle \widetilde{J}^T(f,\varphi)=\int_{G(F)} u(g,T)\int_{G(F)} f(g^{-1}xg) \varphi(x)dx dg.$$
	and
	$$\displaystyle J^T(f,\varphi):=\sum_{M\in \mathcal{L}(M_{min})} \frac{\lvert W^M\rvert}{\lvert W\rvert} \int_{\Gamma_{ell}(M)} J^T_M(\gamma,f,\varphi) d\gamma$$
	where
	$$\displaystyle J^T_M(\gamma,f,\varphi)=D^G(\gamma)\int_{\left(A_M(F)\backslash G(F) \right)^2} f(g_1^{-1}\gamma g_1) \varphi(g_2^{-1}\gamma g_2) v_M(g_1,g_2,T) dg_1dg_2.$$

	\begin{prop}\label{Howe conjecture}
		\begin{enumerate}[(i)]
			\item For every $f\in \mathcal{C}(G(F))$ and $\varphi\in C_c^\infty(G(F))$, the function
			$$\displaystyle T\in \CA_{min,\mathbb{Q}}\mapsto J^T(f,\varphi)$$
			is a unitary polynomial-exponential. Moreover, for any compact-open subgroup $J\subset G(F)$ and subset $\Omega\subset G(F)$ that is bounded modulo conjugation, the subspace of $\cC(J\backslash G(F)/J)^*$ spanned by
			$$\displaystyle \left\{f\in \mathcal{C}(J\backslash G(F)/J)\mapsto J^T(f,\varphi) \mid \varphi\in C_c^\infty(\Omega), T\in \CA_{\min,\mathbb{Q}}\right\}$$
			is finite dimensional.
			
			\item Let $\epsilon>0$. Then, for every $r>0$ there exists a constant $C_r>0$ such that
			$$\displaystyle \left\lvert J^T(f,\varphi)-\widetilde{J}^T(f,\varphi)\right\rvert\leqslant C_r \lVert T\rVert^{-r}$$
			for every $T\in \CA_{P_{min},\mathbb{Q}}^+$ satisfying $\langle \alpha,T\rangle \geqslant \epsilon \lVert T\rVert$ for every $\alpha\in \Delta_{min}$. 
			
			\item For every $T\in \CA_{\min,\mathbb{Q}}$, $\varphi\in C_c^\infty(G(F))$ and strongly cuspidal function $f\in \cC_{scusp}(G(F))$ we have
			\begin{equation*}
				\displaystyle J^T(f,\varphi)=\int_{G(F)} \Theta_f(g)\varphi(g) dg.
			\end{equation*}
		\end{enumerate}
	\end{prop}
	
	\begin{proof}
		\begin{enumerate}[(i)]
			\item is a consequence of the splitting formulas \cite[Equation (5.5)]{Art94} and of Howe's conjecture for weighted orbital integrals \cite[Lemma 8.2]{Art94}. (see Appendix \ref{Appendix Howe} for the twisted case.)
			
			\item follows directly from the proof of the geometric side of the local trace formula in \cite{A91}.
			
			\item follows from the splitting formula \cite[Equation (5.5)]{Art94} and the same argument as in Section \ref{Section TLTF strongly cuspidal}.
		\end{enumerate}
	\end{proof}

	Let $\omega_x\subset \mathfrak{g}_x(F)$ be a sufficiently small invariant neighborhood of $1$ and define, for $n$ sufficiently large, $\phi_n^G\in C_c^\infty(G(F))$ by (here $K_x=K\cap G_x$)
	$$\displaystyle \phi_n^G(g)=\left\{\begin{array}{ll}
		\int_{K_x}\phi_n(a_n^{-1}k_x^{-1}yk_xa_n) dk_x & \mbox{ if } g=kxyk^{-1} \mbox{ for some } (y,k)\in \exp(\omega_x)\times K \\
		0 & \mbox{ otherwise.}
	\end{array} \right.$$
	
	For every $f\in \mathcal{C}(G(F))$, we set
	\begin{equation}
		\displaystyle \widetilde{J}_{x,\xi,n}^T(f)=\int_{B_x(F)\backslash G(F)} \int_{N_x(F)} f(g^{-1}xug) \widetilde{v}_{B_x,\xi}(u, \mathcal{Y}(g,T)+H(a_n)) du dg
	\end{equation}
	where $\mathcal{Y}(g,T)$ denotes the $(G,M(x))$-family defined by 
	$$\displaystyle \mathcal{Y}_P(g,T)=H_P(g)+T_P, \;\; P\in \mathcal{P}(M(x)),$$
	and where we recall that for every $(G,M(x))$-family $\mathcal{X}$ and $u\in N_x(F)$, we have set
	$$\displaystyle \widetilde{v}_{B_x,\xi}(u, \mathcal{X}):=\int_{T_x(F)}^{reg} \xi(t^{-1}ut) \Gamma_{B_x}(H(t), \mathcal{X}) dt.$$

	\begin{prop}\label{prop germ formula}
		\begin{enumerate}[(i)]
			\item Let $f\in \cC_{scusp}(G(F))$. Then, for $n$ large enough we have
			$$\displaystyle J^T(f,\phi_n^G)=c_{f,-\mathcal{O}}(x)$$
			for every $T\in \CA_{P_{min},\mathbb{Q}}$.
			
			\item Let $f\in C_c^\infty(G(F))$. Then, there exists $n_f>0$ and $C=C_f>0$ (both depend on the support and the level of $f$) such that
			$$\displaystyle \widetilde{J}^T(f,\phi_n^G)=\widetilde{J}_{x,\xi,n}^T(f)$$
			for every $n\geqslant n_f$ and $T\in \CA_{P_{min}}^+$ satisfying $\alpha(T)\geqslant C$ for every $\alpha\in \Delta_{min}$.
		\end{enumerate}
	\end{prop}
	
	\begin{proof}
		\begin{enumerate}[(i)]
			\item Applying Proposition \ref{Howe conjecture} (iii) to $\varphi=\phi_n^G$ and by usual arguments of semisimple descent and descent to the Lie algebra, together with the germ expansion of $\Theta_f$, for $n$ large enough we get
			\[\begin{aligned}
				\displaystyle J^T(f,\varphi) & =\int_{G(F)} \Theta_f(g) \phi_n^G(g) dg =\int_{G_x(F)} \Theta_f(xy)\phi_n(y) dy \\
				& =\sum_{\mathcal{O}'\in \Nil(\widehat{\fg}_x)} c_{f,\mathcal{O}'}(x) \int_{\mathcal{O}'}\widehat{\varphi}_n(Y)dY.
			\end{aligned}\]
			The result then follows immediately from Lemma \ref{lem lattice L_n} (iv).
			
			\item We may assume that $f$ is invariant under $K$-conjugation. For every $g\in G(F)$, we define a function ${}^g f_{x,\omega_x}\in C_c^\infty(\fg_x(F))$ by
			$$\displaystyle {}^g f_{x,\omega_x}(X)=\left\{\begin{array}{ll}
				f(g^{-1}x\exp(X)g) & \mbox{ if } X\in \omega_x, \\
				0 & \mbox{ otherwise.}
			\end{array} \right.$$
			Then, by a standard descent argument to the Lie algebra and Lemma \ref{lem lattice L_n} (iii), we have
			\[\begin{aligned}
				\displaystyle \widetilde{J}^T(f,\phi_n^G) & = \int_{G(F)}  u(g,T)\int_{\fg_x} {}^gf_{x,\omega_x}(X) \varphi_n(a_n^{-1}Xa_n) dX dg \\
				&=\int_{G(F)} u(a_n g,T) \int_{\fg_x} {}^gf_{x,\omega_x}(X) \varphi_n(X) dX dg \\
				& =\vol(K_n)^{-1}\vol(K_n\cap N_x(F))^{-1} \int_{G(F)} u(a_ng,T)\int_{K_n} \int_{\Xi+\mathfrak{n}^\perp_x} \widehat{{}^{kg}f}_{x,\omega_x}(Y) \mathbf{1}_{\Xi+L_n^\perp}(Y)dY dk
			\end{aligned}\]
			where the first equality follows from the change of variables $g\mapsto a_n g$ and $X\mapsto a_n Xa_n^{-1}$. Since the function $g\mapsto u(a_ng,T)$ is left-invariant by $K_n$ for $n$ large enough, this gives
			\[\begin{aligned}
				\displaystyle J^T(f,\varphi_n) & =\vol(K_n\cap N_x)^{-1} \int_{G(F)} u(a_ng,T) \int_{\Xi+(L_n^\perp\cap \Fn_{x}^{\perp})} \widehat{{}^gf}_{x,\omega_x}(Y) dY dg \\
				& =\vol(K_n\cap N_x)^{-1} \int_{B_x(F)\backslash G(F)} \int_{B_x(F)}u(a_nbg,T) \int_{\Xi+(L_n^\perp\cap \Fn_{x}^{\perp})} \widehat{{}^gf}_{x,\omega_x}(b^{-1}Yb) dY d_Lbdg.
			\end{aligned}\]
			Thus, it suffices to show that, for every $g\in G(F)$, we have
			\begin{align}\label{eq4}
				\displaystyle & \vol(K_n\cap N_x)^{-1}\int_{B_x(F)}u(a_nbg,T) \int_{\Xi+(L_n^\perp\cap \Fn_{x}^{\perp})} \widehat{{}^gf}_{x,\omega_x}(b^{-1}Yb) dY d_Lb \\
				\nonumber & =\int_{N_x(F)} f(g^{-1}xug) \widetilde{v}_{B_x,\xi}(u, \mathcal{Y}(g,T)+H(a_n)) du.
			\end{align}
			Note that both sides of the above equation are $(B_x,\delta_{B_x})$-equivariant on the left. Moreover, since $f$ is compactly supported and $x$ is semisimple, for $g$ outside of a set that is compact modulo $G_x(F)$, the function ${}^g f$ is zero. As $B_x(F)$ is cocompact in $G_x(F)$, we may restrict to establish the above identity for $g$ in some fixed compact subset $\mathcal{K}\subset G(F)$. We need a lemma.

			\begin{lem}
			Let us set $T_x^+[\geq -C]:= \{ t\in T_x(F)\mid \langle \alpha, H_{T_x}(t)\rangle\leq -C\}$ for every $C>0$. Then, there exist two large enough constants $C>0$ and $n_1>0$ such that:
			\begin{enumerate}[(i)]
			\item For every $g\in \mathcal{K}$, the function
			$$\displaystyle T_x(F)\times (\Xi+\mathfrak{n}_x^\perp(F))\ni (t,Y)\mapsto \widehat{{}^g f}_{x,\omega_x}(tYt^{-1})$$
			is supported on $T_x^+[\geq -C]\times (\Xi+\mathfrak{n}_x^\perp(F))$.
			
			\item For every $g\in \mathcal{K}$ and $n\geq n_1$, the function
			$$\displaystyle B_x(F)\times (\Xi+L_n^\perp\cap \mathfrak{n}_x^\perp(F))\ni (b,Y)\mapsto \widehat{{}^g f}_{x,\omega_x}(bYb^{-1})$$
			is supported on $T_x^+[\geq -C](N_x(F)\cap K_n)\times (\Xi+\mathfrak{n}_x^\perp(F))$.
		\end{enumerate}
			\end{lem}
		
		\begin{proof}
		(i) There exists a compact $\mathcal{K}^L\subset \fg_x(F)$ such that $\Supp(\widehat{{}^gf}_{x,\omega_x})\subset \mathcal{K}^L$ for every $g\in \mathcal{K}$. Thus, it suffices to see that $(t\Xi t^{-1}+\Fn_{x}^{\perp}(F))\cap \mathcal{K}^L\neq \emptyset$ implies $t\in T_x^+[\geq -C]$ for $C>0$ large enough but this is clear as decomposing $\Xi$ along eigenspaces for $A_x(F)$, it has nonzero components in all the root subspaces corresponding the simple negative roots with respect to $B_x$.
		
		(ii) Let $b\in B_x(F)$ and $Y\in \Xi+L_n^\perp\cap \mathfrak{n}_x^\perp(F)$ be such that $\widehat{{}^g f}_{x,\omega_x}(bYb^{-1})\neq 0$ and write $b$ as $b=tu$ where $t\in T_x(F)$ and $u\in N_x(F)$. Since $uYu^{-1}\in \Xi+\mathfrak{n}_x^\perp(F)$, by point (i) we have $t\in T_x^+[\geq -C]$ for some $C>0$ large enough. Moreover, we have $u (\Xi+(L_n^\perp\cap \Fn_{x}^{\perp}(F))) u^{-1} \cap t^{-1}\CK^L t\neq \emptyset$. As
		$$\displaystyle u(\Xi+(L_n^\perp\cap \Fn_{x}^{\perp}(F))) u^{-1}\subset \Xi+\Fn_{x}^{\perp}(F),$$
		and there exists a compact subset $\mathcal{K}'\subset \Xi+ \mathfrak{n}_x^\perp(F)$ such that $t^{-1}\CK^L t\cap(\Xi+\mathfrak{n}_x^\perp(F))\subset \mathcal{K}'$ for every $t\in T_x^+[\geq -C]$, this implies
		$$\displaystyle u (\Xi+(L_n^\perp\cap \Fn_{x}^{\perp}(F))) u^{-1} \cap \mathcal{K}'\neq \emptyset.$$
		Using that $(L_n^\perp\cap \Fn_{x}^{\perp}(F))_n$ is an increasing and exhausting family of compact subsets of $\mathfrak{n}_x^\perp(F)$, for $n$ large enough it follows that
		$$\displaystyle u (\Xi+(L_n^\perp\cap \Fn_{x}^{\perp})) u^{-1} \cap  \Xi+(L_n^\perp\cap \Fn_{x}^{\perp})\neq \emptyset$$
		which implies, by Lemma \ref{lem lattice L_n}, that $u\in N_x(F)\cap K_n$.
		\end{proof}

			Let $g\in \mathcal{K}$. By the lemma above and since $\Xi+(L_n^\perp\cap \Fn_{x}^{\perp})$ is invariant under the conjugation of $N_x(F)\cap K_n$ (Lemma \ref{lem lattice L_n}), we have
			\begin{align}
				\displaystyle & \vol(K_n\cap N_x)^{-1}\int_{B_x(F)}u(a_nbg,T) \int_{\Xi+(L_n^\perp\cap \Fn_{x}^{\perp})} \widehat{{}^gf}_{x,\omega_x}(b^{-1}Yb) dY d_Lb \\
				\nonumber & = \int_{T_x(F)} u(a_nt^{-1}g,T) \int_{\Xi+(L_n^\perp\cap \mathfrak{n}_x^\perp(F))} \widehat{{}^gf}_{x,\omega_x}(tYt^{-1}) dY \delta_{B_x}(t)^{-1}dt.
			\end{align}
			and moreover the integrand is supported in $T_x^+[\geq -C]$ for some large enough $C>0$ (independent of $n$). For $n$ large enough, as seen in the proof of the last lemma, for every $t\in T_x^+[\geq -C]$ the function
			$$\displaystyle Y\in \Xi+\mathfrak{n}_x^\perp(F)\mapsto \widehat{{}^gf}_{x,\omega_x}(tYt^{-1})$$
			is supported in $\Xi+(L_n^\perp\cap \mathfrak{n}_x^\perp(F))$. Hence we can replace the above integral over $\Xi+(L_n^\perp\cap \mathfrak{n}_x^\perp(F))$ by an integral over $\Xi+\mathfrak{n}_x^\perp(F)$ (since the integrand stays supported in $T_{x}^+[\geq -C]$ by the lemma). Then by a simple change of variable, we obtain that the above expression is equal to
			\begin{equation}\label{eq3}
			\displaystyle	\int_{T_x(F)} u(a_nt^{-1}g,T) \int_{t\Xi t^{-1}+\mathfrak{n}_x^\perp(F)} \widehat{{}^gf}_{x,\omega_x}(Y) dY dt.
			\end{equation}
			Next we show that
			
			\begin{num}
				\item\label{eq1} For $n$ sufficiently large and $T$ sufficiently regular, we have
				$$\displaystyle u(a_nt^{-1}g,T)=\Gamma_{B_x}(H(t),\mathcal{Y}(g,T)+H(a_n))$$
				for every $g\in \mathcal{K}$ and $t\in T_{x}^+[\geq -C]$.
			\end{num}
			
			Let $t$ be as above. First, we note that
			$$\displaystyle \Gamma_{B_x}(H(t),\mathcal{Y}(g,T)+H(a_n))=\Gamma_{B_x}(H(a_n^{-1}t),\mathcal{Y}(g,T))$$
			and
			$$\displaystyle  \langle \alpha, H_{T_x}(a_n^{-1}t)\rangle=2n\log(q)+\langle \alpha, H_{T_x}(t)\rangle\geqslant 2n\log(q)-C$$
			for all $\alpha\in \Delta_x$. From this, we see that \eqref{eq1} actually reduces to the following statement:
			
			\begin{num}
				\item\label{eq2} There exist $C'>0$ such that
				$$\displaystyle u(t^{-1}g,T)=\Gamma_{B_x}(H(t),\mathcal{Y}(g,T))$$
				for every $g\in \mathcal{K}$, $t\in T_x(F)$ and $T\in \cA_{P_{min}}^+$ satisfying
				$$\displaystyle \langle \alpha, H(t)\rangle\geqslant 0,\;\; \forall \alpha\in \Delta_x,$$
				and
				$$\displaystyle \langle \alpha, T\rangle\geqslant C',\;\; \forall \alpha\in \Delta_{min}.$$
			\end{num}
			
			Let $t\in T_x(F)$ and $T\in \cA_{P_{min}}^+$ be elements satisfying the above inequalities. Then, provided $C'$ is large enough, the $(G,M(x))$-orthogonal set $\mathcal{Y}(g,T)$ is positive for every $g\in \mathcal{K}$ and therefore by Proposition \ref{prop partition gammaBx functions} and the assumption on $t$ we have
			$$\displaystyle \Gamma_{B_x}(H(t),\mathcal{Y}(g,T))=\Gamma_{M(x)}(H(t),\mathcal{Y}(g,T)).$$
			Furthermore, by the identity at the bottom of p.38 of \cite{A91}, provided again that $C'$ is large enough, we also have
			$$\displaystyle \Gamma_{M(x)}(H(t),\mathcal{Y}(g,T))=u(t^{-1}g,T)$$
			for every $g\in \mathcal{K}$. This shows \eqref{eq2} and ends the proof of \eqref{eq1}.
			
			Now, from \eqref{eq1} and \eqref{eq3}, we deduce that, for $n$ sufficiently large and $T$ sufficiently regular, we have
			\[\begin{aligned}
				\displaystyle & \vol(K_n\cap N_x)^{-1}\int_{B_x(F)}u(a_nbg,T) \int_{\Xi+(L_n^\perp\cap n_x)} \widehat{{}^gf}_{x,\omega_x}(b^{-1}Yb) dY d_Lb \\
				& =\int_{T_x(F)} \Gamma_{B_x}(H(t),\mathcal{Y}(g,T)+H(a_n)) \int_{t\Xi t^{-1}+\mathfrak{n}_x^\perp} \widehat{{}^gf}_{x,\omega_x}(Y) dY dt \\
				& =\int_{T_x(F)} \Gamma_{B_x}(H(t),\mathcal{Y}(g,T)+H(a_n)) \int_{N_x(F)} f(g^{-1}xug) \xi(t^{-1}ut) du dt \\
				& =\int_{N_x(F)} f(g^{-1}xug) \int^{reg}_{T_x(F)} \xi(t^{-1}ut) \Gamma_{B_x}(H(t),\mathcal{Y}(g,T)+H(a_n)) dt du \\
				& =\int_{N_x(F)} f(g^{-1}xug) \tilde{v}_{B_x,\xi}(u, \mathcal{Y}(g,T)+H(a_n)) du
			\end{aligned}\]
			for every $g\in \mathcal{K}$. This gives \eqref{eq4} and therefore closes the proof of the proposition.
		\end{enumerate}
	\end{proof}

	\subsection{End of the proof of Theorem \ref{thm: formula for the germ}}
	In this subsection we will prove the formula of regular germs in Theorem \ref{thm: formula for the germ}. Fix a strongly cuspidal function $f\in \mathcal{C}(G(F))$, we need to show that
	\begin{equation*}
		\displaystyle c_{f,-\CO_{\xi}}(x)=\int_{B_x(F)\backslash G(F)} \int_{N_x(F)} f(g^{-1}x u g) v_{B_x,\xi}(u,g) du dg.
	\end{equation*}
	By Proposition \ref{prop germ formula}(i), there exists $n_f>0$ such that for $n>n_f$ we have
	$$\displaystyle J^T(f,\phi_n^G)=c_{f,-\mathcal{O}_\xi}(x)$$
	for every $T\in \CA_{P_{min},\mathbb{Q}}$. Let $J$ be an open compact subgroup of $G(F)$ by which $f$ is biinvariant. By Proposition \ref{Howe conjecture}(i), once we choose $n_f$ large enough, we can find a sequence of functions $f_N\in C_{c}^{\infty}(G(F))^{J\times J}$ such that
	\begin{itemize}
		\item $f_N\rightarrow f$ (in $\mathcal{C}(G(F))$) as $N\rightarrow \infty$;
		\item $J^T(f,\phi_n^G)=J^T(f_N,\phi_n^G)$ for all $N>0$, $T$ and $n>n_f$.
	\end{itemize}
	Indeed, by Proposition \ref{Howe conjecture}(i), we know that the span of the linear forms 
	$$f\in \CC(G(F))^{J\times J} \mapsto J^T(f,\phi_n^G)$$
	for all $T$ and $n>n_f$ is finite dimensional. Let $J_1,\cdots, J_k$ be a basis of this span. Since these linear forms are continuous, by density we know that $J_1,\cdots, J_k$ are also linearly independent when restricted to $C_c^\infty(G(F))^{J\times J}$. Thus, we can find $g_i\in C_c^\infty(G(F))^{J\times J}$ ($i=1,\ldots,k$) such that $J_i(g_j)=\delta_{i,j}$. Choose now an arbitrary sequence $f'_N\in C_c^\infty(G(F))^{J\times J}$ converging to $f\in \CC(G(F))^{J\times J}$. Then the modified sequence $f_N=f'_N+\sum_{i=1}^k (J_i(f)-J_i(f'_N))g_i$ satisfies the required conditions.

	Let $N>0$ be fixed for the moment. By Proposition \ref{prop germ formula}(ii), we have
\begin{equation}\label{eq10}
\tilde{J}^T(f_N,\phi_n^G)=\tilde{J}_{x,\xi,n}^T(f_N)=\int_{B_x(F)\backslash G(F)} \int_{N_x(F)} f_N(g^{-1}xug) \widetilde{v}_{B_x,\xi}(u, \mathcal{Y}(g,T)+H(a_n)) du dg
\end{equation}
	for $n$ sufficiently large and $T$ sufficiently regular (both with respect to $f_N$). Let $\epsilon>0$ be sufficiently small. Then, by the local trace formula (Proposition \ref{Howe conjecture} (ii)), as $T\in \mathcal{A}_{P_{\min},\mathbb{Q}}^+$ goes to $\infty$ in the cone
	$$\mathcal{C}=\{\langle \alpha, T\rangle\geq \epsilon \lVert T\rVert,\; \forall\; \alpha\in \Delta_{min}\},$$ $\tilde{J}^T(f_N,\phi_n^G)$ is asymptotic to the polynomial-exponential (with unitary exponents) $T\mapsto J^T(f_N,\phi_n^G)=J^T(f,\phi_n^G)$. On the other hand, by Lemma \ref{lem exponential polynomial}, there exists a constant $C>0$ such that the equality of weights
	$$\widetilde{v}_{B_x,\xi}(u, \mathcal{Y}(g,T)+H(a_n))=v_{B_x,\xi}(u, \mathcal{Y}(g,T)+H(a_n))$$
	holds whenever the depth of the $(G, M(x))$-orthogonal set $\mathcal{Y}(g,T)+H(a_n)$ is bigger than $C\sigma(u)$. Thus, as $f_N$ is compactly supported, this holds for $T$ large enough in the cone $\mathcal{C}$ and for every $(g,u)\in \mathcal{K}\times N_x(F)$ such that $f_N(g^{-1}xug)\neq 0$. Thus, as $T\to \infty$ in $\mathcal{C}$, the right hand side of \eqref{eq10} is also asymptotic to the polynomial-exponential
	$$T\mapsto \int_{B_x(F)\backslash G(F)} \int_{N_x(F)} f_N(g^{-1}xug) v_{B_x,\xi}(u, \mathcal{Y}(g,T)+H(a_n)) du dg.$$
	From this, we deduce the equality of polynomial-exponentials
	$$\displaystyle J^T(f,\phi_n^G)=\int_{B_x(F)\backslash G(F)} \int_{N_x(F)} f_N(g^{-1}xug) v_{B_x,\xi}(u, \mathcal{Y}(g,T)+H(a_n)) du dg$$
	for every $T\in \mathcal{A}_{P_{min}, \mathbb{Q}}$ and $n$ large enough. However, we know that the left hand side is identically equal to $c_{f,-\mathcal{O}_\xi}(x)$ whereas the right-hand side is a polynomial-exponential in both $T$ and $H(a_n)$. This polynomial-exponential is therefore constant and the same identity holds for $T=0$ and $H(a_n)=0$ which gives the equality of Theorem \ref{thm: formula for the germ} except with $f_N$ instead of $f$ in the right-hand side. Thus, letting $N\to \infty$, we obtain the desired identity.

	\subsection{A descent formula}\label{secion descent formula}
	In this subsection, we will prove a descent formula that will be used in later section. We keep the notation as in the previous subsections. Moreover, we assume that $\iota(x)=x$, $B_x$ is $\iota$-split and $T_{x}^{\iota}\subset Z_{G_x}$.
	
	The action of $\iota$ naturally descends to $\CA_x$ and this induces a decomposition $\CA_x=\CA_x^\iota\oplus \CA_{x,\iota}$ where $\CA_x^\iota$ (resp. $\CA_{x,\iota}$) denotes the subspace consisting of elements $H\in \CA_x$ satisfying $\iota(H)=H$ (resp. $\iota(H)=-H$). For $H\in \CA_x$ we will denote without further comment by $H^\iota$, $H_\iota$ the respective projections of $X$ with respect to this decomposition. Similarly, if $C$ is a subset of $\CA_x$ (typically the positive cone associated to a parabolic subspace) we will denote by $C^\iota$ the image of its projection to $\CA_x^\iota$.
	
	Let $\mathcal{X}=(X_{\tilde{P}})_{\tP\in \mathcal{P}(\tM(x))}$ be a $(\tG,\tM(x))$-orthogonal set. For every $\tQ\in \mathcal{F}_{B_x,\iota}(\tM(x))$ (resp. $\tQ\in \mathcal{F}_{B_x}(\tM(x))$ such that $\CA_x^\iota\cap \CA^{\tQ}=0$), we define a function $\Gamma_{B_x,\iota}^{\tQ}(.,\mathcal{X})$ (resp. $\Gamma_{B_x}^{\tQ,\iota}(.,\mathcal{X})$) on $\CA_x$ by
	\begin{equation*}
		\displaystyle \Gamma_{B_x,\iota}^{\tQ}(H,\mathcal{X})=\sum_{\tP\in \mathcal{F}_{B_x,\iota}(\tM(x)), \tP\subset \tQ} (-1)^{a_{\tP,\iota}^{\tQ}} \widehat{\tau}^{\tQ}_{\tP,\iota}(H-X_{\tP,\iota}),\;\; H\in \mathfrak{a}_{x},
	\end{equation*}
	\begin{equation*}
		\displaystyle (\mbox{resp. }\Gamma_{B_x}^{\tQ,\iota}(H,\mathcal{X})=\left\{\begin{array}{ll}
			\Gamma_{B_x}^{\tQ}(Y^{\tQ},\mathcal{X}) & \mbox{ if } H\in X_{\tQ}+Y^{\tQ}+\CA_{x}^\iota+\CA_{\tG} \mbox{ for some } Y^{\tQ}\in \CA^{\tQ}; \\
			0 & \mbox{ otherwise.}
		\end{array} \right.,\;\;\; H\in \CA_x.)
	\end{equation*}
	
	\begin{prop}\label{prop iota version}
		For every $\tR\in \mathcal{F}_{B_x,\iota}(\tM(x))$, we have the following identity of functions on $\CA_x$:
		\begin{equation}\label{equation iota version 1}
			\displaystyle \sum_{\tR\supset \tQ\in \mathcal{F}_{B_x,\iota}(\tM(x))} \Gamma_{B_x,\iota}^{\tQ}(.,\mathcal{X}) \tau_{\tQ,\iota}^{\tR}(.-X_{\tQ})=1.
		\end{equation}
		If $\mathcal{X}$ is positive, $\Gamma_{B_x,\iota}^{\tG}(.,\mathcal{X})$ is the characteristic function of either of the two following subsets
		\begin{equation}
			\displaystyle \left\{ H\in \CA_x\mid \varpi_\alpha(H-X_{\tP,\iota})\leqslant 0,\; \forall \tP\in \mathcal{P}_{B_x,\iota}(\tM(x)),\forall \alpha\in \Delta_{\tP,\iota} \right\},
		\end{equation}
		\begin{equation}
			\displaystyle Conv\left\{X_{\tP,\iota}\mid \tP\in \mathcal{P}_{B_x,\iota}(\tM(x)) \right\}+{}^- \CA_{B_x,\iota}+\CA_x^\iota+\CA_{\tG}.
		\end{equation}
	Moreover, if $\CX$ is positive and $\CY=(Y_{\tilde{P}})_{\tP\in \mathcal{P}(\tM(x))}$ is another positive $(\tG,\tM(x))$-orthogonal sets, then, for every $\tQ,\tR\in \mathcal{F}_{B_x,\iota}(\tM(x))$ and $Y\in \CA_x$ we have
	\begin{equation}\label{equation iota version 2}
	\displaystyle \Gamma_{B_x,\iota}^{\tQ}(H,\mathcal{X}) \tau^{\tR}_{\tQ,\iota}(H-X_{\tQ})\Gamma_{B_x,\iota}^{\tR}(H,\mathcal{X}+\mathcal{Y})=\Gamma_{B_x,\iota}^{\tQ}(H,\mathcal{X}) \tau^{\tR}_{\tQ, \iota}(H-X_{\tQ})\phi_{\tQ,\iota}^{\tR}(H-X_{\tQ}-Y_{\tQ})
\end{equation}
	where $\phi_{\tQ,\iota}^{\tR}$ denotes the characteristic function of the set of those $Z\in \CA_{\tQ,\iota}$ such that $\langle \varpi, Z\rangle\leq 0$ for every $\varpi\in \widehat{\Delta}_{\tQ,\iota}^{\tR}$.
	\end{prop}

	\begin{proof}
		The proof is basically the same as for Proposition \ref{prop partition gammaBx functions}, Proposition \ref{prop aleternative description GammaBxX} and Corollary \ref{cor gammaBx functions} adding some $\iota$'s in indices along the way. We skip the details.
	\end{proof}

	\begin{prop}\label{prop descent}
		Assume that $\mathcal{X}$ is positive and let $\varepsilon\in (\CA_{B_x}^+)^\iota$ that is in general position. For every $\tQ\in \CF_{B_x}(\tM(x))$ such that $\varepsilon\in (\CA_{\tQ}^+)^\iota$ we define a number $d_\varepsilon(\tQ)$ inductively by the relation
		\begin{equation}\label{eq0 prop descent}
			\displaystyle \sum_{\substack{\tQ\subset \tR\in \CF_{B_x}(\tM(x)) \\ \varepsilon\in (\CA_{\tR}^+)^\iota}} d_\varepsilon(\tR)=1.
		\end{equation}
		Then, we have the following  equality of functions on $\CA_x$
		\begin{equation}\label{eq1 prop descent}
			\displaystyle \Gamma_{B_x,\iota}^{\tG}(.,\mathcal{X})=\sum_{\substack{\tQ\in \mathcal{F}_{B_x}(\tM(x)) \\ \varepsilon\in (\CA_{\tQ}^{+})^\iota}} d_\varepsilon(\tQ) \Gamma_{B_x}^{\tQ,\iota}(.,\mathcal{X}).
		\end{equation}
	\end{prop}
	
	\begin{proof}
		For $\tQ\in \CF_{B_x}(\tM(x))$, we set
		$$\displaystyle \mathcal{C}^{\tQ}_{B_x}(\CX):=Conv \left\{X_{\tP}\mid \tP\in \CP_{B_x}(\tM(x)), \tP\subset \tQ \right\}+{}^- \CA^{\tQ}_{B_x}.$$
		Then $\Gamma_{B_x}^{\tQ}(.,\CX)$ is the characteristic function of $\CC^{\tQ}_{B_x}(\CX)+\CA_{\tQ}$ (by Proposition \ref{prop aleternative description GammaBxX}) and:
		\begin{num}
			\item\label{eq1 proof prop descent} If $\CA_x^\iota\cap \CA^{\tQ}=0$, $\Gamma_{B_x}^{\tQ,\iota}(.,\CX)$ is the characteristic function of $\CC_{B_x}^{\tQ}(\CX)\oplus (\CA_x^\iota+\CA_{\tG})$.
		\end{num}
		(This follows from the definition of $\Gamma_{B_x}^{\tQ,\iota}(.,\CX)$ and the previous point.) Furthermore, we claim that:
		\begin{num}
			\item\label{eq2 proof prop descent} $\Gamma_{B_x,\iota}^{\tG}(.,\mathcal{X})$ is the characteristic function of $\mathcal{C}_{B_x}^{\tG}(\CX)+\CA_x^\iota+\CA_{\tG}$.
		\end{num} 
		Indeed, by the previous proposition it suffices to check that for every $\tP\in \CP_{B_x}(\tM(x))$, $\tP'\in \CP_{B_x,\iota}(\tM(x))$ and $\alpha\in \Delta_{\tP',\iota}$ we have
		$$\displaystyle \varpi_\alpha(X_{\tP,\iota}-X_{\tP',\iota})\leqslant 0.$$
		But this follows, after projection onto $\CA_{x,\iota}$, from the fact that $X_{\tP,\iota}-X_{\tP',\iota}$ is a linear combination with negative coefficients of elements of $\Delta_{\tP'}$ (by definition of a positive $(\tG,\tM(x))$-family).
		
		With the terminology and notation from Appendix \ref{appendix convex sets}, we also have:
		\begin{num}
			\item\label{eq3 proof prop descent} $\CC^{\tG}_{B_x}(\CX)$ is a finitely generated convex set with faces $F^{\tQ}:=\CC^{\tQ}_{B_x}(\CX)$, $\tQ\in \CF_{B_x}(\tM(x))$, and corresponding (open) cones $\CA_{F^{\tQ}}^{+}=\CA_{\tQ}^+$.
		\end{num}
		Indeed, that $\CC^{\tG}_{B_x}(\CX)$ is a finitely generated convex set is clear from its definition. Let $\lambda\in \CA_x$ and $c\in \BR$ be such that $\langle \lambda,H\rangle\leqslant c$ for every $H\in \CC^{\tG}_{B_x}(\CX)$. Applying this inequality to $H\in {}^-\CA_{B_x}$, we see that $\lambda\in \overline{\CA_{B_x}^{+}}$. As
		$$\displaystyle \overline{\CA_{B_x}^{+}}=\bigcup_{\tQ\in \CF_{B_x}(\tM(x))} \CA_{\tQ}^{+},$$
		we have $\lambda\in \CA_{\tQ}^{+}$ for some $\tQ\in \CF_{B_x}(\tM(x))$. For $H\in {}^-\CA_{B_x}$, we have
		$$\displaystyle \langle \lambda, H\rangle\leqslant 0$$
		with equality if and only if $H\in {}^- \CA_{B_x}^{\tQ}$. Furthermore, as $\CX$ is positive, for every $\tP\in \CP_{B_x}(\tM(x))$ we have
		$$\displaystyle \langle \lambda, X_{\tP}\rangle \leqslant \langle \lambda, X_{\tQ}\rangle$$
		with equality if and only if $\tP\subset \tQ$. Therefore,
		$$\displaystyle \langle \lambda,H\rangle\leqslant \langle \lambda, X_{\tQ}\rangle$$
		for $H\in \CC^{\tG}_{B_x}(\CX)$ with equality if and only if $H\in \CC^{\tQ}_{B_x}(\CX)$ and it follows that the intersection
		$$\displaystyle \CC^{\tG}_{B_x}(\CX)\cap \{H\in \CA_x\mid \langle \lambda,H\rangle=c \}$$
		is either empty or equal to $\CC^{\tQ}_{B_x}(\CX)$. The claim \eqref{eq3 proof prop descent} follows.
		
		From \eqref{eq3 proof prop descent} and Proposition \ref{prop projection convex sets}(i) (applied to $\Fb=\CA_{x,\iota}$), we deduce that
		\begin{equation}\label{eq4 proof prop descent}
			\displaystyle \CC_{B_x}^{\tG}(\CX)+\CA_x^\iota+\CA_{\tG}=\bigcup_{\tQ\in \CF_{B_x,\varepsilon}(\tM(x))} \CC_{B_x}^{\tQ}(\CX)+\CA_x^\iota+\CA_{\tG}
		\end{equation}
		where we have denoted by $\CF_{B_x,\varepsilon}(\tM(x))$ the subset of $\tQ \in \CF_{B_x}(\tM(x))$ such that $\varepsilon\in (\CA_{\tQ}^+)^\iota$.
		
		Thus, by \eqref{eq1 proof prop descent} and \eqref{eq2 proof prop descent}, to get the identity \eqref{eq1 prop descent} it only remains to check that
		\begin{equation}\label{eq5 proof prop descent}
			\displaystyle \sum_{\substack{\tQ\in \CF_{B_x,\varepsilon}(\tM(x)) \\ H\in \CC_{B_x}^{\tQ}(\CX)+\CA_x^\iota+\CA_{\tG}}} d_\varepsilon(\tQ)=1
		\end{equation}
		for every $H\in \CC_{B_x}^{\tG}(\CX)+\CA_x^\iota+\CA_{\tG}$. By Proposition \ref{prop projection convex sets}(ii), there exists a minimal $\tQ\in \CF_{B_x,\varepsilon}(\tM(x))$ such that $H\in \CC_{B_x}^{\tQ}(\CX)+\CA_x^\iota+\CA_{\tG}$ and, by the relation \eqref{eq0 prop descent}, it suffices to show that, for $\tR\in \CF_{B_x,\varepsilon}(\tM(x))$, we have $\CC_{B_x}^{\tQ}(\CX)\subset  \CC_{B_x}^{\tR}(\CX)$ if and only if $\tQ\subset \tR$ but this follows from \eqref{eq3 proof prop descent} (as this shows that both inclusions are equivalent to $\CA_{\tR}^+\subset \overline{\CA_{\tQ}^+}$).
	\end{proof}

	As in the previous subsection, let $N_{x,reg}\subset N_x$ be the subset of regular elements in $N_x$ and $T_{x,c}\subset T_x(F)$ be the maximal compact subgroup. We equip $T_{x,c}$ with the Haar measure of total mass $1$ and we also fix a log-norm $\sigma_{reg}: N_{x,reg}(F)\to \mathbb{R}_{\geqslant 1}$ on $N_{x,reg}(F)$. Set $r=\dim(\mathfrak{a}_x)$. The next two lemmas can be proved by the same argument as in Lemma \ref{lem def vBxpsi}. We will skip the proofs here.
	
	\begin{lem}\label{lem def twisted vBxpsi}
		For any $u\in N_{x,reg}(F)$ and any positive $(\tG,\tM(x))$-orthogonal set $\mathcal{X}$, the iterated integral
		\begin{equation}
			\displaystyle \int_{T_x(F)/A_{\tG}(F)T_{x}^{\iota}(F)} \int_{T_{x,c}} \xi(a^{-1}t^{-1}uta) dt \Gamma_{B_x,\iota}(H_{T_x}(a),\mathcal{X}) da
		\end{equation}
		is absolutely convergent in that order and will be denoted by
		\begin{equation*}
			\displaystyle \tilde{v}_{B_x,\xi,\iota}(u,\mathcal{X}):=\int_{T_x(F)/A_{\tG}(F)T_{x}^{\iota}(F)}^* \xi(a^{-1}ua) \Gamma_{B_x,\iota}(H(a),\mathcal{X}) da.
		\end{equation*}
		Moreover, there exists a constant $C>0$ such that for every $u\in N_{x,reg}(F)$ and every positive $(\tG,\tM(x))$-orthogonal set $\mathcal{X}$, we have
		\begin{equation*}
			\displaystyle \left\lvert \tilde{v}_{B_x,\xi,\iota}(u,\mathcal{X})\right\rvert\leqslant C(\sigma_{reg}(u)+N(\mathcal{X}))^r.
		\end{equation*}
	\end{lem}
	
	\begin{lem}
		For any $u\in N_{x,reg}(F)$, $\tQ\in \mathcal{F}_{B_x}(\tM(x))$ such that $\CA_x^\iota\cap \CA^{\tQ}=0$ and any positive $(\tG,\tM(x))$-orthogonal set $\mathcal{X}$, the iterated integral
		\begin{equation}
			\displaystyle \int_{T_x(F)/A_{\tG}(F)T_{x}^{\iota}(F)} \int_{T_{x,c}} \xi(a^{-1}t^{-1}uta) dt \Gamma_{B_x}^{\tQ,\iota}(H_{T_x}(a),\mathcal{X}) da
		\end{equation}
		is absolutely convergent in that order and will be denoted by
		\begin{equation*}
			\displaystyle \tilde{v}_{B_x,\xi}^{\tQ,\iota}(u,\mathcal{X}):=\int_{T_x(F)/A_{\tG}(F)T_{x}^{\iota}(F)}^* \xi(a^{-1}ua) \Gamma_{B_x}^{\tQ,\iota}(H(a),\mathcal{X}) da.
		\end{equation*}
		Moreover, there exists a constant $C>0$ such that for every $u\in N_{x,reg}(F)$ and every positive $(\tG,\tM(x))$-orthogonal set $\mathcal{X}$, we have
		\begin{equation*}
			\displaystyle \left\lvert \tilde{v}_{B_x,\xi}^{\tQ,\iota}(u,\mathcal{X})\right\rvert\leqslant C(\sigma_{reg}(u)+N(\mathcal{X}))^r.
		\end{equation*}
	\end{lem}
	
	\begin{lem}\label{cor invariant 1}
		There exists $C>0$ such that for every $(\tG,\tM(x))$-orthogonal set $\mathcal{X}$ satisfying $d(\mathcal{X})\geqslant C\sigma(u)$, we have
		$$\tilde{v}_{B_x,\xi}^{\tQ,\iota}(u,\mathcal{X})=\tilde{v}_{B_x,\xi}^{\tQ,\iota}(u^Q,\mathcal{X})$$
		where $u=u^Qu_Q$ is the unique decomposition with $u^Q\in L_{\tQ}(F)$ and $u_Q\in U_{\tQ}(F)$.
	\end{lem}

	\begin{proof}
		The proof is the same as the proof of the second bullet point in the proof of Lemma \ref{lem exponential polynomial} (we just need to use our assumption that $T_{x}^{\iota}\subset Z_{G_x}$).
	\end{proof}

	The next two lemmas can be proved by the same argument as in Lemma \ref{lem exponential polynomial}. We will skip the proof here.

	\begin{lem}\label{exponential polynomial twisted weight}
		There exists $C>0,r>0$ and, for every $u\in N_{x,reg}(F)$, a unique exponential polynomial $v_{B_x,\xi,\iota}(u,.)\in Pol_{\leqslant r}$ whose exponents belong to a finite set independent of $u$ such that for every $(\tG,\tM(x))$-orthogonal set $\mathcal{X}$ satisfying $d(\mathcal{X})\geqslant C\sigma(u)$, we have
		$$\displaystyle v_{B_x,\xi,\iota}(u,\mathcal{X})=\widetilde{v}_{B_x,\xi,\iota}(u,\mathcal{X}).$$
		Moreover, there exists $C'>0$ and $R>0$ such that for every $u\in N_{x,reg}(F)$ and every $(\tG,\tM(x))$-orthogonal set $\mathcal{X}$ we have
		\begin{equation*}
			\displaystyle \left\lvert v_{B_x,\xi,\iota}(u,\mathcal{X})\right\rvert\leqslant C'(\sigma_{reg}(u)+N(\mathcal{X}))^R.
		\end{equation*}
	\end{lem}
	
	\begin{lem}
		For $\tQ\in \mathcal{F}_{B_x}(\tM(x))$ such that $\CA_x^\iota\cap \CA^{\tQ}=0$, there exists $C>0,r>0$ and, for every $u\in N_{x,reg}(F)$, a unique exponential polynomial $v_{B_x,\xi}^{\tQ,\iota}(u,.)\in Pol_{\leqslant r}$ whose exponents belong to a finite set independent of $u$ such that for every $(\tG,\tM(x))$-orthogonal set $\mathcal{X}$ satisfying $d(\mathcal{X})\geqslant C\sigma(u)$, we have
		$$\displaystyle v_{B_x,\xi}^{\tQ,\iota}(u,\mathcal{X})=\widetilde{v}_{B_x,\xi}^{\tQ,\iota}(u,\mathcal{X}).$$
		Moreover, there exists $C'>0$ and $R>0$ such that for every $u\in N_{x,reg}(F)$ and every $(\tG,\tM(x))$-orthogonal set $\mathcal{X}$ we have
		\begin{equation*}
			\displaystyle \left\lvert v_{B_x,\xi}^{\tQ,\iota}(u,\mathcal{X})\right\rvert\leqslant C'(\sigma_{reg}(u)+N(\mathcal{X}))^R.
		\end{equation*}
	\end{lem}
	
	Following the above two lemmas, we define
	\begin{equation*}
		\displaystyle v_{B_x,\xi,\iota}(u,g)=v_{B_x,\xi,\iota}(u,\mathcal{Y}(g)),\;v_{B_x,\xi}^{\tQ,\iota}(u,g)=v_{B_x,\xi}^{\tQ,\iota}(u,\mathcal{Y}(g))
	\end{equation*}
	
	The following corollary is a direct consequence of the two lemmas above.
	
	\begin{cor}
		There exists $d>0$ such that 
		$$v_{B_x,\xi,\iota}(u,g)\ll \sigma_{G}(g)^d \sigma_{N_{x,reg}}(u)^d,\;v_{B_x,\xi}^{\tQ,\iota}(u,g)\ll \sigma_{G}(g)^d \sigma_{N_{x,reg}}(u)^d$$
		for all $u\in N_{x,reg}(F)$ and $g\in G(F)$.
	\end{cor}

	\begin{cor}\label{cor invariant 2}
		The function $v_{B_x,\xi}^{\tQ,\iota}(u,g)$ is left $N_x\cap U_{\tQ}(F)$ on $u$ and left $U_{\tQ}(F)$-invariant on $g$.
	\end{cor}
	
	\begin{proof}
		The left $U_{\tQ}(F)$-invariant on $g$ is clear from the definition. The left $N_x\cap U_{\tQ}(F)$ on $u$ follows from Lemma \ref{cor invariant 1}.
	\end{proof}
	
	\begin{cor}\label{cor descent}
		We have the decent formula
		$$v_{B_x,\xi,\iota}(u,g)=\sum_{\substack{\tQ\in \mathcal{F}_{B_x}(\tM(x)) \\ \varepsilon\in (\CA_{\tQ}^{+})^\iota}} d_\varepsilon(\tQ)v_{B_x,\xi}^{\tQ,\iota}(u,g).$$
	\end{cor}
	
	\begin{proof}
		This is a direct consequence of Proposition \ref{prop descent}.
	\end{proof}

	\section{On the spectral expansion}

	Let $(G,\tG)$ be a connected reductive twisted space over $F$. Let $H$ be a closed unimodular subgroup of $G$ defined over $F$ and $(H,\tH)$ be a twisted space over $F$ equipped with an embedding $\tH\subset \tG$ which is $H\times H$-equivariant. Let $(\chi,\tilde{\chi})$ be a one-dimensional unitary representation of $\tH(F)$ i.e. $\chi: H(F)\to \mathbb{C}^\times$ is a (smooth) unitary character and $\tilde{\chi}:\tH(F)\to \mathbb{C}^\times$ is a map satisfying $\tilde{\chi}(h_1\th h_2)=\chi(h_1h_2)\tilde{\chi}(\th)$ for $(\th,h_1,h_2)\in \tH(F)\times H(F)\times H(F)$. Let $\omega$ be a character of $A_{\tG}(F)$ which coincides with $\chi$ on the intersection $A_{\tG}(F)\cap H(F)$.
	
	Denote by $L^2(H(F)A_{\tG}(F)\backslash G(F),\chi\otimes \omega)$ the Hilbert space of functions $\varphi: G(F)\to \mathbb{C}$ satisfying $\varphi(hag)=\chi(h)\omega(a)\varphi(g)$ for $(h,a,g)\in H(F)\times A_{\tG}(F)\times G(F)$ and such that $g\mapsto \lvert \varphi(g)\rvert^2$ is integrable on $H(F)A_{\tG}(F)\backslash G(F)$. The representation by right translation of $G(F)$ on that space will be denoted by $R$. This extends to a twisted representation $\tilde{R}$ of $\tG(F)$ defined by
	$$\displaystyle (\tilde{R}(\th g)\varphi)(x)=\tilde{\chi}(\th)\varphi(\Ad_{\th}^{-1}(x)g)$$
	for every $(\th,g,x)\in \tH(F)\times G(F)\times G(F)$ and $\varphi\in L^2(H(F)A_{\tG}(F)\backslash G(F),\chi\otimes \omega)$. For $f\in C_c^\infty(\tG(F)/A_{\tG}(F),\omega^{-1})$, the operator $\tilde{R}(f)$ is given by
	\begin{equation*}
		\displaystyle (\tilde{R}(f)\varphi)(x)=\int_{\tG(F)/A_{\tG}(F)} f(\tg) (\tilde{R}(\tg)\varphi)(x) d\tg,\;\; \varphi\in L^2(H(F)A_{\tG}(F)\backslash G(F),\chi\otimes \omega).
	\end{equation*}
	This operator is associated with the kernel function $\nu(\tH)^{-1} K_f(x,y)$ where
	\begin{equation}\label{eq1 tss}
		\displaystyle K_f(x,y)=\int_{\tH(F)/A_{\tG}^{H}(F)} f(x^{-1}\th y) \tilde{\chi}(\th) d\th,\;\; x,y\in G(F)
	\end{equation}
	and $\nu(\tH)=|H(F)\cap A_{\tG}(F):A_{\tG}^{H}(F)|$. Here $A_{\tG}^{H}$ is the maximal split torus of $A_{\tG}\cap H$. We define
	$$\displaystyle I(f)=\int_{H(F)A_{\tG}(F)\backslash G(F)} K_f(x,x) dx, \mbox{ for } f\in C_c^\infty(\tG(F)/A_{\tG}(F),\omega^{-1})$$
	provided the integral is absolutely convergent.

	If the pair $(G,H)$ is tempered (see Section \ref{S tempered varieties} for the definition of tempered), we can define in a similar way operators $\tilde{R}(f)$ for $f\in \mathcal{C}(\tG(F)/A_{\tG}(F),\omega)$ and these operators are associated to kernel functions given by the same expression \eqref{eq1 tss} (which is absolutely convergent) and we also define $I(f)$ by the same formula provided the integral is abstolutely convergent.
	
	Let now $f$ be in $C_c^\infty(\tG(F)/A_{\tG}(F),\omega^{-1})$ or, if $X$ is tempered, in $\mathcal{C}(\tG(F)/A_{\tG}(F),\omega^{-1})$ and assume that it satisfies the following very strong condition:
	\begin{num}
		\item\label{assumption Rf} The operator $\tilde{R}(f)$ is of finite rank.
	\end{num}
	This implies that the integral defining $I(f)$ is convergent and equals $\nu(\tH)\Tr \tilde{R}(f)$:
	$$\displaystyle I(f)=\nu(\tH) \Tr \tilde{R}(f).$$
	
	Let $L^2_{\disc}(H(F)A_{\tG}(F)\backslash G(F),\chi\otimes \omega)$ be the sum of all the irreducible unitary subrepresentations of $L^2(H(F)A_{\tG}(F)\backslash G(F),\chi\otimes \omega)$ and $L^2_{\cont}(H(F)A_{\tG}(F)\backslash G(F),\chi\otimes \omega)$ be its orthogonal complement. The assumption \eqref{assumption Rf} also implies that $\tilde{R}(f)$ acts by zero on $L^2_{\cont}(H(F)A_{\tG}(F)\backslash G(F),\chi\otimes \omega)$, therefore
	$$\displaystyle I(f)=\nu(\tH) \Tr \tilde{R}_{\disc}(f)$$
	where $\tilde{R}_{\disc}(f)$ stands for the restriction of $\tilde{R}(f)$ to $L^2_{\disc}(H(F)A_{\tG}(F)\backslash G(F),\chi\otimes \omega)$.
	
	Let $\Pi_{\disc}(H(F)A_{\tG}(F)\backslash G(F),\chi\otimes \omega)$ be the set of isomorphism classes of irreducible subrepresentations of $L^2(H(F)A_{\tG}(F)\backslash G(F),\chi\otimes \omega)$. Then, we have the isotypic decomposition
	$$\displaystyle L_{\disc}^2(H(F)A_{\tG}(F)\backslash G(F),\chi\otimes \omega)=\bigoplus_{\pi\in \Pi_{disc}(H(F)A_{\tG}(F)\backslash G(F),\chi\otimes \omega)} \pi\otimes M_{L^2}(\pi)$$
	where $M_{L^2}(\pi):=\Hom_G(\pi, L^2(H(F)A_{\tG}(F)\backslash G(F),\chi\otimes \omega))$ are multiplicity spaces. Let $\Pi_{\disc}(H(F)A_{\tG}(F)\backslash G(F),\chi\otimes \omega)^\theta$ be the subset of isomorphism classes fixed by $\theta$ and choose for every $\pi \in\Pi_{\disc}(H(F)A_{\tG}(F)\backslash G(F),\chi\otimes \omega)^\theta$ an extension $\tpi$ of $\pi$ to a representation of the twisted space $\tG(F)$. Then, there is an unique endomorphism $\theta\langle \tpi\rangle$ of $M_{L^2}(\pi)$ such that the restriction of $\tilde{R}(\tg)$ to the isotypic component $\pi\otimes M_{L^2}(\pi)$ is equal to $\tpi(\tg)\otimes \theta\langle \tpi\rangle$ for $\tg\in \tG(F)$. Using these notations, and under the assumption \eqref{assumption Rf}, we have
	$$\displaystyle \Tr(\tilde{R}_{\disc}(f))=\sum_{\pi\in \Pi_{\disc}(H(F)A_{\tG}(F)\backslash G(F),\chi\otimes \omega)} \Tr(\tpi(f))\times \Tr(\theta \langle \tpi\rangle\mid M_{L^2}(\pi))$$
	for $f\in C_{c}^{\infty}(\tG(F)/A_{\tG}(F),\omega^{-1})$ (or $f\in \mathcal{C}(\tG(F)/A_{\tG}(F),\omega^{-1})$ if $(G,H)$ is tempered). Note that a priori we didn't assume the multiplicity spaces $M(\pi)$ to be of finite dimension but, by the assumption \eqref{assumption Rf}, this is automatic whenever $\Tr(\tpi(f))\neq 0$, so that the above expression makes sense.
	
	Summarizing the discussion so far, we have the following proposition:
	
	\begin{prop}\label{prop twisted spectral}
		Let $f$ be in $C_c^\infty(\tG(F)/A_{\tG}(F),\omega^{-1})$ or, if $X$ is tempered, in $\mathcal{C}(\tG(F)/A_{\tG}(F),\omega^{-1})$ and assume that it satisfies \eqref{assumption Rf}. Then, the integral defining $I(f)$ converges and, with the above notation, we have
		$$\displaystyle I(f)=\nu(\tH)\sum_{\pi\in \Pi_{\disc}(H(F)A_{\tG}(F)\backslash G(F),\chi\otimes \omega)} \Tr(\tpi(f))\times \Tr(\theta \langle \tpi\rangle\mid M(\pi)).$$
	\end{prop}
	
	When $X=H\backslash G$ is wavefront spherical and $G$ is split \cite[Theorem 9.2.1]{SV} or when $X$ is symmetric \cite[Theorem 4]{DelPlanch}, we have\footnote{This property is of course expected to hold for all spherical varieties.}:
	\begin{num}\label{twisted spectral 1}
		\item For every compact-open subgroup $J\subset G(F)$, the subspace
		$$\displaystyle L^2_{\disc}(H(F)A_{\tG}(F)\backslash G(F),\chi\otimes \omega)^J$$
		of $J$-invariants in $L_{\disc}^2(H(F)A_{\tG}(F)\backslash G(F),\chi\otimes \omega)$ is finite dimensional.
	\end{num}
	This readily implies that for every $f$ in $C_c^\infty(\tG(F)/A_{\tG}(F),\omega^{-1})$ or, if $X$ is tempered, in $\mathcal{C}(\tG(F)/A_{\tG}(F),\omega^{-1})$, the operator $\tilde{R}_{\disc}(f)$ is of finite rank so that, in those cases, we have
	\begin{num}\label{twisted spectral 2}
		\item the assumption \eqref{assumption Rf} is equivalent to $\tilde{R}(f)=\tilde{R}_{\disc}(f)$. 
	\end{num}
	
	Two other situations where condition \eqref{assumption Rf} is automatically satisfied are as follows:
	\begin{num}
		\item\label{twisted spectral 3} $\bar{f} \in C_c^\infty(\tG(F)/A_{\tG}(F),\omega)$ is a matrix coefficient of a supercuspidal representation $(\pi,\tpi)$ of $\tG(F)$ with
		$$\displaystyle m_{L^2}(\pi):=\dim M_{L^2}(\pi)<\infty.$$
		
		\item\label{twisted spectral 4} The pair $(G,H)$ is tempered and $\bar{f}\in \cC(\tG(F)/A_{\tG}(F),\omega)$ is a matrix coefficient of a discrete series representation $(\pi,\tpi)$ of $\tG(F)$ with
		$$\displaystyle m_{L^2}(\pi):=\dim M_{L^2}(\pi)<\infty.$$
	\end{num}
	By \cite{Delcstterm}, the finite multiplicity assumption in \ref{twisted spectral 3} and \ref{twisted spectral 4} is automatically satisfied when $H=H_0\ltimes N$, with $N$ the unipotent radical of some parabolic subgroup $P=MN$ of $G$, $H_0$ a symmetric subgroup of a Levi factor $M$ (i.e. there exists an involution $\iota$ of $M$ such that $(M^\iota)^0\subset H_0\subset M^\iota$), and the restriction of the character $\chi$ to $N(F)$ is {\em generic} (in the sense that its orbit under the adjoint action of $M(F)$ is open in the group of all continuous characters $\Hom_{\cont}(N(F),\BC^\times)$).

	\section{The geometric expansion}

	\subsection{The setup}
	Let $(\tG,\iota)$ be a twisted symmetric pair (see \S \ref{sect twisted symmetric pairs}) with $G$ connected and reductive, $\tP=\tM N$ be a $\iota$-split parabolic subspace with $\tM=\tP\cap \iota(\tP)$, and $\xi:N(F)\rightarrow \BC^{\times}$ be a generic character of $N(F)$ i.e. a character whose the stabilizer in $M$ under the adjoint action is of minimal dimension. Let $H_0=(M^\iota)^{\circ}$ and $\tH_0$ be a connected component of the subvariety of $\iota$-fixed points $\tM^\iota$ so that $(H_0,\tH_0)$ is a twisted reductive space. We make the following two assumptions:
	
	\begin{itemize}
		\item $\tH_0$ stabilizes the character $\xi$ under the adjoint action. Moreover, if $\xi$ is nontrivial (i.e. if $\tP$ is a proper parabolic subspace), we assume that $H_0$ is the neutral component of the stabilizer of the character $\xi$ in $M$.
		\item The twisted symmetric pair $(\tM,\tH_0)$ is coregular in the sense of Subsection \ref{section coregular varieties}, i.e. $\tH_0(F)\cap \tM_{\rs}(F)\neq \emptyset$ and the function
		$$t\in \tH_0(F)\cap \tM_{\rs}(F)\mapsto \frac{D^{\tH_0} (t)}{D^{\tM} (t)^{1/2}}$$
		is locally bounded on $\tH_0(F)$.
	\end{itemize}
	For every $h\in \tH_{0,ss}(F)$, we have
	$$\displaystyle D^{\tH}(h)=D^{\tH_0}(h)D^{\tG}(h)^{1/2}D^{\tM}(h)^{-1/2}\delta_{\tP}(h)^{-1/2}$$
	and $\delta_{\tP}(h)=\delta_{\tP}(\iota(h))^{-1}=\delta_{\tP}(h)^{-1}$, since $\tP$ is $\iota$-split, hence $\delta_{\tP}(h)=1$. Therefore, the second assumption implies that:
	\begin{num}
		\item the function $h\in \tH_{0}(F)\cap \tM_{rs}(F)\mapsto \frac{D^{\tH} (h)}{D^{\tG} (h)^{1/2}}$ is locally bounded on $\tH_0(F)$.
	\end{num}

	We set
	$$\displaystyle \tH=\tH_0\ltimes N$$
	and we denote again by $\xi: \tH(F)\to \mathbb{C}^\times$ the twisted character that is trivial on $\tH_0(F)$ and coincides with the previous character on $N(F)$. We also fix a unitary twisted character $\chi:\tH_0(F)\to \BC^\times$ and we denote by $\xi\otimes \chi$ the twisted character of $\tH(F)$ given by 
	$$\displaystyle \xi\otimes \chi:h_0u\in \tH(F)=\tH_0(F)\ltimes N(F)\mapsto \chi(h_0)\xi(u).$$

	Let $t\in \tH_{0,rs}(F)$. By the coregular assumption, $t$ is also regular in $\tM$ and this implies that $G_t$ is quasi-split over $F$ with $P_t=M_tN_t$ as a Borel subgroup where $G_t$ (resp. $P_t$, $M_t$ and $N_t$) denotes the neutral component of the centralizer of $t$ in $G$ (resp. in $P$, $M$ and $N$). Let $\xi_t$ be the restriction of $\xi$ to $N_t(F)$. Similarly, if we let $S=H_{0,t}$ and $\tT=St$,  $G_{\tS}$ is quasi-split over $F$, $M_{\tS}N_{\tS}$ is a Borel subgroup of $G_{\tS}$ and we let $\xi_{\tS}$ be the restriction of $\xi$ to $N_{\tS}(F)$ where $G_{\tS}$ (resp. $P_{\tS}$, $M_{\tS}$ and $N_{\tS}$) is the centralizer of $\tS$ in $G$ (resp. $P$, $M$ and $N$). Note that $G_t=G_{\tS}$ for almost all $t\in \tS(F)$ and $M_{\tS}\cap H_0=S$ belongs to the center of $G_{\tS}$.

	\begin{lem}\label{lem generic character}
		With the notation above, $\xi_{\tS}$ is a generic character of $N_{\tS}(F)$.
	\end{lem}
	
	\begin{proof}
	We denote by the same letter the pullbacks of $\xi$ and $\xi_{\tS}$ to $\mathfrak{n}(F)$ and $\mathfrak{n}_{\tS}(F)$ (via the exponential maps). Let $\mathfrak{n}^{\tS}$ be the unique $\Ad(\tS)$-stable complement of $\mathfrak{n}_{\tS}$ in $\mathfrak{n}$. Then, since $\tS$ stabilizes $\xi$, $\xi$ is trivial on $\mathfrak{n}^{\tS}(F)$ and it follows that an element of $M_{\tS}$ stabilizes $\xi_{\tS}$ if and only if it stabilizes $\xi$. However, by our first assumption, $S=M_{\tS}\cap H_0$ is the neutral component of the stabilizer of $\xi$ in $M_{\tS}$. As $S$ is included in the center of $G_{\tS}$, this implies that $\xi_{\tS}$ is generic.
	\end{proof}

	\subsection{Truncations}\label{sec truncation}
	
	Let $\widetilde{X}$, $\widetilde{X}_M$ be the twisted symmetric spaces associated to $(\tG,\iota)$, $(\tM,\iota)$ respectively (see Section \ref{sect twisted symmetric pairs}). More precisely, the underlying varieties are $X=G^\iota\backslash G$, $X_M=M^\iota\backslash M$ and these are equipped with the natural twisted actions of $\tG$, $\tM$ respectively. We fix from now on a special compact subgroup $K\subset G(F)$ in good position relative to $M$ and we set $K_M=K\cap M(F)$. In Section \ref{sect twisted symmetric pairs}, we have defined real affine spaces $\CA_{\tX,K}$ and $\CA_{\tX_M,K_M}$. We claim that there is a natural identification $\CA_{\tX,K}\simeq \CA_{\tX_M,K_M}$. Indeed, for any minimal $\iota$-split parabolic subspace $\tP_0\subset \tP$, we have by definition canonical isomorphisms of real affine spaces
	$$\displaystyle \CA_{\tX,K}\simeq \CA_{\tP_0,\iota}\simeq \CA_{\tP_0\cap \tM,\iota}\simeq \CA_{\tX_M,K_M}$$
	and the resulting isomorphism $\CA_{\tX,K}\simeq \CA_{\tX_M,K_M}$ does not depend on the choice of $\tP_0$. We fix a map $H_{X_M}: X_M(F)/K_M\to \CA_{X_M,K_M}$ satisfying the requirements of Proposition \ref{prop map HX} and, as in Section \ref{sect twisted symmetric pairs}, we let $H_{\tX_M}: X_M(F)/K_M\to \CA_{\tX_M,K_M}$ be the composition of $H_{X_M}$ with the natural projection $\CA_{X_M,K_M}\to \CA_{\tX_M,K_M}$.
	
	Recall also that the vector space associated to $\CA_{\tX,K}$ is the limit $\CA_{\tX}=\varprojlim_{\tP_0} \CA_{\tP_0,\iota}$ where $\tP_0$ runs over all minimal $\iota$-split parabolic subspaces $\tP_0\subset \tG$ and the transition maps are given by conjugation by elements of $G(F)$. As explained in Section \ref{sect twisted symmetric pairs}, there is a characteristic function $\phi_{\tX}: \CA_{\tX}\to \{0,1 \}$ which, upon identifying $\CA_{\tX}$ with $\CA_{\tP_0,\iota}$, is given by $\phi_{\tP_0,\iota}$ for any minimal $\iota$-split parabolic subspace $\tP_0\subset \tG$.
	
	Note that by the Iwasawa decomposition $G(F)=P(F)K$, we have a natural identifications of cosets
	\begin{equation*}
		\displaystyle H(F)\backslash G(F)/K=H_0(F)\backslash M(F)/K_M.
	\end{equation*}
	Moreover, there is a natural map $H_0(F)\backslash M(F)/K_M\to X_M(F)/K_M$ given by the composition of the surjection $H_0(F)\backslash M(F)\twoheadrightarrow M^\iota(F)\backslash M(F)$ with the natural inclusion $M^\iota(F)\backslash M(F)\subset X_M(F)=(M^\iota\backslash M)(F)$.
	
	For $Y\in \CA_{\tX,K}$, we define a characteristic function $\kappa_Y: H(F)\backslash G(F)/K\to \{0,1 \}$ by the following (commutative) diagram:
	$$\displaystyle \xymatrix{ H(F)\backslash G(F)/K \ar@{=}[r] \ar[d]^{\kappa_Y} & H_0(F)\backslash M(F)/K_M \ar[r] & X_M(F)/K_M \ar[d]^{H_{\tX_M}} \\ \{0,1 \} & \CA_{\tX,K} \ar[l]^{\phi_{\tX}(.-Y)} \ar@{=}[r] & \CA_{\tX_M,K_M}.}$$
	In other words, identifying elements in $M(F)$ with their image in $\tX_M(F)$, $\kappa_Y$ is characterized by the following property: for every $(m,u,k)\in M(F)\times N(F)\times K$ we have $\kappa_Y(muk)=\phi_{\tX}(H_{\tX_M}(m)-Y)$.

	\begin{prop}\label{prop truncation}
		
		\begin{enumerate}[(1)]
			\item For every $\iota$-split parabolic subspace $\tQ\subset \widetilde{\overline{P}}$, there is a constant $\epsilon>0$ such that, setting $\tL=\tQ\cap \iota(\tQ)$, the following holds: for every $Y\in \cA_{\tX,K}=\cA_{\tX_M,K_M}$, $x\in L[\leq \epsilon d(Y)]$ and $a\in A_{\tQ,\iota}^+$, we have
			$$\displaystyle \kappa_{Y}(ax)=\phi_{\tQ,\iota}(H_{\tL,\iota}(ax)-Y_{\tQ,\iota}).$$
			

			
			\item There exists a constant $C>0$, such that for every $Y\in \cA_{\tX,K}=\cA_{\tX_M,K_M}$, we have
			$$\displaystyle \Supp(\kappa_Y)\cap X_M(F)\subseteq A_{P,\iota}(F) X_M[\leq C N(Y)].$$
			
			\item For any fixed $x\in G(F)$, there exists $C_2>0$ such that for every $a\in A_{P,\iota}(F)$ satisfying $\kappa_Y(ax)=1$ we have $\langle \varpi, H_{\tM,\iota}(a)\rangle \leq C_2$ for every $\varpi\in \widehat{\Delta}_{\tP,\iota}$.
		\end{enumerate}
	\end{prop}

	\begin{proof}
		Let us fix a minimal $\iota$-split parabolic subspace $\tP_0\subset \widetilde{\overline{P}}$ and identify both $\CA_{\tX}$ and $\CA_{\tX,K}$ with $\CA_{\tP_0,\iota}$ via the natural isomorphisms. Thus, for every $x\in X_M(F)$ we have $\kappa_{Y}(x)=\phi_{\tP_0,\iota}(H_{\tX_M}(x)-Y)$ where we recall that
		$$\displaystyle \phi_{\tP_0,\iota}(H_{\tX_M}(x)-Y)=1\Leftrightarrow \varpi(H_{\tX_M}(x)-Y)\leq 0 \; \forall \varpi\in \widehat{\Delta}_{\tP_0,\iota}.$$
		
		Let $\tQ\subset \widetilde{\overline{P}}$ be a $\iota$-split parabolic subspace $Y\in \CA_{\tX,K}$ and $\epsilon>0$. We are going to prove that point (1) is satisfied provided $\epsilon$ is sufficiently small. Of course, we may and will assume that $\tP_0\subset \tQ$. Let $x\in X_L[\leqslant \epsilon d(Y)]$ and $a\in A_{\tQ,\iota}^+$. By Proposition \ref{prop2 HX}, there exists an absolute constant $c>0$ (depending only on $\tQ$) and a $\iota$-split parabolic subspace $\tQ\subset \tQ'\subset \tilde{\bar{P}}$ satisfying (where we have set $\tL'=\tQ'\cap \iota(\tQ')$)
		$$\displaystyle ||H_{\tX_M}(xa)-H_{\tL',\iota}(xa)||<c\epsilon d(Y),\;proj_{\tQ'}H_{\tX_M}(xa)=H_{\tL',\iota}(xa), $$
		and
		$$\displaystyle ||H_{\tL',\iota}(xa)-H_{\tL,\iota}(xa)||< c \epsilon d(Y).$$
		Furthermore, there exists an (absolute) constant $c_0>0$ such that $\varpi(Y)>c_0 d(Y)$ for every $\varpi\in \widehat{\Delta}^{\tQ'}_{\tP_0,\iota}$. This implies, by the above, that for $\varpi\in \widehat{\Delta}^{\tQ'}_{\tP_0,\iota}$, provided $\epsilon$ is small enough, we have
		$$\displaystyle \varpi(H_{\tX_M}(xa)-Y)=\varpi(H_{\tX_M}(xa)-Y- H_{\tL',\iota}(xa))\leq \varpi(H_{\tX_M}(xa)-H_{\tL',\iota}(xa))-c_0 d(Y)\leq 0.$$
		Similarly, for $\varpi\in \widehat{\Delta}^{\tQ'}_{\tQ,\iota}$ and provided $\epsilon$ is small enough, we have
		$$\displaystyle \varpi(H_{\tL,\iota}(xa)-Y_{\tQ,\iota})=\varpi(H_{\tL,\iota}(xa)-Y- H_{\tL',\iota}(xa))\leq 0.$$
		That is we have
		$$\displaystyle \phi_{\tP_0,\iota}^{\tQ'}(H_{\tX_M}(xa)-Y)=\phi_{\tQ,\iota}^{\tQ'}(H_{\tL,\iota}(xa)-Y_{\tQ,\iota})=1.$$
		Now, as is well-known, for $H\in \CA_{\tP_0,\iota}$ with $\phi_{\tP_0,\iota}^{\tQ'}(H)=1$, we have $\phi_{\tP_0,\iota}(H)=\phi_{\tQ',\iota}(H)$. Therefore,
		$$\displaystyle \phi_{\tP_0,\iota}(H_{\tX_M}(xa)-Y)=\phi_{\tQ',\iota}(H_{\tX_M}(xa)-Y)=\phi_{\tQ',\iota}(H_{\tL',\iota}(xa)-Y_{\tQ',\iota})$$
		and
		$$\displaystyle \phi_{\tQ,\iota}(H_{\tL,\iota}(xa)-Y_{\tQ,\iota})=\phi_{\tQ',\iota}(H_{\tL',\iota}(xa)-Y_{\tQ',\iota})$$
		where we have used that the projections of $H_{\tL,\iota}(xa)-Y_{\tQ,\iota}$ and $H_{\tX_M}(xa)-Y$ to $\CA_{\tQ',\iota}$ are $H_{\tL',\iota}(xa)-Y_{\tQ',\iota}$ and $H_{\tL',\iota}(xa)-Y_{\tQ',\iota}$ respectively. This proves point (1).
		
		(2) follows from Proposition \ref{prop map HX}(3) and (4). Indeed, let $Y^-\in \CA_{P_0,\iota}$ be such that $H_{X_M}$ has image in $Y^-+\CA_{P_0^M,\iota}^+$ where we have set $P_0^M=P_0\cap M$ (a minimal $\iota$-split parabolic subgroup of $M$). Then, the restriction of the projection map $\CA^{P}_{P_0,\iota}\to \CA^{\tP}_{\tP_0,\iota}$ to the image of $Y^-+\CA_{P_0^M,\iota}^+$ is proper and this implies that the intersection of the support of $\phi_{\tP_0,\iota}(.-Y_{\tP_0,\iota})$ (seen as a characteristic function on $\CA_{P_0,\iota}$ by the previous projection) with $Y^-+\CA_{P_0^M,\iota}^+$ is contained in the sum of $\CA_{P,\iota}$ with a ball centered at $0$ of radius $C' N(Y)$ for a certain $C'>0$.

		 For (3), by the Iwasawa decomposition we may assume that $x\in M(F)$. Then $H_{\tX_M}(ax)=H_{\tM,\iota}(a)+H_{\tX_M}(x)$. Then (3) just follows from the definition of $\kappa_Y$.

	\end{proof}

	\subsection{The geometric expansion}\label{sec geometric expansion}
	For $f\in C_{c}^{\infty}(\tG(F))$ and $Y\in \cA_{\tX,K}=\cA_{\tX_M,K_M}$, define
	$$\displaystyle I(f,x)=\int_{\tH(F)}f(x^{-1}hx)(\xi\otimes \chi)(h)dh,\;x\in G(F);$$
	$$\displaystyle I_Y(f)=\int_{H(F)A_{\tG}(F)\back G(F)}I(f,x)\kappa_Y(x) dx.$$
	If the $(G,H)$ is tempered (see Subsection \ref{S tempered varieties}), we can also define $I(f,x)$ and $I_Y(f)$ for $f\in \CC(\tG(F))$. It is then clear that the integral defining $I(f,x)$ is absolutely convergent. 
	
	\begin{lem}\label{lem abs conv}
		The integral defining $I_Y(f)$ is absolutely convergent.
	\end{lem}
	
	\begin{proof}
		It is enough to show that the integral 
		$$\displaystyle \int_{H(F)A_{P}(F)\back G(F)}\int_{(A_P\cap H)(F)A_{\tG}(F)\back A_{P}(F)}I(f,ax)\kappa_Y(ax) \delta_P(a)^{-1}da dx$$
		is absolutely convergent. By Proposition \ref{prop truncation}(2) and the Iwasawa decomposition, the integrand of the outer integral  over $H(F)A_{P}(F)\back G(F)$ is compactly supported, hence it is enough to show that, for each $x\in H(F)A_{P}(F)\back G(F)$, the inner integral is absolutely convergent, which (since $(A_P\cap H_0)(F)A_{P,\iota}(F)$ is of finite index in $A_P(F)$) is equivalent to show that the absolute convergence of the expression  
		\[\begin{aligned}
		\displaystyle & \int_{A_{\tG,\iota}(F)\back A_{P,\iota}(F)}I(f,ax)\kappa_Y(ax)\delta_P(a)^{-1} da= \\
		& \int_{A_{G,\iota}(F)\back A_{P,\iota}(F)}\int_{A_{\tG,\iota}(F)\back A_{G,\iota}(F)}I(f,a_1a_2x)da_1 \kappa_Y(a_2x) \delta_P(a_2)^{-1}da_2.
	\end{aligned}\]
		Since $\xi$ is a generic character and the function $\gamma\in \tG(F)\mapsto f(x^{-1}\gamma x)$ is right-invariant by a compact-open subgroup, there exists a constant $C_{f,x}>0$ (depending on $f$ and $x$) such that, for $a\in A_{P,\iota}(F)$,
		\begin{equation}\label{eq1 lem abs conv}
		\displaystyle I(f,ax)\neq 0\Rightarrow   \langle \alpha, H_{M,\iota}(a)\rangle \geq -C_{f,x}\; \text{for every} \; \alpha\in \Delta_{P,\iota}.
		\end{equation}
		Combining this with Proposition \ref{prop truncation}(3), it follows that the integrand of the outer integral over $A_{G,\iota}(F)\backslash A_{P,\iota}(F)$ is compactly supported. Indeed, the above inequality together with Proposition \ref{prop truncation}(3) imply that when $I(f,ax)\kappa_Y(ax)\neq 0$ the image of $H_{\tM,\iota}(a)$ in $\cA^{\tG}_{\tP,\iota}$ belongs to a fixed compact subset (depending only on $x$ and $f$). However, since the automorphism $\theta$ of $\cA_{P,\iota}$ induced from the twisted space $(P,\tP)$ preserves the set of simple roots $\Delta_{P,\iota}$, the restriction of the natural projection $\cA^{G}_{P,\iota}\to \cA^{\tG}_{\tP,\iota}$ to any translate of the chamber $\cA^{G,+}_{P,\iota}$ is proper which implies (using again the inequality \eqref{eq1 lem abs conv}) that the support of the function $a\in A_{P,\iota}(F)\mapsto I(f,ax)\kappa_Y(ax)$ is compact modulo $A_{G,\iota}(F)$.

		Finally, we are reduced to show the convergence of
		\begin{equation}\label{eq2 lem abs conv}
		\displaystyle \int_{A_{G,\iota}(F)/A_{\tG,\iota}(F)}I(f,ax) da
		\end{equation}
		for any given $x\in G(F)$. Up to replacing $f$ by its conjugate by $x$, we may assume that $x=1$. Then 
		\[\begin{aligned}
		\displaystyle I(f,a)=\int_{\tH(F)}f(a^{-1}ha)(\xi\otimes \chi)(h)dh=\int_{\tH(F)}f(h\theta(a)^{-1}a)(\xi\otimes \chi)(h)dh.
	    \end{aligned}\]
        where we have denoted by $\theta$ the automorphism of $A_{G,\iota}$ induced from the twisted space $\tG$. Consider the regular map
        $$\displaystyle H\times A_{G,\iota}/A_{\tG,\iota}\to G$$
        $$\displaystyle (h,a)\mapsto h \theta(a)^{-1}a.$$
        It is a morphism of linear groups with finite kernel and image $HA_G$, in particular it is finite. This implies that 
        \begin{equation}\label{eq3 lem abs conv}
        \displaystyle \sigma(h\theta(a)^{-1}a)\gg \sigma(h)+\sigma_{A_{G,\iota}/A_{\tG,\iota}}(a) \mbox{ for } (h,a)\in \tH(F)\times A_{G,\iota}(F)/A_{\tG,\iota}(F).
        \end{equation}
        In particular, if $f$ is compactly supported, the function $a\in A_{G,\iota}(F)/A_{\tG,\iota}(F)\mapsto I(f,a)$ is also compactly supported which of course implies the convergence of \eqref{eq2 lem abs conv}.
        
        Assume now that the pair $(G,H)$ is tempered and that $f$ is a Harish-Chandra Schwartz function. Then, by \eqref{eq3 lem abs conv} and since the function $\Xi^{\tG}$ is invariant by $A_G(F)$, for every $d>0$ the integral \eqref{eq2 lem abs conv} is essentially bounded by
		$$\displaystyle \int_{\tH(F)} \Xi^{\tG}(h) \sigma(h)^{-d}dh \times \int_{A_{G,\iota}(F)/A_{\tG,\iota}(F)} \sigma_{A_{G,\iota}/A_{\tG,\iota}}(a)^{-d} da.$$
		However, as $(G,H)$ is tempered the first integral converges for $d$ large enough and similarly for the second integral. This proves the lemma.
	\end{proof}
	
	For $t\in (\tH_0)_{rs}(F)$, let $S=H_{0,t}$, $T=M_{t}$, $\tS=St$ and $\tT=Tt$.
	Let $N_{\tS}$ be the centralizer of $\tS$ in $N$, which is a maximal unipotent subgroup of $G_{\tS}$, the centralizer of $\tS$ in $G$. By Lemma \ref{lem generic character}, we know that $\xi|_{N_{\tS}(F)}$ is generic. Also $B_{\tS}=SN_{\tS}$ is a Borel subgroup of $G_{\tS}$.

	Let $M(t)$ be the centralizer of the maximal split torus of $T$ in $G$ and let $\tM(t)=M(t)t$. Define $Y(g)=(Y_{\tQ}(g))_{\tQ\in \CF_{\iota}(\tM(t))}$ to be
	$$Y_{\tQ}(g)=Y_{\tQ,\iota}-H_{\tilde{\bar{Q}},\iota}(g).$$
	We then define
	$$\tilde{v}_{B_{\tS},\xi,\iota,Y}(x,n_S)=\int_{A_{\tG}(F)S(F)\back T(F)}^{\ast}\Gamma_{\bar{B}_{\tS},\iota} (H_{\tM(t),\iota}(t^M),Y(x))\xi(t^M n_S(t^M)^{-1})dt^M$$
	for $x\in G(F)$ and $n_S\in N_{\tS,reg}(F)$. We refer the reader to Lemma \ref{lem def twisted vBxpsi} for the definition of the normalized integral $\int_{A_{\tG}(F)S(F)\back T(F)}^{\ast}$. By Lemma \ref{exponential polynomial twisted weight}, there exists $C>0, r>0$ and, for any $(x,n_S)$, a unique exponential polynomial $v_{B_{\tS},\xi,\iota,\cdot}(x,n_S)\in Pol_{\leqslant r}$ whose exponents belongs to a finite set independent of $(x,n_S)$, such that 
	$$\tilde{v}_{B_{\tS},\xi,\iota,Y}(x,n_S)=v_{B_{\tS},\xi,\iota,Y}(x,n_S)$$
	for all $Y$ and $(x,n_S)$ such that $1_{<d(Y)/C}(x,n_S)\neq 0$. Here, we recall that $d(Y)$ denotes the depth $Y$ defined in Subsection \ref{orthogonal sets twist} (we also refer the reader to this subsection for the norm $N(Y)$ that appears in the lemma below). Also for $c>0$, $1_{<c}(\cdot, \cdot)$ stands for the characteristic function of the subset
	$$\displaystyle \{(x,n_S)\in A_{\tG}(F)\backslash G(F)\times N_{\tS}(F)|\; \sigma_{A_{\tG}\backslash G}(x),\sigma_{N_{\tS,reg}}(n_S)<c\}.$$	
	From Lemma \ref{lem def twisted vBxpsi}, we have:
	
	\begin{lem}\label{lemma weight estimate 1}
		There exists $d>0$ such that
		$$\tilde{v}_{B_{\tS},\xi,\iota,Y}(x,n_S)\ll N(Y)^d \cdot \sigma_{G/A_{\tG}}(x)^d \sigma_{N_{\tS,reg}}(n_S)^d$$
		and
		$$v_{B_{\tS},\xi,\iota,Y}(x,n_S)\ll N(Y)^d \cdot \sigma_{G/A_{\tG}}(x)^d \sigma_{N_{\tS,reg}}(n_S)^d$$
		for all $Y\in \cA_{\tX,K}=\cA_{\tX_M,K_M}$, $x\in G(F)$ and $n_S\in N_{\tS,reg}(F)$.
	\end{lem}

	Recall that the set of regular semisimple conjugacy classes $\Gamma(\tH_0)$ was equipped with a measure in Subsection \ref{Sect measures}. We define
	$$\tilde{J}_{Y}(f)=\int_{\Gamma(\tH_0)}D^{\tH}(t)\tilde{\chi}(t)\int_{B_{\tS}(F)\back G(F)} \int_{N_{\tS}(F)}f(x^{-1}tn_S x)\tilde{v}_{B_{\tS},\xi,\iota,Y}(x,n_S) dn_Sdxdt$$
	and
	$$J_{Y}(f)=\int_{\Gamma(\tH_0)}D^{\tH}(t)\tilde{\chi}(t)\int_{B_{\tS}(F)\back G(F)} \int_{N_{\tS}(F)}f(x^{-1}tn_S x)v_{B_{\tS},\xi,\iota,Y}(x,n_S) dn_Sdxdt.$$
By Lemma 2.9.3 of \cite{BeuGalP}, Proposition \ref{estimate 1} and the lemma above, we know that the above two integrals are absolutely convergent for all $f\in C_{c}^{\infty}(\tG(F))$ (resp. for all $f\in \CC(\tilde{G}(F))$ if the model $(M,H_0)$ is tempered). Note that since we have assumed that $(\tM,\tH_0)$ is coregular, the function $t\in \tH_{0}(F)\cap \tM_{rs}(F)\mapsto \frac{D^{\tH_0} (t)}{D^{\tM} (t)^{1/2}}=\frac{D^{\tH} (t)}{D^{\tG} (t)^{1/2}}$ is locally bounded on $\tH_0(F)$. If we further assume that $(\tM,\tH_0)$ is tempered, that function is globally bounded. The geometric expansion is the following theorem (see Lemma \ref{lem coregular and tempered}).

	\begin{thm}\label{thm geometric side}
		Let $0<\epsilon<1$ and fix $f\in C_{c}^{\infty}(\tG(F))$. For $k>0$, we have
		$$|I_{Y}(f)-J_{Y}(f)|\ll N(Y)^{-k}$$
		for every $Y\in \cA_{\tX,K}=\cA_{\tX_M,K_M}$ with $d(Y)>\epsilon N(Y)$. Moreover, if the model $(M,H_0)$ is tempered, then the estimates hold for all $f\in \CC(\tilde{G}(F))$.
	\end{thm}
	
\begin{rmk}	
	It is clear that in order to prove the above theorem, we only need to prove 
	$$|I_{Y}(f)-\tilde{J}_{Y}(f)|\ll N(Y)^{-k}$$
	and
	$$|J_{Y}(f)-\tilde{J}_{Y}(f)|\ll N(Y)^{-k}.$$
\end{rmk}

To end this subsection, we will also state an analogue of the trace formula in Theorem \ref{thm geometric side} for functions with a fixed central character. We fix a character $\omega: A_{\tG}(F)\to \mathbb{C}^\times$ whose restriction to $A_{\tG}^{H_0}(F)$ coincides with $\chi$ where $A_{\tG}^{H_0}$ stands for the connected component of $A_{\tG}\cap H_0$. For $f\in C_{c}^{\infty}(\tG(F)/A_{\tG}(F),\omega^{-1})$ and $Y\in \cA_{\tX,K}=\cA_{\tX_M,K_M}$, define
	$$I(f,x)=\int_{N(F)}\int_{\tH_0(F)/A_{\tG}^{H_0}(F)}f(x^{-1}hnx)\xi(n)\chi(h)dhdn,\;x\in G(F);$$
	$$I_Y(f)=\int_{H(F)A_{\tG}(F)\back G(F)}I(f,x)\kappa_Y(x) dx.$$
	If $(M,H_0)$ is tempered, we can also define $I(f,x)$ and $I_Y(f)$ for $f\in \CC(\tG(F)/A_{\tG}(F),\omega^{-1})$.

For the geometric side, we define
$$J_{Y}(f)=\int_{\Gamma(\widetilde{\overline{H_0}})}D^{\tH}(t)\chi(t)\int_{B_{\tS}(F)\back G(F)} \int_{N_{\tS}(F)}f(x^{-1}tn_S x)v_{B_{\tS},\xi,\iota,Y}(x,n_S) dn_Sdxdt.$$
where $\overline{H_0}=H_0/A_{\tG}^{H_0}$ and $\Gamma(\widetilde{\overline{H_0}})$ is the set of regular semisimple conjugacy classes of the twisted space $\widetilde{\overline{H_0}}(F)$ equipped with a measure defined in Subsection \ref{Sect measures}. The following theorem is an analogue of Theorem \ref{thm geometric side} and it can be proved by the exactly same argument.

\begin{thm}\label{thm geometric side center}
		Let $0<\epsilon<1$ and fix $f\in C_{c}^{\infty}(\tG(F)/A_{\tG}(F),\omega^{-1})$. For $k>0$, we have
		$$|I_{Y}(f)-J_{Y}(f)|\ll N(Y)^{-k}$$
		for every $Y\in \cA_{\tX,K}=\cA_{\tX_M,K_M}$ with $d(Y)>\epsilon N(Y)$. Moreover, if the model $(M,H_0)$ is tempered, then the estimates hold for all $f\in \CC(\tG(F)/A_{\tG}(F),\omega^{-1})$.
	\end{thm}

The goal of the rest of this section is to prove Theorem \ref{thm geometric side}.
	
	\subsection{Some reduction}
	Recall that we have defined
	$$I(f,x)=\int_{N(F)}\int_{\tH_0(F)}f(x^{-1}hnx)\xi(n)\tilde{\chi}(h)dhdn,\;x\in G(F);$$
	$$I_Y(f)=\int_{H(F)A_{\tG}(F)\back G(F)}I(f,x)\kappa_Y(x) dx.$$
	By the Weyl's integration formula (applied to $\tH_0$), we have
	$$I_Y(f)=\int_{H(F)A_{\tG}(F)\back G(F)}\kappa_Y(x) \int_{N(F)}\int_{\tH_0(F)}f(x^{-1}hnx)\xi(n)\tilde{\chi}(h)dhdn dx$$
	$$=\int_{H(F)A_{\tG}(F)\back G(F)}\kappa_Y(x) \int_{N(F)}\int_{\Gamma(\tH_0)}  D^{\tH_0}(t)\int_{S(F)\back H_0(F)} f(x^{-1}h^{-1}thnx)\xi(n)\tilde{\chi}(t)dh dtdn dx$$
	$$=\int_{\Gamma(\tH_0)}D^{\tH_0}(t)\tilde{\chi}(t)\int_{A_{\tG}(F)S(F)N(F)\back G(F)}\kappa_Y(x) \int_{N(F)} f(x^{-1}tnx)\xi(n) dn dx dt$$
	$$=\int_{\Gamma(\tH_0)}D^{\tH_0}(t)\tilde{\chi}(t)\int_{T(F)N(F)\back G(F)} \int_{A_{\tG}(F)S(F)\back T(F)}\kappa_Y(t^Mx)$$
	$$\int_{N(F)} f(x^{-1}tnx)\xi(t^M n(t^M)^{-1}) dn dt^Mdx dt.$$
	Here as in the previous subsection, for $t\in (\tH_0)_{rs}(F)$, we let $S=H_{0,t}$, $T=M_{t}$, $\tS=St$ and $\tT=Tt$. It is easy to see that the isomorphism
	
	$$N_{\tS}(F)\times^{N_{\tS}(F)} N(F)\rightarrow N(F):\; (n_S,n)\mapsto (Ad_t)^{-1}(n^{-1})n_Sn$$
	has Jacobian $D^{\tH_0}(t)^{-1}D^{\tH}(t)$, hence the above expression is equal to

	$$\int_{\Gamma(\tH_0)}D^{\tH}(t)\tilde{\chi}(t)\int_{T(F)N(F)\back G(F)}\int_{A_{\tG}(F)S(F)\back T(F)} \kappa_Y(t^Mx)\int_{N_{\tS}(F)\times^{N_{\tS}(F)} N(F)}$$ 
	$$f(x^{-1}t(Ad_t)^{-1}(n^{-1})n_Sn x)\xi(t^M (Ad_t)^{-1}(n^{-1})n_Sn(t^M)^{-1}) d(n_S,n) dt^Mdx dt$$

	$$=\int_{\Gamma(\tH_0)}D^{\tH}(t)\tilde{\chi}(t)\int_{T(F)N(F)\back G(F)}\int_{A_{\tG}(F)S(F)\back T(F)} \kappa_Y(t^Mx)\int_{N_{\tS}(F)\times^{N_{\tS}(F)} N(F)}$$ 
	$$f(x^{-1}n^{-1}tn_Sn x)\xi(t^M n_S(t^M)^{-1}) d(n_S,n) dt^Mdx dt$$
	$$=\int_{\Gamma(\tH_0)}D^{\tH}(t)\tilde{\chi}(t)\int_{T(F)N_{\tS}(F)\back G(F)}$$
	$$\int_{A_{\tG}(F)S(F)\back T(F)}\kappa_Y(t^Mx)\int_{N_{\tS}(F)}f(x^{-1}tn_S x)\xi(t^M n_S(t^M)^{-1}) dn_S dt^Mdx dt.$$

	\begin{defn}
		With the notation above, we define
		$$\kappa_{Y,\tS,\xi}(x,n_S)=\int_{A_{\tG}(F)S(F)\back T(F)}^{\ast}\kappa_Y(t^Mx) \xi(t^M n_S(t^M)^{-1})dt^M$$
		$$:=\int_{A_{\tG}(F)S(F)\back T(F)}\int_{T_c}\kappa_Y(t^Mtx) \xi(t^Mt n_S(t^Mt)^{-1})dtdt^M$$
		and
		$$I_{Y,\tS}(f)=\int_{\tS(F)}D^{\tH}(t)\tilde{\chi}(t)\int_{B_{\tS}(F)\back G(F)} \int_{N_{\tS}(F)}f(x^{-1}tn_S x)\kappa_{Y,\tS,\xi}(x,n_S) dn_Sdxdt.$$
		We also define
		$$\tilde{J}_{Y,\tS}(f)=\int_{\tS(F)}D^{\tH}(t)\tilde{\chi}(t)\int_{B_{\tS}(F)\back G(F)} \int_{N_{\tS}(F)}f(x^{-1}tn_S x)\tilde{v}_{B_{\tS},\xi,\iota,Y}(x,n_S) dn_Sdxdt$$
		and 
		$$J_{Y,\tS}(f)=\int_{\tS(F)}D^{\tH}(t)\tilde{\chi}(t)\int_{B_{\tS}(F)\back G(F)} \int_{N_{\tS}(F)}f(x^{-1}tn_S x)v_{B_{\tS},\xi,\iota,Y}(x,n_S) dn_Sdxdt.$$
	\end{defn}
	
	\begin{lem}\label{lemma weight estimate 2}
		There exists $d>0$ such that
		$$\kappa_{Y,\tS,\xi}(x,n_S)\ll N(Y)^d \cdot \sigma_{G/A_{\tG}}(x)^d \sigma_{N_{\tS,reg}}(n_S)^d$$
		for all $Y\in \cA_{\tX,K}=\cA_{\tX_M,K_M}$, $x\in G(F)$ and $n_S\in N_{\tS,reg}(F)$.
	\end{lem}
	
	\begin{proof}
		This follows from Proposition \ref{prop truncation}.
	\end{proof}
	
	We fix $t,S,T,\tS,\tT$ as above. By the integration formula \eqref{eq1 measures}, in order to prove Theorem \ref{thm geometric side}, we only need to establish the following theorem.
	
	\begin{thm}\label{thm change truncation}
		Let $0<\epsilon<1$ and fix $f\in C_{c}^{\infty}(\tG(F))$. For $k>0$, we have
		$$|I_{Y,\tS}(f)-\tilde{J}_{Y,\tS}(f)|\ll N(Y)^{-k}$$
		and
		$$|J_{Y,\tS}(f)-\tilde{J}_{Y,\tS}(f)|\ll N(Y)^{-k}$$
		for every $Y\in \cA_{\tX,K}=\cA_{\tX_M,K_M}$ with $d(Y)>\epsilon N(Y)$. Moreover, if the model $(M,H_0)$ is tempered, then the estimates hold for all $f\in \CC(\tilde{G}(F))$.
	\end{thm}

	\subsection{Comparison of the weights}
	Recall that for $c>0$ we have defined the function $1_{<c}(\cdot, \cdot)$ to be the characteristic function of
	$$\{(x,n_S)\in A_{\tG}(F)\backslash G(F)\times N_{\tS}(F)|\; \sigma_{A_{\tG}\backslash G}(x)<c ,\sigma_{N_{\tS,reg}}(n_S)<c\}.$$
	The goal of this subsection is to prove the following lemma.
	
	\begin{lem}\label{lem comparison weights}
		There exists $C>0$ such that
		$$\kappa_{Y,\tS,\xi}(x,n_S)=\tilde{v}_{B_{\tT},\xi,\iota,Y}(x,n_S)$$
		for all $Y\in \cA_{\tX,K}=\cA_{\tX_M,K_M}$ and $(x,n_S)$ such that $1_{<N(Y)^\epsilon}(x,n_S)\neq 0$, $d(Y)>\epsilon N(Y)$ and $d(Y)>C$.
	\end{lem}

	Let $M(\tS)=M(t)$ be the centralizer of the maximal split torus of $T$ in $G$ and let $\tM(\tS)=\tM(t)=M(\tS)\tS$. For all $\tQ\in \CF_{\iota}^{\tM}(\tM(\tS))$ (resp. $\tQ\in \CF_{\bar{B}_{\tS},\iota}(\tM(\tS))$), let $\tL$ be the Levi factor containing $\tM(\tS)$ and let $\bar{Q}$ be the opposite parabolic subgroup of $Q$ with respect to $L$. We first need a lemma. 

\begin{lem}\label{lem U bar Q}
		For all $\tQ\in \CF_{\bar{B}_{\tS},\iota}(\tM(\tS))$, we have $U_{\bar{Q}}\subset P$.
	\end{lem}
	
	The above lemma is a direct consequence of the next lemma.

	\begin{lem}\label{lem U bar Q 1}
		Let $A_{T,\iota}$ be the maximal split $\iota$-split torus in $T(F)$. Every root in $\Delta(A_{T,\iota},N)$ can be written as a linear combination of roots in $\Delta(A_{T,\iota},N_{\tS})$ with nonnegative coefficients.
	\end{lem}
	
	\begin{proof}
		For $a\in A_{T,\iota}$ and a character $\xi'$ of $N(F)$, we say $a$ shrinks $\xi'$ if 
		$$\lim_{i\rightarrow \infty} \xi'(a^{-i} na^i)=1$$
		for all $n\in N(F)$. To prove the lemma, we only need to show that for $a\in A_{T,\iota}$, if $a$ shrinks $\xi$, then $a$ shrinks all the characters of $N(F)$. The characters of $N(F)$ can be naturally identified with the vector space $\bar{\mathfrak{n}}(F)/[\bar{\mathfrak{n}}(F),\bar{\mathfrak{n}}(F)]$ and we endow it with the natural topology coming from the vector space. We only need to show that $a$ shrinks an open subset of the characters. There exists a $\iota$-split parabolic subgroup $P_{a,\iota}$ of $M=G_0$ such that $\lim_{i\rightarrow \infty}a^{-i}pa^i$ exists for all $p\in P_{a,\iota}(F)$. Then we know that $a$ shrinks all the characters $\xi'$ of $N(F)$ of the form
		$$\xi'(n)=\xi(m^{-1}nm),\;m\in P_{a,\iota}(F)H_0(F).$$
		Since $\xi$ is a generic character and $P_{a,\iota}H_0$ is Zariski open in $G_0=M$, we know that $a$ shrinks an open subset of the characters. This proves the lemma.
	\end{proof}

Let $C'>0$ be a constant large enough (with respect to $\frac{1}{\epsilon}$ and the inverse of the constant in Prop \ref{prop truncation}(1)). For every positve $(\tG,M(\tS),\iota)$-orthogonal set $Y_0=(Y_{0,\tQ,\iota})_{\tQ\in \CF_{\iota}(\tM(\tS))}$, by using the parabolic subgroup $P=MN$ it induces a $(\tM,M(\tS),\iota)$-orthogonal set $Y_{0,M}=(Y_{0,M,\tQ,\iota})_{\tQ\in \CF_{\iota}^{\tM}(\tM(\tS))}$ where $Y_{0,M,\tQ,\iota}=Y_{0,\tQ\cdot \bar{N},\iota}$. We fix an auxiliary positive $(\tG,M(\tS),\iota)$-orthogonal set $Y_0=(Y_{0,\tQ,\iota})_{\tQ\in \CF_{\iota}(\tM(\tS))}$ satisfies the following conditions.
\begin{itemize}
\item We have
\begin{equation}\label{Y0 condition 1}
\frac{1}{C'^2}N(Y)<N(Y_{0,M})<\frac{1}{C'}N(Y),\;\frac{1}{C'^2}N(Y)<N(Y_0)<\frac{1}{C'}N(Y),\;d(Y_{0,M})>\epsilon N(Y_{0,M}).
\end{equation}
\item For any $\tQ\in \CF_{\iota}^{\tM}(\tM(\tS))$ and  $\alpha\in \Delta(A_{T,\iota},N)$, we have
\begin{equation}\label{Y0 condition 2}
\alpha(Y_{0,M,\tQ,\iota})> \frac{1}{C'^3}N(Y).
\end{equation}
\end{itemize}
It is clear such $Y_0$ exists (for example, we can take it to be of the form $(\frac{Y_{\tQ,\iota}}{C''}+X)_{\tQ}$ for some $C''>0$ and $X\in \CA_{\tM,\iota}$). For $\tQ\in \CF_{\iota}^{\tM}(\tM(\tS))$, define
	$$\kappa_{Y,\tS,\xi}^{Y_0,\tQ}(x,n_S)=\int_{A_{\tG}(F)S(F)\back T(F)}^{\ast}\kappa_Y(t^Mx) \xi(t^M n_S(t^M)^{-1}) $$ $$\Gamma_{\tM(\tS),\iota}^{\tQ}(H_{\tM(\tS),\iota}(t^M),Y_{0,M})\tau_{\tQ,\iota}^{\tM}(H_{\tM(\tS),\iota}(t^M)-Y_{0,M,\tQ,\iota}) dt^M,$$
	$$\tilde{v}_{B_{\tS},\xi,\iota,Y}^{Y_0,\tQ}(x,n_S)=\int_{A_{\tG}(F)S(F)\back T(F)}^{\ast}\Gamma_{\bar{B}_{\tS},\iota} (H_{\tM(\tS),\iota}(t^M),Y(x))\xi(t^M n_S(t^M)^{-1})$$
	$$\Gamma_{\tM(\tS),\iota}^{\tQ}(H_{\tM(\tS),\iota}(t^M),Y_{0,M})\tau_{\tQ,\iota}^{\tM}(H_{\tM(\tS),\iota}(t^M)-Y_{0,M,\tQ,\iota}) dt^M.$$
	By \eqref{equation iota version 1}, in order to prove Lemma \ref{lem comparison weights}, it is enough to compare $\kappa_{Y,\tS,\xi}^{Y_0,\tQ}(x,n_S)$ and $\tilde{v}_{B_{\tS},\xi,\iota,Y}^{Y_0,\tQ}(x,n_S)$ for all $\tQ\in \CF_{\iota}^{\tM}(\tM(\tS))$. Fix $\tQ\in \CF_{\iota}^{\tM}(\tM(\tS))$. We are reduced to prove the following lemma.
	
	\begin{lem}
		There exists $C>0$ such that
		$$\kappa_{Y,\tS,\xi}^{Y_0,\tQ}(x,n_S)=\tilde{v}_{B_{\tS},\xi,\iota,Y}^{Y_0,\tQ}(x,n_S)$$
		for all $Y\in \cA_{\tX,K}=\cA_{\tX_M,K_M}$ and $(x,n_S)$ such that $1_{<N(Y)^\epsilon}(x,n_S)\neq 0$, $d(Y)>\epsilon N(Y)$ and $d(Y)>C$.
	\end{lem}
	
	\begin{proof}
		We can rewrite the weighted functions as ($T_{c}$ is the maximal compact subgroup of $T(F)$)
		$$\kappa_{Y,\tS,\xi}^{Y_0,\tQ}(x,n_S)=\int_{A_{\tG}(F)S(F)\back T(F)}\int_{T_{c}} \kappa_Y(tt^Mx) \xi(tt^M n_S(t^M)^{-1}t^{-1}) $$ 
		$$\Gamma_{\tM(\tS),\iota}^{\tQ}(H_{\tM(\tS),\iota}(t^M),Y_{0,M})\tau_{\tQ,\iota}^{\tM}(H_{\tM(\tS),\iota}(t^M)-Y_{0,M,\tQ,\iota}) dtdt^M,$$
		$$\tilde{v}_{B_{\tS},\xi,\iota,Y}^{Y_0,\tQ}(x,n_S)=\int_{A_{\tG}(F)S(F)\back T(F)} \int_{T_{c}}\Gamma_{\bar{B}_{\tS},\iota} (H_{\tM(\tS),\iota}(t^M),Y(x))\xi(tt^M n_S(t^M)^{-1}t^{-1})$$
		$$\Gamma_{\tM(\tS),\iota}^{\tQ}(H_{\tM(\tS),\iota}(t^M),Y_{0,M})\tau_{\tQ,\iota}^{\tM}(H_{\tM(\tS),\iota}(t^M)-Y_{0,M,\tQ,\iota}) dt^M.$$
		Hence it is enough to show that the two functions
		$$t^M\mapsto \Gamma_{\tM(\tS),\iota}^{\tQ}(H_{\tM(\tS),\iota}(t^M),Y_{0,M})\tau_{\tQ,\iota}^{\tM}(H_{\tM(\tS),\iota}(t^M)-Y_{0,M,\tQ,\iota})\cdot \int_{T_{c}} \kappa_Y(tt^Mx) \xi(tt^M n_T(t^M)^{-1}t^{-1}) dt$$
		$$t^M\mapsto \Gamma_{\bar{B}_{\tS},\iota} (H_{\tM(\tS),\iota}(t^M),Y(x))\Gamma_{\tM(\tS),\iota}^{\tQ}(H_{\tM(\tT),\iota}(t^M),Y_{0,M})\tau_{\tQ,\iota}^{\tM}(H_{\tM(\tS),\iota}(t^M)-Y_{0,M,\tQ,\iota})$$
		$$\cdot \int_{T_{c}}  \xi(tt^M n_S(t^M)^{-1}t^{-1}) dt$$
		on $A_{\tG}(F)S(F)\back T(F)$ are equal to each other. We denote these two functions by $F_{1,x,n_S}$ and $F_{2,x,n_S}$.

		Let $x=mnk$ be the Iwasawa decomposition of $x$ with respect to $P=MN$. Since $\kappa_Y$ is left $N(F)$-invariant and right $K$-invariant, the function $F_{1,x,n_S}$ only depends on $m$. By Lemma \ref{lem U bar Q}, we know that $N\subset Q$ for all $\tQ\in \CP_{\bar{B}_{\tS},\iota}(\tM(\tS))$. This implies that the function $F_{2,x,n_S}$ also only depends on $m$. So we may assume that $x=m\in M(F)$. Let $x=luk$ be the Iwasawa decomposition of $x$ with respect to $\bar{Q}=LU_{\bar{Q}}$. We first prove the following statement. 
		
		\begin{itemize}
			\item[(1)] With the assumption on $Y$ and $(x,n_S)$, once we choose $C$ large enough the above two functions $F_{1,x,n_S}$ and $F_{2,x,n_S}$ only depends on $l$.
		\end{itemize}
		
		For the function $F_{1,x,n_S}$, since $\kappa_Y$ is right $K$-invariant, we know that it only depends on $lu$. It is enough to show that for $t^M\in A_{\tG}(F)S(F)\back T(F)$ with  $$\Gamma_{\tM(\tS),\iota}^{\tQ}(H_{\tM(\tS),\iota}(t^M),Y_{0,M})\tau_{\tQ,\iota}^{\tM}(H_{\tM(\tS),\iota}(t^M)-Y_{0,M,\tQ,\iota})\neq 0,$$ 
		we have
		$$\kappa_Y(t^Mlu)=\kappa_Y(t^M l).$$
		For $t^M\in A_{\tG}(F)S(F)\back T(F)$, we can always choose a representative of $t^M$ in $T(F)$, denoted by $t_2$, such that the projection of $H_{\tM(\tS)}(t_2)\in \CA_{\tM(\tS)}$ to $\CA_{\tM(\tS)}^{\iota}$ belongs to a compact subset that only depends on $\tT$. It is enough to show that 
		$$\kappa_Y(t_2lu)=\kappa_Y(t_2 l).$$
		The argument is the same as the proof of (5.3.14) of \cite{BeuGalP} and we will skip it here. This proves (1) for the function $F_{1,x,n_S}$.

		For the function $F_{2,x,n_S}$, let $\tilde{Q}'=\tilde{Q}\bar{N}\in \CF_{\bar{B}_{\tS},\iota}(\tM(\tS))$ and $\tilde{\bar{Q}}'=\tilde{\bar{Q}}N\in \CF_{B_{\tS},\iota}(\tM(\tS))$. By Lemma \ref{lem U bar Q 1}, we know that there exists $c_1>0$ such that if
		$$\int_{T_{c}}  \xi(tt^M n_S(t^M)^{-1}t^{-1}) dt	\neq 0,$$
		then
\begin{equation}\label{change truncation eqn 2}	
\alpha(H_{\tM(\tS),\iota}(t^M))\leq c_1 N(Y)^{\epsilon},\;\forall \alpha\in \Delta(A_{T,\iota},N).
\end{equation}

Combining with the second condition \eqref{Y0 condition 2} of $Y_0$, we know that once we choose $C,C'>0$ large enough, we have
$$\Gamma_{\tM(\tS),\iota}^{\tQ}(H_{\tM(\tS),\iota}(t^M),Y_{0,M})\tau_{\tQ,\iota}^{\tM}(H_{\tM(\tS),\iota}(t^M)-Y_{0,M,\tQ,\iota})\int_{T_{c}}  \xi(tt^M n_S(t^M)^{-1}t^{-1}) dt	\neq 0$$
$$\Rightarrow \Gamma_{\tM(\tS),\iota}^{\tQ'}(H_{\tM(\tS),\iota}(t^M),Y_0)\tau_{\tQ',\iota}(H_{\tM(\tS),\iota}(t^M)-Y_{0,\tQ',\iota})	\neq 0.$$
Combining with \eqref{equation iota version 2} (applied to the case when $\CX=Y_0$ and $\CY=Y(x)-Y_0$, note that by the first condition \eqref{Y0 condition 1} of $Y_0$ we know that $\CY$ is positive with $d(\CY)>\frac{d(Y)}{2}$) we get the following statement.
		\begin{itemize}
			\item[(2)] If 
			$$\Gamma_{\tM(\tS),\iota}^{\tQ}(H_{\tM(\tS),\iota}(t^M),Y_0)\tau_{\tQ,\iota}^{\tM}(H_{\tM(\tS),\iota}(t^M)-Y_{0,\tQ,\iota})\int_{T_{c}}  \xi(tt^M n_S(t^M)^{-1}t^{-1}) dt	\neq 0,$$ 
			then 
			$$\Gamma_{\bar{B}_{\tS},\iota} (H_{\tM(\tS),\iota}(t^M),Y(x))=\phi_{\tQ',\iota}(H_{\tL,\iota}(t^M)-Y_{\tQ',\iota}+H_{\tilde{\bar{Q}}',\iota}(x)).$$
				\end{itemize}
		In particular, we have proved statement (1) for the function $F_{2,x,n_S}$.

		From now on assume that $x=l\in L(F)$. We only need to prove the following two statements.
		
		\begin{itemize}
			\item[(3)] With the assumption on $Y$ and $(x,n_S)$, once we choose $C>0$ large enough, for $t^M\in A_{\tG}(F)S(F)\back T(F)$ with $F_{1,x,n_S}(t^M)\neq 0$, the following holds.
			\begin{itemize}
				\item $\phi_{\tQ',\iota}(H_{\tL,\iota}(t^M x)-Y_{\tQ',\iota})=1$.
				\item $F_{1,x,n_S}(t^M)=\Gamma_{\tM(\tS),\iota}^{\tQ}(H_{\tM(\tS),\iota}(t^M),Y_0)\tau_{\tQ,\iota}^{\tM}(H_{\tM(\tS),\iota}(t^M)-Y_{0,\tQ,\iota})\cdot \int_{T_{c}} \xi(tt^M n_T(t^M)^{-1}t^{-1}) dt$.
			\end{itemize}
		\end{itemize}

\begin{itemize}
			\item[(4)] With the assumption on $Y$ and $(x,n_S)$, once we choose $C>0$ large enough, for $t^M\in A_{\tG}(F)S(F)\back T(F)$ with $F_{2,x,n_S}(t^M)\neq 0$, the following holds.
			\begin{itemize}
				\item $\phi_{\tQ',\iota}(H_{\tL,\iota}(t^M x)-Y_{\tQ',\iota})=1$.
				\item $F_{2,x,n_S}(t^M)=\Gamma_{\tM(\tS),\iota}^{\tQ}(H_{\tM(\tS),\iota}(t^M),Y_0)\tau_{\tQ,\iota}^{\tM}(H_{\tM(\tS),\iota}(t^M)-Y_{0,\tQ,\iota})\cdot \int_{T_{c}} \xi(tt^M n_T(t^M)^{-1}t^{-1}) dt$.
			\end{itemize}
		\end{itemize}

		Statement (4) follows from (2). It remains to prove (3). As in the proof of (2), by Lemma \ref{lem U bar Q 1} and \eqref{Y0 condition 2}, we know that once we choose $C>0$ large enough, then
		\begin{equation}\label{change truncation eqn 1}
\tau_{\tQ,\iota}^{\tM}(H_{\tM(\tS),\iota}(t^M)-Y_{0,\tQ,\iota})\cdot \int_{T_{c}} \kappa_Y(tt^Mx) \xi(tt^M n_S(t^M)^{-1}t^{-1}) dt	\neq 0
\end{equation}
		$$\Rightarrow \tau_{\tQ',\iota}(H_{\tM(\tS),\iota}(t^M)-Y_{0,\tQ',\iota})\neq 0.$$

We can choose a representative of $t^M$ in $T(F)$ of the form $t't_1 a$ where $t'$ belongs to a compact set, $a\in A_{\tL,\iota}(F)$ and $t_1\in T_{\iota}(F)$ with $H_{\tL,\iota}(t_1)=0$. Here $T_{\iota}$ is the maximal $\iota$-split torus of $T$. 
		
		Since  $\Gamma_{\tM(\tS),\iota}^{\tQ}(H_{\tM(\tS),\iota}(t^M),Y_{0,M})\neq 0$, by Proposition \ref{prop iota version} and \eqref{Y0 condition 1}, we know that once we choose $C>0$ large enough, we have 
		\begin{equation}
			\sigma_{G/A_{\tG}}(t't_1)<\frac{d(Y)}{\sqrt{C'}}.
		\end{equation}
		Combining with \eqref{change truncation eqn 1} and \eqref{Y0 condition 2}, we know that once we choose $C>0$ large enough, we can write $a$ as $a=a_1a_2$ such that $a_1\in A_{\tilde{Q}',\iota}^{+}$ and 	
		\begin{equation}\label{change truncation estimate 1}
			\sigma_{G/A_{\tG}}(t't_1a_2x)<\frac{d(Y)}{\sqrt{C'}}.
		\end{equation}  
		
		Combining \eqref{change truncation estimate 1} with Proposition \ref{prop truncation}, we know that 
		$$\kappa_Y(tt^Mx)=\kappa_{Y}(tt't_1a_2xa_1)$$
		is equal to 1 if and only if $\phi_{\tQ',\iota}(H_{\tL,\iota}(t^M x)-Y_{\tQ',\iota})=1$. This proves (3) and finishes the proof of the proposition.
		
	\end{proof}

	\subsection{The proof of Theorem \ref{thm change truncation}}
	
	We have
	$$|I_{Y,\tS}(f)-\tilde{J}_{Y,\tS}(f)|\leq \int_{\tS(F)}D^{\tH}(t)\int_{B_{\tS}(F)\back G(F)} \int_{N_{\tS}(F)}$$
	$$|f(x^{-1}tn_S x)|\cdot |\kappa_{Y,\tS,\xi}(x,n_S)-\tilde{v}_{B_{\tS},\xi,\iota,Y}(x,n_S)| dn_Sdxdt.$$
	
	Let $N>0$. By Lemma \ref{lemma weight estimate 1}, \ref{lemma weight estimate 2} and \ref{lem comparison weights}, there exists $d_0>0$ such that
	$$|\kappa_{Y,\tS,\xi}(x,n_S)-\tilde{v}_{B_{\tS},\xi,\iota,Y}(x,n_S)|\ll N(Y)^{-N} (\sigma_{G/A_{\tG}}(x)+\sigma_{N_{\tS,reg}}(n_S))^{d_0}$$
	for all $Y\in \cA_{\tX,K}=\cA_{\tX_M,K_M},x\in G(F),n_S\in N_{\tS,reg}(F)$ with $d(Y)>\epsilon N(Y)$.  Since the left hand side is invariant under the transform $(x,n_s)\mapsto (bx,bn_sb^{-1})$ for all $b\in B_{\tS}(F)$, by Lemma \ref{lem norm descent property} it follows that
	$$|\kappa_{Y,\tS,\xi}(x,n_S)-\tilde{v}_{B_{\tS},\xi,\iota,Y}(x,n_S)|\ll N(Y)^{-N} (\sigma_{\tG_{reg}}(x^{-1}tn_Sx)+\sigma_{\tS'}(t))^{d_0}$$
	for all $Y\in \cA_{\tX,K}=\cA_{\tX_M,K_M},x\in G(F),n_S\in N_{\tS,reg}(F),t\in\tS'(F)$ with  $d(Y)>\epsilon N(Y)$.  Combining this with Proposition \ref{estimate 1}, we deduce that, for any $d>0$, the integral
	$$D^{\tH}(t)\int_{B_{\tS}(F)\back G(F)} \int_{N_{\tS}(F)}|f(x^{-1}tn_S x)|\cdot |\kappa_{Y,\tS,\xi}(x,n_S)-\tilde{v}_{B_{\tS},\xi,\iota,Y}(x,n_S)| dn_Sdx$$
	is essentially bounded by 
	$$N(Y)^{-N}\cdot \sigma_{\tS'}(t)^{d_0}\cdot \sigma_{\tG}(t)^{-d}$$
	for $t\in \tS(F)$ and $Y\in \cA_{\tX,K}=\cA_{\tX_M,K_M}$ with  $d(Y)>\epsilon N(Y)$. Moreover, we can choose $d>0$ large enough such that the expression 
	$$\int_{\tS(F)} \sigma_{\tG}(t)^{-d} \sigma_{\tS'}(t)^{d_0} dt$$
	is convergent (see Lemma 2.9.3 of \cite{BeuGalP}, note that $\sigma_{\tS'}(t)\sim \sigma_{\tG}(t)+\log(2+D^{\tG}(t)^{-1})$). This proves the first inequality of Theorem \ref{thm change truncation}. The second inequality follows from the same argument except that we replace Lemma \ref{lem comparison weights} by Lemma \ref{exponential polynomial twisted weight}. This finishes the proof of Theorem \ref{thm geometric side}
	and \ref{thm change truncation}.

	\section{Application of the geometric expansion}
	In this section, we will discuss the application of the geometric expansion in Theorem \ref{thm change truncation}. The first application is a simple local twisted trace formula for strongly cuspidal functions in the coregular case. The second application is a multiplicity formula for Whittaker induction of coregular symmetric pairs. We use the same notation as in the previous section.

	\subsection{A simple local trace formula}\label{sec simple LTF}

	With the same notation as in Theorem \ref{thm geometric side center}, for $f\in C_{c}^{\infty}(\tG(F)/A_{\tG}(F),\omega^{-1})$ (or $f\in \CC(\tG(F)/A_{\tG}(F),\omega^{-1})$ if the model $(M,H_0)$ is tempered) and $Y\in \cA_{\tX,K}=\cA_{\tX_M,K_M}$, we have defined
	$$I(f,x)=\int_{N(F)}\int_{\tH_0(F)/A_{\tG}^{H_0}(F)}f(x^{-1}hnx)\xi(n)\chi(h)dhdn,\;x\in G(F);$$
	$$I_Y(f)=\int_{H(F)A_{\tG}(F)\back G(F)}I(f,x)\kappa_Y(x) dx.$$
	We also define
	$$I(f)=\int_{H(F)A_{\tG}(F)\back G(F)}I(f,x)dx$$
	whenever this integral is convergent. The next proposition has been proved in Proposition \ref{prop twisted spectral}.
	
	\begin{prop}\label{convergence}
		The integrals defining $I(f,x)$ and $I(f)$ are absolutely convergent for all $f\in C_{c,scusp}^{\infty}(\tilde{G}(F)/A_{\tG}(F),\omega^{-1})$ satisfies \eqref{assumption Rf}. If $(M,H_0)$ is tempered, then both integrals are absolutely convergent for all $f\in \CC_{scusp}(\tilde{G}(F)/A_{\tG}(F),\omega^{-1})$ satisfies \eqref{assumption Rf}.
	\end{prop}

	\begin{rmk}
		In fact we can even prove the convergence without the assumption on $\tilde{R}(f)$. But since we will not use it here, we will postpone the proof of the general convergence to our next paper.
	\end{rmk}
	
	For $t\in \Gamma_{rs}(\tH_0)$, let $S,T,\tS,\tT$ be the same as in the previous section and we let $\CO_t=\CO_{\tS}\in Nil_{reg}(\Fg_{\tS}^{\ast})$ be the orbit associated to $\xi_{\tS}$ as explained at the beginning of Section 4. For a quasi-character $\Theta$ on $\tG(F)$ with central character $\omega^{-1}$, we define
	$$m_{geom,\tH}(\Theta):=\int_{\Gamma_{ell}(\widetilde{\overline{H_0}})} D^{\tH}(t) c_{\Theta,-\CO_t}(t)\chi(t)dt.$$
Here $\Gamma_{ell}(\widetilde{\overline{H_0}})$ is the set of regular elliptic semisimple conjugacy classes of the twisted space $\widetilde{\overline{H_0}}(F)$ equipped with a measure defined in Subsection \ref{Sect measures}. By Lemma \ref{lem bounded germ} and coregular assumption, we know that the integral defining $m_{geom,\tH}(\Theta)$ is absolutely convergent.

	For $f\in \CC_{scusp}(\tilde{G}(F)/A_{\tG}(F),\omega)$, we define
	$$I_{geom}(f)=\nu(\tH_0)m_{geom,\tH}(\Theta_f),\;\nu(\tH_0)=|H_0(F)\cap A_{\tG}(F):A_{\tG}^{H_0}(F)|.$$
	The next theorem is the geometric side of a simple local twisted trace formula in the coregular case.
	
	\begin{thm}\label{geometric side}
		\begin{enumerate}
			\item We have $I(f)=I_{geom}(f)$ for all $f\in C_{c,scusp}^{\infty}(\tilde{G}(F)/A_{\tG}(F),\omega^{-1})$ such that the integral defining $I(f)$ is absolutely convergent. 
			\item If the symmetric pair $(M,H_0)$ is tempered, then $I(f)=I_{geom}(f)$ for all $f\in \CC_{scusp}(\tilde{G}(F)/A_{\tG}(F),\omega^{-1})$ such that the integral defining $I(f)$ is absolutely convergent.
		\end{enumerate}
	\end{thm}

	\begin{proof}
		It is enough to show that the limit of $I_{Y}(f)$ is equal to $I_{geom}(f)$ as $N(Y)$ goes to infinity where $Y$ runs over all the elements in $\cA_{\tX,K}=\cA_{\tX_M,K_M}$ with  $d(Y)>\epsilon N(Y)$. In the previous section we have also defined
		$$J_{Y}(f)=\int_{\Gamma(\widetilde{\overline{H_0}})}D^{\tH}(t)\chi(t)\int_{B_{\tS}(F)\back G(F)} \int_{N_{\tS}(F)}f(x^{-1}tn_S x)v_{B_{\tS},\xi,\iota,Y}(x,n_S) dn_Sdxdt$$
		and we proved in Theorem \ref{thm geometric side center} that 
		\begin{equation}
			|I_{Y}(f)-J_{Y}(f)|\ll N(Y)^{-k}
		\end{equation}
		for every $Y\in \cA_{\tX,K}=\cA_{\tX_M,K_M}$ with  $d(Y)>\epsilon N(Y)$. This implies that the limit of $I_Y(f)$ is equal to the limit of $J_Y(f)$ as $N(Y)$ goes to infinity. Since the function 
		$$Y\mapsto v_{B_{\tS},\xi,\iota,Y}(x,n_S)$$
		is an exponential polynomial with bounded degree and with exponents in a fixed finite set (both independent of $x$ and $n_S$), we know that the limit of $J_Y(f)$ is equal to 
		$$J(f)=\int_{\Gamma(\widetilde{\overline{H_0}})}D^{\tH}(t)\chi(t)\int_{B_{\tS}(F)\back G(F)} \int_{N_{\tS}(F)}f(x^{-1}tn_S x)v_{B_{\tS},\xi,\iota}(x,n_S) dn_Sdxdt.$$
		Here $v_{B_{\tS},\xi,\iota}(x,n_S)$ is defined in Section \ref{secion descent formula}. In particular, as $Y$ goes to infinity, $J_Y(f)$ is a constant independent of $Y$.

		Fix $t\in\Gamma_{rs}(\tH_0)$ and let $S,T,\tS,\tT$ be as in the previous section. Define
		$$J_{\tS}(f)=\int_{\tS(F)/A_{\tG}^{H_0}}D^{\tH}(t)\chi(t)\int_{B_{\tS}(F)\back G(F)} \int_{N_{\tS}(F)}f(x^{-1}tn_S x)v_{B_{\tS},\xi,\iota}(x,n_S) dn_Sdxdt.$$
		Fix $\varepsilon\in (\CA_{B_{\tS}}^{+})^\iota$ in general position. By Corollary \ref{cor descent}, we have the descent formula
		$$v_{B_{\tS},\xi,\iota}(x,n_S)=\sum_{\substack{\tQ\in \mathcal{F}_{B_{\tS}}(\tM(\tS)) \\ \varepsilon\in (\CA_{\tQ}^{+})^\iota}} d_\varepsilon(\tQ)v_{B_{\tS},\xi}^{\tQ,\iota}(x,n_S).$$
		This implies that
		$$J_{\tS}(f)=\sum_{\substack{\tQ\in \mathcal{F}_{B_{\tS}}(\tM(\tS)) \\ \varepsilon\in (\CA_{\tQ}^{+})^\iota}} d_\varepsilon(\tQ) J_{\tS}^{\tQ}(f)$$
		where
		$$J_{\tS}^{\tQ}(f)=\int_{\tS(F)/A_{\tG}^{H_0}}D^{\tH}(t)\chi(t)\int_{B_{\tS}(F)\back G(F)} \int_{N_{\tS}(F)}f(x^{-1}tn_S x)v_{B_{\tS},\xi}^{\tQ,\iota}(x,n_S) dn_Sdxdt.$$

		By Corollary \ref{cor invariant 2} and the assumption that $f$ is strongly cuspidal, we know that $J_{\tS}^{\tQ}(f)=0$ if $\tQ\neq \tG$. If $\tQ=\tG$, then $\varepsilon\in (\CA_{\tG})^\iota$. Since $\varepsilon\in (\CA_{B_{\tS}}^{+})^\iota$ is in general position, we must have $(\CA_{B_{\tS}}^{+})^\iota=(\CA_{G_x})^{\iota}=(\CA_{\tG})^\iota$. This implies that $\tS$ is elliptic. If this is the case, the function $\Gamma_{B_{\tS}}^{\tG,\iota}$ in the definition of $v_{B_{\tS},\xi}^{\tG,\iota}$ is just the function $\Gamma_{B_{\tS}}$. Hence we have
		$$v_{B_{\tS},\xi}^{\tG,\iota}(x,n_S)=\nu(H_0)v_{B_{\tS},\xi}(x,n_S).$$
		Here $\nu(H_0)=|H_0(F)\cap A_{\tG}(F):A_{\tG}^{H_0}(F)|$ comes from the volume of $S(F)A_{\tG}(F):A_{\tG}(F)$ which is equal to the volume of $S(F)/A_{\tG}^{H_0}(F)$ (which is equal to 1) times $|S(F)\cap A_{\tG}(F):A_{\tG}^{H_0}(F)|=|H_0(F)\cap A_{\tG}(F):A_{\tG}^{H_0}(F)|=\nu(\tH_0)$. This implies that $J_{\tS}(f)$ is equal to 0 if $\tS$ is not elliptic, and is equal to
		$$\nu(H_0)\int_{\tS(F)/A_{\tG}^{H_0}}D^{\tH}(t)\chi(t)\int_{B_{\tS}(F)\back G(F)} \int_{N_{\tS}(F)}f(x^{-1}tn_S x)v_{B_{\tS},\xi}(x,n_S) dn_Sdxdt$$
		if $\tS$ is elliptic. Then the theorem follows from Theorem \ref{thm: formula for the germ}.
	\end{proof}

	\subsection{The multiplicity formula}
	We use the same notation as the previous subsection and we assume that the twisted space $\tG$ is just $G$ (i.e. the automorphism $\theta$ is the identity map). Let $\pi$ be an irreducible smooth representation of $G$ with central character $\omega^{-1}$. Define
	$$m(\pi)=\Hom_{H(F)}(\pi,\chi^{-1}\otimes \xi^{-1}),\;m_{geom}(\pi):=\int_{\Gamma_{ell}(\overline{H_0}} D^{H}(t) \chi(t)c_{\Theta,-\CO_t}(t)dt.$$
	
	By \cite{Delcstterm}, we know that the multiplicity $m(\pi)$ is finite. The goal of this subsection is to prove the following multiplicity formula.

	\begin{thm}\label{thm multiplicity formula}
		The multiplicity formula $m(\pi)=m_{geom}(\pi)$ holds for all supercuspidal representations of $G(F)$ with central character $\omega^{-1}$. If $(M,H_0)$ is tempered, then the multiplicity formula holds for all discrete series of $G(F)$ with central character $\omega^{-1}$.
	\end{thm}
	
	\begin{proof}
		For $f\in C_{c,scusp}^{\infty}(G(F)/A_G(F),\omega^{-1})$ (or $f\in \CC_{scusp}(G(F)/A_G(F),\omega^{-1})$ if $(M,H_0)$ is tempered), we have defined
		$$I(f,x)=\int_{N(F)}\int_{H_0(F)/A_{G}^{H_0}(F)}f(x^{-1}hnx)\xi(n)\chi(h)dhdn,\;I(f)=\int_{H(F)A_{\tG}(F)\back G(F)}I(f,x)dx.$$
		in the previous subsection. Proposition \ref{prop twisted spectral}, \eqref{twisted spectral 3}, \eqref{twisted spectral 4} and Theorem \ref{geometric side} implies that 
		$$I(f)=\nu(H_0)m_{geom,H}(\Theta_{f})$$
		for all $f\in {}^\circ\CC(G(F)/A_G(F),\omega^{-1})\cap C_{c}^{\infty}(G(F)/A_G(F),\omega^{-1})$ (or $f\in {}^\circ\CC(G(F)/A_G(F),\omega^{-1})$ if $(M,H_0)$ is tempered). Here ${}^\circ\CC(G(F)/A_G(F),\omega^{-1})$ is the span of matrix coefficients of discrete series of $G(F)$ with central character $\omega^{-1}$. For $f\in {}^\circ\CC(G(F)/A_G(F),\omega^{-1})$, define
		$$I_{spec}(f)=\nu(\tH_0)\sum_{\pi\in \Pi_{disc}(G,\omega^{-1})} \tr(\pi(\bar{f}))m(\pi).$$
		By Theorem 4.1.1 of \cite{BeuGalP}, we have the spectral expansion
		\begin{equation}\label{multiplicity formula equation 1}
			\nu(H_0)m_{geom,H}(\Theta_{f})=I(f)=I_{spec}(f)
		\end{equation}
		for all $f\in {}^\circ\CC(G(F)/A_G(F),\omega^{-1})\cap C_{c}^{\infty}(G(F)/A_G(F),\omega^{-1})$ (or $f\in {}^\circ\CC(G(F)/A_G(F),\omega^{-1})$ if $(M,H_0)$ is tempered). Then the multiplicity formula follows from \eqref{multiplicity formula equation 1}. All we need to do is to let $f$ be the matrix coefficient of a supercuspidal representation (or a discrete series if $(M,H_0)$ is tempered). This finishes the proof of the multiplicity formula.
	\end{proof}

	\begin{rmk}
		Some special cases of the multiplicity formula proved in the above theorem are the multiplicity formulas for the Galois models and the generalized Shalika models proved in \cite{BeuGalP} and \cite{BW} .
	\end{rmk}

	\section{The unitary Shalika model}
	In this section we will prove our main theorems for the unitary Shalika model (i.e. Theorem \ref{thm multiplicity constant introduction}, \ref{thm Shalika distinguished introduction}, and \ref{thm Shalika summation introduction}). In Section 8.1 we will recall the defintion of the models and prove a comparison between the unitary Shalika model and the twisted Shalika model (for general linear groups). Then in Section 8.2 we will prove Theorem \ref{thm multiplicity constant introduction}, \ref{thm Shalika distinguished introduction}, and \ref{thm Shalika summation introduction}.
	
	\subsection{Some comparison}
	Let $Z$ be a $E$-vector space of finite dimension $n\geqslant 1$. Let $Z^{*,c}$ be the conjugate-dual of $Z$ that is the space of $c$-linear forms on $Z$ (a similar notation will be applied later to other vector spaces). Set $V=Z\oplus Z^{*,c}$ and we equip with the nondegenerate Hermitian form
	$$\displaystyle h(v+v^*,w+w^*)=\langle v,w^*\rangle +\langle w,w^*\rangle^c,\;\;\; (v,v^*), (w,w^*)\in Z\oplus Z^{*,c}.$$
	Here $\langle .,.\rangle$ stands for the canonical pairing between $Z$ and $Z^{*,c}$. Let $G=U(V,h)$ be the unitary group associated to this Hermitian form. We define two maximal parabolic subgroups $Q$ and $\overline{Q}$ of $G$ as the stabilizers of the maximal isotropic subspaces $Z$ and $Z^{*,c}$ respectively. Then, $L=Q\cap \overline{Q}$ is a Levi component of $Q$ and restriction to $Z$ induces an isomorphism
	\begin{equation}
		\displaystyle L\simeq \Res_{E/F}GL(Z).
	\end{equation}
	Let $N$ be the unipotent radical of $Q$. Thus $Q=LN$ and restriction to $Z^{*,c}$ induces an isomorphism
	\begin{equation}
		\displaystyle N\simeq \left\{X\in \Hom(Z^{c,*},Z)\mid {}^T X^c=-X \right\}
	\end{equation}
	where ${}^T X^c$ denotes the transpose conjugate of $X$ (seen as a linear endomorphism $Z\to Z^{*,c}$ through the canonical identification $(Z^{*,c})^{*,c}=Z$). We will actually identify the right hand side above with the Lie algebra $\mathfrak{n}$ of $N$ in a way such that the above isomorphism becomes the exponential map.
	
	We henceforth choose two isomorphisms $W_+,W_-:Z\to Z^{*,c}$ satisfying ${}^T W_{\pm}^c=-W_{\pm}$ and such that the corresponding antihermitian forms on $Z$ are not equivalent (there are actually only two equivalence classes of antihermitian forms on $Z$). For $\epsilon\in \{ \pm\}$, we let $H_{0,\epsilon}\subset L\simeq \Res_{E/F}GL(Z)$ be the unitary group associated to $W_\epsilon$, that is the stabilizer of $W_\epsilon$ for the obvious action. Then, $H_{0,\epsilon}(F)$ coincides with the stabilizer in $L(F)$ of the character
	$$\displaystyle \xi_\epsilon:N(F)\to \mathbb{C}^\times,$$
	$$\displaystyle \exp(X)\mapsto \psi(\Tr(W_\epsilon X))\;\; (X\in \mathfrak{n}(F)).$$
	We will henceforth assume, as we may, that $W_{\pm}$ have been chosen so that $H_{0,+}$ is quasi-split.
	
	Set $H_\epsilon=H_{0,\epsilon}\ltimes N$. We extend $\xi_\epsilon$ to a character of $H_\epsilon(F)$ trivial on $H_{0,\epsilon}(F)$. We also fix a character $\chi$ of $E^1=\ker(N_{E/F})$ that we will consider as a character of $H_{0,\epsilon}(F)$ through composition with the determinant $\det: H_{0,\epsilon}(F)\to E^1$ . For a smooth irreducible representation $\pi$ of $G(F)$, we define the multiplicity
	$$\displaystyle m_\epsilon(\pi,\chi):=\dim(\Hom_{H_\epsilon(F)} (\pi,\chi\otimes \xi_\epsilon)).$$

	For $x\in H_{0,\epsilon}(F)_{\elli}$, the centralizer $G_x=Z_G(x)$ is quasi-split and the intersection $N_x:=G_x\cap N$ is a maximal unipotent subgroup of it. Moreover, by restriction $\xi_\epsilon$ induces a nondegenerate character of $N_x(F)$. We let $\CO_x$ be the regular coadjoint nilpotent orbit in $\Fg_{x}^\ast$ associated to it. For any quasi-character $\Theta$ on $G(F)$, we set
	\begin{equation*}
		\displaystyle J_{\epsilon,\chi,geom}(\Theta)=\int_{\Gamma_{\elli}(H_{0,\epsilon})} D^{G}(x)^{1/2} c_{\Theta,\CO_x}(x)\chi(x)^{-1}dx,\; J_{\chi,geom}(\Theta)=J_{+,\chi,geom}(\Theta)+J_{-,\chi,geom}(\Theta).
	\end{equation*}
	By Theorem \ref{thm multiplicity formula}, the multiplicity formula
	$$m_\epsilon(\pi,\chi)= J_{\epsilon,\chi,geom}(\Theta_\pi)$$
	holds for all discrete series.

	Recall that two semisimple regular elements $x,y\in G_{\rs}(F)$ are said to be {\em stably conjugated} if they are conjugated in $G(\overline{F})$ and that a quasi-character $\Theta$ on $G(F)$ is called {\em stable} if it is constant on stable conjugacy classes (that is if $x,y\in  G_{\rs}(F)$ are stably conjugated then $\Theta(x)=\Theta(y)$). If $\Theta$ is stable it is clear that we have $c_{\Theta,\CO_x}(x)=c_{\Theta}(x)$. The following comparison between the geometric sides will be used in our applications.
	
	\begin{prop}\label{prop unitary Shalika stable}
		Assume that $\Theta$ is a stable quasi-character on $G(F)$. Then
		$$\displaystyle J_{\chi,geom,+}(\Theta)=J_{\chi,geom,-}(\Theta).$$
	\end{prop}

	\begin{proof}
		This follows from the following two facts 
		\begin{itemize}
			\item there is a natural measure-preserving bijection $x_+\leftrightarrow x_-$ between the regular elliptic stable conjugacy classes of $H_{0,+}(F)$ and the regular elliptic stable conjugacy classes of $H_{0,-}(F)$;
			\item under the above bijection $x_+\leftrightarrow x_-$, the number of rational conjugacy classes in a regular elliptic stable conjugation class $x_+$ of $H_{0,+}(F)$ is equal to the number of rational conjugacy classes in a regular elliptic stable conjugation class $x_-$ of $H_{0,-}(F)$ and we have $c_{\Theta}(x_+)=c_{\Theta}(x_-)$.
		\end{itemize}
	\end{proof}

	Set $G'=Res_{E/F}\GL(V)$ and let $Q'$, $\overline{Q}'$ be the maximal parabolic subgroups of $G'$ stabilizing the subspaces $Z$ and $Z^{*,c}$ respectively. Then, $L':=Q'\cap \overline{Q}'$ is a Levi component of $Q'$ and we have an isomorphism (given by restriction)
	\begin{equation}
		\displaystyle L'\simeq \Res_{E/F}(\GL(Z)\times \GL(Z^{*,c})).
	\end{equation}
	We fix an isomorphism $W:Z\simeq Z^{*,c}$ satisfying ${}^TW^c=-W$ and we let $H'_0\subset L'$ be the subgroup $\{(h,WhW^{-1})\mid h\in \Res_{E/F}\GL(Z) \}$.
	
	Let $N'$ be the unipotent radical of $Q'$ (so that $Q'=L'N'$). We will identify its Lie algebra $\mathfrak{n}'$ with $\Res_{E/F}\Hom(Z^{*,c},Z)$ and we define a character of $N'(F)$ by
	$$\displaystyle \xi': \exp(X)\in N'(F)\mapsto \psi(\tr_{E/F}(\Tr(W X))),\;\; X\in \mathfrak{n}'(F).$$
	We let $H'=H_0'\ltimes N'$ be the {\it Shalika subgroup}. The character $\chi$ of $E^1$ induces a character $\chi'$ of $E^\times$ by $\chi'(x)=\chi(x/x^c)$ and we will identify $\chi'$ with the character of $H'_0(F)$ given by $(h,WhW^{-1})\mapsto \chi'(\det h))$.
	
	For every $g\in G$, let us denote by $g^\star$ the adjoint linear map with respect to the Hermitian form $h$ on $V$. We define $\theta$ to be the automorphism $g\mapsto (g^\star)^{-1}$ of $G$ and we let $\widetilde{G}=G\theta$ be the nonneutral component of the nonconnected group $G\rtimes \{1,\theta \}$. It is a twisted space in the sense of \S \ref{Sect twisted spaces}. We also set
	$\tilde{Q}'=Q'\theta$, $\tilde{L}'=L'\theta$. These are respectively a twisted parabolic subspace of $\tilde{G}'$ and a Levi component of it. The automorphism $\theta$ preserves $H_0'$ and $H'$ and we let $\tilde{H}_0'=H_0'\theta$, $\tilde{H}'=H'\theta$ be the corresponding twisted spaces. The character $\chi'$ of $H_0'(F)$ being conjugate self-dual, it can be extended to the twisted space $\tH_0'$. We fix such an extension whose value at $\theta$ is equal to 1 and we still denote by $\chi'$.

	For every quasi-character $\widetilde{\Theta}$ on $\widetilde{G}'(F)$, we define
	$$\displaystyle \tilde{J}_{\chi',geom}(\widetilde{\Theta})=\int_{\Gamma_{\elli}(\widetilde{H}_0')} D^{\widetilde{G}'}(x)^{1/2} c_{\widetilde{\Theta}}(x)\chi'(x)^{-1}dx.$$

	Let
	$$\displaystyle Nr: \widetilde{G}'(F)\to G'(F),\;\; g\theta\mapsto g\theta(g)$$
	be the norm map. Recall that an element $x\in \widetilde{G}'_{\rs}(F)$ is said to be {\em $G$-regular} if $Nr(x)$ is regular and that if $x\in \widetilde{G}'_{\rs}(F)$ is $G$-regular, an element $y\in G_{\rs}(F)$ is called a {\em norm of} $x$ if it is conjugated to $Nr(x)$ inside $G'(F)$ (note that $G(F)\subset G'(F)$). Remark that if $y\in G_{\rs}(F)$ is a norm of $x$ and $y'\in G_{\rs}(F)$ is stably conjugated to $y$ then $y'$ is also a norm of $x$ (this is because in $G'(F)$ there is no difference between conjugation and stable conjugation).
	Let $\Theta$ be a stable quasi-character on $G(F)$. We also recall that a quasi-character $\widetilde{\Theta}$ on $\widetilde{G}'(F)$ is said to be a {\em transfer} of $\Theta$ if for every $G$-regular element $x\in \widetilde{G}'_{\rs}(F)$ and every $y\in G_{\rs}(F)$ that is a norm of $x$, we have
	$$\displaystyle D^{\widetilde{G}'}_0(x)^{1/2}\widetilde{\Theta}(x)=D^G(y)^{1/2}\Theta(y).$$
	Here $D_{0}^{\tG'}(x)=D^{\tG'}(x)d_{\tG'}(x)^{-1}$ where $d_{\tG'}(x)$ is defined in Section 1.6 of \cite{WalGGPIII} (it is 1 unless the residue characteristic is 2). To end this subsection, we prove a comparison between $J_{\chi,geom}$ and $\tilde{J}_{\chi',geom}$. This will be used in our application.
	
	\begin{prop}\label{prop comparison unitary Shalika}
		Let $\Theta$ be a stable quasi-character on $G(F)$ and $\widetilde{\Theta}$ be a quasi-character on $\widetilde{G}'(F)$. If $\widetilde{\Theta}$ is a transfer of $\Theta$, we have
		$$\displaystyle J_{\chi,geom}(\Theta)=J_{\chi,geom,+}(\Theta)+J_{\chi,geom,-}(\Theta)=\tilde{J}_{\chi',geom}(\widetilde{\Theta}).$$
	\end{prop}
	
	\begin{proof}
		Recall that (note that since $\Theta$ is stable we have $c_{\Theta,\CO_x}(x)=c_{\Theta}(x)$)
		$$\displaystyle J_{\chi,geom}(\Theta)=\int_{\Gamma_{\elli}(H_{0,+})\cup \Gamma_{\elli}(H_{0,+})} D^{G}(x)^{1/2} c_{\Theta}(x)\chi(x)^{-1}dx,$$
		$$\displaystyle \tilde{J}_{\chi',geom}(\tilde{\Theta})=\int_{\Gamma_{\elli}(\widetilde{H}_0')} D^{\widetilde{G}'}(x)^{1/2} c_{\widetilde{\Theta}}(x)\chi'(x)^{-1}dx.$$
		There is a natural bijection (denoted by $t\leftrightarrow t'$) given by the norm map described above between the regular stable elliptic conjugacy classes of $H_{0,+}(F)\cup H_{0,-}(F)$ and the regular stable elliptic twisted conjugacy classes of $\tilde{H}_{0}'(F)$. For each $t\leftrightarrow t'$, the number of conjugacy classes in $t$ is equal to the number of twisted conjugacy classes in $t'$. By the definition of the character $\chi'$ we know that $\chi(t)=\chi'(t')$ for  $t\leftrightarrow t'$. Moreover, by Section 2.2 of \cite{BeuGGPtemper}, we know that under this bijection we have $dt'=d_{\tilde{H}_{0}'}(t')^{-1}dt$ where $d_{\tilde{H}_{0}'}(t')$ is defined in Section 1.6 of \cite{WalGGPIII} (it is equal to 1 unless the residue characteristic of $F$ is 2). Hence it is enough to show that for all $t\leftrightarrow t'$, we have
		$$D^G(t)^{1/2} c_{\Theta}(t)=d_{\tilde{H}_{0}'}(t')^{-1} D^{\tilde{G}'}(t')^{1/2} c_{\tilde{\Theta}}(t').$$

		We fix a representative of $t$ (resp. $t'$) and by abusing of language we still denoted it by $t$ (resp. $t'$). Let $a$ be the natural isomorphism between $E^{\times}$ and $Z_L(F)$. Also let $W$ be the Weyl group of $G_t(F)$ (which is also equal to the Weyl group of $(G')_{t'}(F)$). By Proposition 4.5.1 of \cite{BeuGGP}, we have
		$$D^G(t)^{1/2} c_{\Theta}(t)=\lim_{\lambda\in F^{\times}\rightarrow 1} \frac{D^G(ta(\lambda))^{1/2} \Theta(ta(\lambda))}{|W|},$$
		$$D^{\tilde{G}'}(t')^{1/2} c_{\tilde{\Theta}}(t')=\lim_{\lambda\in F^{\times}\rightarrow 1}\frac{D^{\tilde{G}'}(t'a(\lambda))^{1/2} \tilde{\Theta}(t'a(\lambda))}{|W|}.$$
		Hence it is enough to show that 
		\begin{equation}\label{prop comparison unitary Shalika equation 1}
			D^G(ta(\lambda))^{1/2} \Theta(ta(\lambda))=d_{\tilde{H}_{0}'}(t')^{-1} D^{\tilde{G}'}(t'a(\lambda))^{1/2} \tilde{\Theta}(t'a(\lambda))    
		\end{equation}
		for $1\neq \lambda\in F^{\times}$ that is close to 1. For $\lambda\neq 1$ that is close to 1, we know that $ta(\lambda)$ (resp. $t'a(\lambda)$) is a regular semisimple element of $G$ (resp. $\tilde{G}'$). Since $t\leftrightarrow t'$, we know that the stable conjugacy class of $ta(\lambda)$ corresponds to the stable conjugacy class of $t'a(\lambda)$. Then \eqref{prop comparison unitary Shalika equation 1} follows from the fact that $\Theta$ and $\tilde{\Theta}$ are the transfer of each other (note that by the definition of the constant $d(\cdot)$ we have $d_{\tilde{H}_{0}'}(t')^2=d_{\tilde{G'}} (t'a(\lambda))$). This proves the proposition.
	\end{proof}
	
	\subsection{The proof of the main results for the unitary Shalika model}
	In this section, we will prove Theorem \ref{thm multiplicity constant introduction}, \ref{thm Shalika distinguished introduction}, and \ref{thm Shalika summation introduction}. We start with Theorem \ref{thm multiplicity constant introduction}.
	
	\begin{thm}\label{thm multiplicity constant}
		\begin{enumerate}
			\item Let $\pi$ be a finite length discrete series of $G(F)$ with central character $\chi^n$. If $\Theta_\pi$ is a stable distribution, then $m_+(\pi,\chi)=m_-(\pi,\chi)$.
			\item Let $\Pi_\phi(G)$ be a discrete L-packet of $G(F)$ with central character $\chi^n$. Then we have
			$$\sum_{\pi\in \Pi_\phi(G)} m_+(\pi,\chi)=\sum_{\pi\in \Pi_\phi(G)} m_-(\pi,\chi).$$
		\end{enumerate}
	\end{thm}
	
	\begin{proof}
		The first part is a direct consequence of the multiplicity formula and Proposition \ref{prop unitary Shalika stable}. The second part follows from the first part together with the fact that the distribution character $\Theta_{\Pi_\phi(G)}=\sum_{\pi\in \Pi_\phi(G)} \Theta_\pi$ is stable.
	\end{proof}
	
	Next we will prove a necessary condition for a discrete L-packet to be distinguished and compute the summation of the multiplicity for some special cases. Let $(G,H_\epsilon,\chi\otimes \xi_{\epsilon})$ be the unitary Shalika model defined in the previous subsection. Let $\Pi_\phi(G)$ be a discrete L-packet of $G$ and let $\Pi_\phi(G')$ be its base change to $G'(F)=\GL_{2n}(E)$. Then $\Pi_\phi(G')$ is an irreducible tempered representation and we can extend it to a unitary twisted representation on $\tilde{G}'(F)$ (denoted by $\widetilde{\Pi_\phi(G')}$) so that $\Theta_{\widetilde{\Pi_\phi(G')}}$ is a transfer of $\Theta_{\Pi_\phi(G)}$. Our goal is to prove the following theorem.
	
	\begin{thm}\label{thm Shalika distinguished}
		With the notation above, the packet $\Pi_\phi(G)$ is $(H_+,\chi\otimes \xi_+)$-distinguished (i.e. $m_+(\pi)\neq 0$ for some $\pi\in \Pi_\phi(G)$) only if $\Pi_\phi(G')$ is distinguished by the Shalika model $(H',\chi'\otimes \xi')$. 
	\end{thm}
	
	\begin{rmk}
		By Theorem \ref{thm multiplicity constant}, we know that the packet $\Pi_\phi(G)$ is $(H_+,\chi\otimes \xi_+)$-distinguished if and only if it is $(H_-,\chi\otimes \xi_-)$-distinguished.
	\end{rmk}
	
	\begin{proof}
		Assume that $\Pi_\phi(G')$ is not distinguished by the Shalika model, we need to show that the packet $\Pi_\phi(G)$ is not $(H_+,\chi\otimes \xi_+)$-distinguished. It is enough to show that 
		$$J_{\chi,geom}(\Theta_{\Pi_\phi(G)})=0$$ 
		where $\Theta_{\Pi_\phi(G)}=\sum_{\pi\in \Pi_\phi(G)}\Theta_\pi$. By Proposition \ref{prop comparison unitary Shalika}, we only need to show that 
		$$\tilde{J}_{\chi',geom}(\Theta_{\widetilde{\Pi_\phi(G')}})=0.$$
		
		Since $\Pi_\phi(G')$ is not distinguished by the Shalika model, by Corollary 1.1 of \cite{M}, we can choose a small neighborhood $\omega$ of $\Pi_\phi(G')$ in $Temp_{ind}(\GL_{2n}(E))$ such that every element in $\omega$ is not distinguished by the Shalika model. By Proposition \ref{prop pesudo coefficient}, we can find a strongly cuspidal function $\tilde{f}$ on $\tilde{G}'(F)$ such that $\tilde{f}$ is supported on $\omega$ and $\Theta_{\tilde{f}}=\Theta_{\widetilde{\Pi_\phi(G')}}$. Hence it is enough to show that $\tilde{J}_{\chi',geom}(\Theta_{\tilde{f}})=0$.

		By our assumption on the support of $\tilde{f}$ and Plancherel formula of Shalika model in \cite{D} we known that $\tilde{R}(\tilde{f})=0$ and hence $\tilde{f}$ satisfies \eqref{assumption Rf}. Applying Proposition \ref{prop twisted spectral} and Theorem \ref{geometric side} to the twisted Shalika model, we have
		$$\tilde{J}_{\chi',geom}(\Theta_{\tilde{f}})=\int_{H'(F)\back G'(F)}\int_{N'(F)}\int_{\tH_0'(F)}\tilde{f}(x^{-1}hnx)\xi'(n)^{-1}\chi'(h)^{-1}dhdndx=\nu(\tilde{H}')\tr(\tilde{R}_{disc}(\tilde{f}))=0.$$
		This finishes the proof of the theorem.
	\end{proof}
	
	\begin{rmk}
		The Plancherel decomposition proved in \cite{D} is for the case when $\chi'=1$. However, by our definition of $\chi'$ we know that $\chi'(-1)=1$ which implies that the character $\chi'$ is a square of another character $\chi''$ of $E^{\times}$. Then we just need to twist the Plancherel decomposition in \cite{D} by the character $\chi''\circ\det$.
	\end{rmk}
	
	Now assume that $\Pi_\phi(G')$ is distinguished by the Shalika model $(H',\chi'\otimes \xi')$. By Corollary 1.1 of \cite{M}, $\Pi_\phi(G')$ is of the form (note that $\chi''$ is a character of $E^{\times}$ with $\chi'=(\chi'')^2$)
	$$\Pi_{\phi}(G')\otimes (\chi''\circ \det)^{-1}=(\tau_1\times \cdots \times \tau_l)\times (\sigma_1\times \sigma_{1}^{\vee})\times \cdots \times (\sigma_m\times \sigma_{m}^{\vee})$$
	where 
	\begin{itemize}
		\item $\tau_i$ is a discrete series of $\GL_{2a_i}(E)$ that is conjugate self-dual, self-dual and of symplectic type. In particular, $a_i$ is even.
		\item $\sigma_j$ is a discrete series of $\GL_{b_i}(E)$ that is conjugate self-dual, but NOT self-dual.
		\item $\tau_i,\sigma_j$ are all distinct.
		\item $\sum_{i=1}^{l} a_i +2\sum_{j=1}^{m}b_j=2n$.
	\end{itemize}
	We will consider the special case when $m=0$. The general case will be consider in our future paper. When $m=0$, by the Plancherel decomposition proved in \cite{D}, $\Pi_\phi(G')$ appears discretely in the $L^2$-space of the Shalika model. Our goal is to prove the following theorem.
	
	\begin{thm}\label{thm Shalika summation}
		With the notation above, we have 
		$$\sum_{\pi\in \Pi_\phi(G)} m_+(\pi,\chi)=\sum_{\pi\in \Pi_\phi(G)} m_-(\pi,\chi)=2^{l-1}.$$
	\end{thm}

	\begin{proof}
		It is enough to show that $J_{\chi,geom}(\Theta_{\Pi_\phi(G)})=2^{l}$ where $\Theta_{\Pi_\phi(G)}=\sum_{\pi\in \Pi_\phi(G)}\Theta_\pi$. Since $J_{\chi,geom}(\Theta_{\Pi_\phi(G)})=\sum_{\pi\in \Pi_\phi(G)} m_+(\pi,\chi)+m_-(\pi,\chi)$ is a non-negative integer, we only need to show that $|J_{\chi,geom}(\Theta_{\Pi_\phi(G)})|=2^{l}$.
		By our assumption of $\Pi_\phi(G')$ and the Plancherel formula of Shalika model \cite{D}, it appears discretely in the $L^2$ space of the Shalika model and hence we can choose a small neighborhood $\omega$ of $\Pi_\phi(G')$ in $Temp_{ind}(\GL_{2n}(E))$ such that $\Pi_\phi(G')$ is the only element in $\omega$ distinguished by the Shalika model. By Proposition \ref{prop pesudo coefficient}, we can find a strongly cuspidal function $\tilde{f}$ on $\tilde{G}'(F)$ such that $\tilde{f}$ is supported on $\omega$, $\Theta_{\tilde{f}}=\Theta_{\widetilde{\Pi_\phi(G')}}$ and $\tr(\overline{\widetilde{\Pi_\phi(G')}}(\tilde{f}))=2^l$. Note that the number $\lvert \Stab(i\cA_{\tG,F}^*,\tau)\rvert^{-1}D(\tau)$ in Proposition \ref{prop pesudo coefficient} is equal to $2^{-l}$ for $\widetilde{\Pi_\phi(G')}$. By Proposition \ref{prop comparison unitary Shalika}, we only need to show that $|\tilde{J}_{\chi',geom}(\Theta_{\widetilde{\Pi_\phi(G')}})|=2^l$. 
		
		By our assumption on the support of $\tilde{f}$ and Plancherel formula of Shalika model in \cite{D} we known that $\tilde{R}(\tilde{f})$ satisfies \eqref{assumption Rf}. By Proposition \ref{prop twisted spectral} and Theorem \ref{geometric side}, it is enough to show that 
		\begin{equation}\label{equation unitary Shalika 1}		
			|\tr(\overline{\widetilde{\Pi_\phi(G')}}(\tilde{f}))\cdot \tr(\theta\langle\overline{\widetilde{\Pi_\phi(G')}}\rangle |M(\overline{\widetilde{\Pi_\phi(G')}}))|=2^{l}.
		\end{equation}
		Since $\overline{\widetilde{\Pi_\phi(G')}}$ is unitary, so is $\theta\langle\overline{\widetilde{\Pi_\phi(G')}}\rangle |M(\overline{\widetilde{\Pi_\phi(G')}})$. As the multiplicity space $M(\overline{\widetilde{\Pi_\phi(G')}})$ is one dimensional, this implies that $|\tr(\theta\langle\overline{\widetilde{\Pi_\phi(G')}}\rangle |M(\overline{\widetilde{\Pi_\phi(G')}}))|=1$. Then \eqref{equation unitary Shalika 1} follows from the facts that $\tr(\overline{\widetilde{\Pi_\phi(G')}}(\tilde{f}))=2^l$. This proves the theorem.
	\end{proof}

	\section{Galois model for classical groups}
	In this section we will prove our main theorems for the Galois models (i.e. Theorem \ref{thm Galois distinguished introduction} and \ref{thm Galois summation introduction}). In Section 9.1 we will prove a comparison between the Galois model for classical groups and the twisted Galois model for general linear groups. Then in Section 9.2 we will prove Theorem \ref{thm Galois distinguished introduction} and \ref{thm Galois summation introduction}.
	
	\subsection{Some comparison}
	Let $H$ be a quasi-split special orthogonal group or a symplectic group defined over $F$ and $G=Res_{E/F}H_E$. Let $G'=Res_{E/F}H_E'$ where $H'=\GL_{2n}$ if $H=\SO_{2n}$ or $\SO_{2n+1}$ and $H'=\GL_{2n+1}$ if $H=\Sp_{2n}$. Let $\theta$ be the involution of $G$ given by $\theta(g)=w(g^t)^{-1}w^{-1}$ where $w$ is the longest Weyl element. Let $\tilde{G}'$ be the non-neutral component of $G'\rtimes \{1,\theta\}$ and let $\tilde{H}'=H'\theta$. Then $\tilde{G}'$ (resp. $\tilde{H}'$) is a twisted space of $G'$ (resp. $H'$). Finally, if $H$ is the even special orthogonal group, let $H_0$ be a quasi-split special orthogonal group that is not a pure inner form of $H$ and such that $G=Res_{E/F}H_E=Res_{E/F}H_{0,E}$ (i.e. the determinanats of the quadratic forms defining $H$ and $H_0$ belong to the same square class in $E^\times/(E^\times)^2$ but belong to different square classes in $F^\times/(F^\times)^2$). If $H=\Sp_{2n}$ or $\SO_{2n}$, let $\chi$ be the trivial character on $H$ (and $H_0$ if $H=\SO_{2n}$) and let $\chi'$ be the trivial character on $H'$. If $H=\SO_{2n+1}$, let $\chi\in \{1,\eta_n\}$  where $\eta_n$ is the composition of the Spin norm character of $\SO_{2n+1}$ with the quadratic character $\eta_{E/F}$. In this case, we let $\chi'=1$ if $\chi=1$ and $\chi'=\eta_n':=\eta_{E/F}\circ \det$ if $\chi=\eta_n$. In both cases, we can extend the character $\chi'$ to the twisted space $\tilde{H}'$ by making it equal to 1 on $\theta$.
	
	For a quasi-character $\Theta$ (resp. twisted quasi-character $\tilde{\Theta}$) on $G(F)$ (resp. $\tilde{G}'(F)$), define
	$$J_{geom}(\Theta)=\int_{\Gamma_{ell}(H)} D^G(t)^{1/2} \Theta(t)\chi(t)^{-1}dt,\;\text{if}\; H=\SO_{2n+1},\;\Sp_{2n},$$
	$$J_{geom}(\Theta)=\int_{\Gamma_{ell}(H)\cup \Gamma_{ell}(H_0)} D^G(t)^{1/2} \Theta(t)\chi(t)^{-1}dt,\;\text{if}\; H=\SO_{2n},$$
	$$\tilde{J}_{geom}(\tilde{\Theta})=\int_{\Gamma_{ell}(\tH')} D^{\tilde{G}'}(t)^{1/2}\tilde{\Theta}(t)\chi'(t)^{-1}dt.$$

	\begin{prop}\label{prop comparison Galois}
		Let $\Theta$ be a stable quasi-character on $G(F)$ and $\tilde{\Theta}$ be a twisted quasi-character on $\tilde{G}(F)$. If $H$ is the even orthogonal group, we fix a Whittaker datum in the definition of the transfer factor so that the element $\eta$ in Section 1.6 of \cite{Waltransfer} is equal to 1. Assume that $\Theta$ and $\tilde{\Theta}$ are the transfer of each other (in the sense of Section 1.6 of \cite{WalGGPIII}). Then we have
		$$2\cdot J_{geom}(\Theta)=\tilde{J}_{geom}(\tilde{\Theta}).$$
	\end{prop}

	\begin{proof}
		When $H$ is the odd orthogonal group, the proposition is a direct consequence of the following four facts 
		\begin{itemize}
			\item There is a natural bijection (denoted by $t\leftrightarrow \tilde{t}$) between the stable regular elliptic  conjugacy classes of $H(F)$  and of $\tilde{H}'(F)$. Under this bijection, we have $d\tilde{t}=d_{\tilde{H'}}(\tilde{t})^{-1}dt$ (Section 1.4 of \cite{WalGGPIII}). 
			\item We have 
			$$D^G(t)^{1/2}\Theta(t)=D^{\tilde{G}}(\tilde{t})^{1/2}d_{\tilde{G'}}(\tilde{t})^{-1/2}\tilde{\Theta}(\tilde{t})=D^{\tilde{G}}(\tilde{t})^{1/2}d_{\tilde{H'}}(\tilde{t})^{-1}\tilde{\Theta}(\tilde{t})$$ 
			for all $t\leftrightarrow \tilde{t}$ (note that the transfer factor between $t$ and any rational twisted conjugacy class in $\tilde{t}$ is trivial by Section 1.10 of \cite{Waltransfer}).
			\item For $t\leftrightarrow \tilde{t}$, the number of $H(F)$-conjugacy classes in $t$ is half of the number of $\tilde{H}'(F)$-conjugacy classes in $\tilde{t}$ (the other half belongs to the pure inner form of the odd special orthogonal group).
			\item For all $t\leftrightarrow \tilde{t}$, we have $\chi(t)=\chi'(\tilde{t})$.
		\end{itemize}
		
		The first three facts are straightforward. For the last one, it is trivial when $\chi=1$. It remains to consider the case when $\chi=\eta_n$ and $\chi'=\eta_n'$. In this case, the stable conjugacy class $t$ (resp. $\tilde{t}$) corresponds to (see Section 1.3 of \cite{Waltransfer})
		$$(F_i,F_{\pm i},t_i)_{1\leq i\leq h}$$
		where
		\begin{itemize}
			\item $F_{\pm i}/F$ is a finite extension of degree $d_i$ with $\sum_{i=1}^{h} d_i=n$;
			\item $F_i$ is a quadratic extension of $F_{\pm i}$;
			\item $t_i\in ker(N_{F_i/F_{\pm i}})$.
		\end{itemize}
		It is easy to see from the definition that  
		$$\chi(t)=\chi'(\tilde{t})=\eta_{E/F}(\Pi_{i=1}^{h}N_{F_i/F}(e_i))$$
		where $e_i$ is any element in $F_{i}^{\times}$ such that $\frac{e_i}{\bar{e_i}}=t_i$ ($\bar{e_i}$ is the conjugation of $e_i$ under the nontrivial element of $Gal(F_i/F_{\pm i})$). This proves the last fact.

		For the rest two cases, the characters $\chi$ and $\chi'$ are trivial. When $H$ is the symplectic group, the proposition is a direct consequence of the following three facts (all of them are straightforward)
		\begin{itemize}
			\item There is a natural bijection (denoted by $t\leftrightarrow \tilde{t}$) between the stable regular elliptic  conjugacy classes of $H(F)$  and of $\tilde{H}'(F)$. Under this bijection, we have $d\tilde{t}=\frac{1}{2}\cdot |2|_F\cdot d_{\tilde{H'}}(\tilde{t})^{-1}dt$ (Section 1.6 of \cite{WalGGPIII}, note that in this case $|T(F)^\theta:T_\theta(F)|=2$ for any maximal elliptic twisted torus $\tT$ of $\tH'(F)$). 
			\item We have 
			$$D^G(t)^{1/2}\Theta(t)=D^{\tilde{G}}(\tilde{t})^{1/2}d_{\tilde{G'}}(\tilde{t})^{-1/2}\tilde{\Theta}(\tilde{t})=D^{\tilde{G}}(\tilde{t})^{1/2}d_{\tilde{H'}}(\tilde{t})^{-1}\tilde{\Theta}(\tilde{t})$$ 
			for all $t\leftrightarrow \tilde{t}$ (note that the transfer factor between $t$ and any rational twisted conjugacy class in $\tilde{t}$ is trivial by Section 1.10 of \cite{Waltransfer}).
			\item For $t\leftrightarrow \tilde{t}$, the number of $H(F)$-conjugacy classes in $t$ is equal to the number of $\tilde{H}'(F)$-conjugacy classes in $\tilde{t}$ divided by $|F^{\times}/(F^\times)^2|$. Moreover, $|F^{\times}/(F^\times)^2|=4\cdot |\frac{1}{2}|_F$.
		\end{itemize}

		When $H$ is the even special orthogonal group, the proposition is a direct consequence of the following three facts 
		\begin{itemize}
			\item There is a natural map $\Gamma_{st,ell}(H\cup H_0)\rightarrow \Gamma_{st,ell}(\tilde{H}')$ (denoted by $t\rightarrow \tilde{t}$) from the stable regular elliptic  conjugacy classes of $H(F)$ and $H_0(F)$  to the stable regular elliptic  conjugacy classes of $\tilde{H}'(F)$. The fiber of each element in the image of this map has exactly two elements (differed by the outer automorphism of the even special orthgonal group).  Under this map, we have $d\tilde{t}=d_{\tilde{H'}}(\tilde{t})^{-1}dt$ (Section 1.6 of \cite{WalGGPIII}).
			\item  We have 
			$$D^G(t)^{1/2}\Theta(t)+D^G(t)^{1/2}\Theta(t')=D^{\tilde{G}}(\tilde{t})^{1/2}d_{\tilde{G'}}(\tilde{t})^{-1/2}\tilde{\Theta}(\tilde{t})=D^{\tilde{G}}(\tilde{t})^{1/2}d_{\tilde{H'}}(\tilde{t})^{-1}\tilde{\Theta}(\tilde{t})$$ 
			for all $t\rightarrow \tilde{t}$ where $t'$ is another element in the fiber of $\tilde{t}$. On the other hand, if $\tilde{t}$ is a stable regular elliptic  conjugacy class of $\tilde{H}'(F)$ that does not belong to the image of $\Gamma_{st,ell}(H\cup H_0)\rightarrow \Gamma_{st,ell}(\tilde{H}')$, then $\tilde{\Theta}(\tilde{t})=0$.
			\item For $t\rightarrow \tilde{t}$, the number of $H(F)$-conjugacy classes (or $H_0(F)$-conjugacy classes) in $t$ is half of the number of $\tilde{H}'(F)$-conjugacy classes in $\tilde{t}$ (the other half belongs to the pure inner form of the even special orthogonal group).
		\end{itemize}
		The first and third facts are straightforward. For the second one, we need to show that the trasnfer factor between $t$ and any rational twisted conjugacy class in $\tilde{t}$ is trivial for all $t\in \Gamma_{st,ell}(H\cup H_0) \rightarrow \tilde{t}\in \Gamma_{st,ell}(\tilde{H}')$. We follow the notation in \cite{Waltransfer}. Under the notation in Section 1.3 of \cite{Waltransfer}, the stable conjugacy class $t$ is of the form
		$$(F_i,F_{\pm i},t_i)_{1\leq i\leq h}$$
		where
		\begin{itemize}
			\item $F_{\pm i}/F$ is a finite extension of degree $d_i$ with $\sum_{i=1}^{h} d_i=n$;
			\item $F_i$ is a quadratic extension of $F_{\pm i}$;
			\item $t_i\in ker(N_{F_i/F_{\pm i}})$.
		\end{itemize}
		A rational twisted conjugacy class in $\tilde{t}$ is of the form
		$$(F_i,F_{\pm i},t_i,c_i)_{1\leq i\leq h}$$
		where $c_i\in F_{i}^{\times}/Im(N_{E_i/F_i})$. Next we need to describe how does $(F_i,F_{\pm i},t_i,c_i)_{1\leq i\leq h}$ behave under base change. There are three types:
		
		\begin{itemize}
			\item[Type 1] If $E$ is not contained in $E_i$, then $(F_i/F_{\pm i},t_{i},c_{i})$ becomes $(E_i,E_{\pm i},t_{i},1)$ where $E_i=F_i\otimes_F E$ and $E_{\pm i}=F_{\pm i}\otimes_F E$. Here we view $t_{i}$ as an element of $ker(N_{E_i/E_{\pm i}})$ via the canonical embedding from $ker(N_{F_i/F_{\pm i}})$ to $ker(N_{E_i/E_{\pm i}})$.
			\item[Type 2] If $E$ is not contained in $F_{\pm i}$ and $E\subset F_i$, then $(F_i/F_{\pm i},t_{i},c_{i})$ becomes $(F_i\oplus F_i,F_i,(t_{i},t_{i}^{-1}),1)$. 
			\item[Type 3] If $E\subset F_{\pm i}$, let $F_{\pm i}=E[x]/f(x)$ and we define the field $F_{\pm i}'$ to be $F_{\pm i}'=E[x]/\bar{f}(x)$ where $f\mapsto \bar{f}$ is the conjugation map on $E[x]$ induced by the non-trivial element of $Gal(E/F)$. Similarly we can also define the field $F_i'$ which will be a quadratic extension of $F_{\pm i}'$. Moreover, we have a natural isomorphism (denoted by $x\mapsto \bar{x}$) between $ker(N_{F_i/F_{\pm i}})$ and $ker(N_{F_i'/F_{\pm i}'})$. Then $(F_i/F_{\pm i},t_i,c_{i})$ becomes $(F_i/F_{\pm i},t_i,c_i)\cup (F_i'/F_{\pm i}',\bar{t}_{i},\bar{c}_{i})$.
		\end{itemize}
		We decompose the set $I$ into $I_1\cup I_2\cup I_3$ where $I_j$ consists of those $i\in I$ such that $(F_i,F_{\pm i})$ belongs to Type j above. 
		
		Then if we view $t$ as a stable conjugacy class of $G(F)$, it is of the form 
		$$(E_i,E_{\pm i},t_{i})_{i\in I_1}\cup (F_i\oplus F_i,F_i,(t_{i},t_{i}^{-1}))_{i\in I_2}\cup ((F_i/F_{\pm i},t_i)\cup (F_i'/F_{\pm i}',\bar{t}_{i}))_{i\in I_3}.$$
		Similarly, if we view a rational twisted conjugacy class in $\tilde{t}$ as a rational twisted conjugacy class of $\tilde{G}'(F)$, it is of the form
		$$(E_i,E_{\pm i},t_{i},1)_{i\in I_1}\cup (F_i\oplus F_i,F_i,(t_{i},t_{i}^{-1}),1)_{i\in I_2}\cup ((F_i/F_{\pm i},t_i,c_i)\cup (F_i'/F_{\pm i}',\bar{t}_{i},\bar{c}_{i}))_{i\in I_3}.$$
		
		For $i\in I_1$, the quadratic character $\eta_{E_i/E_{\pm i}}$ is trivial on $F_{\pm i}^{\times}$. For $i\in I_2$, the quadratic character $\eta_{F_i\oplus F_i/F_i}$ is the trivial character. For $i\in I_3$, the natural isomorphism from $F_{\pm i}^{\times}$ to $(F_{\pm i}')^{\times}$ maps the quadratic character $\eta_{F_i/F_{\pm i}}$ to the quadratic character $\eta_{F_i'/F_{\pm i}'}$. Combining these three facts with the definition of transfer factor in Section 1.10 \cite{Waltransfer} (note that we have choosen the Whittaker datum so that the number $\eta$ in loc. cit. is equal to 1), we know that  the transfer factor between $t$ and any rational twisted conjugacy class in $\tilde{t}$ is trivial. This finishes the proof of the proposition.
		
	\end{proof}
	
	\subsection{The proof of the main theorem for Galois model}

	We start with a necessary condition for a discrete L-packet to be distinguished (i.e. Theorem \ref{thm Galois distinguished introduction}).
	
	\begin{thm}\label{thm Galois distinguished}
		Let $H=\Sp_{2n},\SO_{2n}$ or $\SO_{2n+1}$, $G=Res_{E/F}H_E$, $\chi=1$ if $H=\Sp_{2n}$ or $\SO_{2n}$, and $\chi\in \{1,\eta_n\}$ if $H=\SO_{2n+1}$. Let $\Pi_\phi(G)$ be a discrete L-packet of $G(F)$ and $\Pi_\phi(G')$ be the endoscopic transfer of the L-packet to the general linear group $G'=\GL_a(E)$ (here $a=2n$ if $H=\SO_{2n}$ or $\SO_{2n+1}$ and $a=2n+1$ if $H=\Sp_{2n}$). Then the packet $\Pi_\phi(G)$ is distinguished (i.e. $m(\pi,\chi)\neq 0$ for some $\pi\in \Pi_\phi(G)$) only if $\Pi_\phi(G')$ is $(\GL_a(F),\chi')$-distinguished. Here $\chi'=1$ if $\chi=1$ and $\chi'=\eta_n':=\eta_{E/F}\circ \det$ if $\chi=\eta_n$.
	\end{thm}

	For the summation of the multiplicities (i.e. Theorem \ref{thm Galois summation introduction}), assume that $\Pi_{\phi(G')}$ is $(\GL_a(F),\chi')$-distinguished. By Theorem 4.2 of \cite{Mat}, $\Pi_\phi(G')$ is of the form 
	$$\Pi_{\phi}(G')=(\tau_1\times \cdots \times \tau_l)\times (\sigma_1\times \bar{\sigma}_{1})\times \cdots \times (\sigma_m\times \bar{\sigma}_{m})$$
	where 
	\begin{itemize}
		\item $\tau_i$ is a discrete series of $\GL_{a_i}(E)$ that is conjugate self-dual. Moreover, if $(H,\chi)=(\SO_{2n+1},\eta_n)$, $\tau_i$ is self-dual of symplectic type; otherwise, $\tau_i$ is self-dual of orthogonal type.
		\item $\sigma_j$ is a discrete series of $\GL_{b_i}(E)$ that is NOT conjugate self-dual. Moreover, if $(H,\chi)=(\SO_{2n+1},\eta_n)$, $\sigma_j$ is self-dual of symplectic type; otherwise, $\sigma_j$ is self-dual of orthogonal type.
		\item $\tau_i,\sigma_j$ are all distinct.
		\item $\sum_{i=1}^{l} a_i +2\sum_{j=1}^{m}b_j=a$.
	\end{itemize}
	We will consider the special case when $m=0$. The general case will be consider in our future paper. When $m=0$, by the Plancheral formula for the Galois model proved in \cite{Beu}, $\Pi_\phi(G')$ appears discretely in the $L^2$ space of the Galois model $(\GL_a(E),\GL_a(F),\chi')$.
	
	\begin{thm}\label{thm Galois summation}
		With the notation above, if $H$ is the symplectic group or the odd special orthogonal group, we have
		$$\sum_{\pi\in \Pi_\phi(G)} m(\pi,\chi)=2^{l-1}.$$
		If $H$ is the even special orthogonal group, we let $H_0$ be another even special orthogonal group as in the previous subsection. We use $m_0(\pi,\chi)$ to denote the multiplicity for the model $(G,H_0,\chi)$. Then we have
		$$\sum_{\pi\in \Pi_\phi(G)} m(\pi,\chi)+m_0(\pi,\chi)=2^{l-1}.$$
	\end{thm}

	The proof of the above two theorem is almost  the same as the unitary Shalika model case. The only differences is to replace Proposition \ref{prop comparison unitary Shalika} by Proposition \ref{prop comparison Galois}, and to replace the Plancherel formula for Shalika model in \cite{D} by the Plancheral formula for the Galois model in \cite{Beu}. We will skip the details here.

	\appendix
	
	\section{Projections of finitely generated convex sets}\label{appendix convex sets}

	In this appendix, we state and prove a decomposition result for orthogonal projections of finitely generated convex sets that is directly inspired from \cite[Appendix]{Art88}.
	
	Let $\Fa$ be a real Euclidean space with scalar product denoted by $(.,.)$. A subset $\CC\subset \Fa$ is a {\em finitely generated convex} set if it satisfies one of the following equivalent properties:
	\begin{itemize}
		\item $\CC$ is a finite intersection of {\em half-spaces} (by which we mean subsets of the form $\{X\in \Fa\mid (Y,X)\leqslant c \}$ for some $Y\in \Fa$ and $c\in \BR$);
		
		\item There exists finite subsets $\{ X_i\mid i\in I\}$ and $\{v_j\mid j\in J \}$ of $\Fa$ such that
		$$\displaystyle \CC=Conv\{X_i\mid \in I \}+\sum_{j\in J} \BR_+ v_j;$$
		
		\item $\CC\times \{ 1\}$ is the intersection of $\Fa\times \{1 \}$ with a finitely generated cone in $\Fa\oplus \BR$.
	\end{itemize}
	Remark that any finitely generated convex set is automatically closed.
	
	Let $\CC\subset \Fa$ be a finitely generated set. We let
	$$\displaystyle C:=\{X\in \Fa\mid \CC+\BR_+X=\CC \}$$
	be its {\em asymptotic cone} and denote by
	$$\displaystyle C^\vee:=\{Y\in \Fa\mid (Y,X)\leqslant 0 \; \forall X\in C \}$$
	be the corresponding {\em dual cone}.
	
	For $H\in \CC$, we define its {\em tangent cone} and {\em normal cone} respectively by
	$$\displaystyle T_{\CC}(H):= \BR_+ (\CC-H),\;\; N_{\CC}(H):=T_{\CC}(H)^\vee=\{Y\in \Fa\mid (Y,X)\leqslant 0 \; \forall X\in T_{\CC}(H) \}.$$
	Note that both are finitely generated cones.
	
	\begin{lem}\label{lem1 appendix}
		Let $\CC_1,\CC_2\subset \Fa$ be two finitely generated convex sets. Then, for $H\in \CC_1\cap \CC_2$ we have
		$$\displaystyle T_{\CC_1\cap \CC_2}(H)=T_{\CC_1}(H)\cap T_{\CC_2}(H) \mbox{ and } N_{\CC_1\cap \CC_2}(H)=N_{\CC_1}(H)+N_{\CC_2}(H).$$
	\end{lem}
	
	\begin{proof}
		The first equality is obvious by definition. The second follows from the first one and the relation $(C_1\cap C_2)^\vee=C_1^\vee+C_2^\vee$ that holds for every finitely generated cones $C_1,C_2\subset \Fa$.
	\end{proof}
	
	A face of $\mathcal{C}$ is its nonempty intersection with a supporting (affine) hyperplane i.e. a subset of $\CC$ of the form
	$$\displaystyle F=\{H\in \CC \mid (\lambda,H)=c\}$$
	where $\lambda\in \Fa$ and $c\in \BR$ are such that $(\lambda,H)\leqslant c$ for every $H\in \CC$ with equality for at least one such $H$. Note that we allow $\lambda=0$ so that $\CC$ is a face of itself. We let $\CF(\CC)$ be the set of faces of $\CC$.
	
	To every face $F\in \CF(\CC)$, we associate the subspace $\Fa^F$ that is the span of $F-X_F$ for any $X_F\in F$. Moreover, the normal cone $N_{\CC}(X_F)$ is independent of $X_F$ when the latter is chosen in the relative interior $\mathring{F}$ of $F$ (that is its interior relative to $F+\Fa^F$) and we shall denote by $\Fa_F^+$ the relative interior of this normal cone. For any $X_F\in \mathring{F}$ we have
	$$\displaystyle \Fa_F^+=\{Y\in \Fa\mid (Y,H-X_F)<0 \; \forall H\in \CC-F \}.$$
	
	\begin{lem}\label{lem2 appendix}
		We have a partition
		$$\displaystyle C^\vee=\bigsqcup_{F\in \CF(\CC)} \Fa_F^+.$$
	\end{lem}
	
	\begin{proof}
		We lemma reduces to the three following claims:
		\begin{num}
			\item the cones $\Fa_F^+$, $F\in \CF(\CC)$, are mutually disjoint.
		\end{num}
		Indeed, let $F,F'\in \CF(\CC)$ be distinct faces and choose $X_F\in \mathring{F}$, $X_{F'}\in \mathring{F}'$. Without loss of generality we may assume that $F'\not\subset F$ so that $X_{F'}\notin F$. Then, for $H\in \Fa_F^+\cap \Fa_{F'}^+$ we have
		$$\displaystyle (H,X_F-X_{F'})\leqslant 0 \mbox{ and } (H,X_{F'}-X_F)<0.$$
		As these two inequalities are incompatible this shows that $\Fa_F^+\cap \Fa_{F'}^+=\emptyset$.
		
		\begin{num}
			\item For every $F\in \CF(\CC)$ we have $\Fa_F^+\subset C^\vee$.
		\end{num}
		Indeed, for $X_F\in \mathring{F}$ we have $C\subset T_{\CC}(X_F)$ hence $\Fa_F^+\subset N_{\CC}(X_F)\subset C^\vee$.
		
		\begin{num}
			\item For every $X\in C^\vee$, there exists $F\in \CF(\CC)$ such that $X\in \Fa_F^+$.
		\end{num}
		Indeed, the function $Y\mapsto (X,Y)$ attains a maximum on $\CC$ (as follows from the fact that $\CC$ can be written as the sum of a convex hull of finitely many points and $C$) say $c\in \BR$. Then,
		$$\displaystyle F:=\{H\in \CC\mid (X,H)=c \}$$
		is a face of $\CC$ and $X\in \Fa_F^+$.
	\end{proof}
	
	
	Let $\Fb$ be a vector subspace of $\Fa$ and $\Fb^\perp$ be its orthogonal complement. Denote  by $p: \Fa\to \Fb$ and $p^\perp:\Fa\to \Fb^\perp$ the two orthogonal projections. For $\xi\in \Fb^\perp$ we set
	$$\displaystyle \CF(\CC,\xi)=\{F\in \CF(\CC)\mid \xi\in p^\perp(\Fa_F^+) \}.$$
	
	\begin{prop}\label{prop projection convex sets}
		Assume that $\dim(C^\vee+\Fb)=\dim(\Fa)$ and $\xi\in p(\Fa_C^+)$ is in general position. (More precisely, we require that for every face $F\in \CF(\CC)$ with $\dim p^\perp(\Fa_F^+)<\dim \Fb^\perp$ we have $\xi\notin p^\perp(\Fa_F^+)$). Then, we have:
		\begin{enumerate}[(i)]
			\item $p$ induces a bijection between
			$$\displaystyle \bigcup_{F\in \CF(\CC,\xi)} F \mbox{ and } p(\CC).$$
			
			\item For $F_1,F_2\in \CF(\CC,\xi)$ we have
			$$\displaystyle p(F_1)\cap p(F_2)=p(F_1\cap F_2).$$
		\end{enumerate} 
	\end{prop}

	\begin{proof}
		Let $H\in \CC$ and consider the intersection
		$$\displaystyle \CC_{H,\Fb}:=(H+\Fb^\perp)\cap \CC.$$
		It is a finitely generated set with asymptotic cone $C_{\Fb}:=\Fb^\perp\cap C$ and dual cone $C_{\Fb}^\vee=\Fb+C^\vee$. As $\xi\in C_{\Fb}^\vee$, by the decomposition of Lemma \ref{lem2 appendix}, there exists a unique face $F_{H,\Fb}\in \CF(\CC_{H,\Fb})$ such that $\xi\in \Fa_{F_{H,\Fb}}^+$. Take $X\in \mathring{F}_{H,\Fb}$. Then, by Lemma \ref{lem1 appendix}, we have
		$$\displaystyle \Fa_{F_{H,\Fb}}^+=N_{\CC_{H,\Fb}}(X)^\circ=N_{\CC}(X)^\circ+\Fb=\Fa_F^++\Fb$$
		where $F\in \CF(\CC)$ is the unique face such that $X\in \mathring{F}$. In particular, we see that $F\in \CF(\CC,\xi)$. This already shows that $p$ induces a surjection
		$$\displaystyle \bigcup_{F\in \CF(\CC,\xi)} F\to p(\CC).$$
		To prove that this map is also injective, it only remains to check that for the face $F_{H,\Fb}$ is a singleton. But, by the assumption that $\xi$ is in general position, $\xi\in \Fb+\Fa_F^+=\Fa_{F_{H,\Fb}}^+$ implies that $\dim(\Fa_{F_{H,\Fb}}^+)=\dim(\Fa)$ i.e. that $F_{H,\Fb}$ is reduced to one extreme point of $\CC_{H,\Fb}$. This proves (i). Note that (ii) is a direct consequence of (i).
	\end{proof}

	\section{Howe's conjecture for twisted weighted orbital integrals}\label{Appendix Howe}
	
	The purpose of this appendix is to establish an analog of Howe's conjecture \cite{Howe} for weighted orbital integrals on a $p$-adic twisted space. This result is needed for the proof of Theorem \ref{thm: formula for the germ}. A similar extension of Howe's conjecture to weighted orbital integrals was established for honest reductive groups by Arthur \cite{Art94} based on his local trace formula but to the best of our knowledge Arthur's argument hasn't been extended to twisted spaces. The proof presented here is a direct adaptation of the work of Barbasch and Moy \cite{BarM} which has the advantage of allowing non-Archimedean local fields of arbitrary characteristics. Actually, the reasoning in \cite{BarM} extends without much effort to twisted spaces but for the comfort of the reader, as well as for the authors own edification, we reproduce below with some details Barbasch and Moy's beautiful argument.
	
	The first section of this appendix contains the precise statement of the ``Howe conjecture for twisted weighted orbital integrals'' as well as a reduction to a certain property of twisted Hecke modules (Proposition \ref{prop appendix Howe}). The proof of this proposition will be given in Section \ref{Sect 3 appendix Howe} following very closely the paper \cite{BarM}. The intermediate Section \ref{Sect 2 appendix Howe} aims to collect necessary material on Bruhat-Tits buldings and the Moy-Prasad filtrations.
	
	\subsection{The statement}\label{Sect 1 appendix Howe}
	
	We will freely use the basic notations introduced in Chapter \ref{preliminaries} for twisted spaces and their subgroups. The main objects under consideration will a priori depend on the choices of Haar measures. However, the precise normalization of those are completely irrelevant for the main result of this appendix and we will therefore assume that Haar measures have been fixed every time they appear in a formula.
	
	Let $(G,\tG)$ be a twisted reductive space defined over $F$. Let $K$ be a special maximal compact subgroup of $G(F)$ so that for every parabolic subgroup $P\subset G$ we have an Iwasawa decomposition $G(F)=P(F)K$. Let $\tM\subset \tG$ be a Levi subspace. For every parabolic subspace $\tP\in \CP(\tM)$ and $g\in G(F)$ we set $H_{\tP}(g):=H_{\tM}(m_{P}(g))$ where $g=m_P(g)u_P(g)k_P(g)$ is an arbitrarily chosen decomposition with $m_P(g)\in M(F)$, $u_P(g)\in N_P(F)$ and $k_P(g)\in K$. Note that for $g\in G(F)$, the convex hull $Conv\{H_{\tP}(g)\mid \tP\in \CP(\tM) \}$ is contained in a translate of the subspace $\Fa_{\tM}^{\tG}$. We let $v_{\tM}(g)$ be the volume of that convex hull with respect to a given Haar measure on $\Fa_{\tM}^{\tG}$.
	
	For $\gamma\in \tM(F)\cap \tG_{rs}(F)$ and $f\in C_c^\infty(\tG(F))$, we can form the {\em weighted orbital integral}
	$$\displaystyle WO_{\tM}(\gamma,f)=WO_{\tM}^{\tG}(\gamma,f)= \int_{G_\gamma(F)\backslash G(F)} f(g^{-1}\gamma g) v_{\tM}(g) dg$$
	for some choice of invariant measure on $G_\gamma(F)\backslash G(F)$.
	
	For any subset $\Omega\subseteq \tM(F)$, we denote by $WO_{\tM}(\Omega)$ the span of the linear functionals $f\in C_c^\infty(\tG(F))\mapsto WO_{\tM}(\gamma,f)$ for $\gamma\in \Omega\cap \tG_{rs}(F)$. Also, for $J\subset G(F)$ a compact open subgroup we set
	$$\displaystyle \widetilde{\CH}_J=\widetilde{\CH}^{\tG}_J=C_c(\tG(F)/J).$$
	Howe's conjecture for twisted weighted orbital integrals can now be stated as follows:
	
	\begin{thm}\label{theo appendix Howe}
		Assume that $\Omega\subseteq \tM(F)$ is compact modulo conjugation and let $J\subset G(F)$ be a compact-open subgroup. Then, the restriction of $WO_{\tM}(\Omega)$ to $\widetilde{\CH}_J$ is finite dimensional.
	\end{thm}

	We will now reduce the above theorem to a statement about twisted Hecke modules. The space $\widetilde{\CH}=\widetilde{\CH}^{\tG}=C_c^\infty(\tG(F))$ is a bimodule over the Hecke algebra $\CH:= C_c^\infty(G(F))$ for the action by left and right convolution (after fixing a Haar measure on $G(F)$). For $J\in CO(G)$, we denote by $[\CH,\widetilde{\CH}_J]$ the span of the commutators $[\phi,f]=\phi\ast f-f\ast \phi$ for $(\phi,f)\in \CH\times \widetilde{\CH}_J$. Equivalently, $[\CH,\widetilde{\CH}_J]$ is the span of the differences ${}^g f-f$ for $(g,f)\in G(F)\times \widetilde{\CH}_J$ where we have set ${}^gf(\gamma):=f(g^{-1}\gamma g)$, $\gamma\in \tG(F)$.
	
	For any subset $\Omega\subset \tG(F)$, we define let $\widetilde{\CH}_J(\Omega)^c$ be the subspace of functions $f\in \widetilde{\CH}_J$ that are supported in $\tG(F)\setminus \Omega J$. We will establish the theorem through the following more technical statement.
	
	\begin{prop}\label{prop appendix Howe}
		Let $\Omega\subset \tG(F)$ be a subset that is compact modulo conjugation and $J\subset G(F)$ be a compact-open subgroup. Then, there exists an open subgroup $J'\subset J$ such that the quotient space
		$$\displaystyle \widetilde{\CH}_J/([\CH,\widetilde{\CH}_{J'}]\cap \widetilde{\CH}_J+\widetilde{\CH}_J(\Omega)^c)$$
		is of finite dimension.
	\end{prop}
	
	To end this section, we now explain why Proposition \ref{prop appendix Howe} implies Theorem \ref{theo appendix Howe}. Let
	$$\displaystyle \tau_\Omega=\tau^{\tG}_\Omega: \widetilde{\CH}\to C^\infty(\Omega\cap \tG_{rs}(F))$$
	be the linear map sending $f\in \widetilde{\CH}$ to the function
	$$\displaystyle \gamma\in \Omega\cap \tG_{rs}(F)\mapsto WO_{\tM}(\gamma,f).$$
	We need to show that $\displaystyle \tau_\Omega(\widetilde{\CH}_J)$ is finite dimensional. The proof is by induction on the semisimple rank of $G$ and thus we assume that the result already holds for all the proper Levi subspaces of $\tG$.
	
	For $g\in G(F)$, $f\in \widetilde{\CH}$ and $\gamma\in \tM(F)\cap \tG_{rs}(F)$ we have the splitting formula \cite[Proposition 2.9.4 (4)]{LW}
	\begin{equation}\label{splitting formula appendix Howe}
		\displaystyle \frac{D^{\tG}(\gamma)^{1/2}}{D^{\tM}(\gamma)^{1/2}}WO_{\tM}(\gamma,{}^g f)=\sum_{\tQ \in \CF(\tM)} WO_{\tM}^{\tL_Q}(\gamma,f_{g,\tQ})
	\end{equation}
	where $\tL_Q$ stands for the unique Levi factor of $\tQ$ containing $\tQ$, $f_{g,\tQ}\in C_c^\infty(\tL_Q(F))$ is the function given by
	$$\displaystyle f_{g,\tQ}(\tm)=\delta_{\tP}(\tm)^{1/2}\int_{K\times N_Q(F)} f(k^{-1}\tm uk) u_{\tQ}(kg^{-1})dudk,\;\; \tm\in \tL_Q(F)$$
	and
	$$\displaystyle u_{\tQ}(h):=\int_{\Fa_{\tQ}^{\tG}} \Gamma_{\tQ}(H,-H_{\tQ}(h)) dH, \mbox{ for } h\in G(F).$$
	
	Let $\Omega^G\subseteq \tG(F)$ be the union of all $G(F)$-conjugates of $\Omega$. As $\Omega$ is compact modulo $M(F)$-conjugation, $\Omega^G$ is similarly compact modulo $G(F)$-conjugation. Let $J'\subset J$ be as in Proposition \ref{prop appendix Howe} with $\Omega^G$ instead of $\Omega$. Then, for every $\tQ\in \CF(\tM)$ we can find a compact-open subgroup $J_Q\subset L_Q(F)$ such that
	\begin{equation}\label{eq1 appendix Howe}
		\displaystyle f_{g,\tQ}\in \widetilde{\CH}^{\tL_Q}_{J_Q}, \mbox{ for every } f\in \widetilde{\CH}_{J'} \mbox{ and } g\in G(F).
	\end{equation}
	
	From \eqref{splitting formula appendix Howe} and \eqref{eq1 appendix Howe} we deduce that
	\begin{equation*}
		\displaystyle \tau_\Omega([\CH,\widetilde{\CH}_{J'}])\subseteq \sum_{\tG\neq \tQ\in \CF(\tM)} \tau_\Omega^{\tL_Q}(\widetilde{\CH}^{\tL_Q}_{J_Q}).
	\end{equation*}
	By the induction hypothesis, this implies that $\tau_\Omega([\CH,\widetilde{\CH}_{J'}])$ has finite dimension. Furthermore, since the distribution $WO_{\tM}(\gamma,.)$ for $\gamma\in \Omega\cap \tG_{rs}(F)$ is supported in the $G(F)$-conjugacy class of $\gamma$, the image of $\widetilde{\CH}_J(\Omega^G)^c$ by $\tau_\Omega$ is zero. By Proposition \ref{prop appendix Howe}, it follows that $\tau_\Omega(\widetilde{\CH}_J)$ is also of finite dimension Q.E.D.
	
	\subsection{Bruhat-Tits building and the Moy-Prasad filtrations}\label{Sect 2 appendix Howe}
	
	Let $\CB$ be the restricted Bruhat-Tits building of $G$. It is a polysimplicial complex carrying polysimplicial actions of $G(F)$ and $\tG(F)$ that are compatible in the sense that
	$$\displaystyle (g\gamma g')\cdot x=g\cdot(\gamma \cdot (g'\cdot x)), \mbox{ for every } (g,\gamma,g')\in G(F)\times \tG(F)\times G(F) \mbox{ and } x\in \CB.$$
	Moreover, these actions factor through $G(F)/Z_G(F)$ and $\tG(F)/Z_G(F)$ respectively and the resulting actions are proper. Picking, for some minimal Levi subgroup $M_0\subset G$, a scalar product on $\Fa_{M_0}^G$ that is invariant under $Norm_{\tG(F)}(M_0)$ yields a distance function
	$$\displaystyle dist:\CB\times \CB\to \BR_+$$
	on $\CB$ that is uniquely geodesic and invariant under $\tG(F)$. For $x,y\in \CB$ we shall denote by $[x,y]$ the unique geodesic joining $x$ and $y$ and we set $]x,y[=[x,y]\setminus \{x,y \}$.
	
	By a {\em chamber} of $\CB$, we shall mean the closure of a facet of maximal dimension.
	
	For $\gamma\in \tG(F)$, the {\em displacement function} $d_\gamma:\CB\to \BR_+$ is defined by
	$$\displaystyle d_\gamma(x)=dist(x,\gamma\cdot x),\;\; x\in \CB.$$
	
	For any chamber $C\subset \CB$ and $\gamma\in \tG(F)$ we set
	$$\displaystyle d_C(\gamma):=\inf_{x\in C} d_\gamma(x).$$
	We also set
	$$\displaystyle d(\gamma):=\inf_{x\in \CB} d_\gamma(x), \mbox{ for } \gamma\in \tG(F).$$
	Note that, as the set of all chambers cover $\CB$, we have
	\begin{equation}\label{eq5 appendix Howe}
		\displaystyle d(\gamma)=\inf_{C\subset \CB} d_C(\gamma)
	\end{equation}
	where the infimum is taken over the set of all chambers in $\CB$.
	
	\begin{lem}\label{lem2 appendix Howe}
		\begin{enumerate}[(i)]
			\item For any chamber $C\subset \CB$, the value set $d_C(\tG(F))$ is a closed discrete subset of $\BR_+$ and can thus be linearly ordered
			$$\displaystyle d_C(\tG(F))=\{0=r_0<r_1<r_2<\ldots \}.$$
			\item Let $x\in \CB$ and $\gamma\in \tG(F)$ be such that $d_\gamma(x)>d(x)$. Then, for every $y\in ]x,\gamma\cdot x[$ we have
			$$\displaystyle d_\gamma(y)<d_\gamma(x).$$
			\item Let $x\in \CB$ and $\gamma\in \tG(F)$. Then, if $d_\gamma$ attains a local minimum at $x$ we have $d_\gamma(x)=d(x)$.
			\item The function $d:\tG(F)\to \BR_+$ is invariant by $G(F)$-conjugation and locally constant.
		\end{enumerate}
	\end{lem}
	
	\begin{proof}
		(i) The statement is equivalent to $d_C(\tG(F))\cap [0,R]$ being finite for every $R>0$. The set
		$$\displaystyle B(C,R):=\{x\in \CB\mid \inf_{y\in C} dist(x,y)\leqslant R \}$$
		is compact. Thus, by the properness of the action of $\tG(F)/Z_G(F)$, the set
		$$\displaystyle \{\gamma\in \tG(F)\mid d_C(\gamma)\leqslant R \}$$
		is compact modulo $Z_G(F)$. However, $d_C$ is also right invariant by the pointwise stabilizer $G_C$ of $C$ which is an open subgroup of $G(F)$ containing $Z_G(F)$. The claim follows.
		
		(ii) By the triangular inequality, and since $y\in ]x,\gamma x[$, we have
		\[\begin{aligned}
			\displaystyle d_\gamma(y)=dist(y,\gamma y) & \leqslant dist(y,\gamma x)+dist(\gamma x,\gamma y) \\
			& =dist(y,\gamma x)+dist(x,y)=dist(x,\gamma x)=d_\gamma(x).
		\end{aligned}\]
		Moreover, as $\CB$ is uniquely geodesic, equality holds if and only if $\gamma x\in ]y,\gamma y[$ or equivalently $\gamma x\in [x,\gamma^2 x]$. Assume by way of contradiction that $d_\gamma(y)=d_\gamma(x)$. Then, we have $\gamma^n x\in [\gamma^{n-1}x,\gamma^{n+1}x]$ for every $n\geqslant 1$ from which if follows that the geodesics $[x,\gamma x],\ldots, [\gamma^{n-1}x,\gamma^n x]$ piece together to form the geodesic $[x,\gamma^n x]$ and so
		$$\displaystyle dist(x,\gamma^n x)=nd_\gamma(x),\mbox{ for every } n\geqslant 0.$$
		On the other hand, as $d_\gamma(x)>d(\gamma)$, we can find $z\in \CB$ such that $d_\gamma(z)<d_\gamma(x)$. By the triangular inequality again, we have
		\[\begin{aligned}
			\displaystyle nd_\gamma(x)=dist(x,\gamma^n x) & \leqslant dist(x,z)+dist(z,\gamma^n z)+dist(\gamma^n z,\gamma^n x) \\
			& \leqslant 2dist(x,z)+nd_\gamma(z)
		\end{aligned}\]
		for each $n\geqslant 0$. Letting $n$ goes to infinity leads to a contradiction. Therefore, $d_\gamma(y)<d_\gamma(x)$ and we are done.
		
		(iii) This follows from (ii), noting that if $\gamma x\neq x$ every neighborhood of $x$ meets $]x,\gamma x[$.
		
		(iv) It is clear that $d$ is invariant by $G(F)$-conjugation. Let us show that it is also locally constant. Let $\gamma\in \tG(F)$. Then, by (i) and \eqref{eq5 appendix Howe} there exists a chamber $C\subset \CB$ such that $d(\gamma)=d_C(\gamma)$. As $C$ is compact and $x\in C\mapsto d_\gamma(x)$ is continuous, $d_\gamma$ attains its infimum on $C$ and therefore $d(\gamma)=d_\gamma(x)$ for some $x\in C$. Let $J\subset G(F)$ be a compact-open subgroup that fixes pointwise some neighborhood of $x$ in $\CB$. Then, for each $k\in J$ the function $d_{\gamma k}$ attains a local minimum at $x$ from which we deduce, by (iii), that $d(\gamma k)=d_{\gamma k}(x)=d_\gamma(x)=d(\gamma)$ i.e. $d$ is constant on the coset $\gamma J$.
	\end{proof}
	
	Let $x\in \CB$. For every real number $r\geqslant 0$, Moy and Prasad have defined an open-compact subgroup $K_{x,r}\subset G(F)$ with the following properties:
	\begin{num}
		\item For every $s\geqslant r\geqslant 0$ and $x\in \CB$, we have $K_{x,s}\subseteq K_{x,r}$;
		\item For each $x\in \CB$, $\bigcap_{r\geqslant 0} K_{x,r}=\{ 1\}$;
		\item\label{eq6 appendix Howe} For any $x\in \CB$, $r\geqslant 0$ and $\gamma\in \tG(F)$, we have $K_{\gamma\cdot x,r}=\Ad_\gamma(K_{x,r})$;
		\item\label{eq7 appendix Howe} There exists $h>0$ such that for each integer $n\geqslant 0$ and $x\in \CB$, $K_{x,nh}$ only depends on the facet $F$ containing $x$;
		\item\label{eq8 appendix Howe} For $r>0$ and $x,y,z\in \CB$ such that $y\in [x,z]$ we have
		$$\displaystyle K_{y,r}\subset K_{x,r}\cdot K_{z,r}.$$
	\end{num}
	By \eqref{eq7 appendix Howe}, for any chamber $C\subset \CB$ we may define $K_{C,n}$ as $K_{x,nh}$ for any point $x$ in the relative interior of $C$.
	
	\subsection{Proof of Proposition \ref{prop appendix Howe}}\label{Sect 3 appendix Howe}
	
	Fix a chamber $C\subset \CB$. It suffices to prove Proposition \ref{prop appendix Howe} for $J=J'=K_{C,n}$ and $n$ large enough. In particular, we will assume that $n$ is sufficiently large that $J$ fixes pointwise all the chambers $C'\subset \CB$ with $C\cap C'\neq \emptyset$.

	By Lemma \ref{lem2 appendix Howe}(i), we can write
	$$\displaystyle d_C(\tG(F))=\{0=r_0<r_1<r_2<\ldots \}$$
	and for each $i\geqslant 0$ we let $\widetilde{\CH}_{J,\leqslant r_i}$ be the subspace of $f\in \widetilde{\CH}_J$ which are supported in the set of $\gamma\in \tG(F)$ with $d_C(\gamma)\leqslant r_i$. Then, $i\mapsto \widetilde{\CH}_{J,\leqslant r_i}$ is an increasing and exhaustive filtration of $\widetilde{\CH}_J$ and since the action of $\tG(F)/Z_G(F)$ on the building is proper, the quotients $\widetilde{\CH}_{J,\leqslant r_i}/\widetilde{\CH}_{J,\leqslant r_i}\cap \widetilde{\CH}_{J}(\Omega)^c$ are finite dimensional. Therefore, it suffices to check that for $i$ sufficiently large we have
	\begin{equation}\label{eq2 appendix Howe}
		\displaystyle \widetilde{\CH}_{J,\leqslant r_{i}}\subseteq \widetilde{\CH}_{J,\leqslant r_{i-1}}+[\CH,\widetilde{\CH}_J]+\widetilde{\CH}_J(\Omega)^c.
	\end{equation}
	We will actually show that the above inclusion holds as soon as
	\begin{equation}\label{eq3 appendix Howe}
		\displaystyle r_i>\sup_{\gamma\in \Omega} d(\gamma).
	\end{equation}
	We thus assume that the above inequality is satisfied. The quotient $\widetilde{\CH}_{J,\leqslant r_{i}}/\widetilde{\CH}_{J,\leqslant r_{i-1}}$ is spanned by the images of the functions $\mathbf{1}_{\gamma J}$ for $\gamma\in \tG(F)$ with $d_C(\gamma)= r_i$ and it suffices to show that for such $\gamma$, $\mathbf{1}_{\gamma J}\in \widetilde{\CH}_{J,\leqslant r_{i-1}}+[\CH,\widetilde{\CH}_J]+\widetilde{\CH}_J(\Omega)^c$. For this we distinguish two cases:
	
	First we assume that $d_C(\gamma)=d(\gamma)$. Let $x\in C$ be such that $d_C(\gamma)=d_\gamma(x)$. Then, as $J$ fixes pointwise a neighborhood of $x$, for every $k\in J$ the displacement function $d_{\gamma J}$ attains a local minimum at $x$ and therefore, by Lemma \ref{lem2 appendix Howe}(iii), we have
	$$\displaystyle d(\gamma k)=d_{\gamma k}(x)=d_\gamma(x)=d_C(x)=r_i.$$
	By \eqref{eq3 appendix Howe}, this implies $\gamma J\subseteq \tG(F)\setminus \Omega$ and therefore $\mathbf{1}_{\gamma J}\in \widetilde{\CH}_J(\Omega)^c$.
	
	Assume now that $d_C(\gamma)>d(\gamma)$. Let again $x\in C$ be such that $d_C(\gamma)=d_\gamma(x)$. Then by Lemma \ref{lem2 appendix Howe}(ii), we have $[x,\gamma x]\cap C=\{x \}$. Let $y\in ]x,\gamma x[$ be sufficiently close to $x$ so that if $F$ denotes the facet containing $y$ we have $x\in \overline{F}$ (where $\overline{F}$ denotes the closure of $F$). We can find a chamber $D$ containing $y$ and $x'\in \mathring{C}$ (where $\mathring{C}$ denotes the interior of $C$) such that $[x',\gamma x']\cap \mathring{D}\neq \emptyset$. By \eqref{eq6 appendix Howe} and \eqref{eq8 appendix Howe}, for any $y'\in [x',\gamma x']\cap \mathring{D}$, we have
	\begin{equation*}
		\displaystyle K_{D,n}= K_{y',nh}\subset K_{x',nh}K_{\gamma x',nh}=J Ad_{\gamma}(J).
	\end{equation*}
	Let $k_1,\ldots, k_\ell\in Ad_{\gamma}(J)$ be such that
	$$\displaystyle J K_{D,n}=\bigcup_{i=1}^\ell J k_i.$$
	Then, since $k_i^{-1}\gamma J=\gamma J$ for any $i$, we have
	\[\begin{aligned}
		\displaystyle \mathbf{1}_{\gamma J K_{D,n}}=\sum_{i=1}^{\ell} \mathbf{1}_{k_i^{-1}\gamma J k_i}=\sum_{i=1}^{\ell}{}^{k_i^{-1}}\mathbf{1}_{\gamma J}.
	\end{aligned}\]
	This shows that
	\begin{equation}\label{eq4 appendix Howe}
		\displaystyle \ell^{-1} \mathbf{1}_{\gamma J K_{D,n}}-\mathbf{1}_{\gamma J}\in [\CH,\widetilde{\CH}_J].
	\end{equation}
	Furthermore, as $J$ fixes $D$ pointwise (since $D\cap C$ contains $x$ and is therefore nonempty) and $y\in D\cap ]x,\gamma x[$, by Lemma \ref{lem2 appendix Howe}(ii) we have
	\begin{equation*}
		\displaystyle d_D(\gamma k)=d_D(\gamma)\leqslant d_\gamma(y)<d_\gamma(x)=d_C(x)=r_i
	\end{equation*}
	for every $k\in J$. Let $g\in G(F)$ be such that $C=gD$. Then, we have $gK_{D,n}g^{-1}=K_{C,n}=J$ and, by the above,
	$$\displaystyle d_C(g\gamma kg^{-1})=d_D(\gamma k)<r_i$$
	for every $k\in J$. This shows that the function
	$$\displaystyle {}^g \mathbf{1}_{\gamma J K_{D,n}}=\mathbf{1}_{g\gamma J K_{D,n}g^{-1}}=\mathbf{1}_{g\gamma Jg^{-1} J}$$
	belongs to $\widetilde{\CH}_{J,\leqslant r_{i-1}}$. Combining this with \eqref{eq4 appendix Howe}, we deduce that $\mathbf{1}_{\gamma J}\in \widetilde{\CH}_{J,\leqslant r_{i-1}}+[\CH,\widetilde{\CH}_J]$ and the claim follows.

	\flushright Rapha\"el Beuzart-Plessis \\
	Aix Marseille Univ, \\
	CNRS, I2M, \\
	Marseille, France \\
	raphael.beuzart-plessis@univ-amu.fr
	
	\vspace{2mm}
	
	Chen Wan \\
	Department of Mathematics \& Computer Science \\
	Rutgers University – Newark \\
	Smith Hall, 101 Warren St., \\
	Newark, NJ 07102, USA \\
	chen.wan@rutgers.edu
\end{document}